\documentclass[a4paper,11pt]{book}

\usepackage{lmodern}
\usepackage[utf8]{inputenc}
\usepackage[T1]{fontenc}
\usepackage{array}
\usepackage{graphicx}
\usepackage{amssymb, amsfonts, amsthm, amsmath}
\usepackage[right=3.5cm,left=3.5cm]{geometry}
\usepackage[french]{babel}
\usepackage{hyperref}

\usepackage{bibunits}

\newtheoremstyle{thmfrench}%
     {\topsep}%
     {\topsep}%
     {\itshape}
     {}
     {\bfseries}
     {}
     {\newline}
     {\thmname{#1}\thmnumber{ #2}\thmnote{ -- #3}}

\newtheoremstyle{deffrench}%
     {\topsep}%
     {\topsep}%
     {}
     {}
     {\bfseries}
     {}%
     {\newline}
     {\thmname{#1}\thmnumber{ #2}\thmnote{ -- #3}}

\newtheoremstyle{rmkfrench}%
     {\topsep}%
     {\topsep}%
     {}
     {}
     {\itshape}
     {}%
     {.5em}
     {}

\theoremstyle{thmfrench}
\newtheorem{theoreme}{Théorème}[chapter]
\newtheorem{corollaire}[theoreme]{Corollaire}
\newtheorem{lemme}[theoreme]{Lemme}
\newtheorem{proposition}[theoreme]{Proposition}
\newtheorem{propriete}[theoreme]{Propriété}

\theoremstyle{deffrench}
\newtheorem{definition}[theoreme]{Définition}
\newtheorem{exercice}{Exercice}[chapter]

\theoremstyle{rmkfrench}
\newtheorem{remarque}{Remarque}[chapter]

\newcommand{\N}{\mathbb{N}}
\newcommand{\Z}{\mathbb{Z}}
\newcommand{\C}{\mathbb{C}}
\newcommand{\Q}{\mathbb{Q}}
\newcommand{\R}{\mathbb{R}}
\renewcommand{\d}{\mathrm{d}}

\DeclareMathOperator{\E}{\mathbb{E}}
\DeclareMathOperator{\probabilitejolie}{\mathbb{P}}
\renewcommand{\P}{\probabilitejolie}
\newcommand{\F}{\mathcal{F}}

\newcommand{\crochet}[1]{{\langle #1 \rangle}}
\newcommand{\e}{\mathbf{e}}
\newcommand{\n}{\mathbf{n}}
\newcommand{\W}{\mathbb{W}}
\renewcommand{\L}{\mathcal{L}}

\newcommand{\ind}[1]{\mathbf{1}_{\{#1\}}}
\newcommand{\indset}[1]{\mathbf{1}_{#1}}

\newcommand{\floor}[1]{{\left\lfloor #1 \right\rfloor}}

\renewcommand{\tilde}{\widetilde}

\newcommand{\conv}{{\underset{n\to +\infty}{\longrightarrow}}}

\newcommand{\egaldistr}{{\overset{(d)}{=}}}
\DeclareMathOperator{\sgn}{sgn}
\newcommand{\norme}[1]{{\Vert #1 \Vert}}

\newcommand{\sh}{\mathrm{sh}}
\newcommand{\ch}{\mathrm{ch}}
\renewcommand{\coth}{\mathrm{coth}}

\newcommand{\convd}{{\underset{\delta \to 0}{\longrightarrow}}}

\newcommand{\wconvt}{{\underset{t \to +\infty}{\Longrightarrow}}}
\renewcommand{\leq}{\leqslant}

\begin{document}

\bibliographyunit[\part]
\bibliographystyle*{plain}
\renewcommand{\proofname}{Solution}

\title{Exercices sur les temps locaux de semi-martingales continues et les excursions browniennes}
\author{Bastien Mallein et Marc Yor}
\date{\today}

\maketitle

~
\vfill
\begin{flushright}
\textit{Saumons bondissants}

\textit{cherchent la source du torrent ;}

\textit{processus de poissons.}
\end{flushright}
\vfill
\newpage

\bigskip

\begin{center}
 \textbf{Préface}
 
 \textit{de Marc Yor}

 1949-2014
\end{center}

\bigskip

Depuis le tout début du XX${}^\text{e}$ siècle, l'étude des processus stochastiques est un domaine très actif de la recherche en mathématiques. Les motivations peuvent avoir des raisons purement \og internes\fg{} : la description, par exemple, de la complexité de la courbe brownienne plane. D'autres sont plus \og externes \fg{} : de nombreux domaines scientifiques font de plus en plus appel à certains aspects des processus aléatoires qui semblent -- même dans des situations a priori très éloignées des probabilités -- correspondre en un certain sens à des propriétés liées aux processus stochastiques, ou simplement à l'aléa. A posteriori, la raison de cette identification est simple : l'activité scientifique amène très souvent à \og mesurer \fg{} tel ou tel phénomène. De ces \og mesures \fg{} aux probabilités, il n'y a qu'un pas. Ajoutez le facteur temps, et vous voilà dans le domaine des processus stochastiques.

Parmi ces processus, le mouvement brownien --- dont l'étude mathématique a été initiée dès 1900, avec la thèse de Bachelier, entre autres travaux --- a joué, et joue encore, un rôle primordial. Ceci peut s'expliquer par le fait que le mouvement brownien est l'objet limite quasi-universel qui apparaît dans le théorème central limite, lorsqu'on fait agir le temps. Pour étudier le mouvement brownien, et ses processus satellites, il y a bien sûr de nombreux moyens : les grands anciens, Paul Lévy par exemple, mais aussi P. Erdös, A. Dvoretzky, S. Kakutani, J. Taylor, utilisaient pour l'essentiel le caractère gaussien du processus, ainsi que l'indépendance et l'homogénéité des accroissements, à la lumière de leurs éclairs de génie. D'autres pionniers, tout aussi impressionnants, tels que A. Kolmogorov ou M. Kac, utilisaient les connexions entre équation de la chaleur, problème de Dirichlet, équation de Poisson ou de Sturm-Liouville, etc. avec le mouvement brownien, pour exprimer les lois, ou des résultats limites sur certaines fonctionnelles du mouvement brownien.

Enfin, après la seconde guerre mondiale, \og l'ère Itô \fg{}, avec son calcul stochastique directement sur les trajectoires browniennes, a commencé. Le merveilleux petit livre de H.P. Mc Kean : \textit{Stochastic Integrals} (1969) a définitivement popularisé le calcul d'Itô, puis des traités entiers sont apparus sur ce sujet. Le calcul d'Itô -- stricto sensus -- demande de travailler avec des fonctions régulières, i.e. de classe $\mathcal{C}^2$, pour, après composition avec le mouvement brownien, écrire explicitement la décomposition de Doob-Meyer en semi-martingale du processus ainsi obtenu. Avec la formule de Tanaka (1963), il est apparu que la régularité de classe $\mathcal{C}^2$ n'était pas nécessaire, et que les temps locaux browniens (en un niveau donné, le temps passé par un mouvement brownien au voisinage infinitésimal de ce niveau) pouvaient être discutés à l'aide du calcul d'Itô-Tanaka, magistralement unifié et généralisé par Paul-André Meyer (1976) dans son \textit{Cours sur les Intégrales stochastiques}.

C'est encore avec un coup de génie que Kiyoshi Itô, en 1970, jette les bases de la théorie des excursions qui, considérée sous l'angle d'étude du mouvement brownien par exemple, remplace la complexité des trajectoires browniennes par la simplicité des escaliers Poissonniens ; le \og truc \fg{} est qu'en toute généralité, il faudrait considérer une infinité de tels escaliers, mais pour un problème donné, seuls un ou plusieurs de ces escaliers -- processus de Poisson -- suffisent, d'où l'impression de simplicité miraculeuse ressentie au début des applications de cette théorie par les chercheurs tels que D. Williams, J. Pitman... Les notes que l'on trouvera développées dans ce volume, réparties en deux ensembles (Partie 1 : Temps locaux de semi-martingales ; Partie 2 : Excursions browniennes) ont pour objectif de familiariser le lecteur, au travers de la résolution d'exercices, à chacun de ces deux domaines.

Pour que ce travail soit profitable, il est demandé au lecteur de connaître raison-nablement le calcul d'Itô \og régulier \fg{}. Les exercices de ce volume ont été élaborés, année après année, par le second auteur, soit à la suite de lectures d'articles présentant, parfois avec des méthodes très différentes, telle ou telle propriété brownienne, soit simplement pour illustrer le contenu de son cours de DEA (anciennement), de M2 aujourd'hui. Le premier auteur a résolu ces exercices et en a organisé la synthèse, de façon économique et néanmoins --- espérons-le --- très lisible. Les chapitres ont été conçus pour créer un aller-retour permanent entre les principaux résultats du cours et les exercices corrigés, afin que la compréhension des uns renforce celle des autres. C'est ainsi que de nombreuses solutions d'exercices données ici offrent un aperçu de la façon de prouver certains des théorèmes rappelés plus haut.

\medskip
\medskip

Nous souhaitons que ces notes puissent être un palier qui permettra à un étudiant de M2 en probabilités de se hisser aux cimes autrement élevées des recherches actuelles, ainsi qu'en témoignent les articles de J.F Le Gall, P. Biane, J. Bertoin et W. Werner dans le volume \textit{A tribute for K. Itô}, SPA, 2010, où chacun de ces auteurs montre comment utiliser la théorie des excursions dans des domaines aussi variés que la théorie des processus SLE (Schramm-Loewner Evolution), l'étude des processus de coagulation et de fragmentation, les processus \og libres \fg{} et processus matriciels, ainsi que -- pour citer un thème en plein essor -- les cartes planaires.

Notre but sera donc atteint si les connaissances et la maîtrise des objets présentés ici, ainsi acquises par le lecteur, lui permettent d'accéder aux résultats de tout premier plan que nous venons de citer. Enfin, nous lui souhaitons autant de plaisir à résoudre ces exercices qu'ils nous en ont donné à les rédiger.
\tableofcontents
\newpage

\part{Temps locaux de semi-martingales}

Cette partie I est composée de sept chapitres, chacun d'entre eux correspondant à une leçon sur les temps locaux de semi-martingales continues, et/ou de diffusions réelles. Au début de chaque chapitre nous rappelons les théorèmes principaux de la leçon correspondante.

Les temps locaux des semi-martingales permettent une approche \og{}en moyenne\fg{} de l'étude de ces semi-martingales. Au lieu de s'intéresser à la valeur d'une semi-martingale en un instant donné, on s'intéresse au temps qu'elle a passé à un niveau donné. On pourra ainsi, par exemple, déterminer le temps passé par un mouvement brownien au voisinage de 0, ou d'un autre niveau. Grâce à ces nouveaux objets, nous pourrons formuler plusieurs extensions de la formule d'Itô, permettant d'obtenir un certain nombre de résultats probabilistes fins, par exemple concernant l'unicité des solutions de certaines équations différentielles stochastiques à coefficients peu réguliers.

Nous nous attacherons également à étudier les temps locaux comme des objets naturels, et à donner des identités en loi avec d'autres processus stochastiques bien connus. Des exemples de ces résultats sont les théorèmes de Lévy et  de Pitman, qui permettent de voir le temps local en zéro comme le supremum passé d'un mouvement brownien ou l'infimum futur d'un processus de Bessel de dimension 3. Nous parlerons également des théorèmes de Ray-Knight, qui permettent de comprendre les temps locaux comme semi-martingales indexées par l'espace.

Finalement, nous donnerons également une courte introduction aux temps locaux d'intersection, domaine qui se rapproche des recherches actuelles en probabilités. Dans cette introduction, nous étudierons quelques propriétés du temps passé par un mouvement brownien en dimension 2 ou 3 à se recouper, et déterminer de bons ordres de grandeur pour des quantités intéressantes, telles que l'aire de la saucisse de Wiener.

\chapter[Introduction aux temps locaux]{Introduction aux temps locaux de semi-martingales continues}

Nous allons ici nous attacher à la définition et aux premières propriétés des temps locaux\footnote{Paul Lévy utilisait le terme~: mesure de voisinage, pour le temps local en un niveau donné}, qui permettent de mesurer le temps passé au voisinage d'un niveau donné par le processus. Il existe différentes théories des temps locaux, mais nous nous intéresserons ici uniquement aux temps locaux de semi-martingales continues.

Rappelons qu'une semi-martingale continue est un processus stochastique qui peut se décomposer de manière unique en la somme d'un processus à variations finies et d'une martingale locale, toutes deux supposées continues. La formule d'It\^o permet justement de donner cette décomposition pour les images de semi-martingales par des fonctions de classe $\mathcal{C}^2$. En tentant d'étendre cette formule à des fonctions qui ne sont pas nécessairement de classe $\mathcal{C}^2$, on définit toute une famille de processus à variations finies, famille associée à la semi-martingale initiale. Il se trouve que ces processus ont une interprétation immédiate dans le cas du mouvement brownien, ils sont les limites, en un certain sens, du temps passé au voisinage d'un point $x$ avant l'instant $t$, pour tout $x \in \R$. Par conséquent, ces processus ont été appelés temps locaux de semi-martingales, car ils représentent en quelque sorte l'échelle de temps ressentie au voisinage du point $x$. Nous allons ici simplement donner les premières définitions et propriétés de ces processus, et nous nous familiariserons avec ces définitions avec quelques exercices de difficulté croissante. Ceux-ci permettent de démontrer de nombreux résultats qui seront utilisés par la suite.

\section{Théorèmes principaux}

Nous allons commencer par donner une définition de ce qu'est un temps local, qui permet de construire de nombreuses sortes de temps locaux donnant chacune des informations différentes sur le comportement des processus associés au voisinage infinitésimal d'une configuration.

\begin{definition}
Un processus continu $(X_t, t \geq 0)$ admet des temps locaux par rapport au processus croissant $(A_t, t \geq 0)$ et à la mesure de Radon $\mu$, si il existe une fonction mesurable $(\Lambda^x_t, x \in \R,t \geq 0)$  telle que pour toute fonction $f$ continue bornée et pour tout $t \geq 0$,~on a~:
\[\int_0^t \d A_s f(X_s) = \int_\R \mu(\d x) \Lambda^x_t f(x).\]
\end{definition}

Pour toute semi-martingale continue $X$, on note $\crochet{X}$ sa variation quadratique, qui est un processus croissant vérifiant
\[ \forall t \in \R_+, \crochet{X}_t = \lim_{n \to +\infty} \frac{1}{n} \sum_{j=1}^n \left( X_{jt/n} - X_{(j-1)t/n} \right)^2 \text{ en probabilité.}  \]
Dans la plus grande partie de la suite, les temps locaux que nous étudierons seront les temps locaux de semi-martingales, dont le théorème suivant garantit l'existence. 

\begin{theoreme}
Une semi-martingale continue $(X_t, t \geq 0)$ admet des temps locaux $(\Lambda_t^ x ; t \geq 0, x \in \R)$ par rapport au processus croissant $(\crochet{X}_t, t \geq 0)$ et à la mesure de Lebesgue, que l'on appelle simplement les temps locaux (de semi-martingales continues) de $X$.

En d'autre termes, il existe une fonction mesurable $(\Lambda^x_t, x \in \R, t \geq 0)$ telle que pour toute fonction $f$ continue bornée et $t \geq 0$,
\[ \int_0^t f(X_s) \d \crochet{X}_s = \int_\R \Lambda^x_t f(x) \d x. \]
\end{theoreme}

Pour tout $x \in \R$, on notera
\[ x^+ = \begin{cases} x & \text{si } x \geq 0\\ 0 & \text{sinon} \end{cases} \quad \text{ et } \quad
x^- = \begin{cases} 0 & \text{si } x \geq 0 \\ -x & \text{sinon} \end{cases}\]
les parties positive et négative de $x$.

\begin{theoreme}[Formule de Tanaka-Meyer]
Pour toute semi-martingale continue, il existe une modification $(L^x_t)$ de $\Lambda$ qui soit conjointement continue en $t$ et càdlàg en $x$. On a de plus~:
\[\dfrac{1}{2}L^x_t = (X_t-x)^+-(X_0-x)^+ - \int_0^t \ind{X_s>x} \d X_s.\]

Soit $X_t=X_0+M_t+A_t$ une semi-martingale continue avec $M$ une martingale locale et $A$ un processus à variations finies. Les sauts du temps local de $X$ sont donnés par la formule~:
\[L^x_t-L^{x-}_t = 2 \int_0^t \d A_s \ind{X_s=x}.\]

En particulier, si $X$ est une martingale locale continue, la fonction de ses temps locaux admet une version conjointement continue en $t$ et $x$.
\end{theoreme}

\begin{remarque}
Ce théorème est en quelque sorte une \og{}formule d'It\^o\fg{} appliquée à la fonction $x \mapsto x^+$, qui n'est pas de classe $\mathcal{C}^2$. Dans le chapitre suivant, nous donnerons un sens précis à cette remarque. Notons de plus que, à $x$ fixé, nous avons accès au temps local en $x$ d'une manière qui permet, par exemple, l'expression d'intégrales contre le temps local, ce qui n'était pas donné par la définition.
\end{remarque}

\begin{corollaire}
Soit $X$ une semi-martingale et $(L^x_t)$ une version continue en $t$ et càdlàg en $x$ de la fonction de ses temps locaux ; on a pour tout $x \in \R$ et $t \geq 0$~:
\[L^x_t = \lim_{\epsilon \to 0} \dfrac{1}{\epsilon} \int_0^t \d \crochet{X}_s \ind{x\leq X_s <x+\epsilon}.\]
En particulier, pour tout $x \in \R$, $(L^x_t,t \geq 0)$ est constant sur $\{t \geq 0 : X_t \neq x\}$.
\end{corollaire}

\begin{remarque}
Si $B$ est un mouvement brownien, et $(L^x_t)$ la version continue de la fonction de ses temps locaux, pour tout $\epsilon >0$,
$\displaystyle \int_0^t \d s \ind{-\epsilon<B_s<\epsilon} = \int_{-\epsilon}^\epsilon \d x L^x_t.$
Par conséquent, $\displaystyle L^0_t = \lim_{\epsilon\to 0} \frac{1}{2\epsilon} \int_0^t \d s \ind{-\epsilon <B_s<\epsilon},$ et on retrouve ainsi la notion de temps passé au voisinage de 0 par $B$, d'où le terme~: temps local en 0.

Remarquons également que cette définition de temps local en 0, qui nécessite une renormalisation en $\frac{1}{\epsilon}$, est nécessaire. En effet, par formule de Fubini, on a~:
\[ \E \left[ \int_0^t \ind{B_s=0} \d s \right] = \int_0^t \P(B_s=0) \d s = 0,\]
par conséquent $\displaystyle \int_0^t \ind{B_s=0} \d s = 0$ p.s. le mouvement brownien passe presque sûrement un temps de mesure nulle en 0.
\end{remarque}

Nous finissons cette introduction aux temps locaux par une remarque d'intér\^et, appelée propriété fondamentale des temps locaux qui étend la remarque précédente.

\begin{proposition}[Propriété fondamentale des temps locaux]
\label{pro_fondamentale}
Soit $Y$ une semi-martingale et $L$ la fonction de ses temps locaux. Pour tout $y \in \R$, la mesure aléatoire $\d_tL^y_t$ est portée par $\{t\geq 0:Y_t=y\}$.
\end{proposition}

\section{Exercices}

\begin{exercice}
Soit $B$ un mouvement brownien standard issu de 0 et $a,b>0, a \neq b$. On pose $Y_t=|B_t|$ et $Z_t=aB_t^+-bB_t^-$.
\begin{enumerate}
  \item Montrer que $x \mapsto L^x_t(Y)$ et $x \mapsto L^x_t(Z)$ sont discontinues en 0.
  \item Calculer $L^0_t(Y)$ et $L^0_t(Z)$ en fonction de $L^0_t(B)$.
\end{enumerate}
\end{exercice}

\begin{proof}
\textit{1.} Par définition des temps locaux, on a pour toute fonction $f$ continue à support compact~:
\[\int_0^t \d\crochet{Y}_s f(Y_s) = \int_{\R} \d x f(x) L^x_t(Y)\]
On obtient, en remplaçant $Y_t$ par $|B_t|$~:
\[\int_{\R} \d x f(x) L^x_t(Y) = \int_0^t \d s f(|B_s|) = \int_{\R} \d x f(|x|) L^x_t(B).\]

On obtient alors pour tout $x>0$, $L^x_t(Y)=L^x_t(B)+L^{-x}_t(B)$, et pour tout $x<0$, $L^x_t(Y)=0$. Par continuité à droite de $x \mapsto L^x_t(Y)$, on obtient $L^0_t(Y)=2L^0_t(B)$, et $L^{0-}_t(Y)=0$.

\textit{2.} Pour ce qui est de $Z$, on le réécrit de la manière suivante, en utilisant la formule de Tanaka-Meyer, et en écrivant $B^-$ en fonction de $B^+$ et $B$~:
\[ Z_t = a B_t^+ - b B_t^- = a \left( \dfrac{1}{2} L^0_t + \int_0^t \ind{B_s>0} \d B_s\right) - b \left( \dfrac{1}{2} L^0_t - \int_0^t \ind{B_s \leq 0} \d B_s\right). \]

On a donc décomposé $Z_t$ en une martingale locale et le processus $(\frac{a-b}{2}L^0_t, t \geq 0)$, à variations finies (car continu et monotone). On peut donc appliquer la formule donnant les sauts du temps local pour obtenir~: 
\begin{align*}
  L^0_t(Z)-L^{0-}_t(Z) &= (a-b) \int_0^t \d_sL^0_s \ind{Z_s=0} = (a-b)\int_0^t \d_sL^0_s(B) \ind{B_s=0}\\
  &=  (a-b)L^0_t,
\end{align*}
car $L^0_t$ ne croit que sur l'ensemble des instants $t$ tels que $B_t=0$ --par Propriété fondamentale des temps locaux. On en déduit que $L^x_t(Z)$ est discontinu en 0, et le saut réalisé par le temps local en ce point vaut $(a-b)L^0_t$.
\end{proof}

\begin{exercice}[Existence des temps locaux via Fourier]
Soit $Y$ une semi-martingale continue. On définit, pour $\omega \in \Omega$ la mesure aléatoire
\[\mu_t(\omega) : f \mapsto \int_0^t \d\crochet{Y}_s(\omega) f(Y_s(\omega)),\]
qui admet pour transformée de Fourier la fonction~:
\[\xi \mapsto \int_0^t \d \crochet{Y}_s(\omega) e^{i \xi Y_s(\omega)}.\]

\begin{enumerate}
  \item Montrer que si $\E\left[ \int_{\R} \d\xi | \int_0^t \d\crochet{Y}_s e^{i\xi Y_s}|^2 \right] < + \infty$, alors $Y$ admet des temps locaux.
  \item Montrer que cette propriété est vérifiée par le mouvement brownien.
  \item Plus généralement, montrer l'existence de temps locaux pour le mouvement brownien fractionnaire d'ordre $H \in (0;1)$ associés au processus croissant $t$. Ce mouvement brownien fractionnaire $B^{(H)}$ est un processus gaussien centré vérifiant $\E((B^{(H)}_t-B^{(H)}_s)^2)=C|t-s|^{2H}$.
\end{enumerate}
\end{exercice}

\begin{proof}
\textit{1.} Soit $t \in \R_+$ et $Y$ une semi-martingale continue telle que~:
\[\E\left[ \int_{\R} \d\xi | \int_0^t \d\langle Y \rangle_s e^{i\xi Y_s}|^2\right] <+\infty,\]
la transformée de Fourier de cette mesure aléatoire $\mu_t$ est donc une fonction $\hat{\mu}_t$ qui est $L^2$ p.s. Par inversion $L^2$ de la transformée de Fourier, il existe p.s. une fonction (aléatoire) $\phi_t(\omega) \in L^2$ telle que $\hat{\phi}=\hat{\mu}_t$.

On peut donc appliquer la formule de Parseval. Soit $f \in \mathcal{C}_c(\R)$ de transformée de Fourier $\hat{f}$, on a~:
\begin{align*}
  \int_0^t f(Y_s) \d\langle Y \rangle_s &= \int_{\R} \d\xi \hat{f}(\xi) \hat{\mu}_t(\xi) = \int_{\R} \d\xi \hat{f}(\xi) \hat{\phi}_t(\xi)
  = \int_{\R} \d y f(y) \phi_t(y).
\end{align*}
Par conséquent, $\phi_t$ est la densité par rapport à la mesure de Lebesgue de la mesure $\mu_t$. C'est donc un représentant du temps local.

\textit{2.} Passons au mouvement brownien. On a~:
\[\E\left(\int_{\R}\d\xi \left|\int_0^t \d s e^{i\xi B_s}\right|^2 \right)=\int_0^t \d u \int_0^t \d s  \underbrace{\int_{\R}\d\xi \E(e^{i \xi (B_s-B_u)})}_{\sqrt{\dfrac{2\pi}{|s-u|}}}.\]
Comme $\int_0^t \d s \int_0^t \d u|s-u|^{-1/2} < + \infty$, l'existence de temps locaux est démontrée.

\textit{3.} Dans le cas plus général du mouvement brownien fractionnaire de paramètre $H \in (0,1)$, on peut encore montrer l'existence des temps locaux grâce à la transformée de Fourier, car $\displaystyle\int_0^t \d u \int_0^t \d s |s-u|^{-H} <+\infty$.
Le cas $H=\frac{1}{2}$ correspond au cas du mouvement brownien standard, que nous venons de traiter.
\end{proof}

\begin{remarque}
Cet exercice fournit une nouvelle méthode très puissante pour définir des temps locaux assez généraux, cette méthode sera d'ailleurs utilisée au chapitre 7 pour définir les temps locaux d'intersection. Néanmoins, cette méthode a une limite~: elle ne permet pas d'obtenir des résultats tels que la régularité (continuité, etc.) des temps locaux.
\end{remarque}

\begin{exercice}
Soit $p\geq 1$ ; montrer que $(x,t) \mapsto L^x_t(B)$ est continue dans $L^p$.
\end{exercice}

\begin{proof}
Soit $p \geq 1$ et $M,T \in \R_+$. Pour tout $-M \leq x \leq y \leq M$ et $0 \leq s \leq t \leq T$, il existe $C_p>0$ tel que~:
\begin{equation*}
  \E\left( \left| L^x_t(B) - L^y_t(B)\right|^p\right) \leq C_p \underbrace{\E\left(\left|L^x_t(B)-L^y_t(B)\right|^p\right)}_{f_t(x,y)}+C_p\underbrace{\E\left(\left|L^y_t(B)-L^y_s(B)\right|^p\right)}_{g_y(s,t)}
\end{equation*}
Gr\^ace à la formule de Tanaka-Meyer, en utilisant le fait que $x \mapsto x^+$ est Lipschitzienne, il existe $C_p'$ tel que
\[
  f_t(x,y) \leq C_p'\left[ |y-x|^p + \E\left( \left|\int_0^t \ind{x < B_s \leq y} \d B_s\right|^p\right)\right].
\]
De la même façon, il existe $C''_p>0$ tel que
\[ g_y(s,t) \leq C_p''\left[ \E\left( \left| B_t - B_s\right|^p\right) + \E\left(\left|\int_s^t \ind{B_u>y}\d B_u\right|^p\right) \right]. \]

Par conséquent, il existe une constante $K_p>0$ telle que
\begin{align*}
  f_t(x,y) &\leq K_p \left[ |y-x|^p + \E\left(\int_0^t \d u \ind{x<B_u\leq y} \right)^\frac{p}{2} \right]\\
  \mathrm{et} \quad g_t(x,y) &\leq K_p \left[ |t-s|^\frac{p}{2} + \E(|\int_s^t \d u \ind{B_u>y}|^\frac{p}{2})\right]
\end{align*}
En calculant ces intégrales, on en déduit l'existence de $K_{p,M,T}>$ tel que
\[ \E\left[ \left| L^x_t(B) - L^y_t(B)\right|^p\right) \leq K_{p,M,T} \left(|t-s|^\frac{p}{2} + |y-x|^\frac{p}{2}\right). \]
On en conclut que $(L^x_t(B))$ est continue en $(t,x)$ dans $L^p$.
\end{proof}

\begin{exercice}[Propriété fondamentale des temps locaux]
On note $Y$ une semi-martingale positive, et $L^y$ son temps local en $y$. On note $\d_tL^y_t$ la mesure sur $\R_+$ dont $(L^y_t, t \geq 0)$ est la fonction de répartition. 

En décomposant de deux façons différentes $(Y_t-y)^2$ en la somme d'une martingale locale et d'un processus à variations finies, montrer que le support de la mesure $\d_tL^y_t$ est inclus dans $\{t \in \R_+ : Y_t = y\}$.
\end{exercice}

\begin{proof}
On utilise ici la formule de Tanaka-Meyer. On a
\[ (Y_t-y)^+ = (Y_0-y)^+ + \int_0^t \ind{Y_s > y} \d Y_s + \frac{1}{2} L^y_t.\]
De plus pour tout $x \in \R$, $|x|= 2x^+-x$, et $Y_t-y = Y_0-y + \int_0^t 1 \d Y_s$, donc
\[ |Y_t-y| = |Y_0-y| + \int_0^t \sgn(Y_s-y) \d Y_s + L^y_t, \]
où $\sgn$ est la fonction signe définie sur $\R$ par
\[ \sgn(x) =
\begin{cases}
  1 & \mathrm{si} \quad x > 0\\
  -1 & \mathrm{si} \quad x \leq 0.
\end{cases}\]

On décompose $(Y_t-y)^2 = |Y_t-y|^2$, par la formule d'It\^o, on a d'une part
\[(Y_t-y)^2 = (Y_0-y)^2 + 2\int_0^t (Y_s-y)\d Y_s + \crochet{Y}_t\]
et d'autre part
\begin{align*}
  |Y_t-y|^2 &= |Y_0-y|^2 + 2 \int_0^t |Y_s-y|\d|Y_s-y| + \crochet{|Y-y|}_t\\
  &=(Y_0 -y)^2 + 2 \int_0^t |Y_s-y|\sgn(Y_s-y)\d Y_s + 2 \int_0^t |Y_s-y|\d_sL^y_s + \crochet{Y}_t.
\end{align*}
Par conséquent, on obtient $\int_0^t |Y_s-y| \d_sL^y_s = 0$, par conséquent, la mesure $\d_sL^y_s$ ne charge que $\{t \in \R_+ : Y_t=y\}$.
\end{proof}

\begin{exercice}[Généralisation de la formule de densité d'occupation]
\label{exo_generalisationdensite}
Nous allons ici étendre la formule de densité d'occupation à des fonctionnelles aléatoires dépendant du temps. Ce résultat permet en particulier de calculer la mesure de Lévy d'un subordinateur naturellement associé au temps local en 0.
\begin{enumerate}
  \item Montrer que pour toute fonctionnelle $\Phi$ mesurable bornée, on a~:
\[\int_0^t \d s \Phi(s,\omega, B_s(\omega)) = \int_{\R} \d x \int_0^t \d_sL^x_s(\omega) \Phi(s,\omega,x).\]
  \item Soit $l \geq 0$, on pose $\tau_l = \inf\{t \geq 0 : L^0_t > l\}$.
 Montrer que $\tau$ est un subordinateur et déterminer la mesure de Lévy qui lui est associée.
\end{enumerate}
\end{exercice}

\begin{remarque}
Ce résultat est à retenir, car il sera utilisé à de nombreuses reprises par la suite. Il étend la définition fondatrice des temps locaux à une classe bien plus large de fonctions, en les autorisant à dépendre du temps, et de l'aléa.
\end{remarque}

\begin{proof}
\textit{1.} Soit $0 \leq a<b$, et $f,g$ deux fonctions mesurables bornées. On pose $\Phi(s,\omega,x) = \indset{[a,b)}(s)f(\omega) g(x)$, on a~:
\begin{align*}
  \int_0^t \d s \Phi(s,\omega, B_s(\omega)) & = \int_0^t \d s \indset{[a,b)}(s) f(\omega) g(B_s(\omega)) = f(\omega) \int_0^t \d s \indset{[a,b)}(s) g(B_s(\omega))\\
  & = f(\omega) \int_{\R} \d x g(x) (L^x_{b\wedge t}(B)-L^x_{a\wedge t}(B))\\
  & = \int_{\R} \d x \int_0^t \d_sL^x_s \Phi(s, \omega, x).
\end{align*}

On étend ensuite ceci à toute fonctionnelle mesurable positive bornée grâce au théorème de classe monotone, puis on conclut en décomposant $\Phi$ en partie positive et partie négative. On a bien l'égalité demandée.

\begin{remarque}
Notons que $\Phi$ n'est pas nécessairement adaptée à la filtration naturelle de $B$.
\end{remarque}

\textit{2.} On observe pour commencer que, à $l$ fixé, $\tau_l$ est un temps d'arrêt. De plus, en utilisant la propriété de Markov forte au temps $\tau_l$, on observe que $\tau_{l+l'}$ est égal en loi à la somme de deux représentants indépendants de $\tau_l$ et $\tau_{l'}$. Le processus $l \mapsto \tau_l$ est donc un processus à accroissements indépendants stationnaire, positif, et continu à droite. C'est un subordinateur. Calculons alors la transformée de Laplace de ce processus.

Soit $\lambda>0$, on observe que $f_\lambda(l) = \E(e^{- \lambda \tau_l})$ vérifie $f_\lambda(l+l') = f_\lambda(l)f_\lambda(l')$. Pour tout $\lambda > 0$, il existe $\phi(\lambda) \in \R$ tel que pour tout $l \geq 0$, on a $\E(e^{-\lambda \tau_l} ) = \exp(-l \phi(\lambda))$.
Par changement de variables, on a $\int_0^{+\infty} \d_sL^0_s e^{-\lambda s} =  \int_0^{+\infty} \d l e^{-\lambda \tau_l}$. On calcule de deux façons différentes l'espérance de cette variable aléatoire.

D'une part part, en utilisant la forme particulière de la mesure de Lévy de $\tau_l$,
\begin{align*}
  \E \left(  \int_0^{+\infty} \d l e^{-\lambda \tau_l} \right)
  =  \int_0^{+\infty} \d l \E\left( e^{-\lambda \tau_l} \right)
  =  \int_0^{+\infty} \d l e^{-l \phi(\lambda)}
  =  \frac{1}{\phi(\lambda)}.
\end{align*}

D'autre part, la mesure $\d_sL^x_s$ ne chargeant que $\{t \geq 0 : B_t = x\}$, par conséquent, pour toute fonction $f$ positive mesurable~:
\begin{align*}
  \E\left(\int \d x f(x) \int_0^{+\infty} \d_sL^x_s e^{-\lambda s}\right) &= \E\left( \int_0^{+\infty} \d t f(B_t) e^{-\lambda t} \right) = \int_0^{+\infty} \d t e^{-\lambda t} \E(f(B_t))\\
  &= \int_\R \d x f(x) \int_0^{+\infty} \d t \frac{e^{-\frac{x^2}{2t}}}{\sqrt{2\pi t}} e^{-\lambda t}
\end{align*}
En particulier, pour $x=0$, on a
\begin{align*}
  \E\left( \int_0^{+\infty} \d_sL^0_s e^{-\lambda s} \right)
  =  \int_0^{+\infty} \dfrac{\d s}{\sqrt{2\pi s}} e^{-\lambda s}
  =  \dfrac{\Gamma\left(\frac{1}{2}\right)}{\sqrt{2\pi \lambda}}
  =  \frac{1}{\sqrt{2\lambda}}.
\end{align*}

On en déduit que $\phi(\lambda) = \sqrt{2\lambda}$. Or, comme $\displaystyle \sqrt{2\lambda} = \int_{\R^+} \d x \frac{(1-e^{-\lambda x})}{\sqrt{2\pi x}}$, la mesure de Lévy associée à $(\tau_l,l \geq 0)$ est $\dfrac{\d x}{\sqrt{2\pi x}}$.
\end{proof}

\begin{exercice}
Montrer que la mesure aléatoire $\displaystyle \mu_t : f \mapsto \int_0^t \d s \int_0^s \d u f(B_u-B_s)$ est absolument continue par rapport à la mesure de Lebesgue, et identifier la dérivée de Radon-Nikod\'ym, i.e. le processus $\alpha^a_t$ vérifiant, pour toute fonction $f$ mesurable positive
\[\int_0^t \d s \int_0^s \d u f(B_u-B_s) = \int_\R \d a f(a) \alpha^a_t,\]
en fonction des temps locaux de $B$.
\end{exercice}

\begin{proof}
Soit $f$ fonction mesurable positive, on utilise une première fois la formule de densité d'occupation du mouvement brownien
\begin{align*}
  \int_0^t \d s \int_0^s \d u f(B_u-B_s) & = \int_0^t  \d s \int_\R \d a f(a-B_s)L^a_s = \int_\R \d a \int_0^t \d s f(a-B_s)L^a_s,
\end{align*}
par formule de Fubini. On utilise ensuite la formule de densité d'occupation généralisée de l'Exercice \ref{exo_generalisationdensite}, qui nous donne
\[\int_0^t \d s f(a-B_s)L^a_s = \int_\R \d b \int_0^t f(a-b)L^a_s \d_sL^b_s.\]

On obtient finalement~:
\begin{align*}
  \mu_t(f) &= \int_\R da\int_\R \d b f(a-b) \int_0^t \d_sL^b_sL^a_s = \int_\R \d a f(a) \int_\R \d b \int_0^t L^{a+b}_s \d_sL^b_s.
\end{align*}
Par identification, on  en déduit $\displaystyle \alpha^a_t = \int_\R \d b \int_0^t L^{a+b}_s \d_sL^b_s$.
\end{proof}

\begin{exercice}
\label{exo_momenttempslocal}
Montrer que pour tout $m \in \N$, on a 
\[\E(L^0_t(B)^{2m}) = \E((B_t)^{2m}) = t^m \dfrac{(2m)!}{2^m m!}.\]

En déduire que à $t$ fixé, $L^0_t(B) \egaldistr |B_t|$.
\end{exercice}

\begin{remarque}
Ce résultat est une conséquence directe du théorème d'équivalence de Lévy, qui sera énoncé par la suite, mais nous allons ici donner une solution utilisant uniquement un simple calcul d'intégrales.
\end{remarque}

\begin{proof}
Commençons par calculer les moments de $L^0_t$. Pour cela, notons que si $f$ est une fonction de $\R^d$ dans $\R$ mesurable bornée, et $B$ un mouvement brownien, on a par définition des temps locaux~:
\[\int_0^{t_1} \d s_1 f(B_{s_1}, B_{s_2}, \ldots, B_{s_d}) =  \int_{x_1 \in \R} \d x_1 f(x_1, B_{s_2}, \ldots, B_{s_d}) L^{x_1}_{t_1}\text{  p.s.}\]
Dès lors, en utilisant le théorème de Fubini on obtient~:
\[\int_0^{t_1} \d s_1 \cdots \int_0^{t_d} \d s_d f(B_{s_1}, \ldots, B_{s_d}) = \int_{\R^d}\d x_1 \cdots \d x_d f(x_1,\ldots, x_d) L^{x_1}_{t_1}\cdots L^{x_d}_{t_d} \text{  p.s.}\]

Or
\[\E(f(B_{s_1}, \ldots, B_{s_d})) = \int_{\R^d} \d x_1\cdots \d x_d f(x_1, \ldots, x_d) D_{s_1,\ldots, s_d}(x_1,\ldots, x_d),\]
où $D_{s_1,\ldots, s_d}(x_1,\ldots, x_d)$ est la densité de la loi de $(B_{s_1},\ldots B_{s_d})$ par rapport à la mesure de Lebesgue. Dès lors, en utilisant la formule de Fubini, et par identification, on obtient~:
\[\E(L^{x_1}_{t_1}\cdots L^{x_d}_{t_d}) = \int_0^{t_1} \d s_1 \cdots \int_0^{t_d}\d s_d D_{s_1,\ldots, s_d}(x_1,\ldots, x_d).\]

On obtient en particulier, pour tout $n \in \N$,
\[\E((L^0_t)^n) = \int_{[0,t]^n} \d s_1\cdots \d s_n D_{s_1,\ldots, s_n}(0, \ldots, 0).\]
On réalise un réarrangement de $s_1,\ldots s_n$ dans l'ordre croissant, on a $n!$ combinaisons possibles. On obtient donc~:
\[\E((L^0_t)^n) = \dfrac{n!}{(2\pi)^\frac{n}{2}} \int_{0 \leq s_1 \leq s_2 \cdots \leq s_n\leq t} \dfrac{\d s_1\cdots \d s_n}{\sqrt{s_1(s_2-s_1)\cdots (s_n-s_{n-1})}}.\]
Par propriété de scaling, on se ramène à $t=1$~:
\begin{equation}
\label{eqn_tempslocalen0}
\E((L^0_t)^n) = \dfrac{n! t^\frac{n}{2}}{(2\pi)^\frac{n}{2}} \int_{0 \leq s_1 \leq s_2 \cdots \leq s_n \leq 1} \dfrac{\d s_1\cdots \d s_n}{\sqrt{s_1(s_2-s_1)\cdots (s_n-s_{n-1})}}.
\end{equation}

On réalise maintenant le changement de variables suivant, pour $i \geq 2$ on pose $u_i=\frac{s_i-s_{i-1}}{1-s_{i-1}}$, on obtient~:
\begin{align*}
  \int_{0 \leq s_1 \leq s_2 \cdots \leq s_n \leq 1} \dfrac{\d s_1\cdots \d s_n}{\sqrt{s_1(s_2-s_1)\cdots (s_n-s_{n-1})}}
  &  = \int_{[0,1]^n} \dfrac{\d u_1 \cdots \d u_n \prod_{i=1}^n \sqrt{1-u_i}^{n-i}}{\sqrt{u_1\cdots u_n}}\\
  & = \prod_{i=0}^{n-1} \int_0^1 \d u \sqrt{\dfrac{(1-u)^i}{u}}\\
  &  = \prod_{i=1}^{n} \int_0^1 \dfrac{\d v v^\frac{i}{2}}{\sqrt{v(1-v)}}
\end{align*}
par changement de variables $v=\sqrt{1-u}$.

Soit $\beta$ une variable aléatoire de loi $\beta_{\frac{1}{2},\frac{1}{2}}$, on a~:
\[ A_i = \int_0^1 \d u \dfrac{u^\frac{i}{2}}{\sqrt{u(1-u)}} = \pi \E(\beta^\frac{i}{2}) = \sqrt{\pi}\dfrac{\Gamma(\tfrac{i+1}{2})}{\Gamma(\tfrac{i+2}{2})}.\]
Par conséquent
\[  \int_{0 \leq s_1 \leq s_2 \cdots \leq s_n \leq 1} \dfrac{\d s_1\cdots \d s_n}{\sqrt{s_1(s_2-s_1)\cdots (s_n-s_{n-1})}} = \dfrac{\pi^\frac{n}{2}}{\Gamma(\tfrac{n}{2}+1)}.\]

D'après \eqref{eqn_tempslocalen0}, on en déduit~:
\[\E((L^0_t)^n) = n! \left(\dfrac{t}{2}\right)^\frac{n}{2}\dfrac{1}{\Gamma(\tfrac{n}{2}+1)},\]
ce qui se traduit, lorsque $n=2p$ par~:
\[\E((L^0_t)^{2p}) = t^p \dfrac{(2p)!}{2^p p!}\]
qui sont les moments pairs de $|B_t|$. Or $\E(\exp(\lambda B_t^2))<+\infty$ pour $\lambda t<\frac{1}{2}$.

Étant donné que les moments de $(L^0_t)^2$ et ceux de $B_t^2$ sont égaux, on en déduit que les transformées de Laplace sont identiques, donc $(L^0_t)^2 \egaldistr B_t^2$. De plus $L^0_t \geq 0$, par conséquent $L^0_t = \sqrt{(L^0_t)^2} \egaldistr |B_t|$.
\end{proof}

\begin{exercice}
\label{exo_integrationtempslocal}
Soit $L$ le temps local en 0 d'un $(\F_t)$-mouvement brownien $B$. On note $\W$ la mesure de Wiener et $\W^u_{0,0}$ la loi du pont brownien de longueur $u$. On admettra le résultat suivant, démontré dans l'Exercice \ref{exo_radonnikodym}~:
\[\W^u_{0,0}{}_{|\F_s} = \sqrt{\frac{u}{u-s}} \exp\left(-\frac{B_s^2}{2(u-s)}\right) \cdot \W_{|\F_s}\]
donnant la dérivée de Radon-Nykodym de la loi du pont brownien, par rapport à la mesure de Wiener.

Montrer que pour tout $H$ processus prévisible positif, on a~:
\[\E\left[ \int_0^{+\infty} \d_uL_u H_u\right] = \int_0^{+\infty} \dfrac{\d u}{\sqrt{2\pi u}} \E(H_u|B_u=0).\]
\end{exercice}

\begin{proof}
Soit $s \leq t$ et $A$ une variable aléatoire $\F_s$-mesurable positive, on pose $H_u = A \ind{u \in (s,t]}$ et on souhaite calculer~:
\[\E\left[ \int_0^{+\infty} \d_uL_u H_u \right] = \E((L_t-L_s)A).\]

On calcule alors l'espérance conditionnelle de $L_t-L_s$ sachant $\F_s$, grâce à la propriété de Markov et à l'Exercice \ref{exo_momenttempslocal}, on a pour tout $x \in \R$,
\begin{align*}
  \E(L_t-L_s|\F_s,B_s=x)
  & = \E_x(L_{t-s}) = \E(L^{-x}_{t-s}) = \int_0^{t-s} \dfrac{\d u}{\sqrt{2\pi u}} e^{-\frac{x^2}{2u}}.
\end{align*}
Nous pouvons donc écrire~:
\[\E(L_t-L_s|\F_s) \int_s^t \dfrac{\d u}{\sqrt{2\pi (u-s)}} e^{-\frac{B_s^2}{2(u-s)}}.\]

Observons que $\E(H_u|B_u=0) = \W^u_{0,0}(H_u)$, on utilise la densité de la loi du pont brownien par rapport à la mesure de Wiener pour écrire~:
\begin{align*}
  \E\left[ \int_0^{+\infty} \d_uL_u H_u\right]
  & = \E[(L_t-L_s)A] = \E[A \E(L_t-L_s|\F_s)]\\
  & = \E\left[\int_s^t \dfrac{\d u}{\sqrt{2\pi (u-s)}} A e^{-\frac{B_s^2}{2(u-s)}} \right] = \int_s^t \dfrac{\d u}{\sqrt{2\pi u}}\W^u_{0,0}(A)\\
  & = \int_0^{+\infty} \dfrac{\d u}{\sqrt{2\pi u}} \E(H_u|B_u=0).
\end{align*}
On utilise ensuite la linéarité de ces expressions et le théorème des classes monotones pour conclure.
\end{proof}

\begin{remarque}
Le résultat de cet exercice est utile. Il permet notamment de représenter de nombreuses espérances d'intégrales par rapport aux temps locaux grâce à la loi des ponts browniens.
\end{remarque}

\chapter[Formules d'Itô et de Tanaka-Meyer]{Unification des formules d'Itô et de Tanaka-Meyer}

La formule de Tanaka-Meyer, que nous avons développée au chapitre précédent donne une décomposition de la semi-martingale $((Y_t-y)^+, t \geq 0)$, en fonction d'une intégrale stochastique par rapport à $Y$ et du temps local de $Y$ au niveau $y$. Cette formule n'est pas sans rappeler la formule d'Itô, qu'elle étend en un certain sens, puisque la fonction $x \mapsto x^+$ n'est pas de classe $\mathcal{C}^2$. Grâce à cette observation, et à la définition des temps locaux donnée au chapitre précédent, nous allons pouvoir étendre le calcul stochastique d'Itô à des fonctions pour lesquelles la dérivée seconde (au sens des distributions) est une mesure de Radon. En particulier, on en déduira que la classe des semi-martingales continues est invariante par ces fonctions, plus générales que celle des fonctions de classe $\mathcal{C}^2$. Pour finir, nous utiliserons une application des temps locaux qui permet de toucher du doigt une partie de la richesse du mouvement brownien, en particulier le fait que toute loi de probabilité de variance finie peut être modélisée à l'aide d'un mouvement brownien pris en un temps d'arrêt d'espérance finie, ce qui répond au problème de plongement brownien de Skorokhod.

\section{Théorèmes principaux}

\begin{theoreme}[Formule d'Itô généralisée]
Soit $\mu$ une mesure de Radon et $Y$ une semi-martingale continue de temps local en $y$ noté $L^y_t$, on a
\[\int \left((Y_t-y)^+ - (Y_0-y)^+\right) \mu(\d y) = \int_0^t \mu(]-\infty, Y_s[)\d Y_s + \dfrac{1}{2}\int L^y_t \mu(\d y).\]
Soit $F$ fonction continue localement différence de deux fonctions convexes, on a
\[F(Y_t)=F(Y_0)+\int_0^ t F'_g(Y_s) \d Y_s + \dfrac{1}{2} \int L^y_t F''(\d y),\]
où la dérivée seconde de $F$ est prise au sens des distributions : c'est, par hypothèse une mesure de Radon $\mu$.
\end{theoreme}

On va maintenant spécifier quelques résultats relatifs au mouvement brownien, permettant d'interpréter le temps local en 0 comme un objet naturel.

\begin{theoreme}[Théorème d'équivalence de Lévy]
Soit $B$ un mouvement brownien, $L$ son temps local en 0 et $S_t=\sup_{s\leq t}B_s$, on a :
\[(S_t-B_t,S_t)_{t \geq 0} \egaldistr (|B_t|,L_t)_{t\geq 0}.\]
\end{theoreme}

\begin{lemme}[Skorokhod]
\label{lem_skorokhod}
Si $x : \R^+ \to \R$ est une fonction continue avec $x(0)=0$, alors il existe un unique couple $(z,l)$ de fonctions continues issues de 0 vérifiant :
\begin{itemize}
  \item $z(t) = -x(t)+l(t)$
  \item $z(t)\geq 0$
  \item $l(t)$ est croissante et $\d l$ est portée par $\{t \geq 0:z(t)=0\}$.
\end{itemize}

Cette solution est donnée par :
\[l(t) = \sup_{s \leq t} x(s) \text{  et  } z(t) = l(t)-x(t).\]
\end{lemme}

\begin{theoreme}[Extension du théorème de Lévy]
Soit $\tau^{(\mu)}_l = \inf \{t \geq 0 : L^{(\mu)}_t \geq l \}$ et $T^{(\mu)}_l = \inf\{t \geq 0: B^{(\mu)}_t \geq l\}$ ; on a :
\[\forall l \in \R, \left(\left. |B^{(\mu)}_t|, t \leq \tau^{(\mu)}_l \right| \tau^{(\mu)}_l< + \infty \right) \egaldistr \left( S^{(\mu)}_t-B^{(\mu)}_t, t \leq T^{(\mu)}_l \right)\]
\end{theoreme}

Un dernier résultat permet de tester la richesse du mouvement brownien : on peut représenter toute variable aléatoire centrée de carré intégrable grâce à un mouvement brownien pris en un temps d'arrêt d'espérance finie. Ce problème, dit du plongement de Skorokhod a été résolu par de très nombreux moyens différents (l'article de Obloj \cite{Obl2004} en dénombrait au moins 21). Nous allons en exprimer un grâce à la fonction de Hardy-Littlewood.

\begin{proposition}
Soit $\mu$ une mesure de probabilité, centrée possédant un deuxième moment fini et $X$ une variable aléatoire de loi $\mu$, on pose :
\[\Psi_\mu (x) = \dfrac{1}{\mu([x,+\infty))}\int_x^{+\infty} y\mu(\d y) = \E(X|X\geq x)\]
la fonction de Hardy-Littlewood associée à $\mu$.

Le temps d'arrêt $T_\mu = \inf\{t>0:S_t \geq \Psi_\mu(B_t)\}$ est un temps d'arrêt intégrable vérifiant $B_{T_\mu} \sim \mu$.
\end{proposition}

\section{Exercices}

\begin{exercice}[Une réciproque à la Proposition \ref{pro_fondamentale}]
Soit $B$ un mouvement brownien et $L^y$ son temps local en $y$. Montrer que le support de $\d_tL^y_t$ est égal à $\{t \geq 0 : Y_t=y\}$ p.s.
\end{exercice}

\begin{proof}
On s'intéresse tout d'abord au cas $y=0$. Pour tout $q \in \Q$, on pose :
\[\sigma_q = \inf\{t \geq q : B_t = 0\},\]
par propriété de Markov forte, $B^{\sigma_q} = (B_{\sigma_q+t},t \geq 0)$ est un mouvement brownien indépendant de $\F_{\sigma_q}$ issu de 0.

Observons que pour tout $t\geq 0$, $L^0_t>0$ p.s., car $\lim_{t \to 0} \{L^0_t>0\}$ est un événement de $\F_{0+}$, et de plus $\P(L^0_t>0) = \P(L^0_1>0)>0$ par propriété de scaling du mouvement brownien. Par conséquent,
\[L^0_{\sigma_q+r} = L^0_{\sigma_q} + L^0_r(B^{\sigma_q}),\]
$\sigma_q$ est donc un instant de croissance de $L$.

Dès lors, avec probabilité 1, $\{\sigma_q, q \in \Q\}$ est dense dans $\{t\geq 0 : B_t=0\}$. On peut conclure sur l'égalité de cet ensemble avec le support de $\d_t L^0_t$.

Soit $T_y = \inf\{t \geq 0 : B_t=y\}$. On observe que
\[L^y_t = L^0_{(t-T_y)^+}(B^{T_y}-y),\]
ce qui permet de se ramener au cas $y=0$.
\end{proof}

Dans le prochain exercice, on discute du choix de $(\crochet{X}_t, t \geq 0)$ comme processus croissant pour définir les temps locaux de la semi-martingale $X$. Plus exactement, nous allons identifier les processus croissants associés à $X$ pour lesquels il existe une identité \og naturelle \fg{} au sens des temps locaux.

\begin{exercice}
Soit $X$ une semi-martingale.
\begin{enumerate}
  \item On note $(L^a_t(X), a \in \R, t \geq 0)$ la famille des temps locaux de $X$ et $h : \R \to \R$ une fonction de classe $\mathcal{C}^2$, strictement croissante et bijective. Exprimer le temps local $L^{h(a)}_t(h(X))$ en fonction de $L^a_t(X)$.
  \item Soit $(A_s(X),s \geq 0)$ une fonctionnelle associant un processus croissant à la semi-martingale $X$. On pose $(\Lambda^a_t(X), a \in \R, t \geq 0)$ les temps locaux définis par la formule suivante : pour tout $f \geq 0$ mesurable,
\[ \int_0^t \d A_s(X) f(X_s) = \int_\R \d x f(x) \Lambda^x_t(X). \]
Déterminer pour quelles fonctionnelles $A$ la formule précédente est satisfaite.
  \item Soit $\mu$ une mesure $\sigma$-finie positive, on note $A^\mu$ le processus croissant défini par
\[A^{\mu}_s(X) = \int_\R \d \mu(x) L^x_s(X).\]
Pour quelles mesures $\mu$ l'identité précédente a-t-elle lieu ?
\end{enumerate}
\end{exercice}

\begin{proof}
\textit{1.} Soit $h$ une fonction de classe $\mathcal{C}^2$ bijective sur $\R$. Par formule d'Itô, on a :
\[h(X_t) = h(X_0) + \int_0^t h'(X_s)\d X_s + \frac{1}{2}\int_0^t h''(X_s) \d\crochet{X}_s.\]

Par conséquent, toute fonction mesurable positive $f$, le temps local de $h(X)$ vérifie
\begin{align*}
  \int_\R \d x f(x) L^x_t(h(X)) & = \int_0^t f(h(X_s)) \d\crochet{h(X)}_s\\
  & = \int_0^t f(h(X_s)) h'(X_s)^2 \d\crochet{X}_s\\
  & = \int_\R \d x f(h(x)) h'(x)^2 L^x_t(X),
\end{align*}
en utilisant la définition des temps locaux de $X$.

Appliquons un changement de variables $x=h(y)$, on obtient :
\[\int_\R \d y f(h(y)) h'(y) L^{h(y)}_t(h(X)) = \int_\R \d x f(h(x)) h'(x)^2 L^x_t(X),\]
d'où l'on tire
\begin{equation}
 \label{eqn_tempslocal}
 L^{h(a)}_t(h(X)) = h'(a)L^a_t(X).
\end{equation}

\begin{remarque}
Cette relation semble raisonnable étant donné la notion de temps local, dans un certain sens, le temps local de $h(X)$ en $h(a)$ est donné par le temps local de $X$ en $a$, multiplié par la \og vitesse \fg{} à laquelle $h(X)$ quitte $h(a)$.
\end{remarque}

\textit{2.} Nous pouvons maintenant réaliser le même raisonnement pour ces nouveaux temps locaux $(\Lambda^x_t)$, on a par changement de variables :
\begin{align*}
  \int_0^t f(h(X_s)) \d A_s(h(X)) & = \int_\R \d x f(x) \Lambda^x_t(h(X))\\
  & = \int_\R \d x h'(x) f(h(x)) \Lambda^{h(x)}_t(h(X)).
\end{align*}

Par conséquent, si $\Lambda$ satisfait l'équation \eqref{eqn_tempslocal}, l'égalité suivante est vérifiée :
\[\int_0^t f(h(X_s)) \d A_s(h(X)) = \int_\R \d x h'(x)^2 f(h(x)) \Lambda^x_t(X) = \int_0^t \d A_s(X) h'(X_s)^2 f(h(X_s)).\]
La relation suivante doit donc être vérifiée par la fonctionnelle $A$ :
\begin{equation}
\label{eqn_fonctionnelle}
A_t(h(X)) = \int_0^t h'(X_s)^2 \d A_s(X).
\end{equation}

\textit{3.} Intéressons-nous maintenant aux fonctionnelles croissantes $A$ définies à partir des temps locaux de semi-martingales. L'équation \eqref{eqn_fonctionnelle} devient :
\begin{align*}
  A^\mu_t(h(X)) = & \int_0^t h'(X_s)^2 \d A^\mu_s(X)\\
  \int_\R \d\mu(x) L^x_t(h(X)) = & \int_\R \d \mu(x) \int_0^t h'(X_s)^2 \d_sL^x_s
\end{align*}
or $\int_0^t h'(X_s)^2 \d_sL^x_s = h'(x)^2 L^x_t$, car le support de la mesure $\d_sL^x_s$ est inclus dans $\{s \geq 0 : X_s=x\}$. On en déduit :
\begin{align*}
  \int_\R \d \mu(x) L^x_t(h(X))
  = & \int_\R \d \mu(x) h'(x)^2L^x_t(X) \\
  = & \int_\R \d \mu(x) h'(x) L^{h(x)}_t(h(X)).
\end{align*}

En particulier, pour $h(x)=x+a$, on a :
\[\int_\R \d \mu(x) L^x_t(X) = \int_\R \d \mu(x) L^{x+a}_t(X),\]
par conséquent, $\mu$ est une mesure $\sigma$-finie, invariante par translations sur $\R$ , elle donc est proportionnelle à la mesure de Lebesgue.

En définitive, le seul temps local \og raisonnable \fg{} de semi-martingale est proportionnel à $L^x_t$.
\end{proof}

\begin{exercice}
\label{exo_martingaleKennedy}
Soit $B$ un mouvement brownien standard issu de 0 et $L$ son temps local en 0.

\begin{enumerate}
  \item Soit $\phi(l,x,t)$ une fonction régulière de $\R^3$ dans $\R$. Déterminer la relation vérifiée par $\phi$ quand $\phi(L_t,|B_t|,t)$ est une martingale locale.
  \item En déduire la transformée de Laplace de $\tau_l=\inf\{t>0:L_t\geq l\}$.
\end{enumerate}
\end{exercice}

\begin{proof}
\textit{1.} On applique la formule d'Itô au processus $\phi(L_t,|B_t|,t)$, on a :
\begin{align*}
 \phi(L_t,|B_t|,t) = \phi(0,0,0) &+ \int_0^t \partial_l \phi (L_s,|B_s|,s) \d L_s + \int_0^t \partial_x \phi(L_s,|B_s|,s) \d|B_s|\\
  &+ \int_0^t \partial_t \phi(L_s,|B_s|,s)\d s+\dfrac{1}{2}\int_0^t \partial_{x,x} \phi(L_s,|B_s|,s) \d s.
\end{align*}

De plus :
\[\int_0^t\partial_x \phi(L_s,|B_s|,s) \d|B_s| = \int_0^t \partial_x \phi(L_s,|B_s|,s) \sgn(B_s) \d B_s+ \int_0^t \partial_x \phi(L_s,|B_s|,s) \d L_s,\]
par conséquent $\phi(L_t,|B_t|,t)$ est une martingale locale si et seulement si les parties à variations bornées sont nulles. En d'autres termes, si $\phi$ vérifie pour tout $t\geq 0$ :
\[\left\{ 
\begin{array}{l}
 \int_0^t \partial_l\phi(L_s,|B_s|,s) + \partial_x \phi(L_s,|B_s|,s) \d L_s = 0\\
 \int_0^t \frac{1}{2}\partial_{x,x} \phi(L_s,|B_s|,s) + \partial_t \phi(L_s,|B_s|,s) \d s =0.
\end{array}
\right.\]
Cela revient donc à demander que pour tout $l,x,t \geq 0$ (rappelons que $B$ est nul sur le support de $\d_sL_s$) $\phi$ vérifie les équations aux dérivées partielles suivantes :
\[ \left\{
\begin{array}{l}
 \partial_l \phi (l,0,t) + \partial_x \phi(l,0,t) = 0\\
 \frac{1}{2}\partial_{x,x} \phi(l,x,t) + \partial_t \phi(l,x,t) = 0.
\end{array}
\right.\]

\textit{2.} Calculons maintenant la transformée de Laplace de $\tau_l=\inf\{t>0:L_t=l\}$. On cherche alors une fonction $\phi$ sous la forme suivante $\phi(l,x,t) = f(l)g(x) e^{-\lambda t}$, avec $\lambda>0$ donné. La solution suivante est obtenue par intégration successives~:
\[\phi(l,x,t) = \exp(\sqrt{2\lambda} (l-x))e^{-\lambda t}.\]

La fonction $\phi$ est bornée sur $[0,l]\times \R^+\times \R^+$, et $L_{t \wedge \tau_l} \leq l$, on peut donc appliquer le théorème d'arrêt à la martingale $\phi(L_{t \wedge \tau_l},|B_{t \wedge \tau_l}|,{t \wedge \tau_l})$, on a alors :
\[1 = \E(\exp (\sqrt{2\lambda} (L_{\tau_l}-|B_{\tau_l}|))e^{-\lambda \tau_l}).\]

On utilise maintenant le fait que $L_{\tau_l}=l$, de plus comme le temps d'arrêt $\tau_l$ est un instant de croissance de $L$ on a $B_{\tau_l}=0$. On obtient donc :
\[\E(e^{-\lambda \tau_l} ) = \exp(-l\sqrt{2\lambda}).\]

Le processus $(\tau_l,l \geq 0)$ est un processus à accroissements indépendants stationnaires, par propriété de Markov appliquée à $B$, donc un subordinateur stable de paramètre $\frac{1}{2}$.
\end{proof}

\begin{exercice}[Calcul stochastique du premier ordre]
\label{exo_mart}
Soit $f : \R^+ \to \R^+$ une fonction $L^1_\text{loc}$ et $F(x) = \int_0^x f(s)\d s$. On pose $B$ un mouvement brownien et $L$ son temps local en 0. Montrer que 
\[(f(L_t)|B_t|-F(L_t), t \geq 0)\]
est une martingale locale.
\end{exercice}

\begin{proof}
Rappelons pour commencer que la formule de Tanaka donne :
\[|B_t| = \int_0^t \sgn(B_s) \d B_s + L_t.\]

Observons que $L$ est un processus continu et croissant, donc à variations finies, et que $\beta_t = \int_0^t \sgn(B_s) \d B_s$ est un mouvement brownien. Soit $f$ une fonction de classe $\mathcal{C}^2$, en appliquant la formule d'Itô, on obtient
\[f(L_t)(\beta_t+L_t) = \int_0^t (f'(L_s)(\beta_s+L_s)+f(L_s))\d L_s + \int_0^t f(L_s) \d \beta_s.\]
Nous utilisons maintenant le fait que $\d L_s$ est portée par $\{s\geq 0: B_s=0\}$,
\[\int_0^t f'(L_s)(\beta_s+L_s)\d L_s=\int_0^t f'(L_s)|B_s|\d L_s = 0.\]

Par changement de variables $u=L_s$, on a
\[F(L_t) = \int_0^{L_t} f(u) \d u = \int_0^t f(L_s)\d L_s.\]
On note $\tau_l = \inf\{ t \geq 0 : L_t \geq \lambda\}$, observons que
\[ \left| \int_0^{\tau_l} f(L_s)^2 \d s \right| \leq \tau_l \sup_{s \in [0,l]}  f(s)^2.  \]
Par conséquent, $f(L_t)|B_t|-F(L_t) = \int_0^t f(L_s) \d \beta_s$ est une martingale locale.

\begin{remarque}
Cette formule peut être généralisée par la formule de balayage, qui sera développée dans un prochain chapitre (c.f. Exercice \ref{exo:balayage}), ou en utilisant la théorie des excursions.
\end{remarque}

Passons maintenant au résultat général. Soit $f \in L^1_\text{loc}$, on approche $f$ par une suite de fonctions $f_n \in \mathcal{C}^2$ majorées par $2f+1$. On a alors la suite de martingales locales $M^n_t = f_n(L_t)|B_t|-F_n(L_t)$ qui converge presque sûrement vers $M_t=f(L_t)|B_t|-F(L_t)$. On note maintenant
\[T_p = \inf\{t>0:L_t>p\},\]
comme $L_t<+\infty$ p.s., cette suite de temps d'arrêts croit bien vers $+\infty$. Grâce au théorème de Lévy, on a :
\[E(f(L_{t \wedge T_p})) = \E(f(S_{t\wedge \tilde{T}_p})) = f(p) \P(\tilde{T}_p<t) + \dfrac{2}{\sqrt{2\pi t}}\int_0^p f(x) \exp\left(- \frac{x^2}{2t}\right)\d x <+\infty\]
pour des valeurs de $p$ bien choisies en fonction de $f$.

Comme $F$ est continue, on obtient de la même manière $\E(F(L_{t \wedge T_p}))<+\infty$. Par conséquent, la suite $(M^n_{t \wedge T_{p}})_{n \in \N}$ est une suite de martingales uniformément intégrables qui converge p.s. vers $M_{t \wedge T_p}$ pour des valeurs de $p$ bien choisies. De plus, à $t$ fixé, ces variables aléatoires sont dominées par une variable aléatoire dans $L^1$, la convergence a donc également lieu dans $L^1$, ce qui nous permet d'obtenir la propriété de martingale pour $(M_{t \wedge T_p})_{t \geq 0}$.

$M$ est donc bien une martingale locale, réduite par $(T_p)_{p \geq 0}$.
\end{proof}

\begin{remarque}
On peut remarquer que pour toute fonction $f \in L^1_\text{loc}$, $\int_0^t f(L_s) \d \beta_s$ est une martingale locale. Montrons que pour tout $l \geq 0, \int_0^{\tau_l} f^2(L_s)\d s <+\infty$ p.s.
\begin{align*}
 \int_0^{\tau_l} f^2 (L_s)\d s & = \sum_{\lambda \leq l} \int_{\tau_{\lambda-}}^{\tau_\lambda} f^2(L_s)\d s  = \sum_{\lambda \leq l} f^2(\lambda)(\tau_{\lambda}-\tau_{\lambda-}).
\end{align*}

On calcule la transformée de Laplace de cette variable aléatoire, en utilisant le fait que $(\tau_l)_{l\geq 0}$ est un subordinateur d'exposant caractéristique associé $\mu \mapsto \sqrt{2\mu}$ et de mesure de Lévy associée $\frac{\d x}{\sqrt{2\pi x^3}}$ :
\begin{align*}
 \E\left(\exp\left( - \mu\int_0^{\tau_l} f^2 (L_s)\d s  \right)\right)
 & = \E\left(\exp\left( - \mu \sum_{\lambda \leq l} f^2(\lambda)(\tau_{\lambda}-\tau_{\lambda-})  \right)\right)\\
 & = \exp \left( - \int_0^l \d \lambda \int_\R \dfrac{\d x}{\sqrt{2\pi x^3}} \left( 1-e^{-\mu f^2(\lambda) x} \right) \right)\\
 & = \exp\left( - \sqrt{2 \mu} \int_0^l |f(\lambda)| \d \lambda \right).
\end{align*}

Par conséquent, comme $f \in L^1_\text{loc}$, on en déduit que $\int_0^{\tau_l} f^2 (L_s)\d s$ est fini presque sûrement, garantissant l'intégrabilité de $f(L_s)$ par rapport à $\d \beta_s$.
\end{remarque}

\begin{exercice}
\label{exo_supmartingale}
Soit $(N_t)$ une martingale locale continue positive issue de $a>0$. On suppose que $\lim_{t \to +\infty} N_t =0$.

Montrer que $\sup_{t \geq 0} N_t \egaldistr \dfrac{a}{U}$ où $U$ est une variable aléatoire uniforme sur $[0,1]$. 
\end{exercice}

\begin{proof}
On note $N^*=\sup_{t \geq 0} N_t$. On pose également pour $x\geq 0$ :
\[T_x = \inf\{t\geq 0 : N_t \geq x\} \text{ le premier temps d'atteinte de }x.\]

Notons que comme $\lim_{t \to +\infty} N_t =0$, pour tout $x \leq a$, on a $\P(T_x <+\infty) =1$. Soit $x \geq a$, comme~$N_{t \wedge T_x}$ est une martingale bornée, elle est uniformément intégrable, et on a :
\[a = \E(N_0) = \E(N_{T_x}) = x \P(T_x<+\infty) +0 \P(T_x=+\infty).\]
Par conséquent, pour tout $x>0$, on a $\P(T_x<+\infty)=\left(\!\frac{a}{x}\!\right)\wedge 1$. On obtient ainsi :
\[\P(N^*\geq x) = \P(T_x<+\infty) = \dfrac{a}{x}\wedge 1,\]
et on observe que $\P(\frac{a}{U}\geq x) = \P(U\leq \frac{a}{x}) = \frac{a}{x}\wedge 1.$ Par conséquent,~$N^* \egaldistr \frac{a}{U}$.
\end{proof}

\begin{remarque}
Le caractère universel de ce résultat peut être retrouvé en utilisant la représentation de Dubins-Schwarz d'une martingale, $N_t = \beta_{\langle N \rangle_t}$, où $\beta$ est un mouvement brownien issu de $a$. On a alors $\sup_{t \geq 0} N_t \egaldistr \sup_{t \leq T_0} \beta_t$ où on a posé $T_0 = \inf \{u \geq 0 : \beta_u = 0\}$.
\end{remarque}

\begin{exercice}
Soit $f : \R^+ \to \R^+$ une fonction de $L^1_\text{loc}$ et $F(x) = \int_0^x f(s)\d s$. Soit $B$ mouvement brownien et $S_t = \sup_{s \in [0,t]} B_s$. Montrer que 
\[(f(S_t)(S_t-B_t)-F(S_t))_{t \geq 0}\]
est une martingale locale.
\end{exercice}

\begin{proof}
Une preuve directe est d'utiliser le théorème de Lévy pour conclure en utilisant l'exercice précédent.

Une autre preuve est simplement d'utiliser la formule d'Itô comme précédemment pour $f$ de classe $\mathcal{C}^2$ et une suite d'approximations régularisantes.
\end{proof}

\begin{remarque}
Ce résultat peut être généralisé à toute martingale continue.
\end{remarque}

\begin{exercice}
Soit $\alpha \geq 0$, $(X_t)$ un mouvement brownien avec dérive négative $-\alpha$ issu de $x$ et $I_t = \inf_{s \leq t} X_s$, on pose :
\[M^{(-\alpha)}_t = I_{t \wedge T_0} + \frac{1}{2\alpha} \left[\exp\left( 2\alpha(X_t-I_t) \right)-1\right] \ind{t\leq T_0},\]
montrer que $M^{(-\alpha)}$ est une martingale.
\end{exercice}

\begin{proof}
On observe que la loi du mouvement brownien avec dérive $(X_t)$ est à densité par rapport à la mesure de Wiener, de densité 
\[\mathcal{E}^{(-\alpha)}_t = \exp\left( - \alpha (X_t-x) -\frac{\alpha^2}{2}t \right).\]

Montrer que $(M^{(-\alpha)}_t)$ est une martingale revient donc à montrer que 
\[\Phi^{(-\alpha)}_t = M^{(-\alpha)}_t \mathcal{E}^{(-\alpha)}_t,\]
en est une lorsque $X$ est un mouvement brownien standard issu de $x$.
Gr\^ace à la formule d'Itô, on écrit
\[ \Phi^{(-\alpha)}_t = \Phi^{(-\alpha)}_0 + \int_0^t \mathcal{E}^{(-\alpha)}_s \d M^{(-\alpha)}_s + \int_0^t M^{(-\alpha)}_s \d \mathcal{E}^{(-\alpha)}_s + \crochet{M^{(-\alpha)},\mathcal{E}^{(-\alpha)}}_t. \]

On note :
\[F^\alpha_t = \frac{1}{2\alpha} \left[\exp\left( 2\alpha(X_t-I_t) \right)-1\right] \ind{t\leq T_0}, \]
on observe que
\begin{align*}
  \crochet{M^{(-\alpha)},\mathcal{E}^{(-\alpha)}}
  &= - \alpha \int_0^t \d \crochet{X,F^\alpha}_s \mathcal{E}^{(-\alpha)}_s \\
  &= - \alpha \int_0^{t \wedge T_0} \mathcal{E}^{(-\alpha)}_s\exp(2\alpha(X_s-I_s)) \d s,
\end{align*} 
et que
\[ \int_0^t \mathcal{E}^{(-\alpha)}_s \d M^{(-\alpha)}_s = \int_0^{t \wedge T_0} \mathcal{E}^{(-\alpha)}_s \d I_s + \int_0^{t \wedge T_0} \mathcal{E}^{(-\alpha)}_s \d F^\alpha_s.  \]
Or,
\begin{multline*}
  \int_0^{t \wedge T_0} \mathcal{E}^{(-\alpha)}_s \d F^\alpha_s\\
  = \int_0^{t \wedge T_0} \exp\left( 2\alpha(X_s-I_s) \right) \mathcal{E}^{(-\alpha)}_s (\d X_s-\d I_s)
    + \alpha \int_0^t \exp\left( 2\alpha(X_s-I_s) \right) \mathcal{E}^{(-\alpha)}_s \d s.
\end{multline*}

Par conséquent, on a
\begin{multline*}
  \Phi^{(-\alpha)}_t = \Phi^{(-\alpha)}_0 + \int_0^t M^{(-\alpha)}_s \d \mathcal{E}^{(-\alpha)}_s+\int_0^{t \wedge T_0} \exp\left( 2\alpha(X_s-I_s) \right) \d X_s\\
   + \int_0^{t \wedge T_0} \mathcal{E}^{(-\alpha)}_s\left[1 - \exp\left(2\alpha(X_s - I_s)\right) \right] \d I_s.
\end{multline*}
En utilisant le fait que sur tout intervalle de croissance de $I_t$, on a $X_t=I_t$, on conclut que $\Phi^{(-\alpha)}$ est une martingale locale, car la partie à variations finies est nulle.

Une autre preuve de ce résultat est de vérifier que $f(t,i,b)$ vérifie les équations différentielles obtenues dans l'Exercice \ref{exo_martingaleKennedy}, puis d'utiliser le Théorème d'équivalence de Lévy.
\end{proof}

\begin{exercice}
Soit $B$ un mouvement brownien, $L$ son temps local en 0 et $S$ le processus de son supremum courant.
\begin{enumerate}
  \item Montrer que :
  \begin{itemize}
   \item à $t$ fixé $S_t \egaldistr L_t \egaldistr |B_t| \egaldistr S_t-B_t$,
   \item pour tout $a \geq 0$ et tout $b \leq a$, on a $\P(S_t\geq a,B_t\leq b) = \P(B_t\geq 2a -b)$.
  \end{itemize}
  \item En déduire qu'à $t$ fixé, $2S_t-B_t$ est égal en loi à la norme d'un mouvement brownien dans $\R^3$.
  \item Montrer que conditionnellement à $2S_t-B_t=r,$ $S_t-B_t$ et $S_t$ sont uniformément distribuées sur $[0,r]$.
  \item En déduire la loi jointe de $L_t$ et $|B_t|$ et que conditionnellement à $L_t+|B_t|=r$, $|B_t|$ et $L_t$ sont uniformément distribuées sur $[0,r]$.
\end{enumerate}
\end{exercice}

\begin{proof}
\textit{1.} Soit $a\geq 0$ et $b \leq a$, on pose $T_a =\inf\{t\geq 0:B_t=a\}$, on a :
\begin{align*}
  \P(S_t\geq a, B_t\leq b) = & \P(T_a \leq t, B_t\leq b)\\
  = & \P(T_a\leq t, B_{t-T_a+T_a}-B_{T_a}\leq b-a)\\
  = & \P(T_a \leq t, B^{(T_a)}_{t-T_a}\leq b-a),
\end{align*} 
où on a posé, $B^{(T_a)}_u=B_{T_a+u}-B_{T_a}$. Or le mouvement brownien$B^{(T_a)}$ est en particulier indépendant de $T_a$, par conséquent, on a $(B^{(T_a)},T_a) \egaldistr (-B^{(T_a)},T_a)$. Ainsi on obtient :
\begin{align*}
  \P(S_t\geq a, B_t\leq b) = & \P(T_a \leq t, B^{(T_a)}_{t-T_a}\geq a-b)\\
  = & \P(T_a\leq t, B_{t-T_a+T_a}\geq 2a-b)\\
  = & \P(B_t \geq 2a-b)
\end{align*} 
car on a $2a-b\geq a$.

En particulier, il s'ensuit :
\[\P(S_t \geq x) = \P(S_t \geq x, B_t\geq x) + \P(S_t\geq x , B_t\leq x) = 2\P(B_t \geq x) = \P(|B_t|\geq x).\]
Le théorème d'équivalence de Lévy implique que $S_t \egaldistr L_t$ et $|B_t| = \egaldistr S_t-B_t$. Le fait que $S_t \egaldistr |B_t|$ est simplement dû au Lemme de réflexion, que nous venons de prouver.

\textit{2.} Ce lemme nous donne également la densité du couple $(S_t,B_t)$ par rapport à la mesure de Lebesgue :
\[\P(S_t \in \d a ,B_t \in \d b) = 2\dfrac{2a-b}{\sqrt{2\pi t^3}}e^{-\dfrac{(2a-b)^2}{2t}}\ind{a\geq 0, b\leq a} \d a \d b.\]

Par conséquent, en utilisant un bon changement de variable et intégrant, on obtient :
\[\P(2S_t-B_t \in \d r) = \int_{b=0}^r \dfrac{r}{\sqrt{2\pi t^3}} e^{-\dfrac{r^2}{2t}} \ind{r\geq 0}\d b\d r =\dfrac{r^2}{\sqrt{2\pi t^3}}e^{-\dfrac{r^2}{2t}}\ind{r\geq 0}\d r,\]
qui est la densité d'un processus de Bessel de dimension 3 à l'instant $t$.

\textit{3.} Il suffit maintenant de réaliser un autre changement de variables pour obtenir :
\[\P(S_t -B_t \in \d x, S_t \in \d y ) = 2\dfrac{x+y}{\sqrt{2\pi t^3}}e^{-\dfrac{(x+y)^2}{2t}}\ind{y\geq 0, x\geq 0} \d x \d y.\]
On observe alors simplement que conditionnellement à $S_t-B_t+S_t=r$, $S_t-B_t$ et $S_t$ sont uniformes sur $[0,r]$.

\textit{4.} Finalement, par théorème d'équivalence de Lévy, ce résultat est aussi vrai pour $|B_t|$ et $L_t$, et leur loi jointe est donc donnée par :
\[\P(|B_t| \in \d x, L_t \in \d y ) = 2\dfrac{x+y}{\sqrt{2\pi t^3}}e^{-\dfrac{(x+y)^2}{2t}}\ind{y\geq 0, x\geq 0} \d x \d y.\]
\end{proof}

\begin{remarque}
Les résultats obtenus dans l'exercice précédent seront par la suite généralisés par le théorème de Pitman, on verra alors que $R_t=2S_t-B_t$ peut être vu comme un processus de Bessel 3 et $S_t$ comme l'infimum de $R_u$ pour $u \geq t$ ; or $\frac{1}{R_t}$ est une martingale locale positive, donc l'exercice \ref{exo_supmartingale} peut \^etre appliqué.
\end{remarque}

\begin{exercice}
\label{exo_representationSkorokhod}
Soit $B$ mouvement brownien standard issu de 0 et $S_t = \sup_{s\leq t} B_s$.

Soit $0\leq \alpha \leq 1$, on pose $T_a^{(\alpha)}=\inf\{t\geq 0:S_t=\alpha B_t+a\}$.

Déterminer la loi de $B_{T_a^{(\alpha)}}$ et $S_{T_a^{\alpha}}$.
\end{exercice}

\begin{proof}
Supposons pour commencer $\alpha \in (0,1)$. On a $S_{T_a^{(\alpha)}}=\alpha B_{T_a^{(\alpha)}} +a$. Notons également que $B_{T_a^{(\alpha)}} \in \left[ -\frac{a}{\alpha}, \frac{a}{1-\alpha} \right].$

On rappelle que pour tout $f \in C_c(\R^+)$, le processus : 
\[(f(S_t)(S_t-B_t)+\int_{S_t}^{+\infty} f(u)\d u, t\geq 0)\]
est une martingale bornée dans $L^2$, et est donc uniformément intégrable. En appliquant le théorème d'arrêt en $T=T^{(\alpha)}_a$, on obtient
\[\int_0^{+\infty} f(u)\d u = \E\left[ f(S_T)(S_T-B_T)+\int_{S_T}^{+\infty} f(u)\d u \right],\]
que l'on peut réécrire ainsi :
\[\E\left(\int_0^{\alpha B_T +a} f(x)\d x\right) = \E(f(\alpha B_T +a)((\alpha-1)B_T+a).\]

On pose alors $g(x)= \P(B_T\geq x)$, on a :
\[\E\left(\int_0^{\alpha B_T +a} f(x)\d x\right)= \alpha \int_{\frac{-a}{\alpha}}^{+\infty} f(\alpha x+a)g(x)\d x\]
ainsi que :
\[\E(f(\alpha B_T +a)((\alpha-1)B_T+a)=-\int_{\frac{-a}{\alpha}}^{+\infty} f(\alpha u+a)((\alpha-1)u+a)\d g(u).\]

La fonction $g$ satisfait donc l'équation différentielle suivante :
\[\alpha g(x) = -((\alpha-1)x+a)g'(x),\]
avec condition initiale $\P(B_T \geq \frac{-a}{\alpha})=1$, dont la solution est :
\[g(x) = \left(\frac{\alpha}{a}\right)^\frac{\alpha}{1-\alpha} ((\alpha-1)x+a)^\frac{\alpha}{1-\alpha}.\]
La densité de $B_T$ par rapport à la mesure de Lebesgue est donc égale à :
\[-g'(x) = \left(\frac{\alpha}{a}\right)^\frac{\alpha}{1-\alpha}  \alpha \left((\alpha-1)x+a\right)^\frac{2\alpha-1}{1-\alpha}.\]

Dans le cas où $\alpha=1$, on peut raisonner de la même manière ou en passant à la limite. En gardant les mêmes notations, on obtient :
\[g = -a g',\]
par conséquent en intégrant, on observe que $B_{T_a^{(1)}}$ est une variable aléatoire de loi exponentielle de paramètre $a$ centrée (à valeurs dans $[-a,+\infty)$).

On peut s'intéresser au cas où $\alpha=0$. On entre dans une autre classe de variables aléatoires pour lesquelles $B_T$ n'est plus centrée. En effet $T_a^0$ est simplement le temps d'atteinte de $a$ par $S_t$. On a par conséquent $S_{T_a^0}=a$, et $B_{T_a^0}=a$.
\end{proof}

\begin{exercice}[Lemme de Skorokhod et changement de temps]
Soit $B$ un mouvement brownien standard issu de 0 et $L$ son temps local en 0.

\begin{enumerate}
 \item Montrer l'égalité entre processus :
\[\left( \log\left(1+\frac{|B_t|}{h}\right),\frac{L_t}{h} \right)_{t \geq 0} = \left( S_{H_t}-X_{H_t},S_{H_t} \right)_{t\geq 0},\]
où $\beta$ désigne un autre mouvement brownien issu de 0 bien choisi, et pour $u \geq 0$ 
\[X_u=\beta_u+\frac{u}{2}, \text{  } S_u=\sup_{s \leq u} X_s \text{  et  }H_t=\int_0^t\frac{\d s}{(h+|B_s|)^2}.\]

 \item Soit $\tau_l=\inf\{t>0:L_t\geq l\}$ ; montrer que $H_{\tau_l}=\inf\{u>0:X_u=\frac{l}{h}\}$ est un subordinateur ; calculer son exposant caractéristique et sa mesure de Lévy.

 \item Soit $T^*_a=\inf\{t>0: |B_t|=a\}$, montrer que $H_{T^*_a} = \inf\{u>0:S_u-X_u = \log (1+\frac{a}{h})\}$. $(H_{T^*_a},a\geq 0)$ est-il un subordinateur ? un processus à accroissements indépendants ?
\end{enumerate}
\end{exercice}

\begin{proof}
\textit{1.} Rappelons pour commencer la formule de Tanaka :
\[|B_t|=\int_0^t \sgn(B_s)\d B_s + L_t.\]
Nous pouvons donc écrire, grâce à la formule d'Itô :
\[\log \left(1+\frac{ |B_t| }{h}\right) = \int_0^t \dfrac{1}{h+|B_s|} \sgn(B_s)\d B_s + \int_0^t \dfrac{1}{h+|B_s|} \d L_s - \dfrac{1}{2} \int_0^t \dfrac{1}{(h+|B_s|)^2}\d s.\]

On utilise alors que $\d_sL_s$ est portée par $\{ s \geq 0: |B_s| = 0 \}$
\[\ln \left(1+\frac{|B_t|}{h}\right) = \int_0^t \dfrac{1}{h+|B_s|} \sgn(B_s)\d B_s + \dfrac{L_t}{h} - \dfrac{1}{2} \int_0^t \dfrac{1}{(h+|B_s|)^2}\d s.\]

Par conséquent, $\log (1+\frac{|B_t|}{h})-\frac{L_t}{h}$ est la somme d'une martingale (car $\frac{1}{h+|B_s|} \leq \frac{1}{h}$) de variation quadratique $\langle M\rangle_t=H_t$ et de $-\frac{1}{2}H_t$. Il existe donc un mouvement brownien $\beta$ tel qu'on ait :
\[\log \left(1+\frac{|B_t|}{h}\right)-\frac{L_t}{h} = - \beta_{H_t} - \frac{1}{2}H_t.\]
On pose alors $X_u=\beta_u + \frac{u}{2}$, on a, d'après l'égalité précédente :
\[X_{H_t} = \frac{L_t}{h} - \log\left(1+\frac{|B_t|}{h}\right).\]

On utilise alors le Lemme de Skorokhod \ref{lem_skorokhod} : on peut décomposer une fonction continue $x$, en un unique couple $(z,l)$ de fonctions continues telles que $l$ soit croissante, $z(t)=-x(t)+l(t)$ et $\d l$ est portée par $\{t\geq 0:z(t)=0\}$.

Cette décomposition se retrouve pour $(X_{H_t})_{t\geq 0}$, le couple associé est $(S_{H_t}-X_{H_t},S_{H_t})_{t\geq 0}$. De plus, la décomposition de $(-\log (1+\frac{|B_t|}{h})+\frac{L_t}{h})_{t\geq 0}$ est bien $(\log (1+\frac{|B_t|}{h}),\frac{L_t}{h})_{t\geq 0}$. En effet, $L_t$ est croissante, et ne croit que lorsque $|B_t|$ vaut 0.

Grâce au Lemme de Skorokhod, on a bien obtenu l'égalité des deux processus : 
\[\left( \log(1+\frac{|B_t|}{h}),\frac{L_t}{h} \right)_{t \geq 0} = \left( S_{H_t}-X_{H_t},S_{H_t} \right)_{t\geq 0}.\]

\textit{2.} Cette égalité nous permet d'étudier $H_{\tau_l}$. On a
\[\tau_l = \inf\{t\geq 0: L_t \geq l\} = \inf\{t \geq 0:hS_{H_t} \geq l\},\]
par conséquent $H_{\tau_l}$ est le premier temps d'atteinte de $\frac{l}{h}$ par $S$, et donc, $H_{\tau_l}$ est également le premier temps d'atteinte de $\frac{l}{h}$ par $X$, d'où :
\[H_{\tau_l} = \inf\{u>0:\beta_u + \frac{u}{2}=\frac{l}{h}\}.\]

On calcule alors la transformée de Laplace de $H_{\tau_l} = \tilde{\tau}_l$ en appliquant le théorème d'arrêt à la martingale exponentielle de $\beta$,
\[1 = \E\left[\exp\left(\lambda \beta_{\tilde{\tau}_l} -\frac{\lambda^2}{2}\tilde{\tau}_l\right)\right] = \E\left[\exp\left(\lambda \left(\frac{l}{h} - \frac{\tilde{\tau}_l}{2}\right) -\frac{\lambda^2}{2} \tilde{\tau}_l\right)\right].\]

On obtient donc :
\[\E\left(\exp\left(-\tilde{\tau}_l\frac{\lambda+\lambda^2}{2}\right)\right) = \exp\left(-\lambda \frac{l}{h}\right).\]

Par changement de variables, on a :
\[\E(\exp(-\mu \tilde{\tau}_l)) = \exp \left(-l \frac{-1+\sqrt{1+8\mu}}{2h}\right).\]
Par conséquent $\tilde{\tau}_l$ est un subordinateur d'exposant caractéristique $\frac{-1+\sqrt{1+8\lambda}}{2h}$. Cherchons la mesure de Lévy $\nu$ associée à $\tilde{\tau}_l$, vérifiant
\[\int_\R \nu(\d t) \left(1-e^{-\lambda t}\right) = \frac{-1+\sqrt{1+8\lambda}}{2h}.\]

Pour ce faire, on dérive l'expression précédente par rapport à $\lambda$, en utilisant la représentation des puissances négatives par une intégrale :
\[\int_\R \nu(\d t) t e^{-\lambda t} = \frac{2}{h\sqrt{1+8\lambda}} = \int_{\R^+} \frac{2\d t}{h\sqrt{\pi t}} \exp(-(1+8\lambda)t).\]
Par unicité de la transformée de Laplace, et en réalisant le bon changement de variables, on obtient :
\[t \nu(\d t) = \d t\frac{e^{-\frac{t}{8}}}{\sqrt{2\pi t}}.\]

La mesure de Lévy associée à $(\tilde{\tau}_l,l\geq 0)$ est donc $\nu(\d t) = \dfrac{1}{h\sqrt{2\pi t^3}}e^{-\frac{t}{8}}\d t$.

\textit{3.} En raisonnant de la même manière avec $H_{T^*_a}$, on a :
\begin{align*}
  T^*_a & = \inf\{t\geq 0: |B_t| = a\}\\
  & = \inf\left\{t \geq 0:\log\left(1+\frac{|B_t|}{h} \right) = \log \left(1+\frac{a}{h}\right)\right\}\\
  & = \inf\left\{t \geq 0:S_{H_t}-X_{H_t} = \log\left(1+\frac{a}{h}\right)\right\}.
\end{align*}

Par conséquent on a :
\[H_{T^*_a} = \inf\left\{u\geq 0: S_u-X_u = \log \left(1+\frac{a}{h}\right)\right\}.\]
Le processus $(H_{T^*_a},a \geq 0)$ n'est pas un subordinateur, car bien qu'il soit à accroissements indépendants, ces derniers ne sont pas stationnaires. En effet, si l'on pose :
\[U_a = \inf \{t \geq 0: X_t \geq \log\left(1+\frac{a}{h}\right)\}\]
\[\text{et  } V_{a,b}=\inf\{t \geq 0 :X_t \leq \log\left(1+\frac{a+b}{h}\right) - \log\left(1+\frac{a}{h}\right)\},\]
alors par application de la propriété de Markov forte en $H_{T^*_a}$, qui est bien un temps d'arrêt de $S-X$ on a :
\[H_{T^*_{a+b}}-H_{T^*_a} \egaldistr (U_a + H_{T^*_{a+b}}) \ind{U_a \leq V_{a,b}} + V_{a,b} \ind{U_a \leq V_{a,b}}.\]
\end{proof}

\chapter[Théorème de Pitman]{Théorème de Pitman}

Nous nous intéressons ici à une égalité en loi entre processus, qui complète celle donnée par le théorème de Lévy (vu au chapitre 2). Cette égalité permet une identification du processus joint du mouvement brownien et de son supremum courant avec un processus de Bessel de dimension 3 et son infimum futur. Ce résultat s'obtient en utilisant le théorème général de retournement du temps de Nagasawa pour de \og bons \fg{} processus de Markov, théorème que nous citons ici, bien qu'une démonstration détaillée de celui-ci nous amènerait à des considérations éloignées de notre sujet principal. Les résultats que nous obtenons ici sont d'une importance certaine par la suite, le processus de Bessel de dimension 3 jouant un rôle fondamental dans l'étude du mouvement brownien linéaire, par exemple pour certains calculs sur les lois des temps locaux.

\section{Théorèmes principaux}

On introduit ici pour commencer le processus de Bessel de dimension 3, objet central de ce chapitre.

\begin{definition}
Les trois définitions suivantes sont équivalentes.

Soit $(B_t,t \geq 0)$ un mouvement brownien dans $\R^3$. Le processus $(\norme{B_t},t \geq 0)$ est un processus de Bessel de dimension 3.

Soit $(B_t, t \geq 0)$ un mouvement brownien dans $\R$. La solution de l'équation différentielle stochastique $\d X_t = \d B_t + \frac{\d t}{X_t}$ est un processus de Bessel de dimension 3.

Soit $(B_t, t \geq 0)$ un mouvement brownien réel issu de $a \in \R$ sous la mesure de Wiener $\W_a$. On note $T_0=\inf\{t \geq 0 : B_t = 0\}$. On définit la loi
\[ {\P^{(3)}_a}_{|\F_t} = \frac{B_{t\wedge T_0}}{a} . {\W_a}_{|\F_t}.\]
Sous la loi $\P^{(3)}_a$, le processus $B$ est un processus de Bessel de dimension 3 issu de $a$.
\end{definition}

On note $\P^{(3)}_a$ la loi d'un processus de Bessel de dimension 3 issu de $a$, et $\E^{(3)}_a$ l'espérance correspondance. Par la suite, nous parlerons également de processus de Bessel 3, voire de BES(3), pour simplifier l'énoncé. Citons dès maintenant le théorème central de ce chapitre.

\begin{theoreme}[Pitman]
Soit $B$ mouvement brownien issu de 0, et $S_t = \sup_{s \leq t} B_s$ ; soit $(R_t)_{t \geq 0}$ un processus de Bessel 3 et $J_t = \inf_{s \geq t} R_s$ son infimum futur. On a :
\[(2S_t-B_t, S_t)_{t \geq 0} \egaldistr (R_t, J_t)_{t \geq 0}.\] 
\end{theoreme}

La preuve de ce théorème repose sur deux identités en loi, relatives au mouvement brownien retourné en temps. Cette égalité en loi découlent du théorème de retournement du temps de Nagasawa.

\begin{theoreme}[Nagasawa]
Soit $X_t$ une diffusion réelle transiente et $\gamma_a = \sup\{t \geq 0 ; X_t=a\}$ le dernier temps d'atteinte de $a$. On pose :
\[\tilde{X}_t = \left\{
\begin{array}{ll}
  X_{\gamma_a-t} & \text{ si  }  0\leq t \leq \gamma_a <+\infty\\
  \partial \text{  (point \og cimetière\fg)}& \text{  sinon.}
\end{array}
\right.\]
le processus $X$ retourné au dernier temps d'atteinte de $a$.

Soit $\mu$ une mesure de probabilité sur $\R$, on note $\nu$ la mesure définie par,
\[ \forall f \in \mathcal{C}_c, \int f \d \nu = \int U(f) \d \mu, \]
où $U$ est l'opérateur potentiel de $X$, lui-m\^eme défini par
\[ \forall f \in \mathcal{C}_c, y \in \R, U(f)(y) = \int_0^{+\infty} \E_y(f(X_t)) \d t.  \]

Soit $(P_t)$ le semi-groupe de $X$ ; on suppose qu'il existe $\widehat{P}_t$ un semi-groupe tel que
\[ \forall f, g \in \mathcal{C}_c, \int P_t (f) g \d \nu = \int f \widehat{P}_t (g) \d \nu.\]
Dans ce cas, sous $\P_\mu$, $\tilde{X}$ est un processus de Markov de semi-groupe associé $\widehat{P}_t$.
\end{theoreme}

\begin{corollaire}
Soit $B$ mouvement brownien issu de 0 et $R$ processus de Bessel 3. On introduit $\gamma_a = \sup\{t \geq 0:R_t=a\}$, on a alors :
 \[(B_u, u \leq \tau_l) \egaldistr (B_{\tau_l-u}, u \leq \tau_l)\text{  et  } (R_u, u \leq \gamma_a) \egaldistr (a-B_{T_a-u}, u \leq T_a).\]
\end{corollaire}

La démonstration de ce corollaire figure en exercice, et est en grande partie une conséquence du lemme suivant.

\begin{lemme}[$P^{(3)}_a$ comme transformée de Doob de $W_a$]
Soit $P^{(3)}_a$ la loi du processus de Bessel en dimension 3 issu de $a$, et $W_a$ la mesure de Wiener vérifiant $W_a(X_0=a)=1$. On note (classiquement) $\F_t = \sigma(X_s, s \leq t)$. On a :
\[P^{(3)}_a{}_{|\F_t} = \frac{X_{t \wedge T_0}}{a} \cdot W_a{}_{|\F_t} \text{  pour tout } t \geq 0.\]
\end{lemme}

Nous pouvons maintenant étendre le théorème de Pitman à des mouvements browniens avec dérive, comme nous l'avions fait avec le théorème de Lévy dans le chapitre précédent.

\begin{theoreme}[Extension du théorème de Pitman]
Soit $\mu \in \R$ ; on pose 
\[B^{(\mu)}_t= B_t +\mu t,\]
où $B$ est un mouvement brownien standard. Les quantités se rapportant à $B^{(\mu)}$ sont également notées avec $(\mu)$ en exposant.

Les processus $(L^{(\mu)}_t + |B^{(\mu)}_t|)_{t \geq 0}$ et $(2S^{(\mu)}_t-B^{(\mu)}_t)_{t\geq 0}$ ont même loi, celle de la diffusion de générateur étendu $\frac{1}{2}\frac{\d^2}{\d x^2} +\mu \coth(\mu x) \frac{\d}{\d x}.$
\end{theoreme}

\section{Exercices}

\begin{exercice}
Soit $(B_t)$ un mouvement brownien réel et $S_t= \sup_{s \leq t} B_s$. On pose :
\[\Gamma_a = \inf\{t \geq 0 : 2S_t-B_t=a\},\]
déterminer la loi de $B_{\Gamma_a}$.
\end{exercice}

\begin{proof}
Le théorème de Pitman donne immédiatement~:
\[(2S_t-B_t,S_t ; t \geq 0) \egaldistr (R_t,J_t ; t\geq 0).\]
En particulier, on a
\[(\Gamma_a,B_{\Gamma_a} ) \egaldistr (U_a,2J_{U_a} - R_{U_a}),\]
où on a posé $U_a = \inf\{t \geq 0 : R_t=a\}$.

On observe pour commencer que $R_{U_a}=a$. De plus on peut réécrire
\[J_{U_a} = \inf_{s \geq U_a} R_s = \frac{1}{\sup_{s\geq 0} \frac{1}{R_{U_a+s}}},\]
or $(\frac{1}{R_{U_a+s}}, s \geq 0)$ est une martingale locale positive qui tend vers 0, issue de $\frac{1}{a}$, en appliquant le résultat de l'Exercice \ref{exo_supmartingale}, $J_{U_a} \egaldistr Ua$, où $U$ est une variable aléatoire uniforme sur $[0,1]$.

Par conséquent, $B_{\Gamma_a} \egaldistr 2aU-a$, donc $B_{\Gamma_a}$ est uniformément distribuée sur $[-a,a]$.
\end{proof}

\begin{exercice}[Démonstration du Théorème 3.4.]
Soit $\mu \in \R$, on pose $B^{(\mu)}_t=B_t+\mu t$ et $S^{(\mu)}_t = \sup_{s\leq t} B^{(\mu)}_s$.

Montrer que $(2S^{(\mu)}_t-B^{(\mu)}_t)_{t \geq 0}$ est un processus de diffusion de générateur étendu $\frac{1}{2}\frac{\d^2}{\d x^2}+\mu \coth (\mu x) \frac{\d}{\d x}.$ 
\end{exercice}

\begin{remarque}
Le générateur étendu $\mathcal{G}$ d'un processus de Markov $X$ est défini pour toute fonction $f$ telle qu'il existe $\psi$ fonction mesurable vérifiant :
\[\left(f(X_t)-\int_0^t \psi(X_s) \d s, t\geq 0\right) \text{   est une martingale}.\]

Dans ce cas, on note $\psi = \mathcal{G} (f)$. Cette définition étendue du générateur d'un processus de Markov est due à Kunita.
\end{remarque}

\begin{proof}
Soit $F : \mathcal{C}([0,t], \R) \to \R^+$. On a, en utilisant le théorème de Girsanov
\[\E(F(2S^{(\mu)}_s-B^{(\mu)}_s,s \leq t)) = \E\left[ F(2S_s-B_s, s \leq t) \exp\left(\mu B_t - \frac{\mu^2}{2}t\right)\right].\]

Soit $R$ un processus de Bessel de dimension 3 et $J_t = \inf_{s \geq t} R_s$. En appliquant le théorème de Pitman, on a :
\[
  \E( F(2S_s-B_s,s \leq t) \exp(\mu B_t - \frac{\mu^2}{2}t)) = \E(F(R_s,s \leq t) \exp(\mu (2J_t-R_t)-\frac{\mu^2}{2}t)).
\]
Conditionnellement à $R_t$, la variable $J_t$ est distribuée uniformément sur $[0,R_t]$, donc~:
\begin{align*}
  & \E(F(2S^{(\mu)}_s-B^{(\mu)}_s,s \leq t))\\
  = & \E \left( F(R_s, s \leq t) \exp(-\frac{\mu^2}{2}t)  \frac{1}{R_t}\int_0^{R_t} \exp(\mu (2u-R_t))\d u \right)\\
  = & \E \left( F(R_s,s \leq t) \exp(-\frac{\mu^2}{2}t) \frac{1}{2R_t} \int_{-R_t}^{R_t}  \exp(\mu u)\d u \right)\\
  = & \E \left( F(R_s,s \leq t) \exp(-\frac{\mu^2}{2}t) \frac{\sh (\mu R_t) }{\mu R_t} \right).
\end{align*}

Pour utiliser le théorème de Girsanov, on pose $\mathbb{Q}_{|\F_t} = D_t \P_{|\F_t}$, où on a noté 
\[D_t = \exp(-\frac{\mu^2}{2}t) \frac{\sh (\mu R_t) }{\mu R_t}.\]
Sous $\P$, $R$ est solution de l'équation différentielle stochastique :
\[ X_t = B'_t +\int_0^t \frac{\d s}{X_s},\]
avec $B'$ un mouvement brownien. Par conséquent, sous $\Q$,
\[R_t = \tilde{B}_t + \int_0^t \frac{\d r}{R_r} + \int_0^t \frac{\d \langle D,B \rangle_r}{D_r},\]
où $\tilde{B}$ est un mouvement brownien (par le théorème de Lévy). De plus :
\begin{align*}
 & \int_0^t \frac{\d r}{R_r} + \int_0^t \frac{\d \langle D,R \rangle_r}{D_r}\\
 = & \int_0^t \frac{1}{R_r} + \frac{\mu R_r}{\sh(\mu R_r)}\left(\frac{\ch(\mu R_r)}{R_r} - \frac{\sh(\mu R_r)}{\mu R_r^2}\right)\d r\\
 = & \int_0^t \mu \coth( \mu R_r) \d r.
\end{align*}

Par conséquent la loi de $(2S^{(\mu)}_s-B^{(\mu)}_s, s \leq t)$ sous $\P$ est égale à la loi de $(R_s)_{s \leq t}$ sous $\Q$, c'est donc une diffusion qui vérifie l'équation différentielle stochastique suivante :
\[ \d X_t = \d B_t + \mu \coth (\mu X_t) \d t.\]
C'est donc bien la diffusion espérée, de générateur étendu $\frac{1}{2}\frac{\d^2}{\d x^2} +\mu \coth(\mu x) \frac{\d}{\d x}$.
\end{proof}

\begin{exercice}[Preuve du Corollaire 3.3.]
Soit $Q_t$ le semigroupe de $(B_{t\wedge T_0})_{t \geq 0}$.

\begin{enumerate}
 \item Montrer que $Q_t(a,\d y) = q_t(a,y)\d y$, où on a posé :
\[q_t(a,y) = \dfrac{1}{\sqrt{2\pi t}} \left[ \exp\left( -\dfrac{(a-y)^2}{2t} \right) - \exp\left( -\dfrac{(a+y)^2}{2t} \right) \right].\]

 \item Expliciter la mesure potentielle en 0 de $X_t$ sous $\P^{(3)}$ donnée par :
\[U(f) : a \mapsto \int_0^{+\infty} \d t \E^{(3)}_a(f(X_t)).\]
En déduire la loi du processus défini comme le retourné (en temps) du processus de Bessel en son dernier temps de passage en $a$.
\end{enumerate}
\end{exercice}

\begin{proof}
\textit{1.} Soit $B$ mouvement brownien issu de $a>0$, et $T_0 = \inf \{t>0: B_t= 0\}$. Calculons, pour $f \in \mathcal{C}_b$ et $a \in \R$ la valeur de $Q_t f (a) = \E(f(B_t) \ind{t \leq T_0})$. Pour cela, on définit $g$ une fonction mesurable bornée définie sur $\R$ de la manière suivante :
\[g(x) = \left\{
\begin{array}{ll}
 f(x) & \text{ si } x > 0,\\
 -f(-x) & \text{ si } x < 0,\\
 0 & \text{ si } x = 0.
\end{array}
\right.\] 
On remarque alors que
\[\E(g(B_t)) = \E(g(B_t) \ind{t \leq T_0} + g(B_{t-T_0+T_0}) \ind{t > T_0}).\]
On utilise la propriété de Markov forte, $\beta$ désignant un mouvement brownien indépendant de $B$ issu de 0, on a
\[
  \E(g(B_t)) = \E ( f(B_t) \ind{t \leq T_0} ) + \E(g(\beta_{t - T_0})\ind{t \geq T_0}).
\]
On observe que $\beta$ a une distribution symétrique et est indépendant de $T_0$, et que $g$ est impaire. Par conséquent on a :
\[\E(g(\beta_{t - T_0})\ind{t \geq T_0}) = \E(g(-\beta_{t - T_0})\ind{t \geq T_0}) = -\E(g(\beta_{t - T_0})\ind{t \geq T_0}),\]
d'où $\E(g(\beta_{t - T_0})\ind{t \geq T_0})=0$.

On en tire
\[Q_tf(a) = \E(g(B_t)) = \frac{1}{\sqrt{2\pi t}}\int_0^{+\infty} \d y f(y) e^{-\frac{(y-a)^2}{2t}} - \frac{1}{\sqrt{2\pi t}}\int_{-\infty}^0\d u f(-u) e^{-\frac{(y-u)^2}{2t}}\]
ce qui implique
\[q_t(a,y) = \dfrac{1}{\sqrt{2\pi t}} \left(e^{-\frac{(y-a)^2}{2t}} - e^{-\frac{(y+a)^2}{2t}} \right).\]

\textit{2.} Soit $X$ processus de Bessel 3 issu de $a$, nous allons calculer $\int_0^{+\infty} \E^{(3)}_a(f(X_t))\d t$.
Pour cela, on utilise la densité de la loi de ce processus par rapport à la mesure de Wiener
\[\E^{(3)}_a(f(X_t)) = \E\left(\frac{B_{t \wedge T_0}}{a} f(B_t)\right) = \frac{1}{a}\E(B_tf(B_t) \ind{t \leq T_0}).\]
En appliquant la formule de la densité obtenue précédemment pour le mouvement brownien tué en 0, on obtient
\begin{align*}
  Uf(a) = & \dfrac{1}{a} \int_0^{+\infty} \d t \int_0^{+\infty} \d y f(y)y q_t(a,y)\\
  = & \dfrac{1}{a} \int_0^{+\infty} \d y f(y)y \int_0^{+\infty} \d t q_t(a,y).
\end{align*}

Pour calculer $Uf(0)$, on observe que $\lim_{a \to 0} \frac{q_t(a,y)}{a} = \frac{y}{\sqrt{2 \pi t^3}} e^{-\frac{y^2}{2t}}$. On a donc
\begin{align*}
  Uf(0) = & \int_0^{+\infty} \d y f(y)y \int_0^{+\infty} \dfrac{y \d t}{\sqrt{2\pi t^3}} e^{-\frac{y^2}{2t}}\\
  = & \int_0^{+\infty} \d y f(y) y.
\end{align*}
En effet, $\frac{y}{\sqrt{2\pi t^3}}\exp(-\frac{y^2}{2t})\ind{t>0}$ est la densité de la loi de $T_y$, et est donc d'intégrale 1. Nous pouvons alors vérifier la relation de dualité entre le processus de Bessel 3 $X_t$ et le mouvement brownien tué en 0
\begin{align*}
  \int_{\R^+} \d a \E_a^{(3)}(f(X_t)) g(a) a = & \int_{\R^+} \d a \frac{a}{a} \E^{(3)}_a(B_t f(B_t) \ind{t \leq T_0}) g(a)\\
  = & \int_{\R^+} \d a \int_{\R^+}\d y f(y)y q_t(a,y)g(a)\\
  = & \int_{\R^+}\d y f(y) y \int_{\R^+}\d a q_t(y,a) g(a) \\
  = & \int_{\R^+}\d y \E_y(g(B_t) \ind{t \leq T_0}) f(y) y .
\end{align*}
en observant que $q_t(a,y)=q_t(y,a)$.

On termine en utilisant le théorème de Nagasawa. Les hypothèses du théorème sont bien vérifiées, par conséquent
\[(R_t, t \leq \gamma_a)  \egaldistr (a-B_{T_a-t}, t \leq T_a).\]
\end{proof}

\begin{exercice}
\label{exo_scalingaleatoire}
Soit $(b_u)_{u \leq 1}$ un pont brownien standard indexé par $[0,1]$ ; on pose $B^{\tau_l}_u = \frac{1}{\sqrt{\tau_l}} B_{\tau_l u}$.

En étudiant, pour $F$ mesurable positive et $\phi$ continue positive à support compact, l'intégrale
\[\E\left(\int_0^{+\infty} \d L_s \phi(s) F\left(\frac{1}{\sqrt{s}}B_{s u}, u \leq 1\right) \right),\]
montrer que
\[\E(F(b_u,u \leq 1)) = \E\left(\sqrt{\dfrac{\pi}{2\tau_1}} F(B^{\tau_1}_u,u \leq 1)\right)\]
ou de manière équivalente que
\[\E\left(\sqrt{\dfrac{2}{\pi}}\dfrac{1}{l_1}F(b_u,u \leq 1)\right) = \E\left(F(B^{\tau_1}_u,u \leq 1)\right).\]
\end{exercice}

\begin{proof}
Observons pour commencer que pour tout $\lambda > 0$, par propriété de scaling du mouvement brownien, on a
\begin{equation}
  \label{eqn:scaling}
  (B^{\tau_1}_u,u \leq 1) \egaldistr (B^{\tau_\lambda}_u,u \leq 1).
\end{equation}

Soit $F : \mathcal{C}([0,1]) \to \R^+$ mesurable positive et $\phi : \R^+ \to \R^+$ continue à support compact, par changement de variable $\d t = \d L_s$
\[ \E\left(\int_0^{+\infty} \d L_s \phi(s) F\left(\frac{1}{\sqrt{s}}B_{s u}, u \leq 1\right)\right) = \E\left(\int_0^{+\infty} \d t \phi(\tau_t) F\left(\frac{1}{\sqrt{\tau_t}}B_{\tau_t u}, u \leq 1\right)\right).\]
On applique le théorème de Fubini, on obtient d'une part
\begin{align*}
  \E \left( \int_0^{+\infty} \d s \phi(\tau_s) F(B^{\tau_s}_u,u \leq 1) \right)
  = & \int_0^{+\infty} \d s \E(\phi(\tau_s)F(B^{\tau_s}_u,u\leq 1))\\
  = & \int_0^{+\infty} \d s \E(\phi(s^2\tau_1) F(B^{\tau_1}_u,u \leq 1))
\end{align*}
grâce à l'Équation \eqref{eqn:scaling}. En définitive, par changement de variables $s^2\tau_1=t$
\[
  \E \left( \int_0^{+\infty} \d s \phi(\tau_s) F(B^{\tau_s}_u,u \leq 1) \right)
  = \int_0^{+\infty} \d t \frac{\phi(t)}{\sqrt{t}} \E\left( \dfrac{1}{2\sqrt{\tau_1}} F(B^{\tau_1}_u,u\leq 1) \right).
\]

D'autre part, en utilisant l'Exercice \ref{exo_integrationtempslocal}, on obtient
\begin{align*}
  &\E\left( \int_0^{+\infty} \d L_s \phi(s) F\left(\frac{1}{\sqrt{s}}B_{su},u \leq 1\right) \right)\\
  & \qquad \qquad \qquad \qquad = \int_0^{+\infty} \dfrac{\d u}{\sqrt{2\pi u}} \phi(u) \E\left[ \left. F\left(\frac{1}{\sqrt{s}}B_{su},u \leq 1\right)\right|B_u=0\right]\\
  & \qquad \qquad \qquad \qquad =\int_0^{+\infty} \dfrac{\d u}{\sqrt{2\pi u}} \phi(u) \E\left[ F(b_u,u \leq 1)\right]
\end{align*}
où $(b_u,0\leq u \leq 1)$ est un pont brownien de longueur 1.

On en conclut l'égalité
\[ \E \left( \int_0^{+\infty} \phi(t) F\left(\frac{1}{\sqrt{t}}B_{t u},u \leq 1\right) \d t \right) =  \E(F(b_u,u \leq 1)) \int_0^{+\infty} \frac{\d t}{\sqrt{2\pi t}}  \phi(t).\]
En particulier, on obtient
\[ \E \left(\sqrt{\frac{\pi}{2\tau_1}}F( B^{\tau_1}_u)_{u \leq 1})\right) =\E \left(F\left((b_u)_{u \leq 1}\right)\right)\]
ce qui est le résultat attendu.

On souhaite obtenir l'identité réciproque. On observe que pour tout $\epsilon >0$,
\begin{align*}
  \dfrac{1}{\epsilon} \int_0^1 \ind{0<B^{\tau_1}_u<\epsilon}\d u
  = & \dfrac{1}{\epsilon \tau_1} \int_0^{\tau_1} \ind{0<B_v<\epsilon \sqrt{\tau_1}} \d v\\
  = & \dfrac{1}{\sqrt{\tau_1}} \left[ \dfrac{1}{\epsilon \sqrt{\tau_1}} \int_0^{\tau_1} \ind{0<B_v<\epsilon \sqrt{\tau_1}} \d v \right].
\end{align*}
Par conséquent, comme $L^0_{\tau_1}=1$, le temps local en 0 de $(B^{\tau_1}_u)_{u \leq 1}$ est égal à $\frac{1}{\sqrt{\tau_1}}$. Ainsi si l'on note $l_1$ le temps local en 0 du pont brownien, on a également :
\[\E\left( \sqrt{\frac{2}{\pi}} \frac{1}{l_1} f((b_u)_{u\leq 1}) \right) = \E\left( f(B^{\tau_1}_u)_{u \leq 1})\right).\]
\end{proof}

\begin{remarque}
Ce résultat, obtenu ici par ce qui semble \^etre une astuce de calcul, peut en réalité \^etre généralisé à de nombreuses autres fonctionnelles que celles relatives aux temps locaux. Cette construction sera étudiée plus en détails dans le Chapitre~5.
\end{remarque}

\chapter[Unicité trajectorielle de solutions d'EDS]{Unicité trajectorielle de solutions d'EDS à coefficients non réguliers}

Nous allons maintenant donner un exemple d'application de la théorie des temps locaux. Les temps locaux permettent en effet de démontrer l'unicité trajectorielle de solutions d'équations différentielles stochastiques à coefficients peu réguliers, pour peu que l'on ait déjà l'unicité en loi. Or on sait que lorsqu'une équation différentielle stochastique vérifie la propriété d'unicité trajectorielle, les solutions s'expriment en terme d'une fonction mesurable du mouvement brownien directeur. C'est le théorème de Yamada-Watanabe. 

Prouver de telles résultats est donc d'une importance cruciale, tant au niveau de la théorie que des applications. Nous verrons ainsi de nombreux exemples pour lesquels la méthode de temps local, due à Jean-François Le Gall et développée ici s'applique, ou non. Pour finir, nous montrerons également un résultat, la \og formule de balayage \fg{} prélude à la théorie des excursions.

\section{Théorèmes principaux}

Soit $\sigma$ et $b$ deux fonctions mesurables, et $B$ un mouvement Brownien. On rappelle que $X$ est une solution de l'équation différentielle stochastique
\[ \d X_t = \sigma(X_t) \d B_t + b(X_t) \d t \]
si pour tout $t \geq 0$, on a 
\[ X_t - X_0 = \int_0^t \sigma(X_s) \d B_s + \int_0^t b(X_s) \d s. \]

\begin{theoreme}[Argument de temps local nul]
Soit $\sigma,b$ deux fonctions mesurables, qui sont les coefficients de l'équation différentielle stochastique :
\[\d X_t = \sigma(X_t) \d B_t +b(X_t)\d t.\]
Si les solutions de cette équation différentielle stochastique sont uniques en loi, et si, pour toute paire $X^{(1)},X^{(2)}$ de solutions définies sur un même espace, on a $L^0_t(X^{(1)}-X^{(2)})=0$, alors on a unicité trajectorielle.
\end{theoreme}

Nous allons maintenant déterminer une condition suffisante sous laquelle le temps local en 0 d'une semi-martingale est bien nul. Soit $\rho:\R^+ \to \R^+$ une fonction continue strictement croissante telle que pour tout $\epsilon > 0$, on a $\int_0^\epsilon \frac{\d x}{\rho(x)}=+\infty$. Une telle fonction $\rho$ une mauvaise fonction ; on notera en particulier que $\rho(x)=x$ est une mauvaise fonction. Le lemme suivant donne une condition suffisante d'annulation du temps local en zéro de $X$.

\begin{lemme}
Si pour une mauvaise fonction $\rho$, il existe $\epsilon>0$ tel que 
\[\int_0^t \frac{\d \crochet{X}_s}{\rho(X_s)}\ind{0<X_s<\epsilon}<+\infty,\]
alors $L^0_t(X)=0$ p.s.
\end{lemme}

En utilisant ce lemme, on peut trouver plusieurs classes d'équations pour lesquelles on a bien unicité trajectorielle.

\begin{theoreme}
L'équation différentielle stochastique associée à $(\sigma,b)$ jouit de l'unicité trajectorielle si l'un des jeux de conditions suivantes est réalisée :
\begin{enumerate}
 \item 
 \begin{itemize}
   \item[\textbullet] $(\sigma(x)-\sigma(y))^2 \leq \rho(|y-x|)$ et $|\sigma| \geq \epsilon >0$,
   \item[\textbullet] $b$ et $\sigma$ bornées.
 \end{itemize}
 \item
 \begin{itemize}
   \item[\textbullet] $(\sigma(x)-\sigma(y))^2 \leq \rho(|y-x|),$
   \item[\textbullet] $b$ lipschitizienne.
 \end{itemize}
 \item
 \begin{itemize}
  \item[\textbullet] $(\sigma(x)-\sigma(y))^2\leq |f(x)-f(y)|$ avec $f$ fonction croissante bornée,
  \item[\textbullet] $\sigma\geq \epsilon>0$ et $b$ bornée.
 \end{itemize}
\end{enumerate}
\end{theoreme}

\section{Exercices}

\begin{exercice}[Formule de balayage]
\label{exo:balayage}
Soit $(Y_t, t\geq 0)$ une semi-martingale continue et $(k_u,u\geq 0)$ un processus $(\F_u)$-prévisible borné ; on pose $g_t = \sup\{s \leq t : Y_s=0\}$ le dernier zéro de $Y$ avant l'instant $t$.
\begin{enumerate}
  \item Montrer que $k_{g_t} Y_t = \int_0^t k_{g_s} \d Y_s$.
  \item Calculer le temps local en 0 de $k_{g_t}|B_t|$, où $B$ est un mouvement brownien.
  \item Soit $T$ un temps d'arrêt du mouvement brownien tel que $\P(B_T=0)=0$ et $B_T$ uniformément intégrable, montrer qu'il existe un unique processus prévisible croissant $A$ tel que pour tout processus prévisible positif $h$ on a :
\[ \E(h_{g_T}) = \E\left( \int_0^{+ \infty} h_u \d A_u \right).\]
  \item Montrer que $A_T$ est une variable aléatoire exponentielle de paramètre 1. En particulier, montrer que si $\tilde{T}_a = \inf\{t \geq 0 : |B_t|=a\}$, alors $L_{\tilde{T}_a}$ est une variable aléatoire exponentielle de paramètre $a$.
\end{enumerate}
\end{exercice}

\begin{proof}
\textit{1.} On pose $d_t = \inf\{s \geq t : Y_s = 0\}$ le premier zéro de $Y$ après l'instant $t$. C'est temps d'arrêt pour la filtration $(\F_u)$. De plus, pour tout $(a,b) \in {\R^+}^2$
\[ \ind{a \leq g_t < b}  = \ind{d_a \leq t < d_b}.\]

Considérons $(k_u,u\geq 0)$ processus prévisible simple : $k = H\ind{(a,b]},$ où $H$ est une variable aléatoire $\F_a$-mesurable bornée. On a alors :
\[ \int_0^t k_{g_s} \d Y_s = H (Y_{t\wedge d_b}-Y_{t \wedge d_a}).\]

De plus, on sait que pour tout $t \geq 0, Y_{d_t}=0$ ; par conséquent $k_{g_t} Y_t = \int_0^t k_{g_s} \d Y_s$. Par linéarité, et en utilisant le lemme des classes monotones, cette identité est donc vérifiée pour tout processus prévisible borné $(k_u, u \geq 0)$.

\textit{2.} Appliquons maintenant la formule de Tanaka à l'équation :
\[|k_{g_t}||B_t| = \int_0^t |k_{g_s}| \d|B_s|,\]
on a
\begin{align*}
  |k_{g_t} |B_t|| = & \int_0^t |k_{g_s}| \sgn(B_s) \d B_s + \int_0^t |k_{g_s}| \d L_s\\
  = & \int_0^t |k_{g_s}| \sgn(B_s) \d B_s + \int_0^t |k_s| \d L_s,
\end{align*}
où on a utilisé que $L$ ne croit que sur l'ensemble des zéros de $B$. Par conséquent, le temps local en 0 de $(k_{g_t}B_t, t \geq 0)$ est $\int_0^t |k_s|\d L_s$.

\textit{3.} Soit $h$ un processus prévisible positif, grâce au calcul précédent, on a, par théorème d'arrêt : 
\[\E[ |h_{g_T} |B_T|| ] = \E\left[\int_0^T |h_s| \d L_s\right], \]
or, on peut également écrire, en conditionnant par rapport à $\F_{g_t}$ la tribu engendrée par les variables aléatoires $k_{g_t}$, pour tout $k$ processus prévisible :
\[\E[ h_{g_T} |B_T| ] = \E[ \E[B_T|\F_{g_T}] h_{g_T}],\]

Posons alors $\xi_t = \E[|B_t||\F_{g_t}]$, qui est également un processus prévisible, on a :
\begin{align*}
  \E(h_{g_T}) = & \E\left( \frac{h_{g_T}}{\xi_T} |B_T|\right)\\
  = & \E\left[\int_0^t \frac{h_s}{\xi_s} \d L_s\right].
\end{align*}
Par conséquent, $A_t = \int_0^{t \wedge T} \frac{\d L_s}{\xi_s}$ est un processus croissant correspondant aux attentes. L'unicité de ce processus croissant s'obtient sans difficulté.

\textit{4.} Afin de calculer la loi de $A_T$, on calcule sa transformée de Laplace. Pour cela on utilise le fait que $A_T=A_{g_T}$. On a donc :
\begin{align*}
  \E(\exp(-\lambda A_T))
  = & \E(\exp(-\lambda A_{g_T})) = \E\left[ \int_0^T \exp(-\lambda A_s)\d A_s \right] \\
  = & \E \left[ \frac{1-\exp(-\lambda A_T)}{\lambda} \right].
\end{align*}
On en déduis $\E(\exp(-\lambda A_{g_T})) = \frac{1}{1+\lambda}$, donc $A_T$ est une variable aléatoire exponentielle de paramètre 1.

Dans le cas particulier où $T=\tilde{T}_a$, on a $\E(|B_{\tilde{T}_a}||\F_s)=a$, donc $A_t = \frac{L_t}{a}$, $L_{\tilde{T}_a}$ est bien une variable aléatoire exponentielle de paramètre $a$.
\end{proof}

\begin{remarque}
Grâce à la formule de balayage, on peut retrouver de façon simple le résultat de l'Exercice \ref{exo_mart} :
\[\left(f(L_t)B_t - F(L_t),t \geq 0 \right) \text{  est une martingale.}\]

En effet, on a $F(L_t)= \int_0^{L_t} f(l)\d l = \int_0^t f(L_u)\d L_u$, or comme $f$ est supposée positive, $\int_0^t f(L_u)\d L_u$ est le temps local en 0 de $f(L_t)|B_t|$. Par conséquent, le processus étudié est bien une martingale.
\end{remarque}

\begin{exercice}[Équation de Tanaka]
\label{exo_tanaka}
Soit $(\F_t)_{t\geq 0}$ une filtration et $B$ un $(\F_t)$-mouvement brownien, on étudie la structure des solutions de l'équation différentielle stochastique de Tanaka :
\[X_t = \int_0^t \sgn(X_s) \d B_s.\]

\begin{enumerate}
 \item Montrer que si $X$ existe, alors à $t$ fixé, $\sgn(X_t)$ est une variable aléatoire de Bernoulli symétrique indépendante de $B$.
 
 \item En déduire que si $\F_t = \sigma(B_s,s\leq t)$, cette équation différentielle stochastique n'admet pas de solution.
 
 \item Montrer que s'il existe une solution $X^{(0)}$ à cette équation différentielle, alors en posant $g^{(0)}_t = \sup\{s \leq t : X^{(0)}_s =0\}$, pour tout processus prévisible $k$ à valeurs dans $\{-1,1\}$ le processus $X^{(k)}_t = k_{g^{(0)}_t}X^{(0)}_t$ est une autre solution.
 
 \item Montrer que toutes les solutions de l'équation, définies sur un même espace de probabilité avec un même mouvement brownien peuvent être représentées à partir de $X^{(0)}$ sous la forme $X^{(k)}$, pour $k$ processus prévisible bien choisi.
\end{enumerate}
\end{exercice}

\begin{remarque}
Classiquement, concernant l'équation de Tanaka, on dit qu'il y a unicité en loi, mais pas d'unicité trajectorielle, car si $X$ est solution, alors $-X$ est une autre solution définie sur le même espace de probabilité. Malgré tout, ces deux solutions sont des mouvements browniens. Ici nous allons expliciter toutes les solutions de cette équation à partir de l'une d'entre elles.
\end{remarque}

\begin{proof}
\textit{1.} Soit $X$ un processus vérifiant l'équation différentielle stochastique de Tanaka. On a $\crochet{X}_t=t$, donc, grâce au théorème de Lévy, $X$ est un $(\F_t)$-mouvement brownien. Il est par conséquent immédiat d'observer que $\sgn(X_t)$ est une variable aléatoire de Bernoulli symétrique.

Afin de vérifier que $\sgn(X_t)$ est indépendant de $B$, il suffit d'observer que $-X$ est une autre solution de l'équation de Tanaka, construite sur le même espace de probabilité avec le même mouvement brownien. On a donc
\[\P(\sgn(X_t)=1|B) = \P(\sgn(-X_t)=1|B) = \P(\sgn(X_t)=-1|B)=\frac{1}{2}.\]
On a donc bien l'indépendance.

\textit{2.} Comme $\sgn(X_t)$ est indépendant de $B$, le processus $X$ ne peut être mesurable par rapport à la tribu canonique du mouvement brownien. Par conséquent il n'existe pas de solution adaptée à la filtration de $B$.

\textit{3.} On utilise maintenant la formule de balayage de l'Exercice \ref{exo:balayage}. Soit $X$ une solution de l'équation différentielle stochastique et $k$ un processus prévisible à valeurs dans $\{-1,1\}$. On a :
\[X^{(k)}_t = k_{g^{(0)}_t}X^{(0)}_t = \int_0^t k_{g^{(0)}_s}\d X^{(0)}_s = \int_0^t \sgn(X^{(k)}_s) \d B_s.\]
Par conséquent $X^{(k)}$ est bien une solution de cette équation différentielle stochastique.

\textit{4.} On pose maintenant $Y$ une solution quelconque de l'équation différentielle stochastique de Tanaka associée au mouvement brownien $B$. Observons pour commencer que l'on a :
\[|Y_t| = \int_0^t \sgn{Y_s}\d Y_s + L^0_t(Y) = B_t + L^0_t(Y).\]

Par conséquent, $(|Y|,L^0(Y))$ est la seule solution de l'équation de Skorokhod associée au mouvement brownien $B$. On obtient par conséquent, pour tout $t \geq 0$ :
\[(|B_t|,L^0_t(B))=(|Y_t|,L^0_t(Y)),\]
ces quantités sont donc les mêmes pour toutes les solutions de l'équation différentielle stochastique construites sur un même espace de probabilités et associées au même mouvement brownien. 

\begin{remarque}
Notons en particulier que la filtration de $B$ est égale à celle de $|Y|$.
\end{remarque}

Les zéros de $Y$ et $X^{(0)}$ sont donc en particulier les mêmes, et le rapport $\frac{Y_t}{X^{(0)}_t}$ reste constant égal à $\pm 1$ sur tout intervalle du type $(g^{(0)}_t,d^{(0)}_t)$. On note $k_u=\frac{Y_u}{X^{(0)}_u}$ pour tout $u$ tel que $Y_u=0$, que l'on prolonge par continuité à droite.

On observe enfin que pour presque tout $t\geq 0$, $B_t \neq 0$, donc on a $Y_t=X^{(k)}_t$, ce qui permet de conclure.
\end{proof}

\begin{remarque}
Il existe bien des solutions faibles à l'équation de Tanaka. En effet, prenons un mouvement brownien $\beta$ muni de sa filtration canonique $(\F_t)$, on pose :
\[B_t = \int_0^t \sgn(\beta_s) \d \beta_s,\]
par théorème de Lévy, $B$ est bien un mouvement brownien. De plus, il est immédiat de constater que $\beta$ est une solution de l'équation :
\[X_t = \int_0^t \sgn(X_s)\d B_s.\]
Dans ce cas, on vérifie bien que la filtration canonique associée à $B$ est strictement plus petite que celle associée à $\beta$.
\end{remarque}

\begin{exercice}
Soit $B$ un mouvement brownien.

\begin{enumerate}
 \item Soit $\alpha>1$, on pose $X_t=|B_t|^\alpha$. Montrer que $L^0_t(X)=0$.
 \item Soit $\phi$ une fonction meusrable bijective. Sous quelle condition suffisante sur $\phi$, le processus $(\phi(B_t))_{t \geq 0}$ est-il une semi-martingale de temps local en 0 nul ?
\end{enumerate}
\end{exercice}

\begin{proof}
\textit{1.} On commence par calculer le crochet de $X$, et on cherche une fonction $\rho$ strictement croissante telle que $\frac{1}{\rho}$ soit non-intégrable au voisinage de 0 et vérifie : 
\[\exists \epsilon >0 \, : \, \int_0^t \dfrac{\d \crochet{X}_s}{\rho(X_s)}\ind{0<X_s<\epsilon} <+\infty.\]
Pour cela, on utilise la formule d'Itô-Tanaka à la fonction convexe $x \mapsto |x|^\alpha$ est convexe. On a
\[X_t = \alpha \int_0^t \sgn(B_s)|B_s|^{\alpha-1} \d B_s + \dfrac{\alpha(\alpha-1)}{2} \int_\R L^x_t |x|^{\alpha-2} \d x.\]

On en déduit que le crochet de $X$ vaut
\[\crochet{X}_t= \alpha^2 \int_0^t |B_s|^{2\alpha-2}\d s.\]
On pose alors $\rho(x) = x$, dont l'inverse est non-intégrable au voisinage de 0. On a :
\begin{align*}
  \int_0^t \dfrac{\d \crochet{X}_s}{\rho(X_s)}\ind{0<X_s<1}
  = & \alpha^2 \int_0^t |B_s|^{\alpha-2} \ind{0<|B_s|<1}\d s\\
  = & \alpha^2 \int_{[-1,1]} L^x_t(B) |x|^{\alpha-2}\d x <+\infty,
\end{align*}
grâce à l'expression du temps local du mouvement brownien en fonction de la densité de la loi gaussienne, et l'intégrabilité de $x \mapsto |x|^{\alpha-2}$ au voisinage de 0.

On en tire immédiatement que le temps local en 0 de $X$ est nul, ce qui est contre-intuitif. En effet, $X_t$ et $B_t$ possèdent exactement les mêmes zéros, mais au voisinage de 0, les oscillations sont \og écrasées \fg par la puissance.

\textit{2.} On cherche maintenant une condition portant sur les fonctions $\phi$ telles que $\phi(B_t)$ est une semi-martingale de temps local en 0 nul. Par Itô-Tanaka, pour que $\phi(B_t)$ soit une semi-martingale, il suffit que $\phi$ soit localement différence de deux fonctions convexes. Par conséquent, on considère désormais que $\phi$ est continue et admet une dérivée à gauche et à droite en tout point. Sans perdre de généralités, on peut supposer que $\phi$ est croissante et ne s'annule qu'en zéro.

On cherche maintenant des conditions suffisantes sur $\phi$ pour que le temps local en 0 de $\phi(B_t)$ soit nul. On a, en prenant $\rho(x) = x$ :
\begin{align*}
  \int_0^t \dfrac{\d \crochet{X}_s}{\rho(X_s)}\ind{0<X_s<\epsilon}
  = & \int_0^t \dfrac{\phi'(B_s)^2}{\phi(B_s)} \ind{0<\phi(B_s)<\epsilon} \d s \\
  = & \int_\R \dfrac{\phi'(y)^2}{\phi(y)}L^y_t \ind{0\leq \phi(y) \leq \epsilon}\d y.
\end{align*}

Par changement de variables, la formule devient donc :
\[ \int_0^{\phi^{-1}(\epsilon)} \dfrac{ \phi'(y)^2}{\phi(y)}\d y = \int_0^{\epsilon} \dfrac{1}{(\phi^{-1})'(u)u}\d u<+\infty.\]
La fonction $\phi$ est positive à droite de $0$, donc si $\phi'_d(0) \neq 0$, la condition suffisante n'est pas obtenue, car à droite de 0, $\phi$ est équivalent à $\lambda y$, qui n'est pas intégrable. Il faut par conséquent imposer $\phi'_d(0)=0$, et qu'il existe $\epsilon >0$ tel que $|\phi'_d(y)| = O(|y|^\epsilon)$. Il faut également imposer le même type de condition sur $\phi'_g$.

Par conséquent, pour que $\phi(B)$ soit une semi-martingale dont le temps local en $0$ est nul, il suffit d'imposer que $\phi'$ soit à croissance au plus de type puissance au voisinage de $0$. Dans le cas où $\phi$ est de classe $\mathcal{C}^2$, on demande $\phi'(x) = 0$ dès lors que $\phi(x)=0$.
\end{proof}

\begin{remarque}
Nous développerons plus tard une théorie des temps locaux de diffusions pour lesquelles le temps local en 0 de l'image de $B$ par une fonction $\phi$ bijective est donné par $L^0_t(\phi(B)) = c L^0_t(B)$.
\end{remarque}

\begin{exercice}[Autour du résultat de Zvonkin]
On étudie le processus solution de l'équation différentielle stochastique :
\[X_t = x + B_t - \lambda \int_0^t  \sgn(X_s) \d s,\]
appelé processus bang-bang de paramètre $\lambda$.

Déterminer le semi-groupe de ce processus de Markov.
\end{exercice}

\begin{proof}
Soit $B$ un mouvement brownien, et $f \in \mathcal{C}_b$, on calcule
\[\E_x(f(X_t)) = \E\left( f\left(x + B_t - \lambda \int_0^t \sgn(X_s) \d s\right) \right).\]
Grâce à la formule de Girsanov, on a
\[\E_x(f(X_t)) =  \E\left(f(x + B_t) \exp\left( -\lambda \int_0^t \sgn(B_s)\d B_s -\dfrac{\lambda^2}{2} t \right)\right).\]

\begin{remarque}
\label{rem_girsanov}
C'est un fait qui peut être généralisé, soit $b$ une fonction mesurable bornée et $Y$ un processus vérifiant, pour tout $t \geq 0$, l'équation
\[Y_t = y + B_t + \int_0^t b(Y_s) \d s,\]
par formule de Girsanov :
\[\E_y(f(Y_t)) = \E\left(f(y + B_t) \exp\left( -\int_0^t b(B_s)\d B_s -\dfrac{1}{2} \int_0^t b^2(B_s)\d s \right)\right).\]
\end{remarque}

On utilise ensuite la formule de Tanaka pour exprimer $\int_0^t \sgn(B_s)\d B_s = |B_t|-L_t$. De plus, conditionnellement à $(|B_t|, L_t)$ on a $\sgn(B_t)=1$ avec probabilité $\frac{1}{2}$. On pose $g(u) = \frac{f(x+u)+f(x-u)}{2}$, par conditionnement, l'espérance est égale à :
\[ \E_x(f(X_t)) =\E\left[g(|B_t|)\exp\left( -\lambda (|B_t|-L_t) -\dfrac{\lambda^2}{2} t \right)\right].\]

Il suffit enfin d'utiliser la formule donnant la loi jointe de $(L_t,|B_t|)$. On obtient :
\begin{align*}
\int_{\R^+} \d u \int_{\R^+} &  \d v g(u) \exp(\lambda (v-u) - \frac{\lambda^2}{2} t) \dfrac{2(u+v)}{\sqrt{2\pi t^3}}e^\frac{-(u+v)^2}{2t}\\
 = &  2\int_{\R^+} \d u \dfrac{g(u)}{\sqrt{2\pi t^3}} e^{-2\lambda u- \frac{\lambda^2}{2} t} \int_{\R^+} \d v(u+v)e^{\lambda (u+v)} e^\frac{-(u+v)^2}{2t}\\
 = & 2\int_{\R^+} \d u \dfrac{g(u)}{\sqrt{2\pi t}} e^{-2\lambda u- \frac{\lambda^2}{2} t} \int_u^{+\infty} \d w \dfrac{w}{t}e^{\lambda w} e^\frac{-w^2}{2t}.
\end{align*}

Or, par intégration par partie,
\[\int_u^{+\infty} \d w \dfrac{we^{\lambda w}}{t}e^{\frac{-w^2}{2t}} = e^\frac{-u^2}{2t}e^{\lambda u} + \lambda e^{\frac{\lambda^2}{2}t} \int_{u-\lambda t}^{+\infty}\d w e^\frac{-w^2}{2t}.\]
Par conséquent,
\[\E_x(f(X_t)) = \int_{\R^+}\d u \dfrac{f(x+u)+f(x-u)}{\sqrt{2\pi t}} \left[ e^{\frac{-(u+\lambda t)^2}{2t}} + \lambda e^{-2\lambda u} \left(1-\Phi\left(\frac{u-\lambda t}{\sqrt{t}}\right)\right) \right],\]
où $\Phi$ est la fonction de répartition de la gaussienne centrée réduite. On en déduit que le semi-groupe du processus $X$ est :
\[P_t(x,\d y) =  \dfrac{\d y}{\sqrt{2\pi t}}\left[ \exp\left(\frac{-(|y-x|+\lambda t)^2}{2t} \right) + \lambda e^{-2\lambda |y-x|} \left(1-\Phi\left(\frac{|y-x|-\lambda t}{\sqrt{t}}\right)\right) \right].\]
\end{proof}

\begin{exercice}[Équation de Tsirel'son]
Posons, pour $x \in \R$, $\{x\}=x - \floor{x}$ la partie fractionnaire de $x$. On étudie la solution (unique en loi) de l'équation différentielle stochastique
\[X_t = B_t + \int_0^t T(X)_s \d s,\]
où on a posé :
\[T(X)_s = \sum_{k \in -\N} \left\{ \dfrac{X_{t_k}-X_{t_{k-1}}}{t_k-t_{k-1}} \right\} \ind{(t_k,t_{k+1}]}(s),\]
avec $(t_k)_{k \in -\N}$ une suite croissante de réels positifs avec $\lim_{k \to -\infty} t_k = 0$. Cette fonctionnelle est appelée dérive de Tsirel'son

\begin{enumerate}
 \item Montrer que $\forall k \in -\N, \eta_k = \left\{\frac{X_{t_k}-X_{t_{k-1}}}{t_k-t_{k-1}}\right\}$ suit une loi uniforme sur $[0,1]$.
 
 \item Montrer de plus que cette variable aléatoire est indépendante de $B$.
\end{enumerate}
\end{exercice}

\begin{proof}
\textit{1.} Remarquons pour commencer que pour tout $t \in [t_k,t_{k+1}]$, on a :
\[X_t = X_{t_k} + B_t-B_{t_k} + (t-t_k) \eta_{k}.\]
En particulier $\frac{X_{t_k}-X_{t_{k-1}}}{t_k-t_{k-1}} = \frac{B_{t_k}-B_{t_{k-1}}}{t_k-t_{k-1}} + \eta_{k-1}$, et ces deux dernières variables aléatoires sont indépendantes, car $B_{t_k}-B_{t_{k-1}}$ est indépendant de $\F_{t_{k-1}}$. Par conséquent
\begin{equation}
  \label{eqn:recurrence}
  \eta_k = \left\{\frac{B_{t_k}-B_{t_{k-1}}}{t_k-t_{k-1}} + \eta_{k-1}\right\}.
\end{equation}

On calcule la transformée de Fourier discrète de $\eta_k$ ; pour $p \in \Z\backslash\{0\}$, on a
\[ \E\left(e^{i2\pi p \eta_k}\right)  =  \E\left(\exp\left(2 i \pi p \left({\eta_{k-1}} +\frac{B_{t_k}-B_{t_{k-1}}}{t_k-t_{k-1}}\right)\right)\right),\]
car $e^{2i \pi n}=1$ pour tout $n \in  \N$. Par indépendance, on a alors
\begin{align*}
 \E(e^{2i\pi p \eta_k})
 = & \E(e^{2i\pi p \eta_{k-1}}) \E\left( \exp \left( 2 i \pi p \frac{B_{t_k}-B_{t_{k-1}}}{t_k-t_{k-1}}\right) \right)\\
 = & \E(e^{i2\pi p \eta_{k-1}})\exp\left(- \frac{2\pi^2 p^2}{t_k-t_{k-1}}\right).
\end{align*}

Or $\frac{1}{t_k-t_{k-1}}\geq \frac{1}{t_0}$ donc pour tout entier $n$,
\begin{align*}
  \left|\E\left(e^{i2\pi p \eta_k}\right)\right| \leq & \left|\E\left(e^{i2\pi p \eta_{k-1}}\right)\right| \exp\left( - \frac{2\pi^2 p^2}{t_0}\right)\\
  \leq & \left|\E\left(e^{i2\pi p \eta_{k-n}}\right)\right| \exp\left(- n\frac{2\pi^2 p^2}{t_0}\right)
  \leq  \exp\left(- n\frac{2\pi^2 p^2}{t_0}\right).
\end{align*}
Dès lors, pour tout $p \in \Z\backslash\{0\}$, $\E(e^{i2\pi p \eta_k}) = 0$. De plus, la variable aléatoire $\eta_k$ est à valeurs dans $[0,1[$. Par conséquent, $\eta_k$ est une variable aléatoire de loi uniforme sur $[0,1[$.

\textit{2.} Notons $\widehat{\mathcal{B}}_t=\sigma(B_{t+s}-B_t,s\geq 0)$, et calculons, pour $n \in -\N$, la quantité $\E(e^{i2\pi p \eta_k}|\widehat{\mathcal{B}}_{t_n})$. En utilisant la relation de récurrence \eqref{eqn:recurrence}, on observe que pour tout $n \in \N$, on a $\E(e^{i2\pi p \eta_k}|\widehat{\mathcal{B}}_{t_n}) = 0$. Grâce aux propriétés de continuité de la tribu brownienne, on en déduit $\E(e^{i2\pi p \eta_k}|\widehat{\mathcal{B}}_{0})=0$, donc $\eta_k$ est indépendant du mouvement brownien.

De la même façon que pour l'Exercice \ref{exo_tanaka}, on peut conclure que l'équation de Tsirel'son ne vérifie pas la propriété d'unicité trajectorielle. En fait, cette équation différentielle stochastique donne un exemple de dérive bornée $T$, qui dépende uniquement du passé de $X$, et qui est telle qu'on n'ait pas d'unicité trajectorielle.

L'existence et l'unicité en loi des solutions de cette équation différentielle stochastique se prouvent grâce à la formule de Girsanov (cf Remarque \ref{rem_girsanov}).
\end{proof}

\begin{remarque}
Si on considère l'équation différentielle obtenue en n'imposant pas de partie fractionnaire, i.e. 
\[X_t=X_0+B_t +\int_0^t \tilde{T}(X)_s\d s,\]
où on a posé :
\[\tilde{T}_k = \sum_{k \in -\N} \frac{X_{t_k}-X_{t_{k-1}}}{t_k-t_{k-1}} \ind{(t_k,t_{k+1}]}(s)\]
en posant $\tilde{\eta}_k=\frac{X_{t_k}-X_{t_{k-1}}}{t_k-t_{k-1}},$ on peut appliquer le même raisonnement que précédemment. La transformée de Fourier $\E(e^{i\lambda \tilde{\eta}_k})$ est nulle pour tout $\lambda \neq 0$, ce qui ne se peut. Par conséquent, l'équation de Tsirel'son modifiée dirigée par $\tilde{T}$ n'admet pas de solution.
\end{remarque}

\chapter[Lois de l'arcsinus et scaling aléatoire]{Lois de l'arcsinus et scaling aléatoire}

Les développements de la théorie des temps locaux ont permis de découvrir un certain nombre d'identités en lois relatives aux fonctionnelles du mouvement brownien. Parmi celles-ci, on compte de nombreux développements autour de la loi de l'arcsinus. En particulier, on doit à Paul Lévy le résultat suivant : le dernier zéro avant l'instant 1, et le temps passé au dessus de l'axe des abscisses avant l'instant 1 suivent une même loi appelée loi de l'arcsinus. De façon générale, lorsqu'on s'intéresse à la loi jointe de plusieurs fonctionnelles browniennes, les temps locaux apparaissent de façon naturelle, et s'avèrent être des éléments importants pour la compréhension de ces lois jointes.

Dans un second temps, nous nous intéresserons à la notion de scaling aléatoire, une relation déjà entrevue dans l'Exercice \ref{exo_scalingaleatoire}, qui permet souvent d'obtenir la dérivée de Radon-Nikod\'ym entre deux processus construits à partir du mouvement brownien, et conditionnés d'une certaine façon que nous détaillerons. Pour finir et revenir à la loi de l'arcsinus, nous nous attarderons sur la décomposition du mouvement brownien en deux mouvements browniens réfléchis indépendants changés de temps.

\section{Théorèmes principaux}

\begin{definition}
La loi de l'arcsinus est la loi de probabilité sur $[0,1]$, qui admet comme densité par rappport à la mesure de Lebesgue $\dfrac{1}{\pi\sqrt{x(1-x)}}\ind{0 < x < 1}$. C'est un cas particulier de loi $\beta$, la loi $\beta_{\frac{1}{2},\frac{1}{2}}$.
\end{definition}

\begin{theoreme}[Lois de l'arcsinus dans le cadre brownien ; Paul Lévy]
On note $g_t=\sup\{s \leq t : B_s=0\}$ le dernier temps de passage en zéro avant l'instant $t$ d'un mouvement brownien, et $A^{\scriptscriptstyle (+)}_t=\int_0^t \ind{B_s>0} \d s$ le temps durant lequel le mouvement brownien est resté positif avant l'instant $t$.

Les variables aléatoires $\frac{g_t}{t}$ et $\frac{A^{\scriptscriptstyle (+)}_t}{t}$ vérifient toutes deux la loi de l'arcsinus.
\end{theoreme}

La loi de l'arcsinus favorise surtout les valeurs proches de 0 et de 1, en d'autres termes le mouvement brownien réalise un choix. Soit il s'éloigne relativement beaucoup de 0, auquel cas $g_t$ est petit, soit il reste proche de $0$ et $g_t$ est très proche de $t$. Le m\^eme type d'interprétation s'applique très bien à $A^{\scriptscriptstyle (+)}_t$

Une des nombreuses preuves du Théorème 5.1 pour la loi de $A^{\scriptscriptstyle (+)}_t$ s'appuie sur la proposition suivante.

\begin{proposition}
Soit $\alpha^{\scriptscriptstyle (+)}_u = \inf\{t \geq 0 : A^{\scriptscriptstyle (+)}_t > u \}$ l'inverse continu à droite de $t \mapsto A^{\scriptscriptstyle (+)}_t$ et $(B^+_t,t \geq 0)$ la partie positive du mouvement brownien. On utilisera également, avec des notations évidentes, $A^{\scriptscriptstyle (-)}_t$, $B^-$ et $\alpha^{\scriptscriptstyle (-)}_u$.

Les processus $(B^+_{\alpha^{\scriptscriptstyle (+)}_u},u \geq 0)$ et $(B^-_{\alpha^{\scriptscriptstyle (-)}_u},u \geq 0)$ sont deux mouvements browniens réfléchis indépendants.
\end{proposition}

Les processus $(A^{\scriptscriptstyle (+)}_t, t \geq 0)$ et $(g_t,t \geq 0)$ jouissent à l'évidence d'une propriété de scaling, héritée de celle du mouvement brownien. Il est donc utile d'énoncer --et d'utiliser-- le théorème suivant.

\begin{theoreme}[Scaling aléatoire]
Soit $(A_t,t\geq 0)$ un processus continu croissant tel qu'il existe $r >0$ vérifiant :
\[\forall c >0, \left(A_{ct},B_{ct},t\geq 0\right) \egaldistr \left(c^{r+1} A_t, \sqrt{c}B_t,t\geq 0\right).\]

On pose $\alpha_u = \inf\{t \geq 0  : A_t >u\}$, pour toute fonction mesurable bornée $F$ sur $\mathcal{C}([0,1])$, pour toute fonction $\phi : \R^+ \to \R^+$, on a :
\begin{multline*}
  \E\left[\int_0^{+\infty} F\left( \frac{1}{\sqrt{t}}B_{ut} ,u \leq 1\right) \phi(t) \d A_t \right]\\
  = \E\left[\frac{r+1}{\alpha_1^{r+1}} F\left( \frac{1}{\sqrt{\alpha_1}}B_{\alpha_1u} , u \leq 1 \right) \right] \int_0^{+\infty} \phi(t) t^r \d t.
\end{multline*}
\end{theoreme}

Lorsque le processus croissant continu admet une dérivée, on peut appliquer des techniques de preuves similaires à celle utilisée dans l'Exercice \ref{exo_scalingaleatoire}, pour obtenir le résultat suivant.

\begin{corollaire}
\label{cor_radonnikodym}
Si $A_t = \int_0^t \theta_s \d s$ avec $\theta$ processus mesurable vérifiant 
\[\forall c >0, \left(\theta_{ct},B_{ct};t\geq 0\right) \egaldistr \left(c^{r} \theta_t, \sqrt{c}B_t; t\geq 0\right),\]
alors
\[\E(\theta_1 F(B_v,v\leq 1)) = \E\left( \dfrac{r+1}{\alpha_1^{r+1}}F\left( \frac{1}{\sqrt{\alpha_1}} B_{\alpha_1u} ,u \leq 1 \right) \right).\]
\end{corollaire}

Voici maintenant deux exemples d'application de ce corollaire ; tout d'abord une relation entre le mouvement brownien arrêté en $\tau_1$ et le pont brownien.
\begin{proposition}
\label{pro_pontbrownien}
Soit $(b_u,u\leq 1)$ un pont brownien, et $l$ son temps local en 0. Pour toute fonction $F : \mathcal{C}([0,1]) \to \R^+$, on a :
\[\E(F(b_u,u\leq 1)) = \E\left( \sqrt{\dfrac{\pi}{2\tau_1}}F\left( \frac{1}{\sqrt{\tau_1}} B_{\tau_1u},u \leq 1\right) \right) \mathrm{  et }\]
\[\E\left(\sqrt{\dfrac{2}{\pi}} \dfrac{1}{l_1}F(b_u,u\leq 1)\right) = \E\left( F\left( \frac{1}{\sqrt{\tau_1}} B_{\tau_1u}, u \leq 1 \right) \right).\]
\end{proposition}

On s'intéresse maintenant au processus associé au temps passé par le mouvement brownien au dessus de 0.
\begin{proposition}
\label{pro_pontbizarre}
Pour toute fonction $F : \mathcal{C}([0,1]) \to \R^+$, on a :
\[\E(\ind{B_1>0}F(B_u,u\leq 1)) = \E\left( \dfrac{1}{\alpha^{\scriptscriptstyle (+)}_1} F\left( \frac{1}{\sqrt{\alpha^{\scriptscriptstyle (+)}_1}} B_{\alpha^{\scriptscriptstyle (+)}_1u}, u \leq 1 \right) \right) \mathrm{  et }\]
\[\E\left( \dfrac{\ind{B_1>0}}{A^{\scriptscriptstyle (+)}_1}F(B_u,u\leq 1)\right) = \E\left( F\left( \frac{1}{\sqrt{\alpha^{\scriptscriptstyle (+)}_1}} B_{\alpha^{\scriptscriptstyle (+)}_1u},u \leq 1\right)\right).\]
\end{proposition}

On s'intéresse maintenant à d'autres types de résultats, précisant la loi jointe de plusieurs variables aléatoires apparaissant dans la décomposition du mouvement brownien en partie positive et partie négative.
\begin{theoreme}
Soit $Z_t = \frac{1}{t}(A^{\scriptscriptstyle (+) }_t,A^{\scriptscriptstyle (-) }_t,(L_t)^2)$, et $T$ un temps aléatoire.

Pour $T = t \in \R$, ou $T = \tau_l$ ou $T=\alpha^{\scriptscriptstyle (+) }_s$, $Z_T$ obéit à la même loi. 
\end{theoreme}

\begin{theoreme}
\label{thm_identiteenloi}
Soit $\tilde{T}_1 = \inf\{t>0 : |B_t|=1\}$ et $\widehat{T}_1$ une copie indépendante de $T_1$.

On a $(A^{\scriptscriptstyle (+) }_{T_1},A^{\scriptscriptstyle (-) }_{T_1},\frac{1}{2}L_{T_1}) \egaldistr (\tilde{T}_1, (L_{\tilde{T}_1})^2\widehat{T}_1, L_{\tilde{T}_1})$.

La transformée de Laplace de $(A^{\scriptscriptstyle (+) }_{T_1},A^{\scriptscriptstyle (-) }_{T_1},\frac{1}{2}L_{T_1})$ est donnée par :
\[\E\left[\exp\left(-\left(\dfrac{\lambda^2}{2}A^{\scriptscriptstyle (+) }_{T_1}) + \dfrac{\mu^2}{2}A^{\scriptscriptstyle (-) }_{T_1} + \alpha L_{T_1}\right) \right) \right] = \left[ \ch(\lambda) + \dfrac{\mu +2\alpha}{\lambda}\sh(\lambda) \right]^{-1}.\]
\end{theoreme}

\section{Exercices}

\begin{exercice}
Montrer que $\P(B_1>0|A^{\scriptscriptstyle (+) }_1)=A^{\scriptscriptstyle (+) }_1$.
\end{exercice}

\begin{proof}
On applique la Proposition \ref{pro_pontbizarre}, pour toute fonction $\Phi$ mesurable positive
\[\E(\ind{B_1>0} \Phi(B_u,u\leq 1)) = \E\left( \dfrac{1}{\alpha^{\scriptscriptstyle (+)}_1} \Phi\left( \frac{1}{\sqrt{\alpha^{\scriptscriptstyle (+)}_1}} B_{\alpha^{\scriptscriptstyle (+)}_1u}, u \leq 1 \right) \right).\]
En particulier, si $\Phi(X) = \int_0^1 \ind{X_s > 0} \d s$, on a $\Phi(B)=A^{\scriptscriptstyle (+)}_1$, et
\begin{align*}
  \Phi\left( \frac{1}{\sqrt{\alpha^{\scriptscriptstyle (+)}_1}} B_{\alpha^{\scriptscriptstyle (+)}_1u}, u \leq 1 \right)
  &=\int_0^1 \ind{B_{u\alpha^{\scriptscriptstyle (+) }_1}>0}\d u
   = \frac{1}{\alpha^{\scriptscriptstyle (+) }_1} \int_0^{\alpha^{\scriptscriptstyle (+)}_1} \ind{B_r>0} \d r = \frac{1}{\alpha^{\scriptscriptstyle (+) }_1}.
\end{align*}

Par conséquent, pour toute fonction $f$ mesurable positive, on a
\begin{equation}
  \label{eqn:pontbizarre} \E (\ind{B_1>0} f(A^{\scriptscriptstyle (+) }_1)) =  \E\left( \dfrac{1}{\alpha^{\scriptscriptstyle (+) }_1} f\left(\dfrac{1}{\alpha^{\scriptscriptstyle (+) }_1}\right)\right).
\end{equation}
On calcule maintenant la loi de $\frac{1}{\alpha^{\scriptscriptstyle (+) }_1}$. Par propriété de scaling brownienne, pour tout $c>0$, on a $A^{\scriptscriptstyle (+) }_{c}\egaldistr cA^{\scriptscriptstyle (+) }_1$. Dès lors, pour tout $u>0$,
\begin{align*}
  \P\left(\frac{1}{\alpha^{\scriptscriptstyle (+) }_1}>u\right)
  = \P\left(\alpha^{\scriptscriptstyle (+) }_1<\frac{1}{u}\right)
  = \P(A^{\scriptscriptstyle (+) }_\frac{1}{u}>1)
  = \P(A^{\scriptscriptstyle (+) }_1>u),
\end{align*}
donc $\frac{1}{\alpha^{\scriptscriptstyle (+) }_1} \egaldistr A^{\scriptscriptstyle (+) }_1$ suit la loi de l'arcsinus.

On en conclut que pour toute fonction $f$ mesurable bornée
\[\E(\ind{B_1>0}f(A^{\scriptscriptstyle (+) }_1)) = \E(A^{\scriptscriptstyle (+) }_1 f(A^{\scriptscriptstyle (+) }_1)),\]
d'où on tire $\P(B_1>0|A^{\scriptscriptstyle (+) }_1) = A^{\scriptscriptstyle (+) }_1$.
\end{proof}

\begin{exercice}
Calculer la loi de $l_1$ le temps local en 0 du pont brownien.
\end{exercice}

\begin{proof}
On utilise la Proposition \ref{pro_pontbrownien}, afin de calculer $\E(f(l_1))$ pour toute fonction fonction mesurable positive $f$. Pour cela on utilise le fait suivant :
\[l_1 = \lim_{\epsilon \to 0}\dfrac{1}{\epsilon}\int_0^1 \ind{0<b_s<\epsilon} \d s,\]
donc $l_1$ est une fonctionnelle mesurable du pont brownien. On calcule alors la valeur de cette fonctionnelle appliquée au mouvement brownien stoppé en $\tau_1$
\[\lim_{\epsilon \to 0} \dfrac{1}{\epsilon}\int_0^1 \ind{0<B_{u\tau_1}<\epsilon \sqrt{\tau_1}} \d u = \dfrac{1}{\sqrt{\tau_1}} L_{\tau_1} = \frac{1}{\sqrt{\tau_1}}.\]

On obtient donc $\E(f(l_1)) = \E\left[\sqrt{\frac{\pi}{2\tau_1}} f\left(\frac{1}{\sqrt{\tau_1}}\right)\right]$.
En appliquant le Théorème d'équivalence de Lévy, la loi de $\tau_1$ est égale à celle de $\frac{1}{N^2}$ où $N$ est une variable aléatoire gaussienne centrée réduite. Par conséquent
\[\E(f(l_1)) = \E(|N|f(|N|)) = \int_0^{+\infty} xe^{-\frac{x^2}{2}} f(x)\d x.\]

\begin{remarque}
En particulier, $l_1 \egaldistr \sqrt{2\e}$ avec $\e$ variable aléatoire exponentielle de paramètre 1. Cette loi, appelée loi de Rayleigh, est la loi de la norme d'un vecteur gaussien bidimensionnel.
\end{remarque}
\end{proof}

\begin{exercice}
\label{exo_transformeelaplace}
Montrer, en étudiant $f$ telle que $f(|B_t|,L_t)\exp\left(-\frac{\lambda^2}{2}t - \nu L_t \right)$ est une martingale que l'on a :
\[\E\left[\exp\left(-\left(\dfrac{\lambda^2}{2}A^{\scriptscriptstyle (+) }_{T_1}) + \dfrac{\mu^2}{2}A^{\scriptscriptstyle (-) }_{T_1} + \alpha L_{T_1}\right) \right) \right] = \left[ \ch(\lambda) + \dfrac{\mu +2\alpha}{\lambda}\sh(\lambda) \right]^{-1}.\]
\end{exercice}

\begin{proof}
Soit $f$ une fonction de classe $\mathcal{C}^2$ sur $\R^3$, en appliquant la formule d'It\^o, on a
\begin{multline*}
 f(t,|B_t|,L_t) = f(0,0,0) + \int_0^t \partial_2 f(s,|B_s|,L_s) \sgn(B_s)\d B_s\\
  + \int_0^t \partial_1 f(s,|B_s|L_s) \d s + \dfrac{1}{2} \int_0^t \partial^2_{2} f(s,|B_s|,L_s) \d s\\
 + \int_0^t \partial_2 f(s,|B_s|,L_s) \d L_s + \int_0^t \partial_3 f(s,|B_s|,L_s)\d L_s.
\end{multline*}
Par conséquent, ce processus est une martingale locale dès lors que $f$ vérifie :
\[\left\{
\begin{array}{l}
  \dfrac{1}{2} \partial^2_{2} f(t,x,l) + \partial_1 f(t,x,l) = 0 \\
  \partial_2 f(t,0,l) +  \partial_3 f(t,0,l) = 0.
\end{array}
\right.\]

On cherche $f$ sous la forme $\phi(x)\psi(l)\exp(-\frac{\lambda^2}{2}t)$. La fonction $\phi$ s'obtient grâce à la première équation différentielle : $\phi(x) = A\ch(\lambda x) + B \sh(\lambda x)$. La fonction $\psi$ s'en déduit alors, $\psi(y) = C\exp(-\lambda B y)$. Le processus suivant est donc une martingale locale
\[\left( \ch(\lambda |B_t|) + \dfrac{\nu}{\lambda}\sh(\lambda |B_t|) \right) \exp \left( - \dfrac{\lambda^2}{2} t - \nu L_t \right). \]

On peut appliquer le théorème d'arrêt en $\tilde{T}_a$, car la martingale locale arrêtée est bornée ; on obtient
\[\E\left[ \left( \ch(\lambda a) + \dfrac{\nu}{\lambda} \sh(\lambda a) \right) \exp \left( - \dfrac{\lambda^2}{2} \tilde{T}_a - \nu L_{\tilde{T}_a}\right)\right] = 1. \]
En particulier, pour $a=1$, on a
\[\E\left[\exp\left(-\dfrac{\lambda^2}{2}\tilde{T}_1 - \nu L_{\tilde{T}_1}\right)\right] = \left[\ch(\lambda)+\dfrac{\nu}{\lambda}\sh(\lambda)\right]^{-1}.\]

Gr\^ace au Théorème \ref{thm_identiteenloi}, on en déduit
\begin{align*}
  &\E\left[\exp\left(-\left(\dfrac{\lambda^2}{2}A^{\scriptscriptstyle (+) }_{T_1} + \dfrac{\mu^2}{2}A^{\scriptscriptstyle (-) }_{T_1} + \alpha L_{T_1}\right) \right)\right]\\
  & \qquad \qquad\qquad \qquad = \E\left[ \exp\left(-\left(\dfrac{\lambda^2}{2}\tilde{T}_1 + \dfrac{(\mu L_{\tilde{T}_1})^2}{2}\widehat{T}_1 + 2\alpha L_{\tilde{T}_1}\right) \right) \right]\\
  &\qquad \qquad \qquad \qquad = \E\left[ \exp\left(-\left(\dfrac{\lambda^2}{2}\tilde{T}_1 + (\mu + 2\alpha) L_{\tilde{T}_1}\right) \right) \right]\\
  &\qquad \qquad \qquad \qquad = \left[ \ch(\lambda) + \dfrac{\mu + 2\alpha}{\lambda}\sh(\lambda) \right]^{-1},
\end{align*}
ce qui conclut cet exercice.
\end{proof}

\begin{exercice}
Soit $B$ un mouvement brownien complexe, et $z_1, \ldots, z_n$ des complexes distincts et non-nuls. Pour tout $t \geq 0$, on définit le nombre de tours réalisés par le mouvement brownien autour de chacun de ces points :
\[ \theta^{(i)}_t = \text{Im}\left[ \log (B_t-z_i) \right] = \text{Im}\left[ \int_0^t \dfrac{\d B_s}{B_s-z_i} \right]. \]
Gr\^ace au théorème de Spitzer, on sait en particulier que
\[\frac{2}{\log t} \theta^{(1)}_t \wconvt C_1,\]
où $C_1$ suit une loi de Cauchy standard. Une généralisation de ce théorème permet d'obtenir la convergence jointe en loi pour les nombres de tours autour de chacun des points considérés, on a en fait :
\[ \frac{2}{\log t }(\theta^{(1)}_t,\ldots, \theta^{(n)}_t) \wconvt (W_1,\ldots W_n),\]
avec 
\[\E\left[ \exp\left( i \sum_{j=1}^n \lambda_j W_j\right) \right] = \left[ \ch\left(\sum_{j=1}^n \lambda_j\right) + \frac{\sum_{j=1}^n |\lambda_j|}{\sum_{j=1}^n \lambda_j} \sh\left(\sum_{i=1}^n \lambda_j\right) \right]^{-1}. \]

Gr\^ace à ce résultat, exprimer la loi jointe du $n-1$-uplet $(W_1-W_n, \ldots W_{n-1}-W_n)$ comme $(C^{(1)}_\e - C^{(n)}_\e, \ldots, C^{(n-1)}_\e - C^{(n)}_\e)$ où $(C^{(j)}_t, j \leq n, t \geq 0)$ est une collection de processus de Cauchy standard indépendants, et $\e$ est une variable aléatoire exponentielle indépendante. Ces variables aléatoires mesurent en quelque sorte les tours réalisés autour de chaque point, qui n'est pas un tour réalisé autour de l'ensemble des points.
\end{exercice}

\begin{proof}
On calcule la transformée de Fourier de ce $n-1$-uplet, on obtient alors, en utilisant la formule donnée dans l'énoncé :
\[\E\left[ \exp\left( i \sum_{j=1}^{n-1} \lambda_j (W_j-W_n)\right) \right] = \left[ 1 + \sum_{j=1}^{n-1} |\lambda_j| + \left|\sum_{j=1}^{n-1} \lambda_j\right| \right]^{-1}.\]

On peut réécrire cette transformée de Fourier de la manière suivante :
\begin{align*}
  & \E\left[ \exp\left( i \sum_{j=1}^{n-1} \lambda_j (W_j-W_n)\right) \right]\\
  &\qquad \qquad\qquad \qquad =  \int_0^{+\infty} \d t e^{-t} \exp\left( -t \sum_{j=1}^{n-1} |\lambda_j| - t \left|\sum_{j=1}^{n-1} \lambda_j\right| \right)\\
  &\qquad \qquad\qquad \qquad = \int_0^{+\infty} \d t e^{-t} \E\left[ \exp\left( i \sum_{j=1}^{n-1} \lambda_j C^{(j)}_t - i \left(\sum_{j=1}^{n-1} \lambda_j\right) C^{(n)}_t \right)  \right]\\
  &\qquad \qquad\qquad \qquad = \E\left[ \exp\left( i \sum_{j=1}^{n-1} \lambda_j (C^{(j)}_\e - C^{(n)}_\e) \right)\right],
\end{align*}
où les $(C^{(j)}_t, j \leq n, t \geq 0)$ sont des processus de Cauchy standard indépendants, et $\e$ est une variable aléatoire exponentielle indépendante.

On tire de cette égalité
\[(W_1-W_n, \ldots, W_{n-1}-W_n) \egaldistr (C^{(1)}_\e - C^{(n)}_\e, \ldots, C^{(n-1)}_\e - C^{(n)}_\e). \]
\end{proof}

\begin{remarque}
En réalité, on peut décomposer les variables aléatoires $W_i$ en deux parties : une partie $W_i^-$ correspondant au nombre de tours réalisés au voisinage du point $z_i$, et une partie commune $W^+$ correspondant aux nombres de tours réalisés loin des points. On a en particulier :
\[ (W_1^-, \ldots, W^-_n) \egaldistr  (C^{(1)}_\e, \ldots, C^{(n)}_\e), \]
qui est bien entendu un résultat plus fort que celui que nous venons de démontrer.
\end{remarque}

\begin{exercice}[Sur la décomposition de $B$ en partie positive et partie négative]
Soit $B$ mouvement brownien standard, on pose $S_t = \sup_{s \leq t} B_s$ et $I_t=-\inf_{s\leq t} B_s$.

\begin{enumerate}
  \item Montrer que $I_{T_1}$ a pour densité $\frac{1}{(1+x)^2}\ind{x>0}$  par rapport à la mesure de Lebesgue.
  \item Montrer que le processus $(S_{\tau_l})_{l\geq 0}$ est un processus de Feller dont on explicitera le semi-groupe. Ce processus est appelé processus de Watanabe.
  \item Montrer que $(I_{T_t})_{t\geq 0} \egaldistr (U_{V^{-1}_t})_{t\geq 0}$, où $U$ et $V$ sont deux processus de Watanabe indépendants.
\end{enumerate}
\end{exercice}

\begin{proof}
\textit{1.} On observe pour commencer que $1-B_{t\wedge T_1}$ est une martingale positive qui converge presque sûrement vers 0. L'Exercice \ref{exo_supmartingale} prouve que $\sup_{t \geq 0} (1-B_{t\wedge T_1})$ est distribuée comme $\frac{1}{U}$, où $U$ est une variable aléatoire uniforme sur $[0,1]$. En particulier
\[\sup_{t \geq 0} (1-B_{t\wedge T_1}) = 1- \inf_{t \geq 0} B_{t\wedge T_1} = 1 - \inf_{t \leq T_1}B_t = 1 +I_{T_1}.\]
Par conséquent, $I_{T_1} \egaldistr \frac{1}{U}-1$, et a donc pour densité $\frac{\d x}{(1+x)^2}\ind{x>0}$.

\textit{2.} Étudions maintenant le processus $(S_{\tau_l})_{l\geq 0}$. On peut calculer $S_{\tau_{l+l'}}$ de la manière suivante :
\[S_{\tau_{l+l'}} = \max\left(\sup_{s\leq \tau_l} B_s ; \sup_{\tau_l\leq s \leq \tau_{l'}} (B_s-B_{\tau_l})\right).\]
En appliquant la propriété de Markov forte du mouvement brownien et en observant que $B_{\tau_l}=0$, on remarque que $S_{\tau_{l+l'}}$ est le maximum de deux variables aléatoires indépendantes de loi respectivement égales à celles de $S_{\tau_l}$ et $S_{\tau_{l'}}$. Par conséquent
\[\E(f(S_{\tau_{l+l'}})|\F_{\tau_l}) = \E(f(S_{\tau_l} \vee \widehat{S}_{l'}|\F_{\tau_l}),\]
où $\widehat{S}_{l'}$ est une copie de $S_{\tau_{l'}}$ indépendante de $\F_{\tau_l}$. Le processus $(S_{\tau_l})_{l\geq 0}$ est donc un processus de Markov, et il suffit pour calculer son semi-groupe de déterminer la loi de $S_{\tau_l}$.

Or
\[\P(S_{\tau_l}\leq t) = \P(T_t \geq  \tau_l) = \P(L_{T_t} \geq l).\]
D'après l'Exercice \ref{exo_transformeelaplace}, la transformée de Laplace de $L_{T_1}$ est
\begin{align*}
  \E(\exp(-\alpha L_{T_1})) &= \lim_{\mu,\lambda \to 0} \dfrac{1}{\ch(\lambda)+\dfrac{\mu+2\alpha}{\lambda}\sh(\lambda)}\\
  & = \dfrac{1}{1 + 2\alpha},
\end{align*}
donc $L_{T_1}$ est une variable aléatoire de loi exponentielle de paramètre 2. Par propriété de scaling, $L_{T_t}$ suit la loi exponentielle de paramètre $\frac{2}{t}$. Par conséquent
\[P(S_{\tau_l}\leq t) = \P(L_{T_t}\geq l) = \exp\left(-\dfrac{2l}{t}\right).\]

\begin{remarque}
On peut également directement prouver que $L_{T_t}$ suit une loi exponentielle de paramètre  $\frac{2}{t}$. On a
\[\P(L_{T_t}>l+l'|L_{T_t}>l)=\P(\tau_{l+l'} \leq T_t | \tau_l \leq T_t)= \P(\tau_{l'} \leq T_t) = \P(L_{T_t}>l'),\]
par propriété de Markov, donc $L_{T_t}$ suit une loi exponentielle dont le paramètre reste à déterminer. Or, par symétrie, on a~:
\[\E\left[\int_0^{T_t} \ind{B_s>0}\d B_s\right]  =\E\left[\int_0^{T_t} \ind{B_s\leq 0}\d B_s\right] = \frac{t}{2},\]
donc en écrivant $B_{T_t}^-$ avec la formule de Tanaka, on obtient :
\[\E(L_{T_t})=\frac{2}{t},\]
ce qui permet de conclure.
\end{remarque}

On calcule le semi-groupe $(P_l, l \geq 0)$ du processus de Markov $(S_{\tau_l},l\geq 0)$~:
\begin{align*}
  P_{l'}f(x) = & \E(f(S_{\tau_{l+l'}}|S_{\tau_l}=x)
  =  \E(f(S_{\tau_{l'}}\vee x))\\
  = & f(x) \P(S_{\tau_{l'}} \leq x) + \E(f(S_{\tau_{l'}}) \ind{S_{\tau_{l'}}\geq x}),
\end{align*}
par conséquent
$P_{l}(x,\d y) = \exp\left(-\frac{2l}{x}\right)\delta_x(\d y) + \frac{2l}{y^2}\exp\left(-\frac{2l}{y}\right)\ind{y>x}\d y$.

On lit directement sur ce semi-groupe que le processus associé est de Feller. En effet pour toute fonction $f \in \C_b$ on a
\[e^{-\frac{2l}{x}}f(x)- \Vert f \Vert_\infty(1-e^{-\frac{2l}{x}}) \leq P_lf(x) \leq e^{-\frac{2l}{x}}f(x)+ \Vert f \Vert_\infty(1-e^{-\frac{2l}{x}}),\]
donc $\Vert P_lf-f\Vert_{\infty} \underset{l \to 0}{\longrightarrow} 0$.

\textit{3.} On montre maintenant que si $U$ et $V$ sont deux processus de Watanabe indépendants, alors $(I_{T_t}, t \geq 0) \egaldistr (U_{V^{-1}_t},t \geq 0)$. Soit $B$ et $\widehat{B}$ deux mouvements browniens indépendants et
\[ U_t = \widehat{S}_{\widehat{\tau}_t}, \quad V_t = S_{\tau_t}. \]
On observe pour commencer que $V^{-1}_x = L_{T_x}$. Par conséquent
\begin{align*}
  \P(U_{V^{-1}_t} \geq x) &= \P(U^{-1}_x \leq V^{-1}_t) = \P(\widehat{L}_{\widehat{T}_x} \leq L_{T_t})
  = \P\left(\frac{x}{2}\e \leq \frac{t}{2}\widehat{\e}\right) = \dfrac{t}{t+x}.
\end{align*}
Grâce à la propriété de scaling du mouvement brownien, $I_{T_t} \egaldistr tI_{T_1}$, donc
\[ \P(I_{T_t} \geq x ) = \P(U_{V^{-1}_t} \geq x), \]
les marginales unidimensionnelles des deux processus sont bien égales.

Pour passer aux marginales finies-dimensionnelles, on observe que
\[ I_{T_{t+t'}} = \max\left\{ I_{T_t}, -t - \inf_{0 \leq s \leq T_{t+t'}-T_t} (B_{T_t+s}-B_{T_t})\right\}. \]
Par propriété de Markov forte du mouvement brownien, les deux termes de ce maximum sont indépendants, et le second est égal en loi à $-t+I_{T_{t'}}$. Les marginales finies-dimensionnelles de $I_{T_t}$ se déduisent donc des marginales unidimensionnelles.

On montre maintenant le même type de propriété pour $U_{V^{-1}}$. Par propriété du processus de Watanabe, on a
\begin{align*}
  U_{V^{-1}_{t+t'}} = U_{L_{T_{t+t'}}} = \max\left\{ U_{L_{T_t}} , \bar{U}_{L_{T_{t+t'}}-L_{T_t}} \right\},
\end{align*}
où $\bar{U}$ est un processus de Watanabe indépendant de $U$. On pose $d_{T_t} = \inf\{u \geq T_t:B_u=0\}$ et on étudie la loi de $L_{T_{t+t'}}-L_{T_t}$ sous l'alternative suivante.
\begin{itemize}
  \item Conditionnellement à $d_{T_t} \leq T_{t+t'}$, par propriété de Markov, $L_{T_{t+t'}}-L_{T_t}$ est indépendant de $\F_{T_t}$ et est distribué comme $L_{T_{t+t'}}$.
  \item Conditionnellement à  $T_{t+t'} \leq d_{T_t}$, $L_{T_{t+t'}}-L_{T_t} = 0$ car $L$ ne croit que sur l'ensemble des zéros du mouvement brownien.
\end{itemize} 
On applique la propriété de Markov forte en $T_t$, puis un théorème d'arrêt à $B$,
\begin{align*}
  (t+t') \P(T_{t+t'} \leq d_{T_t}) = \E\left[ B_{T_{t+t'} \wedge d_{T_t}} \right]
  = \E_t\left[ B_{T_0 \wedge T_{t'}} \right]
  = t,
\end{align*}
donc $\P(T_{t+t'} \leq d_{T_t}) = \frac{t}{t+t'}$.

On en déduit le résultat suivant
\begin{align*}
  \P(U_{V^{-1}_{t+t'}} \leq x|\F_{V^{-1}_t})
  = & \ind{U_{V^{-1}_t}\leq x} \P(U_{V^{-1}_{t+t'}-V^{-1}_t} \leq x)\\
  = & \ind{U_{V^{-1}_t}\leq x} \left[ \frac{t}{t+t'} + \frac{t'}{t+t'} \P(U_{V^{-1}_{t+t'}} \leq x)\right]\\
  = & \ind{U_{V^{-1}_t}\leq x} \left[ \frac{t}{t+t'} + \frac{t'}{t+t'} \frac{x}{t+t'+x}\right]\\
  = & \ind{U_{V^{-1}_t}\leq x} \frac{t+x}{t+t'+x}\\
  = & \ind{U_{V^{-1}_t}\leq x} \P(U_{V^{-1}_{t'}}\leq x+t)\\
  = & \ind{U_{V^{-1}_t}\leq x} \P(U_{V^{-1}_{t'}}-t\leq x)\\
  = & \P\left(\max\left\{ U_{V^{-1}_{t}} , (\bar{U}_{\bar{V}^{-1}_{t'}}-t)\right\} \leq x|\F_{V^{-1}_t}\right).
\end{align*}
Les deux processus ont donc le même semi-groupe et les mêmes marginales unidimensionnelles, par conséquent
$(U_{V^{-1}_t})_{t \geq 0} \egaldistr (I_{T_t})_{t \geq 0}$.
\end{proof}

\begin{exercice}
\label{exo_raykdemo}
Soit $g$ une fonction continue ; pour quelles fonctions $f$ le processus
\[f(B^+_t,L_t) \exp\left(-\frac{1}{2}\int_0^t g(B_s)\ind{B_s>0}\d s\right)\]
est-il une martingale locale ?
\end{exercice}

\begin{proof}
On pose $Z_t = \int_0^t g(B_s)\ind{B_s>0}\d s$ et $h(x,y,z) = f(x,y)\exp\left(-\frac{1}{2}z\right)$, on applique la formule d'Itô à $h\left(B^+_t,L_t,Z_t\right)$, on a :
\begin{multline*}
  h(B^+_t,L_t,Z_t) = f(0) + \int_0^t \partial_1 h(B_s^+,L_s,Z_s)\ind{B_s>0}\d B_s\\
  + \dfrac{1}{2}\int_0^t \partial_1 h(B_s^+,L_s,Z_s)\d L_s + \int_0^t \partial_2 h(B^+_s,L_s,Z_s)\d L_s\\
  -\dfrac{1}{2} \int_0^t \left(h(B^+_s,L_s,Z_s) g(B_s) - \partial^2_{1} h(B_s^+,L_s,Z_s)\right)\ind{B_s>0}\d s.
\end{multline*}

Par conséquent, $f$ doit satisfaire les deux équations suivantes
\[\left\{
\begin{array}{l}
  \dfrac{1}{2} \partial_1 f(0,y) +\partial_2 f(0,y) = 0 \\
  \partial^2_{1}f(x,y) = g(x)f(x,y).
\end{array}
\right.\]

En cherchant $f$ sous la forme $\phi(x)\psi(y)$, on a $\phi=\Phi_\mu$ la solution de l'équation de Sturm-Liouville associée à $\mu(\d x)=g(x)\d x$ et $\psi(y) = \exp(-\frac{1}{2}y\Phi'_\mu(0))$.
Par conséquent
\[\left(\Phi_\mu(B^+_t) \exp\left(-\frac{1}{2}\left( L_t \Phi'_\mu(0^+) + \int_0^t g(B_s)\ind{B_s>0}\d s\right)\right)\right)_{t\geq 0}\]
est une martingale locale.

\begin{remarque}
On observe que ce processus est borné lorsqu'arrêté en $\tau_l$, on peut donc appliquer le théorème d'arrêt. Or
\[\int_0^t g(B_s)\ind{B_s>0}\d s = \int_0^{+\infty} g(x)L^x_t\d x,\]
par conséquent
\[\E\left[\exp\left(-\dfrac{1}{2}\int_0^{+\infty}L^x_{\tau_l}g(x)\d x\right)\right] = \exp\left(\dfrac{l}{2}\Phi'_\mu(0^+)\right).\]
Par densité, ce résultat s'étend à toute mesure $\mu$ sur $[0,+\infty)$. Par conséquent la loi de $(L^x_{\tau_l})_{x \geq 0}$ est $Q^0_l$, carré de Bessel de dimension 0 issu de $l$.
\end{remarque}
\end{proof}

\chapter[Théorèmes de Ray-Knight]{Sur les théorèmes de Ray-Knight et les identités de Ciesielski-Taylor}

Dans ce chapitre, nous étudions les temps locaux du mouvement brownien linéaire, non pas en temps que processus temporel à niveau $x$ fixé, mais en temps que processus indexé par la variable d'espace $x \in \R$, pris en certains temps aléatoires. Nous verrons que ces processus suivent les lois de certains carrés de Bessel. Ce type de théorème peut s'étendre à de nombreuses semi-martingales continues vérifiant en outre une propriété de type Markov, permet de démontrer certaines égalités en loi classiques, par exemple celle de Ciesielski-Taylor.

\section{Théorèmes principaux}

On définit tout d'abord le processus stochastique appelé carré de processus de Bessel, parfois abrégé en carré de Bessel.

\begin{definition}
Soit $\delta, l$ deux réels positifs et $B$ un mouvement brownien. Un carré de processus de Bessel de dimension $\delta$ issu de $l$ est la seule solution de l'équation différentielle stochastique
\[ \forall a \geq 0, Z_a = l + 2 \int_0^a \sqrt{Z_b} \d B_b + \delta a. \]
On note $Q^\delta_l$ la loi de ce processus.
\end{definition}

Ces processus stochastiques vérifient une propriété d'additivité remarquable.

\begin{propriete}
Soit $\delta,\delta',l,l' \in \R^+$, si $Z$ et $Z'$ sont deux processus stochastiques indépendants tels que $Z$ a pour loi $Q^\delta_l$ et $Z'$ a pour loi $Q^\delta_{l'}$, alors $Z+Z'$ a pour loi $Q^{\delta + \delta'}_{l + l'}$.

On a $Q^{\delta + \delta'}_{l + l'} = Q^\delta_l \ast Q^{\delta'}_{l'}$.
\end{propriete}

On introduit maintenant le théorème principal du chapitre, qui identifie la loi des temps locaux d'un mouvement Browien arrêtés en certain temps aléatoires.

\begin{theoreme}[Ray-Knight]
Soit $B$ un mouvement Brownien et $L$ ses temps locaux. Pour tout $a \in \R$ et $l \geq 0$, on a
\begin{enumerate}
  \item $\left(L^{a-x}_{T_a}\right)_{0\leq x \leq a}$ a pour loi $Q^2_0$,
  \item $(L^x_{\tau_l})_{x\geq 0}$ et $(L^{-x}_{\tau_l})_{x\geq 0}$ sont indépendants et ont pour loi $Q^0_l$.
\end{enumerate}
\end{theoreme}

\begin{remarque}
La seconde partie de ce théorème a été démontrée dans l'Exercice \ref{exo_raykdemo}, mais nous allons ici détailler un autre type de preuves, basées sur l'identification des processus comme solutions d'équations différentielles stochastiques.
\end{remarque}

L'une des nombreuses façons de démontrer ce théorème consiste à s'appuyer sur le lemme suivant, qui permet de donner une représentation des variables aléatoires mesurables dans la tribu des temps locaux, et donc par la suite d'étudier des martingales par rapport à ces temps locaux.

\begin{lemme}
\label{lem_rep-marttempslocal}
Pour tout $0\leq b \leq 1$ on pose $\mathcal{Z}_b = \sigma(L^{1-a}_{T_1},0\leq a \leq b)$ la filtration associée au processus $(L^{1-a}_{T_1})_{0\leq a \leq 1}$.

Toute variable aléatoire $H$ dans $L^2(\mathcal{Z}_b)$ peut s'écrire sous la forme :
\[H = \E(H) + \int_0^{T_1} h_s \ind{B_s>1-b} \d B_s,\]
avec $h$ processus prévisible vérifiant $\E\left[\int_0^{T_1}h_s^2 \ind{B_s>1-b}\d s\right] <+\infty$
\end{lemme}

Pour tout $\delta\geq 0$, on note $R^{(\delta)}$ un processus de Bessel de dimension $\delta$ issu de 0.

\begin{theoreme}[Extension du théorème de Ray-Knight aux processus de Bessel]
\label{thm_extensionrk}
On note $B$ un mouvement brownien dans $\R^2$ et $\beta$ un pont brownien dans $\R^2$. Soit $\delta >  0$ et $0\leq \gamma < 2$, on a les identités en loi suivantes
 \[\left(L^a_{\infty}(R^{(\delta+2)}),a \geq 0\right) \egaldistr \left(\frac{1}{\delta a^{\delta-1}}|B_{a^\delta}|^2 ,a \leq 1\right)\]
 \[\text{et  } \left(L^a_{T_1}(R^{(\delta+2)}),a \leq 1\right) \egaldistr \left(\frac{1}{\delta a^{\delta-1}}|\beta_{a^\delta}|^2,a \leq 1\right).\]
De plus, on a également :
\[\left(L^a_{T_1}(R^{(2)}),a \leq 1\right) \egaldistr \left(a|B_{-\log a}|^2,a \leq 1\right)\]
\[\text{et  } \left(L^a_{T_1}(R^{(2-\gamma)}),a \leq 1\right) \egaldistr \left(\frac{1}{\gamma a^{\gamma-1}}|B_{1-a^\gamma}|^2,a \leq 1\right).\]
\end{theoreme}

Ce théorème permet également d'obtenir l'identité de Ciesielski-Taylor.

\begin{theoreme}[Identité de Ciesielski-Taylor]
Pour tout $\delta > 0$, on a
\[\int_0^{+\infty} \ind{R^{(\delta)}\leq 1}\d s \egaldistr T_1(R^{(\delta)}).\]
\end{theoreme}

Les résultats suivants, démontrés en exercices, permettent de lier entre eux des processus de Bessel de dimension différentes par des changements de variables espaces-temps, ce qui sera utile pour démontrer le Théorème \ref{thm_extensionrk}

\begin{theoreme}
Soit $\alpha>0$, on pose $h(x) = x^\alpha$ et, pour $d >0$, on note $R^{(d)}$ un processus de Bessel de dimension $d$. On a les identités suivantes :
\begin{itemize}
 \item si $d>2(1-\alpha)$, soit $d_\alpha = \frac{d}{\alpha}+2(1-\frac{1}{\alpha})$, on a 
 \[h(R^{(d)}_t) = R^{(d_\alpha)}\left( \int_0^t h'(R^{(d)}_u)^2 \d u \right),\]
 \item si $\alpha>0$, on a 
 \[h(R^{(2+\alpha)}_t) = R^{(3)}\left( \int_0^t h'(R^{(2+\alpha)}_u)^2 \d u \right),\]
 \item si $0<\alpha<2$, on a 
 \[(h(R^{(2-\alpha)}_t))_{t \geq 0} = R^{(1)}\left( \int_0^t h'(R^{(2-\alpha)}_u)^2 \d u \right).\] 
\end{itemize}
\end{theoreme}

\section{Exercices}

\begin{exercice}[Tribu horizontale et tribu verticale]
Montrer que la tribu $\L_t = \sigma\left\{L^a_s, a \in \R, s \leq t \right\}$ est égale à $\mathcal{B}_t = \sigma\{B_s,s \leq t\}$, aux ensembles négligeables près.
\end{exercice}

\begin{proof}
Observons pour commencer que, grâce à la généralisation de la formule de densité d'occupation (Exercice \ref{exo_generalisationdensite}), pour toute fonction $\phi$ mesurable positive et tout $t \geq 0$, on a :
\[\int_0^t \phi(s,B_s)\d s = \int_\R \d x \int_0^t \d_sL^x_s \phi(s,x) \text{   p.s.}\]

On observe que l'expression de gauche est $\mathcal{B}_t$-mesurable, et celle de droite est $\L_t$-mesurable. Il reste donc à montrer que la tribu 
\[\mathcal{T} = \sigma\left( \int_0^t \phi(s,B_s) \d s, \phi \text{ mesurable positive} \right)\]
contient $\mathcal{B}_t$ et $\L_t$.

Pour commencer, soit $0\leq t_1<\cdots<t_p\leq t$, $\lambda_1,\ldots, \lambda_p$ des réels positifs et $\Gamma_1,\ldots, \Gamma_p$ des ouverts de $\R$. On pose :
\[\phi_n(t,x) = \sum_{i=1}^p \lambda_i 2n\ind{t \in (t_i-\frac{1}{n} , t_i + \frac{1}{n}} \ind{x \in \Gamma_i},\]
on observe que
\[ \int_0^t \d s \phi_n(s,B_s) \underset{n \to +\infty}{\longrightarrow} \sum_{i=1}^p \lambda_i\ind{B_{t_i} \in \Gamma_i} \text{   p.s.}\]
par conséquent $\mathcal{T}=\mathcal{B}_t$ aux ensembles négligeables près.

On étudie maintenant la fonction suivante :
\[\phi_n(t,x) = \sum_{i=1}^p \lambda_i 2n \ind{t \leq t_i}\ind{x \in (x_i-\frac{1}{n},x_i+\frac{1}{n})},\]
on obtient la limité suivante
\[ \int_\R \d x \int_0^t \d_sL^x_s \phi_n(s,x) \conv \sum_{i=1}^n \lambda_i L^x_{t_i} \text{   p.s.}\]
par conséquent on a également $\mathcal{T}=\L_t$.
\end{proof}

\begin{exercice}
\label{exo_RKbang}
Soit $k \in \R$, on considère l'équation différentielle stochastique :
\[X_t = B_t + k \int_0^t \ind{X_s>0} \d s.\]
On note $P^{(k)}$ la loi de cette solution.

\begin{enumerate}
 \item Écrire la relation d'absolue continuité entre $P^{(k)}$ et la mesure de Wiener.
 \item Exprimer la loi de $(L^{1-a}_{T_1}(X))_{0 \leq a \leq 1}$ et celle de $\left((L^x_{\tau_l})_{x \geq 0},L^{-x}_{\tau_l})_{x \geq 0}\right)$ sous $P^{(k)}$ en fonction des lois ${}^{(\beta)}Q^\delta_l$ de processus de diffusion solutions de la famille d'équations différentielles stochastiques :
 \[Z_a = l + 2\int_0^a \sqrt{Z_b}\d B_b + \int_0^a (2 \beta Z_b +\delta) \d b.\]
 \item Étudier de même le temps local de la solution de l'équation différentielle stochastique :
\[X_t = B_t + k \int_0^t \sgn(X_s) \d s.\]
\end{enumerate}
\end{exercice}

\begin{proof}
\textit{1.}
Pour toute martingale locale $X$ on pose
\begin{equation}
  \label{eqn:martingaleExponentielle}
  \forall t \geq 0, \mathcal{E}(X)_t = \exp\left( X_t - \frac{1}{2} \crochet{X}_t \right)
\end{equation}
la martingale exponentielle associée à $X$. En particulier, on remarquera que $\mathcal{E}(X)$ est une martingale locale positive.

On note $(M_t)$ la martingale définie, pour tout $t \geq 0$, par
\[M_t = \mathcal{E}\left(k\int_0^. \ind{X_s>0} \d X_s\right)_t = \exp\left( k\int_0^t \ind{X_s>0} \d X_s -\dfrac{k^2}{2}\int_0^t \ind{X_s>0}\d s \right).\]
On note $\P^{(k)}$ la loi vérifiant, pour tout $t \geq 0$, $P^{(k)}_{|\F_t} = M_t . W_{|\F_t}$.

Soit $X$ un processus tel que sous $W$, $X$ est un mouvement brownien. En appliquant les théorèmes de Girsanov et de Lévy, le processus $B$ défini par :
\[B_t = X_t - k\int_0^t \ind{X_s>0} \d s\]
est un mouvement brownien sous $\P^{(k)}$. Par conséquent, sous $P^{(k)}$, $X$ est solution de l'équation différentielle stochastique :
\[X_t = B_t + k \int_0^t \ind{X_s>0} \d s.\]

\textit{2.} Posons pour $0 \leq a \leq 1, Z_a = L^{1-a}_{T_1}(X)$, on souhaite étudier la loi de $(Z_a)_{0 \leq a \leq 1}$ sous $P^{(k)}$. Sous $W$, ce processus est une semi-martingale de loi $Q^2_0$. Par conséquent, sous $P^{(k)}$, la loi de cette fonctionnelle de $X$ est $M_{T_1} \cdot Q^2_0$ (en effet, $M_{T_1}$ est bien mesurable par rapport à $\sigma(Z_a,0\leq a \leq 1)$), donc est encore une semi-martingale, par la formule de Girsanov.

On détermine maintenant une équation différentielle stochastique satisfaite par le processus $(Z_a)_{0\leq a \leq 1}.$ Par la formule de Tanaka, on a
\begin{equation}
\label{eqn_Z}
Z_a = 2 (X_{T_1}-(1-a))^+ - 2 \int_0^{T_1} \ind{X_s>1-a} \d X_s = 2a -2 \int_0^{T_1} \ind{X_s>1-a} \d X_s.
\end{equation}

On utilise maintenant l'équation différentielle satisfaite par $X$, pour réécrire :
\[Z_a = 2 a - 2 \int_0^{T_1} \ind{X_s>1-a} \d B_s - 2k \int_0^{T_1} \ind{X_s>1-a} \d s.\]
Par conséquent, par définition des temps locaux :
\[Z_a -2a + 2 k \int_0^a Z_b \d b = -2 \int_0^{T_1} \ind{X_s>1-a} \d B_s.\]
En utilisant le Lemme \ref{lem_rep-marttempslocal} de représentation des variables aléatoires dans la tribu $(\mathcal{Z}_a)$, on en déduit que le processus $\left(Z_a - 2a+2k \int_0^a Z_b \d b, a \geq 0 \right)$ est une martingale. De plus, cette martingale a pour crochet $4\int_0^a Z_b \d b$, le crochet n'étant pas modifié par un changement de probabilité entre probabilités équivalentes. Par conséquent, $Z$ satisfait l'équation différentielle stochastique
\[ Z_a = 2 \int_0^a \sqrt{Z_b}\d \beta_b + \int_0^a \left( 2 - 2k Z_b\right) \d b,\]
où $\beta$ est un mouvement brownien adapté à la filtration $(\mathcal{Z}_a)$ ; $Z$ a pour loi ${}^{-k}Q^2_0$.

On s'intéresse maintenant à $Y_x = L^x_{\tau_l}$, on étudie la loi de $(Y_x)_{x \geq 0}$ et $(Y_{-x})_{x\geq 0}$. De la même manière que précédemment, on identifie les équations différentielles satisfaites par ces processus gr\^ace à la formule de Tanaka.
\[f(X_{\tau_l})-f(X_0) - \int_0^{\tau_l} \ind{0<X_s\leq x} \d X_s = \dfrac{1}{2}(L^0_{\tau_l}(X)-L^x_{\tau_l}(X)). \]
Par conséquent, pour tout $x \geq 0$
\begin{equation}
\label{eqn_Y+}
Y_x = l + 2 \int_0^{\tau_l} \ind{0<X_s\leq x} \d X_s
\end{equation}
et de même pour $x \leq 0$, on a
\begin{equation}
\label{eqn_Y-}
Y_x = l + 2 \int_0^{\tau_l} \ind{x < X_s\leq 0} \d X_s.
\end{equation}

On s'intéresse tout d'abord au cas $x \geq 0$ ; on a :
\begin{align*}
  Y_x & = l + 2 \int_0^{\tau_l} \ind{0<X_s\leq x} \d X_s \\
  & = l + 2 \int_0^{\tau_l} \ind{0<X_s\leq x} \d B_s + 2 k \int_0^{\tau_l} \ind{0<X_s\leq x}\d s\\
  & = l + 2 \int_0^{\tau_l} \ind{0<X_s\leq x} \d B_s + 2 k \int_0^x Y_b \d b.
\end{align*}
Par conséquent, étant donné que, sous $W$, $Y$ est une semi-martingale de crochet $4\int_0^a Y_b \d b$, sous $P^{(k)}$ ce processus suit l'équation différentielle stochastique suivante :
\[Y_x = l + 2 k \int_0^x Y_b\d b + 2 \int_0^x \sqrt{Y_b}\d \beta_b,\]
donc suit la loi ${}^{(k)}Q^0_l$. Ces calculs peuvent être réalisés de la même manière pour $L^{-x}_{\tau_l}$, on observe ainsi que ces deux processus sont indépendants et que $(L^{-x}_{\tau_l})_{x \geq 0}$ a pour loi ${}^{(0)}Q^0_l$.

\textit{3.} On peut réécrire le même type d'équations dans le cas où la dérive est donnée par $k \text{ } \sgn(X_t)$, et grâce aux équations \eqref{eqn_Z}, \eqref{eqn_Y+} et \eqref{eqn_Y-}, on obtient: $(L^{1-a}_{T_1})_{0\leq a \leq 1}$ suit la loi ${}^{(-k)}Q^2_0$, et $(L^{x}_{\tau_l})_{x \geq 0}$ et $(L^{-x}_{\tau_l})_{x \geq 0}$ sont indépendants et de même loi ${}^{(k)}Q^0_l$.
\end{proof}

\begin{exercice}
Soit $\mu$ mesure de Radon ; on note $\Phi_\mu$ la solution de l'équation de Sturm-Liouville 
\[\Phi'' = \mu \Phi, \text{  } \Phi \text{ décroissante issue de 1} \]

On note $\Pi^{(\mu)}$ la probabilité définie par :
\[\Pi^{(\mu)}_{|\F_t} = \exp\left( - \dfrac{L_t}{2}\Phi'_\mu(0^+)\right) \Phi_\mu (B^+_t) \exp\left( -\dfrac{1}{2}\int_{\R^+} L^x_t\mu(\d x) \right) . W_{|\F_t}.\]
Énoncer un théorème de Ray-Knight pour $L$ sous $\Pi^{(\mu)}$, pour cela nous utiliserons la loi :
\[Q^{(\delta,\mu)}_x = \exp\left(\dfrac{1}{2}\int f(s) \d M_s - \dfrac{1}{2} \int_0^t f(s)^2 X_s \d s\right) . Q^\delta_x,\]
où on a posé $f = \dfrac{\Phi'_\mu}{\Phi_\mu}$ et $M_s = X_s-x-\delta s$.
\end{exercice}

\begin{proof}
Posons 
\[M_t = \exp\left( - \dfrac{L_t}{2}\Phi'_\mu(0^+)\right) \Phi_\mu (B^+_t) \exp\left( -\dfrac{1}{2}\int_{\R^+} L^x_t\mu(\d x) \right).\]

Dans un premier temps, on cherche à caractériser la loi de $Z_a = L^{(1-a)}_{T_1}$ sous $\Pi^{(\mu)}$ ; ce processus est une semi-martingale, on cherche une équation différentielle stochastique vérifiée par celui-ci. On utilise le fait que sous $\Pi^{(\mu)}$, $B$ s'écrit :
\[B_t = \tilde{B}_t + \int_0^t \dfrac{\d \crochet{B,M}_s}{M_s}.\]

Or la martingale $M_t$ s'écrit, au moyen de la formule d'Itô-Tanaka :
\[1 + \int_0^t \dfrac{\Phi'_\mu(B_s^+)}{\Phi_\mu(B_s^+)}M_s \ind{B_s>0} \d B_s.\]

Par conséquent, sous $\Pi^{(\mu)}$,
\[B_t = \tilde{B}_t + \int_0^t \dfrac{\Phi'_\mu(B_s^+)}{\Phi_\mu(B_s^+)}\ind{B_s>0} \d s.\]

On sait que $Z_a$ satisfait l'équation \eqref{eqn_Z}, donc :
\begin{align*}
  Z_a  = & 2a - 2 \int_0^{T_1} \ind{B_s>1-a} \d B_s\\
   = & 2a - 2\int_0^{T_1} \ind{B_s>1-a} \d \tilde{B}_s - 2 \int_0^{T_1} \dfrac{\Phi'_\mu(B_s^+)}{\Phi_\mu(B_s^+)}\ind{B_s>1-a}\d s.
\end{align*}
En particulier on a $\crochet{Z}_a = 4 \int_0^a Z_b \d b$, donc $Z$ satisfait l'équation différentielle stochastique :
\[Z_a = 2a - 2 \int_0^{T_1} \sqrt{Z_b} \d \beta_b - 2 \int_0^a Z_b \dfrac{\Phi'_\mu(1-b)}{\Phi_\mu(1-b)}\d b.\]

De la même manière on peut s'intéresser à la loi de $Y_x = L^x_{\tau_l}$. Sous $W$ ce processus vérifie l'équation \eqref{eqn_Y+}:
\[Y_x = l + 2 \int_0^{\tau_l} \ind{0<B_s\leq x} \d B_s = l + 2\int_0^{T_1} \ind{0<B_s\leq x} \d \tilde{B}_s + 2 \int_0^{T_1} \dfrac{\Phi'_\mu(B_s^+)}{\Phi_\mu(B_s^+)}\ind{0<B_s\leq x}\d s.\]
Donc $Y$ satisfait l'équation différentielle stochastique
\[Y_x = l + 2\int_0^x \sqrt{Z_y} \d \beta_y + 2 \int_0^x \dfrac{\Phi'_\mu(y)}{\Phi_\mu(y)} Y_y \d y.\]
On utilise la formule de Girsanov, $(Y_x)_{x \geq 0}$ suit donc la loi $Q^{2,\mu}_l$. Mutatis mutandis, $(Y_{-x})_{x \geq 0}$ suit la loi $Q^0_l$.
\end{proof}

\begin{exercice}
Soit $f$ décroissante et $g$ croissante, deux fonctions continues sur $[a,b]$ et $B$ un mouvement brownien.

\begin{enumerate}
 \item Montrer que $-\int_a^b B^2_{g(x)}\d f(x) +f(b)B^2_{g(b)} \egaldistr \int_a^b B^2_{f(x)}\d g(x) +g(a)B^2_{f(a)}.$
 \item En déduire une extension des identités de Ciesielski-Taylor.
\end{enumerate}
\end{exercice}

\begin{proof}
\textit{1.} On calcule la transformée de Laplace de ces variables aléatoires, soit $f$ et $g$ deux fonctions continues croissantes et $B$ un mouvement brownien, on a
\begin{multline*}
  \E\left[\exp\left( \dfrac{\lambda^2}{2} \int_a^b B^2_{g(x)}\d f(x) - \dfrac{\lambda^2}{2}f(b)B^2_{g(b)}\right) \right]\\
 =  \E\left[\exp\left( i \lambda \left(\int_a^b B_{g(x)}\d C_{f(x)} + B_{g(b)}C_{f(b)}\right)\right)\right]
\end{multline*}
où $(C_t,t \geq 0)$ est un mouvement brownien indépendant de $B$. Par intégration par parties stochastique, on a
\begin{align*}
  &\E\left[\exp\left( \dfrac{\lambda^2}{2} \int_a^b B^2_{g(x)}\d f(x) - \dfrac{\lambda^2}{2}f(b)B^2_{g(b)}\right) \right]\\
 & \qquad \qquad =  \E\left[\exp\left( i \lambda \left(-\int_a^b C_{f(x)}\d B_{g(x)} + C_{f(a)}B_{g(a)}\right)\right)\right]\\
 & \qquad \qquad =  \E\left[\exp\left( - \dfrac{\lambda^2}{2} \int_a^b C^2_{f(x)}\d g(x) - \dfrac{\lambda^2}{2}g(b)C^2_{f(b)}\right) \right].
\end{align*}
Par égalité des transformées de Laplace, on en déduit l'égalité en loi des variables aléatoires.

\textit{2.} On utilise alors les extensions du théorème de Ray-Knight (le Théorème \ref{thm_extensionrk}) pour traduire les équations de Cieselski-Taylor, on a, pour $\delta < 2$, $0\leq x \leq y \leq 1$
\[\int_x^y (2-\delta)a^{1-\delta} B_{a^\delta}^2 \d a + (1-y^{2-\delta})B^2_{y^\delta} \egaldistr \int_x^y \delta a^{\delta-1} B_{1-a^{2-\delta}}^2 \d a + x^\delta B^2_{1-x^{2-\delta}},\]

d'où on tire :
\[\int_0^{+\infty} \ind{x\leq R^{(2+\delta)}_t\leq y} \d t + \dfrac{(y^{\delta-1}-y)}{2-\delta}L^y_{\infty}(R^{(\delta+2)}) \egaldistr \int_0^{T_1} \ind{x \leq R^{(\delta)}_t\leq y} \d t + \dfrac{x}{\delta} L^x_{T_1}(R^{(\delta)}).\]

%
%
\end{proof}

\begin{exercice}
\label{exo_representationbessel}
Soit $\delta>0$, trouver $c >0$ tel que $\left(L^x_{+\infty}(|B|+cL) \right)_{x \geq 0}$ a pour loi $Q^\delta_0$.
\end{exercice}

\begin{proof}
Soit $\delta >0$ et $c>0$, on pose $U_t = |B_t|+cL_t$ et  $X_x = L^x_{+\infty}(U)$. En appliquant la formule d'Itô-Tanaka à $f(y) = (y\vee (-1))\wedge x$,  on a
\[ f(U_\infty)-f(U_0) = \int_0^{+\infty} \ind{-1<U_s\leq x} \d U_s +\dfrac{1}{2}(X_{-1}-X_x).\]
 On remarque que $U_t >0$ sur un ensemble de probabilité 1, et que $\lim_{t \to +\infty} U_t = +\infty$. Par conséquent $X_{-1} = 0$, et on a
\[X_x = - 2 x + 2\int_0^{+\infty} \ind{0 \leq U_s\leq x} \d U_s.\]

On utilise ensuite l'expression de $U_t$ ; on obtient successivement
\begin{align*}
  X_x & =  - 2x + 2\int_0^{+\infty} \ind{0<U_s \leq x} \d|B|_s + 2 c \int_0^{+\infty} \ind{0<U_s\leq x} \d L_s\\
  & =  - 2 x + 2\int_0^{+\infty} \ind{0<U_s \leq x} \sgn(B_s) \d B_s + 2(c+1)\int_0^{+\infty} \ind{0<U_s\leq x} \d L_s.
\end{align*}
Or, sur le support de $\d L_s$, on a $|B_s|=0$, par conséquent
\begin{align*}
  X_x  & =  -2 x + 2\int_0^{+\infty} \ind{0<U_s \leq x} \sgn(B_s) \d B_s + 2(c+1)\int_0^{+\infty} \ind{0<L_s \leq \frac{x}{c}} \d L_s.\\
  & =  -2x + 2 \int_0^{+\infty} \ind{0<U_s\leq x} \sgn(B_s) \d B_s +  2(c+1)\dfrac{x}{c}.
\end{align*}
On note que $X_x$ est une semi-martingale de crochet $4\int_0^x X_y \d y$, donc $X$ est solution de l'équation différentielle stochastique suivante
\[X_x = \int_0^x 2\sqrt{X_y}\d \beta_y +\frac{2x}{c},\]
et a pour loi $Q^{\frac{2}{c}}_0$. Il faut donc choisir $c=\frac{2}{\delta}$ pour que $X$ soit de loi $Q^\delta_0$.

\begin{remarque}
Ce résultat peut être --et sera-- obtenu de manière plus simple en utilisant la théorie des excursions.
\end{remarque}
\end{proof}

\begin{exercice}
\label{exo_rayk-general}
Enoncer le théorème de Ray-Knight pour les temps locaux de $X_t = B_t +\int_0^t b(X_s) \d s$.
\end{exercice}

\begin{proof}
Ce résultat est une généralisation de l'exercice \ref{exo_RKbang}, on le traitera donc de la même manière. On sait en particulier que $\crochet{X}_t = t$. On étudie pour commencer l'équation différentielle stochastique vérifiée par le processus $(Z_a)_{0\leq a \leq 1}$. Grâce à l'équation \eqref{eqn_Z}, on a :
\begin{align*}
  Z_a & =  2a - 2 \int_0^{T_1} \ind{X_s>1-a} \d X_s\\
  & =  2a - 2 \int_0^{T_1} \ind{X_s>1-a} \d B_s - 2 \int_0^{T_1} b(X_s) \ind{X_s>1-a} \d s\\
  & =  2a - 2 \int_0^a \sqrt{Z_c} \d \beta_c - 2 \int_0^a b(1-c)Z_c \d c,
\end{align*}
où on utilise $\int_0^{T_1} \ind{X_s>1-a} \d B_s = \int_0^a \sqrt{Z_c} \d \beta_c$, grâce au théorème de représentation de martingales, on a déjà vu que $\crochet{\int_0^{T_1} \ind{X_s>1-.} \d B_s}_a = Z_a$. Par conséquent, $Z$ est la solution de l'équation différentielle stochastique :
\[Z_a = 2\int_0^a  \sqrt{Z_c} \d \beta_c +  \int_0^a 2-2 b(1-c) Z_c \d c.\]

De la même manière, on traite le cas de $(Y_x)_{x\geq 0}$ et de $(Y_{-x})_{x \geq 0}$. On sait que ces deux processus sont indépendants car déterminés à partir de processus réfléchis indépendants. On sait de plus que pour $x \geq 0$, on a :
\begin{align*}
  Y_x & =  l + 2 \int_0^{\tau_l} \ind{0<X_s\leq x} \d X_s\\
  & =  l + 2 \int_0^{\tau_l} \ind{0<X_s\leq x} \d B_s + \int_0^{\tau_l} 2b(X_s)\ind{0<X_s\leq x} \d s\\
  & =  l + 2 \int_0^x \sqrt{Y_y}\d \beta_y + \int_0^x 2 b(y) Y_y \d y.
\end{align*}
De la m\^eme façon, on obtient l'équation différentielle stochastique satisfaite par $(Y_{-x})_{x \geq 0}$, que l'on résume en
\[\left\{
\begin{array}{l}
  Y_x = l +\int_0^x 2 \sqrt{Y_y}d \beta_y + \int_0^x 2 b(y) Y_y \d y\\
  Y_{-x} = l +\int_0^x 2 \sqrt{Y_{-y}}d \tilde{\beta}_y + \int_0^x b(-y)Y_{-y}\d y.
\end{array}
\right.\]
\end{proof}

\begin{exercice}
On pose $h(x)=x^\alpha$ et $R^{(d)}$ un processus de Bessel de dimension $d$. Montrer les égalités suivantes entre processus :
\begin{enumerate}
 \item si $d>2(1-\alpha)$, soit $d_\alpha = \frac{d}{\alpha}+2(1-\frac{1}{\alpha})$, on a 
 \[h(R^{(d)}_t) = R^{(d_\alpha)}\left( \int_0^t h'(R^{(d)}_u)^2 \d u \right),\]
 \item si $\alpha>0$, on a 
 \[h(R^{(2+\alpha)}_t) = R^{(3)}\left( \int_0^t h'(R^{(2+\alpha)}_u)^2 \d u \right),\]
 \item si $0<\alpha<2$, on a 
 \[h(R^{(2-\alpha)}_t) = R^{(1)}\left( \int_0^t h'(R^{(2-\alpha)}_u)^2 \d u \right).\]
\end{enumerate}

En déduire les extensions du théorème de Ray-Knight pour les processus de Bessel.
\end{exercice}

\begin{proof}
Plutôt que de raisonner avec $R^{(d)}$, on va étudier le carré de ce processus $X$, de loi $Q^d_l$. On va commencer par montrer que $h(X^{(d)})= h({R^{(d)}}^2)$ est un carré de Bessel changé de temps. La formule d'Itô nous donne
\[X^\alpha_t = l^\alpha + \alpha \int_0^t X^{\alpha-1}_s \d X_s +\dfrac{\alpha(\alpha-1)}{2}\int_0^t X^{\alpha-2}_s \d \crochet{X}_s.\]

On utilise l'équation différentielle satisfaite par $X$ pour réécrire :
\[X^\alpha_t = l^\alpha + 2 \alpha \int_0^t \sqrt{X_s} X^{\alpha-1}_s \d B_s + \alpha( d + 2(\alpha-1)) \int_0^t X^{\alpha-1}_s\d s.\]

On pose alors $\sigma_u = \inf\{t \geq 0 : \int_0^t \alpha^2 X^{\alpha-1}_s \d s >u\}$, on a directement
\[\alpha( d + 2(\alpha-1)) \int_0^{\sigma_u} X^{\alpha-1}_s\d s = \left(\frac{d}{\alpha}+2\left(1-\frac{1}{\alpha}\right)\right)u=: d_\alpha u.\]
Il suffit maintenant d'observer que $M_t = 2 \alpha \int_0^t \sqrt{X_s} X^{\alpha-1}_s \d B_s$ est une martingale de crochet $4 \alpha^2 \int_0^t X^{2\alpha-1}_s\d s$, par conséquent, en changeant le temps $t$ en $\sigma_u$ une famille continue de temps d'arrêts, on obtient encore une martingale (locale) $M_{\sigma_u}$, de crochet
\[ 4  \int_0^{\sigma_u} X^\alpha_s \alpha^2 X^{\alpha-1}_s \d s = 4 \int_0^u X^\alpha_{\sigma_v}\d v.\]

Par conséquent, $X^\alpha_{\sigma_u}$ est solution de l'équation différentielle stochastique
\[Y_u = l^\alpha + \int_0^u 2\sqrt{Y_v} \d \beta_v + d_\alpha u,\]
donc est un carré de Bessel de dimension $d_\alpha$ issu de $l^\alpha$. Par conséquent, en prenant la racine carrée de cette égalité, on a :
\[h(R^{(d)}_{\sigma_u}) = R^{(d_\alpha)}_u,\]
d'où l'expression souhaitée du théorème en utilisant que $\sigma_u$ est la fonction inverse de $t \mapsto \int_0^t h'(R^{(d)}_s)^2 \d s$.

On s'intéresse maintenant au cas $\alpha>0$ et $d = 2+\alpha$, on a alors $d_\alpha = 3$, donc $X^\alpha_{\sigma_u}$ est solution de l'équation différentielle stochastique :
\[Y_u = l^\alpha + \int_0^u 2\sqrt{Y_v} \d \beta_v + 3 u,\]
dès lors $X^\alpha_{\sigma_u}$ est un carré de Bessel de dimension 3, par conséquent :
\[h(R^{(2+\alpha)}_t) = R^{(3)}\left( \int_0^t h'(R^{(2+\alpha)}_u)^2 \d u \right).\]

Pour finir, on étudie le cas où $0<\alpha<2$ et $d=2-\alpha$. Dans ce cas, $d_\alpha = 1$ et $X^\alpha_{\sigma_u}$ est alors solution de l'équation différentielle :
\[Y_u = l^\alpha + \int_0^u 2\sqrt{Y_v}\d \beta_v + t,\]
donc est un mouvement brownien réfléchi en dimension 1 issu de $l^\alpha$.
\end{proof}

\begin{exercice}
Soit $R$ un processus de Bessel de dimension 2, on pose $M_s=R_s^2-2s$.

Notons $T^*_t = \inf\{t \geq 0 : |B_s|=t\}$, montrer que pour tout $t \geq 0$ on a
\[T^*_t \egaldistr \int_0^t \d s R^2_s.\]

En déduire que $M'_t=T^*_t-t^2$ est une martingale vérifiant pour tout $t \geq 0$ :
\[M'_t \egaldistr \int_0^t \d s M_s.\]
\end{exercice}

\begin{proof}
On exprime un théorème de Ray-Knight pour $(Z_a = L^{t-a}_{T^*_t}(|B|),0\leq a \leq t)$, de façon analogue au cas du mouvement brownien, on a
\[Z_a = 2 a - 2 \int_0^{T^*_t} \ind{|B_s|>t-a} \d|B|_s = 2a - 2 \int_0^{T^*_t} \sgn(B_s)\ind{|B_s|>t-a} \d B_s,\]
donc $Z$ est solution de l'équation différentielle stochastique
\[Z_a = 2a + 2\int_0^a \sqrt{Z_s}\d \beta_s,\]
c'est un carré de Bessel de dimension 2 issu de 0.

On en déduit l'égalité en loi suivante
\[\int_0^t \d s R_s^2 \egaldistr \int_0^t \d a L^{t-a}_{T^*_t}(|B|) = \int_0^{T^*_t} \d s \ind{|B|_s \in [0,t]} = T^*_t.\]

Cette égalité en loi entraîne en particulier
\[M'_t \egaldistr \int_0^t \d s M_s,\]
il reste néanmoins à prouver que $M'_t$ est bien une martingale.

On utilise le fait que $(B^2_t-t)$ est une martingale pour calculer :
\begin{align*}
  \E((t+s)^2-T^*_{t+s}|\F_{T^*_t}) = \E\left[\left. B_{T^*_{t+s}}^2 - T^*_{t+s} \right| \F_{T^*_t} \right] = B^2_{T^*_t}-T^*_t = t^2-T^*_t,
\end{align*}
on a donc bien le résultat escompté.
\end{proof}

\begin{exercice}
Montrer que pour toute martingale locale continue $M$, positive issue de 1 tendant vers 0 p.s. en $+\infty$, le processus $(L^a_\infty(M), a \geq 0)$ admet la même loi, indépendante du choix de $M$.
\end{exercice}

\begin{proof}
On utilise la représentation de Dubins-Schwarz d'une martingale locale, on a
$M_t = \beta_{\crochet{M}_t},$
avec $\beta$ un mouvement brownien issu de 1 et $\crochet{M}_\infty = T_0(\beta)$.

On a alors pour toute fonction $f$ positive mesurable, par changement de variables :
\[\int_{\R^+} \d_a L^a_\infty f(a) = \int_0^{+\infty} f(M_s) \d \crochet{M}_s = \int_0^{T_0(\beta)} f(\beta_u) \d u = \int_{\R^+} \d_a L^a_{T_0}(\beta) f(a).\]

Par conséquent, 
\[(L^a_\infty(M), a \geq 0) \egaldistr (L^a_{T_1}(1-B), a \geq 0) = (L^{1-a}_{T_1}(B), a \geq 0),\]
on utilise alors le théorème de Ray-Knight pour conclure. On a :
\[(B_{T_1}-(1-a))^+ - (a-1)^+ = \int_0^{T_1} \ind{B_s>1-a}\d B_s + \frac{1}{2} L^{1-a}_{T_1}.\]

Par conséquent, $L^{1-a}_{T_1}$ est solution de l'équation différentielle stochastique :
\[Z_a = 2(a^+ - (a-1)^+) + \int_0^a 2 \sqrt{Z_b}\d \beta_b,\]
donc $(Z_a)$ suit la loi d'un carré de Bessel de dimension 2 issu de 0 sur $[0,1]$, puis celle d'un carré de Bessel de dimension 0 issu de $Z_1$ ensuite.
\end{proof}

\begin{exercice}[Processus de populations]
On dit qu'un processus de Markov vérifie la propriété de branchement lorsque, pour tout couple $x$ et $y$, si $(X^x_t)$ et $(\widehat{X}^y_t)$ sont deux versions indépendantes du processus de Markov respectivement issues de $x$ et de $y$, alors $(X^x_t+ \widehat{X}^y_t)$ suit la loi du processus issu de $x+y$.

\begin{enumerate}
  \item Caractériser les solutions d'équations différentielles stochastiques à coefficients continus vérifiant la propriété de branchement.
  \item En utilisant un théorème de Ray-Knight, identifier ces processus comme les temps locaux, en un temps d'arr\^et convenablement choisis d'un processus stochastique.
\end{enumerate}
\end{exercice}

\begin{proof}
\textit{1.} Soit $B$ et $\widehat{B}$ deux mouvements browniens indépendants. On note $X^x_t$ la solution de l'équation différentielle stochastique :
\[X^x_t = x + \int_0^t \sigma(X_s)\d B_s + \int_0^t b(X_s)\d s \]
et $\widehat{X}^y_t$ celle de :
\[\widehat{X}^y_t = y + \int_0^t \sigma(\widehat{X}^y_s)\d \widehat{B}_s + \int_0^t b(\widehat{X}^y_s)\d s.\]

Si $b$ et $\sigma$ sont telles que le processus vérifie la propriété de branchement, $X^x_t+\widehat{X}^y_t$ est solution de l'équation différentielle stochastique :
\[Y_t = x + y + \int_0^t \sigma(Y_s) \d B_s + \int_0^t b(Y_s) \d s.\]

Par identification des parties à variation finies, on a :
\[\int_0^t b(X^x_s+\widehat{X}^y_s)\d s = \int_0^t b(X^x_s) + b(\widehat{X}^y_s) \d s,\]
donc pour tout $x,y$, on a $b(x+y)=b(x)+b(y)$. On en déduit $b(x) = \lambda x$.

De la m\^eme manière, on a par identification des martingales continues :
\[\sigma^2(x+y) = \sigma^2(x)+\sigma^2(y),\]
d'où $\sigma(x) = \delta \sqrt{x}$.

On en déduit que les diffusions vérifiant la propriété de branchement sont solutions d'équations différentielles stochastiques du type :
\[ X_t = X_0 + \int_0^t \delta \sqrt{X_s}\d B_s + \lambda \int_0^t X_s \d s.\]

\textit{2.} Soit $Y$ la solution de l'équation différentielle stochastique :
\[Y_t = \lambda B_t +  \int_0^t f(Y_s)\d s.\]
Grâce à l'Exercice \ref{exo_rayk-general}, on sait que $(U_x = L^x_{\tau_l}(X), x \geq 0)$ est solution de l'équation différentielle stochastique :
\[U_x = l + \int_0^x 2\lambda \sqrt{U_y}\d \beta_y + \int_0^x 2f(y)U_y\d y.\]

Par conséquent, un processus de branchement issu de $l$ peut être vu comme les temps locaux en $\tau_l$ d'un mouvement brownien avec dérive. Ainsi le processus associé à $b(x)=\lambda x$ et $\sigma(x) = \delta \sqrt{x}$ peut être vu comme les temps locaux en $\tau_l$ de $\frac{\delta}{2} \beta_t+ \frac{\lambda}{2}t$.
\end{proof}

\chapter[Temps locaux d'intersection]{Introduction aux temps locaux d'intersection du mouvement brownien en dimension 2 et 3}

Pendant la décennie 1950-1960, 4 articles fondamentaux sur les points multiples du mouvement brownien $d$-dimensionnel ont été publiés, par A. Dvoretzky, P. Erdös et S. Kakutani pour les trois premiers, et S. Taylor pour le $4^\text{e}$ ; ils ont ainsi montré que :
\begin{itemize}
 \item si $d=2$, la trajectoire brownienne possède des points multiples de tous ordres ;
 \item si $d=3$, la trajectoire brownienne possède des points doubles, mais pas de points triples ;
 \item pour tout $d\geq 4$, la trajectoire brownienne ne présente pas de points doubles.
\end{itemize}

Ces résultats ont été obtenus à l'aide seulement des propriétés fondamentales d'accroissements indépendants et de stationnarité du mouvement brownien. Il a fallu attendre le début des années 80 (Wolpert, Rosen) pour observer les premières constructions de temps locaux d'intersection, permettant une compréhension plus fine des points multiples de ces courbes.

Il apparaît à la lecture des articles correspondants que ces constructions des temps locaux d'intersection, permettant de mesurer le temps que la courbe brownienne passe à proximité des points déjà visités par elle-même, n'a été permise que grâce au mûrissement de la compréhension des propriétés essentielles des temps locaux du mouvement brownien linéaire, des diffusions, des semi-martingales dans les années 1970. Nous nous intéresserons principalement ici aux versions càdlàg de ces temps locaux d'intersection.

\section{Théorèmes principaux}

La théorie des temps locaux que nous avons développée jusqu'ici ne s'intéresse qu'au temps passé par un processus au voisinage d'un niveau. Mais il existe de nombreuses autres notions de temps locaux, particulièrement les temps locaux d'intersection de la courbe brownienne en dimension 2 et 3, qui comptent le temps passé par la courbe à se recouper sur elle-même. Nous ne donnerons ici qu'une très rapide introduction à ce domaine important.

\begin{theoreme}
Soit $B$ un mouvement brownien $d$-dimensionnel, pour $d \leq 3$ il existe une famille mesurable $(\alpha(y,t) ; y \in \R^d, t \geq 0\}$ qui permet d'écrire, pour toute fonction $f$ continue bornée :
\[ \int_0^t \d s \int_s^t \d u f(B_u-B_s) = \int_{\R^d} \d y f(y) \alpha(y,t).\]
\end{theoreme}

Nous allons maintenant introduire une expression pour $\alpha$, l'équivalent de la formule de Tanaka-Meyer qui donnait l'expression des temps locaux de semi-martingale au Chapitre 1. Cette formule dépend de la dimension du mouvement brownien considéré. On se placera pour commencer en dimension $d=2$.

\subsection{Formule de Tanaka-Rosen en dimension 2}

\begin{theoreme}[Formule de Tanaka-Rosen]
Pour tout $y \neq 0$ et $t \geq 0$, il existe $a(y,t)$ vérifiant l'équation suivante :
\[\int_0^t \d s \left[ \log|B_t-B_s-y| - \log |y| \right] = \int_0^t \d B_u .\int_0^u \d s \dfrac{B_u-B_s-y}{|B_u-B_s-y|^2} + \pi a(y,t).\]

De plus, le processus $(a^y_t ; y \in {\R^2}^*,t \geq 0)$ est une version continue de la famille $(\alpha(y,t), y \in {\R^2}^*, t \geq 0)$ en dimension 2 définie dans le Théorème 7.1. 
\end{theoreme}

Ces temps locaux ne sont pas définis en zéro. Cependant, on connait le comportement des temps locaux de $B$ au voisinage de $0$ --c'est-à-dire le temps local d'intersection proprement dit-- grâce au corollaire suivant.

\begin{corollaire}
\label{cor:renormalisation}
On pose $\gamma(y,t) = a(y,t)-t \log \frac{1}{|y|}$.

Dans ce cas, $\gamma(y,t)$ converge dans $L^2$ quand $y \to 0$ vers une limite notée $\gamma(0,t)$. De plus la famille $(\gamma(y,t), y \in \R^2, t \geq 0)$ admet une version continue.
\end{corollaire}

\begin{remarque}
En dimension 2, la notion de temps locaux d'intersection s'étend de la façon suivante : il existe des temps locaux de tous ordres $k$ vérifiant :
\[ \int_{0 \leq s_1 \leq \cdots \leq s_k} \d s_1 \cdots \d s_k f(B_{s_2}-B_{s_1}, \ldots, B_{s_k}-B_{s_{k-1}}) = \int \d z A^{(k)}(t,z) f(z). \]

De plus, $A^{(k)}(t,z)$, convenablement recentrée, converge p.s. et dans $L^2$ quand $z \to 0$, on a donc bien une certaine notion de temps local d'intersection $k^\text{ième}$.
\end{remarque}

On s'intéresse maintenant au mouvement brownien dans l'espace.

\subsection{Formule de Tanaka-Rosen en dimension 3}

\begin{theoreme}[Formule de Tanaka-Rosen]
Pour tout $y \neq 0$ et $t \geq 0$, il existe $a(y,t)$ vérifiant l'équation :
\[ \int_0^t \d s \left[ \dfrac{1}{|B_t-B_s-y|} - \dfrac{1}{|y|} \right] = - \int_0^t \d B_u. \int_0^u \d s \dfrac{B_u-B_s-y}{|B_u-B_s-y|^3} - 2 \pi \alpha(y,t).\]

De plus, $(a(y,t), y \in{\R^3}^*,t \geq 0)$ est une version continue de la famille $(\alpha^y_t ; y \in {\R^3}^*,t \geq 0)$ en dimension 3. 
\end{theoreme}

Le résultat suivant nous donne une version faible de la valeur du temps local d'intersection du mouvement brownien plan.

\begin{theoreme}
La suite de processus 
\[\left[ B_t, \dfrac{1}{\sqrt{- \log |y|}}\left( 2 \pi a(y,t)-\dfrac{t}{|y|} \right);t \geq 0\right]\]
converge en loi, quand $y \to 0$ vers $(B_t, 2 \beta_t ; t \geq 0)$, avec $\beta$ un mouvement brownien indépendant de $B$.
\end{theoreme}

\section{Exercices}

\begin{exercice}
\begin{enumerate}
  \item Montrer l'existence des temps locaux d'intersection en dimension 2 et 3 en utilisant la transformée de Fourier.
  \item Étudier l'existence de temps locaux d'intersections de la courbe brownienne avec un passé distant d'au moins $\epsilon>0$.
\end{enumerate}
\end{exercice}

\begin{proof}
\textit{1.} On définit la transformée de Fourier de la mesure aléatoire $\mu_t$ par :
\[\widehat{\mu_t}(\xi) = \int_0^t\d u \exp(i \xi.B_u) \int_u^t \d s \exp(-i\xi. B_s).\]

Il suffit de vérifier que cette transformée de Fourier est dans $L^2$ presque sûrement pour conclure à l'existence d'une densité pour cette mesure, par formule de Parseval. Or, la formule d'Itô nous donne
\[\exp(i \xi . B_t) = 1 + i \int_0^t \exp(i \xi.B_s) \d(\xi.B_s)- \dfrac{|\xi|^2}{2}\int_0^t \d s \exp(i \xi.B_s),\]
par conséquent, on a
\[ \left| \int_0^t \d s \exp(i \xi.B_s)\right| \leq C\left(\dfrac{1}{|\xi|^2} + \dfrac{1 }{|\xi|} I_t(\xi) \right),\]
où on a posé $I_t(\xi) =  \int_0^t \exp(i \xi.B_s) \d(\frac{\xi}{|\xi|}.B_s)$, qui est une intégrale relative à un mouvement brownien.

On en déduit ainsi 
\[\E\left[ \left(\int_0^t \d s \exp(i \xi.B_s)\right)^2 \right] \leq \dfrac{C}{|\xi|^2} \wedge t^2,\]
on obtient donc sans difficulté
\[\E \left[ \widehat{ \mu_t}(\xi)^2 \right] \leq \dfrac{C}{|\xi|^4}\wedge t^4 \in L^2,\]
dès lors que $d \leq 3$. Par conséquent, $\widehat{\mu_t}$ est dans $L^2$ p.s. par Fubini, ce qui nous permet de conclure.

\textit{2.} On peut de la même manière s'intéresser à la transformée de Fourier de la mesure aléatoire $\mu^\epsilon_t$ définie par :
\[ \int f \d \mu^\epsilon_t = \int_0^{t-\epsilon} \d u \int_{u+\epsilon}^t \d s f(B_u-B_s). \]
Sa transformée de Fourier s'exprime de la façon suivante :
\[\widehat{\mu^\epsilon_t}(\xi) = \int_0^{t-\epsilon} \d u \exp(i \xi.B_u) \int_{s+\epsilon}^t \d s \exp(-i\xi. B_s).\]

On en déduit donc, de la même façon que précédemment, que 
\[\E \left[ \widehat{ \mu_t^\epsilon}(\xi)^2 \right] \leq \dfrac{C}{|\xi|^4}\wedge t^4 \in L^2,\]
dès lors que $d \leq 3$. On en déduit l'existence des temps locaux d'intersections de la courbe brownienne avec son passé \og lointain \fg.
\end{proof}


\begin{exercice}
Montrer la formule de Tanaka-Rosen, en dimension 2 et 3.
\end{exercice}

\begin{proof}
On se place pour commencer en dimension $d=2$. Soit $f$ fonction de classe $\mathcal{C}^\infty$ à support compact, on définit la fonction :
\[F(x) = \int_{\R^2} \d y \log |y-x| f(y),\]
solution de l'équation différentielle $\Delta F = 2\pi f$. En appliquant la formule d'Itô entre les instants $s$ et $t$, on a :
\[F(B_t-B_s) = F(0) + \int_s^t \d B_u. \nabla F(B_u-B_s) + \pi \int_s^t f(B_u-B_s)\d u.\]
En intégrant ensuite par rapport à $s$, on obtient :
\[\pi \int_0^t \d s \int_s^t \d u f(B_u-B_s) = \int_0^t \d B_u . \int_0^u \d s \nabla F(B_u-B_s) - \int_0^t \d s F(B_t-B_s)-F(0).\]

On utilise alors la formule de Fubini pour obtenir
\begin{multline*}
   \pi \int_0^t \d s \int_s^t \d u f(B_u-B_s)\\
  = \int \d y f(y)\left[\int_0^t \d s \log |B_t-B_s-y|-\log |y| + \int_0^t \d B_u. \int_0^u \d s \dfrac{B_u-B_s-y}{|B_u-B_s-y|^2}\right],
\end{multline*}
d'où la formule annoncée en dimension 2, par identification.

Le résultat s'obtient de la même manière en dimension 3 en utilisant la solution de $\Delta u = 4\pi \delta_0$ donnée par :
\[u(x) = \frac{1}{|x|}.\]
\end{proof}

\begin{exercice}[Renormalisation de Varadhan]
Soit $f : \R^2 \to \R^+$ continue à support compact vérifiant $\int \d x f(x) = 1$. On pose $f_n(x) = n^2 f(nx)$ et 
\[Z_n = \int_0^t \d s \int_s^t \d u f_n(B_u-B_s).\]
\begin{enumerate}
  \item Soit $y \in {\R^2}^*$, calculer $\E(\alpha(y,t))$
  \item Montrer que $(Z_n-\E(Z_n))$ converge, quand $n \to +\infty$ vers $\bar{\gamma}(t)$ défini par
\[\bar{\gamma}(t) = \lim_{y \to 0} \alpha(y,t)-\E(\alpha(y,t)).\]
\end{enumerate}
\end{exercice}

\begin{proof}
\textit{1.} Soit $f$ une fonction positive mesurable, on a par formule de Fubini :
\[\E \left[ \int_0^t \d s \int_s^t \d u f(B_u-B_s) \right] = \int_{\R^2} \d y f(y) \E\left( \alpha(y,t) \right).\]
Or, en utilisant le fait que $B_u-B_s$ est un vecteur aléatoire gaussien de dimension 2, de variance $u-s$, on a
\[ \E \left[ \int_0^t \d s \int_s^t \d u f(B_u-B_s) \right] = \int_{\R^2} \d y f(y)\int_0^t \d s \int_s^t \dfrac{\d u}{2\pi (u-s)} \exp\left(-\frac{|y|^2}{2(u-s)}\right). \]
On en déduit
\[\E(\alpha(y,t)) = \int_0^t \dfrac{(t-u)}{2\pi u} \exp\left( - \frac{|y|^2}{2u} \right)\d u.\]

\begin{remarque}
L'équivalent au voisinage de 0 de $\E(\alpha(y,t))$ est bien $t \log \frac{1}{|y|}$ comme dans le Corollaire \ref{cor:renormalisation}.
\end{remarque}

\textit{2.} On réécrit :
\[\int_0^t \d s \int_s^t \d u f_n(B_u-B_s) = \int \d y f(y) \alpha\left(\frac{y}{n},t\right).\]
On a alors, en recentrant les variables aléatoires :
\[\alpha(y,t)-\E(\alpha(y,t)) = \gamma(y,t)-\E(\gamma(y,t)) \underset{y \to 0}{\longrightarrow} \gamma(0,t)-\E(\gamma(0,t)),\]
qui existe donc bien, par Corollaire \ref{cor:renormalisation}.

On peut alors bien passer à la limite pour les variables aléatoires recentrées, on obtient :
\[Z_n-\E(Z_n) \conv \int \d y f(y) \bar{\gamma}(0,t) = \bar{\gamma}(0,t).\]
\end{proof}

\part{Excursions browniennes}

Cette partie II est composée de six chapitres, chacun d'entre eux correspondant à une leçon sur les excursions du mouvement brownien réel, que nous étendrons ensuite aux diffusions réelles. Au début de chaque chapitre, nous rappelons les théorèmes principaux de la leçon correspondante.

Nous nous attachons ici à introduire la théorie des excursions du mouvement brownien, laquelle remplace ce processus par le processus de sauts de ses excursions. En d'autres termes, on considère que le mouvement brownien est une trajectoire, que nous décomposons selon les \og{}arches de ponts\fg{} qu'elle réalise hors de 0. Cette description du processus permet de nombreux calculs, et simplifie l'expression de certaines quantités. C'est un regard entièrement différent de celui porté par le calcul stochastique, qui s'avère ainsi être complété et renforcé. Nous verrons de cette façon que de nombreux résultats obtenus dans la première partie peuvent être retrouvés en utilisant cette théorie.

Nous nous intéresserons en particulier à donner la loi de certaines de ces \og{}arches de pont\fg{}, considérées comme des processus, ce qui nous permettra d'expliciter la mesure de Lévy du processus de ces \og{}excursions hors de zéro\fg{} du mouvement brownien ; on la décomposera selon la longueur (ou la hauteur) des excursions, une fois que nous aurons donné un sens à tout ceci. Tous ces résultats nous permettront de donner une preuve de la formule de Feynman-Kac ainsi que de nombreux autres résultats.

Pour finir, nous étudierons une courte extension aux diffusions linéaires des résultats que nous aurons obtenus pour le mouvement brownien. La théorie des temps locaux et celle des excursions s'adaptent en effet très bien à ces processus, de plusieurs manières différentes (selon qu'on les regarde comme des semi-martingales ou des processus de Markov), ce qu'il sera intéressant d'explorer.

\chapter[Introduction à la théorie des excursions]{Introduction à la théorie des excursions}

Certains des travaux de Paul Lévy relatifs aux \og excursions \fg{}, permettent d'étudier sous un jour nouveau le mouvement brownien. Dans cette optique, il est vu comme un processus de Lévy à valeurs dans un espace de trajectoires, un \og saut \fg{} de ce processus de Lévy correspondant à l'allongement de la trajectoire du mouvement brownien entre deux zéros successifs. La nouvelle échelle de temps qu'il convient de considérer ici est alors l'inverse du temps local au niveau 0. Pour l'heure, nous traduirons ceci, de façon assez libre, comme l'étude du processus $l \mapsto (B_t)_{0 \leq t \leq \tau_l}$.

Ce dernier peut être considéré comme un processus de Lévy, lorsque l'addition est remplacée par la concaténation de trajectoires. En effet, il vérifie les deux propriétés clé des processus de Lévy, l'indépendance entre la trajectoire du mouvement brownien jusqu'en l'instant $\tau_l$ et celle de son \og accroissement \fg{} entre les instants $\tau_l$ et $\tau_{l +l'}$ ; ainsi que le caractère stationnaire de la loi de cet accroissement (la trajectoire entre $\tau_l$ et $\tau_{l+l'}$ a même loi que celle entre 0 et $\tau_{l'}$). Un des objectifs de la théorie des excursions est de déterminer la mesure de Lévy de ce processus, que l'on appeller mesure d'Itô, souvent notée $\n$. Une fois cette mesure $\n$ bien comprise, le lecteur sera capable de traiter, au moyen cette fois-ci de la théorie des excursions, la plupart des questions abordées dans la partie I. Grâce à un mélange de théorie des excursions et de théorie des temps locaux on pourra même obtenir de façon assez simple des résultats à première vue compliqués.

\section{Théorèmes principaux}

L'ensemble des trajectoires continues issues de $0$ et absorbées en 0 est noté $\Omega^*$. Pour toute fonction $f \in \Omega^*$, on note $V(f) = \inf\{t>0:f(t)=0\}$. En particulier, si $f \in \Omega^*$, pour tout $s \geq V(f)$ on a $f(s)=0$. On appelle $f$ une excursion, de longueur $V(f)$. On remarque en particulier que $\Omega^*$ peut se décomposer en deux sous-ensembles symétriques, que nous noterons $\Omega^+$ pour les trajectoires positives et $\Omega^-$ pour les négatives.

Pour tout $l \geq 0$, on note $\tau_l = \inf\{t \geq 0 : L_t > l\}$, et $\tau_{l-} = \inf\{t \geq 0 : L_t \geq l\}$. Les excursions de $B$ sont notées $e_l = \left(B_{(\tau_{l-}+s)\wedge \tau_l}\right)_{s \geq 0}$.

\begin{definition}
Soit $U,\mathcal{U}$ un ensemble mesuré. Un processus de Poisson ponctuel sur $U$ de mesure d'intensité $\mu$ est un processus $(e_t,t \geq 0)$ à valeurs dans $U \cup \{\partial\}$, tel que pour tout choix d'ensembles mesurables disjoints $\Gamma_1, \ldots, \Gamma_k$ de $U$, les processus 
\[N^{\Gamma_j}_l = \sum_{s \leq t} \ind{e_s \in \Gamma_j}\]
sont soit infinis, soit des processus de Poisson indépendants de paramètre $\mu(\Gamma_j)$.
\end{definition}

\begin{theoreme}[Théorème d'Itô]
Le processus $(e_l)_{l \geq 0}$ est un processus de Poisson ponctuel de mesure d'intensité $\n$. $\n$ est une mesure sigma-finie sur $\Omega^*$ appelée la mesure d'Itô du mouvement brownien.

Le processus de Poisson pour la concaténation $X_l = (B_t)_{t \leq \tau_l}$ a pour mesure de Lévy $\n$.
\end{theoreme}

On présente maintenant deux formules clés pour le calcul de certaines fonctionnelles de processus de Lévy.

\begin{proposition}[Formule additive]
Pour tout $H : \R^+ \times \Omega \times \Omega^* \to \R^+$, processus prévisible par rapport à $(\F_{\tau_l} \otimes \mathcal{E})_{l \geq 0}$, où $\mathcal{E}$ est la tribu Borélienne associée à $\Omega^*$, on a :
\[\E\left(\sum_{\lambda \leq l} H(\lambda, \omega ; e_\lambda) \ind{\tau_\lambda >\tau_{\lambda-}} \right) = \E\left[\int_0^l \d \lambda \int_{\Omega^*} \n(\d \epsilon) H(\lambda, \omega ; \epsilon)\right].\]
\end{proposition}

\begin{remarque}
Une autre version de cette formule clé, pour laquelle les excursions sont indexées par l'instant de début, et non le temps local en cet instant, est donnée en exercice.
\end{remarque}

\begin{proposition}[Formule multiplicative]
Pour toute fonction positive mesurable $f :\R^+ \times \Omega^* \to \R^+$, on a :
\[\E\left[ \exp\left( - \sum_{\lambda \leq l} f(\lambda, e_\lambda) \ind{\tau_\lambda>\tau_{\lambda-}} \right) \right] = \exp\left( -\int_0^l \d \lambda \int_{\Omega^*} \n(\d \epsilon)  \left(1 - e^{-f(\lambda, \epsilon)}\right) \right)\]
\end{proposition}

\begin{remarque}
Cette proposition peut s'étendre à des fonctionnelles positives, mais possiblement infinies, en utilisant le théorème de convergence monotone.
\end{remarque}

\section{Exercices}

\begin{exercice}
Montrer, en utilisant un processus de Poisson réel que la formule multiplicative ne peut pas s'étendre avec des fonctionnelles aléatoires, i.e. que les expressions 
\[\E\left[\exp\left( -\sum_{s \leq t} f(s, \omega ; N_s-N_{s-}) \right) \right] \text{  et  } \E\left[\exp\left( -\int_0^t \d s \left(1-e^{-f(s, \omega ; 1)}\right)\right) \right]\]
sont en général différentes.
\end{exercice}

\begin{proof}
Il suffit d'utiliser l'identité suivante $\sum_{s \leq t, \Delta N_s = 1} \ind{N_s \geq 1}  = N_t$, on a alors :
\[\E\left[\exp\left( -\sum_{s \leq t, \Delta N_s = 1} \ind{N_s \geq 1}\right)\right] = \E\left[e^{-N_t} \right] = \exp(-t(1-e^{-1})).\]

D'autre part, on a
\[\E \left[ \exp \left( - \int_0^t \d s  \left(1 - e^{-\ind{N_s \geq 1}}\right) \right) \right] = \E\left[ \exp\left( - (1-e^{-1}) (t-T_1)^+\right)\right],\]
et cette quantité est strictement plus petite que l'autre, car $T_1>0$ p.s.
\end{proof}

\begin{exercice}[Formule additive modifiée]
\label{exo_formuleclemodifiee}
Soit $F : \R^+ \times \Omega \times \Omega^* \to \R$ processus prévisible par rapport à la filtration $\left(\F_t \times \mathcal{E}\right)_{t \geq 0}$. On note $G = \{ \tau_{l-}, l \geq 0\}$ l'ensemble des débuts d'excursions du mouvement brownien et $i_s=e_\lambda$ l'excursion débutant à l'instant $s=\tau_{\lambda-}$.
\begin{enumerate}
\item Montrer que :
\[\E \left[ \sum_{s \leq t, s \in G} F(s, \omega ; i_s) \right] =  \E\left[ \int_0^t \d L_s \int_{\Omega^*} \n(\d \epsilon) F(s, \omega ; \epsilon) \right].\]
\item Soit $t \geq 0$, et $\Gamma \subset \Omega^*$ tel que $\n(\Gamma)<+\infty$, on pose 
\[C^\Gamma_t = \sum_{\lambda \geq 0 : \tau_{\lambda-} \leq t} \ind{e_\lambda \in \Gamma},\]
montrer que $C^\Gamma_t- \n(\Gamma)L_t$ est une martingale.
\end{enumerate}
\end{exercice}

\begin{proof}
\textit{1.} On réécrit la somme étudiée de façon à pouvoir utiliser la formule clé additive, on a
\begin{align*}
  \E \left[ \sum_{s \leq t, s \in G} F(s, \omega ; i_s) \right]
   = & \E\left[ \sum_{\lambda \geq 0 : \tau_\lambda \leq t} F(\tau_{\lambda-}, \omega ; e_\lambda)  \right]\\
   = & \E \left[ \int_0^{+\infty} \d \lambda \ind{\lambda \leq L_t} F(\tau_{\lambda-}, \omega ; e_\lambda) \right]\\
   = & \E \left[ \int_0^t \d_sL_s F(s, \omega ; i_s) \right]
\end{align*}
par changement de variable $\lambda = L_s$.

\textit{2.} On réécrit $C^\Gamma$ de la façon suivante
\[C^\Gamma_t = \sum_{s \in G, s < t} \ind{i_s \in \Gamma},\]
par conséquent, $\E(C^\Gamma_t- \n(\Gamma)L_t) = 0$. On utilise alors la propriété de Markov forte en $d_t = \inf\{s \geq t : B_s = 0\}$. On observe que $C^\Gamma_{t+s} - C^\Gamma_t$ est indépendant et de même loi que $C^\Gamma_s$. Cette propriété étant partagée par $L$, on en déduit que $(C^\Gamma_t - \n(\Gamma) L_t$ est une martingale, et plus particulièrement un processus de Lévy.
\end{proof}

\begin{exercice}
Soit $B$ un mouvement brownien, et $e_\lambda(t)=B_{(t+\tau_{\lambda-})\wedge \tau_\lambda}$.
\begin{enumerate}
  \item Montrer que le processus défini par 
\[ X_t = |e_\lambda(t-\tau_{\lambda-})|-\lambda, \text{  pour  } \tau_{\lambda-} \leq t \leq \tau_\lambda,\]
est le mouvement brownien $\int_0^t \sgn(B_s) \d B_s$.
  \item De la même manière montrer que le processus
\[ Y_t = |e_\lambda(t-\tau_{\lambda-})|+\lambda, \text{  pour  } \tau_{\lambda-} \leq t \leq \tau_\lambda,\]
est un processus de Bessel de dimension 3.
\end{enumerate}
\end{exercice}

\begin{proof}
\textit{1.} Rappelons dans un premier temps la formule de Tanaka :
\[|B_t| = \int_0^t \sgn{B_s} \d B_s + L_t,\]
donc pour $0 \leq t \leq \tau_\lambda - \tau_{\lambda-}$, on a :
\[|e_\lambda(t)| = \int_0^{\tau_{\lambda-} + t} \sgn(B_s) \d B_s + \lambda,\]
on en déduit $\int_0^t \sgn(B_s) \d B_s = X_t$.

\textit{2.} Le deuxième résultat s'obtient de la même manière, on observe que $Y_t=|B_t|+L_t$, donc par théorèmes de Lévy puis de Pitman, ce processus est un processus de Bessel de dimension 3.
\end{proof}

\begin{exercice}[\og Loi \fg{} de la longueur d'une excursion]
\label{exo_integrationlongueur}

Soit $\epsilon \in \Omega^*$, on rappelle que $V(\epsilon) = \inf\{s >0 : \epsilon(s)=0\}$ est la longueur de l'excursion, déterminer la mesure $n_V$ vérifiant pour toute fonction mesurable $f:\R^+ \to \R$ :
\[\int_{\Omega^*} \n(\d \epsilon) f(V(\epsilon))  = \int_{\R^+} n_V(\d v) f(v).\]
\end{exercice}

\begin{proof}
Soit $\mu \geq 0$ donné, on utilise la formule multiplicative avec la fonction $f(\lambda,\epsilon) = \mu V(\epsilon)$, on a
\[\E\left[\exp\left(-\mu \sum_{\lambda \leq l} V(e_\lambda)\right) \right] = \exp \left(-l \int_{\Omega^*} \left(1-e^{-\mu V(\epsilon)}\right) \n(\d \epsilon)\right).\]

D'autre part, $\sum_{\lambda \leq l} V(e_\lambda) = \tau_l$, et on connaît la transformée de Laplace de $\tau_l$ :
\[\E(\exp(-\mu \tau_l)) = \exp(-l \sqrt{2 \mu}).\]
Dès lors, par identification
\[\sqrt{2\mu} = \int_{\Omega^*} \n(\d \epsilon) \left(1-e^{-\mu V(\epsilon)}\right) = \int_{\R^+} n_V(\d v) \left( 1-e^{-\mu v} \right).\]

Identifions maintenant la mesure $n_V$. En dérivant par rapport à $\mu$, il vient :
\[\int_{\R^+} n_V(\d v) v e^{-\mu v} = \dfrac{1}{\sqrt{2\mu}}\]
Or, par définition de la fonction $\Gamma$ d'Euler et un changement de variable adapté, on a :
\begin{align*}
  \dfrac{1}{\sqrt{2\mu}} &= \dfrac{1}{\sqrt{2}\Gamma(\frac{1}{2})}\int_0^{+\infty} h^{-\frac{1}{2}} e^{-2\mu h}dh\\
  &= \dfrac{1}{\sqrt{2\pi}} \int_0^{+\infty} e^{-\mu v} \dfrac{\d v}{\sqrt{v}}.
\end{align*}
Par injectivité de la transformée de Laplace, on en déduit
\[n_V(\d v) = \dfrac{\d v}{\sqrt{2\pi v^3}}.\]
\end{proof}

\begin{exercice}[\og Loi \fg{} de la hauteur d'une excursion]
\label{exo_integrationhauteur} 
\begin{enumerate}
  \item Soit $\epsilon \in \Omega^*$, on note $M(\epsilon) = \sup_{s \leq V(\epsilon)} \epsilon(s)$ le maximum atteint au cours de l'excursion, déterminer la mesure $n_M$ vérifiant pour toute fonction mesurable $f:\R^+ \to \R$ :
\[\int_{\Omega^*} \n(\d \epsilon) f(M(\epsilon)) = \int_{\R^+} n_M(\d m) f(m).\]
  \item Retrouver, grâce à la théorie des excursions, l'indépendance des processus $(S_{\tau_l}, l \geq 0)$ et $(I_{\tau_l}, l \geq 0)$, où on a posé $I_t = \inf_{s \leq t} B_s$.
\end{enumerate}
\end{exercice}

\begin{proof}
\textit{1.} On étudie la loi de $S_{\tau_l}$, en utilisant la théorie des excursions. En particulier, on note que $S_{\tau_l} = \max_{\lambda \leq l} M(e_\lambda)$. On peut ainsi écrire, 
\[ \P(S_{\tau_l} \leq a) = \E(\ind{S_{\tau_l} \leq a})= \E\left(\prod_{\lambda \leq l} \ind{M(e_\lambda) \leq a}\right). \]
Dès lors, en utilisant la formule multiplicative, on obtient
\begin{align*}
   \E\left[ \exp \left( - \sum_{\lambda \leq l} -\ln (\ind{M(e_\lambda) \leq a}) \right)\right]
   = & \exp\left(-l\int\n(\d \epsilon)(1-\ind{M(\epsilon)\leq a})\right)\\
   = & \exp\left(-l \text{ } n_M((a,+\infty))\right).
\end{align*}

D'autre part, on a également
\[
  \P(S_{\tau_l} \leq a) = \P(\tau_l \leq T_a) = \P(L_{T_a} \geq l) = \exp\left(-\frac{l}{2a}\right),\]
car $L_{T_a}$ est une variable aléatoire exponentielle de paramètre $\frac{1}{2a}$ (voir Exercice \ref{exo_transformeelaplace}). On en déduit $n_M((a,+\infty)) = \frac{1}{2a}$, ce qui permet de conclure que $n_M(\d m) = \dfrac{\d m}{2m^2}$.

\textit{2.} Nous allons maintenant étudier, à $l$ fixé, l'indépendance de $S_{\tau_l}$ et $I_{\tau_l}$. Posons $M^*(\epsilon) = \sup_{s \leq V(\epsilon)} |\epsilon(s)|$ le plus grand écart à 0 de l'excursion, positive ou négative. On a alors de la même manière que précédemment, pour tous $a, b \geq 0$ :
\begin{align*}
  \P(S_{\tau_l} \leq a, I_{\tau_l} \geq -b)
   = & \E \left[ \prod_{\lambda \leq l} \left(\ind{M^*(e_\lambda) \leq a, e_\lambda >0} + \ind{M^*(e_\lambda) \leq b, e_\lambda<0} \right)  \right]\\
   = & \exp\left( - l \int_{\Omega^*} \n(\d \epsilon) 1 - \left(\ind{M^*(\epsilon) \leq a, \epsilon >0} + \ind{M^*(\epsilon) \leq b, \epsilon<0} \right) \right)\\
   = & \exp \left( - l \left( \n^+(M(\epsilon) \geq a) + \n^+(M(\epsilon) \geq b) \right)  \right)\\
   = & \P(S_{\tau_l} \leq a) \P(I_{\tau_l} \geq -b),
\end{align*}
en utilisant la symétrie de la mesure d'Itô sur les ensembles $\Omega^+$ et $\Omega^-$. Les deux variables aléatoires sont donc indépendantes à $l$ fixé.

On utilise alors la propriété de Markov forte du mouvement brownien : le couple $(S_{\tau_{l+l'}}, I_{\tau_{l+l'}}$ est obtenu à partir de $(S_{\tau_l}, I_{\tau_l})$ et d'une copie $(\widehat{S}, \widehat{I})$ de $(S_{\tau_{l'}}, I_{\tau_{l'}})$, indépendante de $\F_{\tau_l}$ en prenant le maximum des deux premières coordonnées et minimum des secondes. Par conséquent, les deux processus étudiés sont indépendants.
\end{proof}

\begin{exercice}
Soit $k$ une fonction mesurable positive et $b \geq 0$.
\begin{enumerate}
  \item Expliciter une fonction $\phi : \R^+ \to \R^+$ vérifiant pour tout $l \geq 0$ :
\[ \P\left( \sup_{s \leq \tau_l} (|B_s| - k(L_s)) \leq b \right) =  \exp\left( - \int_0^l\d \lambda \phi(\lambda)\right).\]
  \item En déduire l'expression de la fonction 
\[\Phi_k(b) = \P\left( \sup_{s \geq 0} (|B_s| - k(L_s)) \leq b \right).\]
\end{enumerate}
\end{exercice}

\begin{proof}
\textit{1.} On explicite la probabilité $\P\left( \sup_{s \leq \tau_l} (|B_s| - k(L_s)) \leq b \right)$ à l'aide de la formule multiplicative
\begin{align*}
  \P\left( \sup_{s \leq \tau_l} (|B_s| - k(L_s)) \leq b \right)
   = & \E\left( \prod_{\lambda \leq l} \ind{ \sup_{u \leq V(e_\lambda)} |\epsilon(u)|-k(\lambda) \leq b}\right)\\
   = & \exp\left( - \int_0^l \d \lambda \n( \sup_{u \leq V(\epsilon)} |\epsilon(u)|-k(\lambda) > b ) \right).
\end{align*}
On note $M^*(\epsilon)= \sup_{u \leq V(\epsilon)} |\epsilon(u)|$. Grâce à l'Exercice \ref{exo_integrationhauteur}, on a
\begin{align*}
  \P\left( \sup_{s \leq \tau_l} (|B_s| - k(L_s)) \leq b \right)
   = & \exp\left( - \int_0^l \d \lambda \n(M^*(\epsilon)>b+k(\lambda)\right)\\
   = & \exp \left( - \int_0^l \d \lambda \frac{1}{2(b+k(\lambda))} \right).
\end{align*}
On en déduit
\[\phi( \lambda) = \frac{1}{2(b+k(\lambda))}.\]

\textit{2.} On fait tendre $l$ vers $+\infty$, on obtient par convergence monotone
\[\Phi_k(b) = \exp\left( - \int_0^{+\infty} \frac{\d \lambda}{2(b+k(\lambda))} \right).\]
\end{proof}

\begin{exercice}
Montrer que l'on peut indexer les excursions par $\Q^+$ de telle manière que si $q \leq q'$ alors $e_q$ arrive avant $e_{q'}$.
\end{exercice}

\begin{proof}
On construit la bijection entre $\Q^+$ et $G=\{\tau_{\lambda-}, \lambda \geq l\}$ comme suit. On commence par indexer par $1$ l'excursion qui enjambe l'instant 1, et pour tout $n$ entier, on note $e_n$ l'excursion enjambant l'instant de fin de l'excursion $e_{n-1}$ augmenté de 1.

On se place alors sur un intervalle entre deux excursions marquées par $p$ et $p+1$. On associe à $\frac{2p+1}{2}$ l'excursion qui enjambe le milieu de cet intervalle. On raisonne alors par récurrence sur $q$. Si $\frac{p}{q+1}$ est un rationnel différent des rationnels auxquels on a déjà associé une excursion, on associe ce nombre à l'excursion qui enjambe le point au milieu de l'intervalle entre les excursions associées aux deux rationnels qui l'encadrent. Ainsi $\frac{3}{5}$ est associé à l'excursion qui enjambe le point situé au milieu entre la fin de l'excursion $e_\frac{1}{2}$ et le début de l'excursion $e_\frac{2}{3}$.

Comme l'ensemble des excursions est dénombrable, et que chaque excursion enjambe un intervalle de mesure strictement positive, après un certain temps un rationnel lui est attribué presque sûrement.
\end{proof}

\chapter{Sur la loi de l'excursion standard}

Dans ce chapitre, nous allons, suivant l'exemple de Chung, décrire la loi de certaines excursions individuelles. Les excursions que nous choisirons enjamberont soit l'instant $t$ --autrement dit le début de l'excursion est le dernier zéro du mouvement brownien avant l'instant $t$, et la fin de l'excursion est le premier zéro après l'instant $t$ ; soit l'instant $T_c$ --ce sera alors la première excursion dont la hauteur dépasse $c$. Les lois de ces deux excursions individuelles sont différentes ; elles seront toutes deux utilisées pour décomposer la mesure d'It\^o, soit en fonction de la longueur de l'excursion générique, soit en fonction de sa hauteur. Certains de ces résultats peuvent être démontrés en utilisant la méthode du grossissement de filtrations, dont nous parlerons dans le prochain chapitre.

\section{Théorèmes principaux}

Soient $a \leq b$ deux instants aléatoires ; on note pour $0 \leq u \leq 1$, $B^{[a,b]}_u = \frac{1}{\sqrt{b-a}}B_{a+u(b-a)}$. On étudie maintenant de manière fine la décomposition de la loi du mouvement brownien entre les instants 0, $g_t$ et $d_t$.

\begin{proposition}
\label{pro_decomposition-MBt}
Notons $\P^{(3)}_0$ la loi du processus de Bessel de dimension 3 indexé par $[0,1]$, et issu de 0.

\begin{itemize}
  \item Le processus $B^{[0,g_t]}$ est un pont brownien standard de longueur 1.
  \item Le processus $|B^{[g_t,t]}|$ est un méandre brownien, un processus dont la loi est $\mathbb{M} = \frac{c}{X_1} \cdot  \P^{(3)}_0$.
  \item Le processus $|B^{[g_t,d_t]}|$ un pont de Bessel de dimension $3$ de longueur 1 appelé excursion brownienne normalisée.
\end{itemize}
\end{proposition}

On peut réaliser une description semblable du mouvement brownien entre les instants $g_{T_c}, T_c$ et $d_{T_c}$.

\begin{proposition}
\label{pro_decomposition-MBTc}
Soient $R$ et $R'$ deux processus de Bessel de dimension 3 indépendants, $\theta_c = \inf\{t \geq 0 : R_t =c\}$ et $\theta'_c = \sup\{t \geq 0 :R'_t = c\}$, on a :
\[(B_{(g_{T_c} + u)\wedge T_c}, u \geq 0) \egaldistr (R_{u\wedge \theta_c}, u \geq 0)\]
\[(B_{(d_{T_c}-u)\vee T_c}, u \geq 0) \egaldistr (R'_{u\wedge \theta'_c}, u \geq 0).\]

Si on note $M^{(c)} = \sup\{B_u, u \in[g_{T_c},d_{T_c}]\}$, on peut préciser ce résultat :
\begin{itemize}
  \item $M^{(c)} \egaldistr \frac{c}{U}$,
  \item Conditionnellement à $M^{(c)}=m$, on pose $\rho_c = \inf\{t \leq d_{T_c} : B_t = M^{(c)}\}$, les processus $(B_{(g_{T_c}+u) \wedge \rho_c}, u \geq 0)$ et $(B_{(d_{T_c}-u) \vee \rho_c}, u \geq 0)$ sont des processus de Bessel de dimension 3 indépendants issus de 0, considérés jusqu'en leur premier temps d'atteinte de $m$ respectif.
\end{itemize}
\end{proposition}

Soit $(\F_t)_{t \geq 0}$ la filtration brownienne, on définit, pour tout $\Gamma$ instant aléatoire, la tribu du passé jusqu'en $\Gamma$,
\[\F_{\Gamma} = \sigma(H_{\Gamma}, H \text{  processus prévisible}),\]
ainsi que celle du futur infinitésimal jusqu'en $\Gamma$ :
\[\F_{\Gamma+} = \cap_{\epsilon>0} \F_{\Gamma+\epsilon}.\]

Pour tout réel $s \geq 0$, toute fonction $\gamma \in \mathcal{C}(\R^+,\R)$ et tout $u \geq 0$, on pose $i_s(\gamma)(u) = \gamma(s+u)\ind{u \leq d_s(\gamma)}$. En particulier, $i_{g_t}(B)$ représente l'excursion brownienne qui enjambe l'instant $t$. On peut décrire de la façon suivante la loi de ces excursions.

\begin{lemme}
\label{lem_decomposition_longueur}
Pour toute fonction $F$ mesurable bornée, on a :
\[\E[ F(i_{g_t}) | \F_{g_t}] = \mathbf{q}(A_t, F),\]
où on a posé $A_t = t -g_t$ et $\mathbf{q}(v,F) = \frac{1}{\n(V>v)}\int_{\Omega^*} F(\epsilon) \ind{V(\epsilon)>v} \n(\d \epsilon).$

Ce résultat peut être raffiné de la manière suivante :
\[ \E[F(i_{g_t})|\F_{g_t},d_t-g_t=v] = \mathbf{\nu}(v,F),\]
où on a posé $\mathbf{\nu}(v, F) = \n(F|V=v),$ et cette dernière quantité est déterminée grâce à la décomposition de la mesure d'Itô en fonction de la longueur des excursions :
\[ \n(\d \epsilon)= \int_0^{+\infty} n_V(\d v) \n(\d \epsilon|V=v). \]
\end{lemme}

On appelle $T$ un temps d'arrêt terminal si pour tout $t \geq 0$, sur l'événement $\{T>t\}$, on a
\[ T(B) = t + T(B_{.+t}-B_t). \]
Les résultats prouvés dans le lemme précédent s'étendent ainsi.

\begin{lemme}
Soit $T$ un temps d'arrêt terminal, on a :
\[\E[ F(i_{g_T}) | \F_{g_T}] = \frac{1}{\n(V>T)} \int_{\Omega^*} F(\epsilon) \ind{V(\epsilon)>T} \n(\d \epsilon),\]
en particulier, $i_{g_T}$ est indépendant de $\F_{g_T}$.

Ce résultat peut également être raffiné ainsi :
\[ \E[F(i_{g_T})|\F_{g_T},d_T-g_T=v] = \dfrac{\mathbf{\nu}(v,F\ind{V>T})}{\mathbf{\nu}(v,\ind{V>T})},\]
où on rappelle $\mathbf{\nu}(v, F) = \n(F|V=v).$
\end{lemme}

\section{Exercices}

\begin{exercice}
Soit $t > 0$, calculer la loi de $d_t$ conditionnellement à $g_t$.
\end{exercice}

\begin{proof}
On calcule $\E(f(d_t)|\F_{g_t}) = \E(f(V(i_{g_t}) + g_t) | \F_{g_t})$, grâce au Lemme \ref{lem_decomposition_longueur}
\begin{align*}
  \E(f(d_t)|\F_{g_t}) = & \frac{1}{\n(V>t -g_t)}\int_{\Omega^*} \n(\d \epsilon) f(V(\epsilon)+g_t) \ind{V(\epsilon)>t-g_t} \\
  = & \frac{1}{\n(V>t-g_t)} \int_{t-g_t}^{+\infty} \dfrac{\d v}{\sqrt{2\pi v^3}} f(v+g_t).
\end{align*}

On en déduit 
\[\P(d_t \geq a|g_t) = \left(\frac{t-g_t}{a-g_t}\right)^\frac{1}{2} \wedge 1.\]
\end{proof}

\begin{exercice}
\begin{enumerate}
  \item Montrer que $\E(f(|B_t|)|\F_{g_t}) = \frac{1}{t-g_t} \int_0^{+\infty} e^{-\frac{y^2}{2(t-g_t)}} y f(y) \d y$.
  \item En déduire que $t-g_t$ et $\sqrt{\frac{\pi}{2}(t-g_t)}$ sont des $(\F_{g_t})$-martingales.
  \item Montrer que $f(L_t)-\sqrt{\frac{\pi}{2}(t-g_t)}f'(L_t)$ est une $(\F_{g_t})$-martingale
\end{enumerate}
\end{exercice}

\begin{proof}
\textit{1.} En utilisant la Proposition \ref{pro_decomposition-MBt}, on sait que $B^{[g_t,t]}$ est un méandre brownien, indépendant de $\F_{g_t}$ par conséquent :
\begin{align*}
  \E(f(|B_t|)|\F_{g_t}) &= \E\left(\left.f\left(\frac{B_t}{\sqrt{t-g_t}}\sqrt{t-g_t}\right)\right|\F_{g_t}\right)\\
  &= \frac{1}{t-g_t}\int_0^{+\infty} \d y e^{-\frac{y^2}{2(t-g_t)}} y f(y),
\end{align*}
par changement de variables.

\textit{2.} On applique la formule précédente à $f(y) = y$ et $f(y)=y^2$, on obtient :
\[\E(|B_t||\F_{g_t}) = \dfrac{1}{t-g_t} \int_0^{+\infty} \d y e^{-\frac{y^2}{2(t-g_t)}} y^2 = \sqrt{\frac{\pi (t-g_t)}{2}} \text{  et}\]
\[\E(B_t^2|\F_{g_t}) = \dfrac{1}{t-g_t} \int_0^{+\infty} \d y e^{-\frac{y^2}{2(t-g_t)}} y^3 = t-g_t.\]

On a alors 
\begin{align*}
 \E((t+s)-2g_{t+s}|\F_{g_t}) & = & \E\left(\E\left(\E(2(B_{t+s}^2-(t+s))|\F_{g_{t+s}})|\F_t\right)| \F_{g_t}\right)\\
  = & 2\E(B_t^2-t|\F_{g_t})\\
  = & t-2g_t,
\end{align*}
donc $(t-2g_t, t \geq 0)$ est une $(\F_{g_t})$-martingale.

Un même raisonnement donne
\[\left(\sqrt{\frac{\pi(t-g_t)}{2}}, t \geq 0\right) \text{  est une  } (\F_{g_t}) \text{-martingale.}\]
En effet, $(L_t = L_{g_t})$ en est une, et $(|B_t|-L_t)$ est une $(\F_t)$-martingale.

\textit{3.} La troisième martingale s'obtient en utilisant que $(f(L_t) - |B_t|f'(L_t), t \geq 0)$ est une $(\F_t)$-martingale --cf. Exercice \ref{exo_mart}-- donc par le même raisonnement que précédemment, on obtient :
\[f(L_t) - \sqrt{\frac{\pi}{2} (t-g_t)}f'(L_t) \text{  est une  } (\F_{g_t})\text{-martingale}.\]
\end{proof}

\begin{exercice}
Pour tout instant $t$ tel que $B_t=0$, on note $D_t$ la longueur de la plus longue excursion avant l'instant $t$. Soit $l \geq 0$ et $x \geq 0$, on pose
\[L_\beta(x) = \E\left[\int_0^{+\infty} \d t e^{-\beta t} \ind{D(g_t)\geq x}\right],\]
\[f(x) = \int_x^{+\infty} \dfrac{\d v}{\sqrt{2\pi v^3}} e^{-\beta v} \text{  et}\]
\[\Phi_l(x,\beta) = \E \left(e^{-\beta \tau_l} \ind{D(\tau_l)\geq x}\right).\]

\begin{enumerate}
  \item Pour $l \geq 0$, déterminer la loi de $D_{\tau_l}$.
  \item Montrer que $\beta L_\beta(x) = \sqrt{2\beta}\int_0^{+\infty} \d l\Phi_l(x,\beta)$.
  \item Trouver une équation différentielle satisfaite par $\Phi$.
  \item En déduite $L_\beta$ en fonction de $f(x)$.
  \item Résoudre le même problème en considérant $D_{d_t}$ au lieu de $D_{g_t}$.
  \item Par scaling brownien, calculer les transformées de Laplace de $\frac{1}{D(g_1)}$ et $\frac{1}{D(d_1)}$.
\end{enumerate}
\end{exercice}

\begin{proof}
\textit{1.} Soit $l \geq 0$, on utilise la formule multiplicative des excursions pour obtenir :
\begin{align*}
  \P(D_{\tau_l} \leq a) = & \E \left[ \prod_{\lambda \leq l} \ind{V(e_\lambda) \leq a} \right]\\
  = & \exp\left[ - l \int_{\Omega^*} \n(\d \epsilon) \left(1-\ind{V(\epsilon) \leq a} \right) \right]\\
  = & \exp\left( - l \n(V(\epsilon) \geq a)\right).
\end{align*}
Par conséquent, pour tout $a\geq 0$, on a 
\[ \P(D_{\tau_l} \leq a)  =  \exp\left( - l \sqrt{\frac{2}{\pi a}}\right).  \]
Par conséquent, $D_{\tau_l} \egaldistr \sqrt{\frac{2}{\pi}}\frac{l}{\e}$, où $\e$ est une variable aléatoire de loi exponentielle de paramètre 1.

\textit{2.} Calculons maintenant $L_\beta$, par formule additive, on a :
\begin{align*}
  \beta L_\beta(x) = & \beta \E\left[ \sum_{\lambda \geq 0}  \int_{\tau_{\lambda -}}^{\tau_\lambda} \d t e^{-\beta t} \ind{D_{g_t} \geq x}\right]\\
  = & \beta \E\left[ \sum_{\lambda \geq 0}  e^{-\beta \tau_{\lambda-}} \ind{D(\tau_{\lambda-})\geq x} \int_0^{V(e_\lambda)} \d t e^{-\beta t} \right]\\
  = & \E\left[ \int_0^{+\infty} \d l e^{-\beta \tau_l} \ind{D_{\tau_l} \geq x} \right] \int_{\Omega^*} \n(\d \epsilon) (1-e^{-\beta V(\epsilon)})
\end{align*}
Cette dernière équation peut alors se réécrire, en fonction de $f$ et $\Phi_l$ comme :
\[ \beta L_\beta(x) = \sqrt{2\beta} \int_0^{+\infty} \d l \Phi_l(x, \beta).\]

\textit{3.} On utilise à nouveau la formule additive pour déterminer une équation différentielle satisfaite par $\Phi_l$ :
\begin{align*}
 \Phi_l(x,\beta)
    = & \E \left[\sum_{\lambda \leq l} e^{-\beta \tau_\lambda}\ind{D(\tau_\lambda)\geq x} - e^{- \beta \tau_{\lambda-}} \ind{D(\tau_{\lambda-})\geq x}\right]\\
   = & \E \Big[\sum_{\lambda \leq l} e^{- \beta \tau_{\lambda-}}\left( \ind{D(\tau_{\lambda-})\geq x} \left( \exp(-\beta V(e_\lambda)) - 1 \right)\right. \\
  &  \qquad \qquad \left. + \ind{D(\tau_{\lambda-} \leq x} \ind{V(e_\lambda) \geq x} \exp(-\beta V(e_\lambda)) \right)\Big]\\
   = & - \int \n(d\epsilon) \left(1 - e^{-\beta V(\epsilon)}\right) \int_0^l \d s \Phi_s(x, \beta)\\
   & \qquad \qquad + \int \n(\d\epsilon) e^{-\beta V(\epsilon)} \ind{V(\epsilon) \geq x} \int_0^l \d s (\Phi_s(0, \beta)-\Phi_s(x, \beta)).
\end{align*}

On rappelle maintenant les deux résultats suivants :
\[\Phi_l(0,\beta) = \E(e^{-\beta \tau_l} ) = \exp(-l \sqrt{2 \beta}),\]
\[\int \n(\d \epsilon) \left(1 - e^{-\beta V(\epsilon)}\right) = \int_0^{+\infty} \dfrac{\d v}{\sqrt{2\pi v^3}} (1 - e^{-\beta v}) = \sqrt{2\beta},\]
grâce à l'Exercice \ref{exo_integrationlongueur}. On en déduit
\[\Phi_l(x, \beta) = - \sqrt{2\beta} \int_0^l \d s \Phi_s(x, \beta) + \int_x^{+\infty} \dfrac{\d v}{\sqrt{2\pi v^3}} e^{-\beta v} \int_0^l \d s (e^{-s\sqrt{2\beta}} - \Phi_s(x, \beta)).\]

En utilisant $f(x) = \int_x^{+\infty} \frac{\d v}{\sqrt{2\pi v^3}} e^{-\beta v}$, on obtient finalement l'équation différentielle satisfaite par $\Phi$ :
\[ \Phi_l (x, \beta) = -(\sqrt{2\beta} + f(x) ) \int_0^l \d s \Phi_s(x, \beta) + \frac{f(x)(1-e^{-l \sqrt{2\beta}})}{\sqrt{2\beta}}.\]

\textit{4.} Cette équation différentielle est résolue ainsi
\[ \Phi_l(x,\beta) =  -f(x) \exp\left( -l(\sqrt{2\beta} + f(x) ) \right) + f(x) \exp\left( -l \sqrt{2\beta} \right), \]
d'où on tire immédiatement :
\[\beta L_\beta(x) = \dfrac{f(x)}{f(x)+\sqrt{2\beta}}.\]

\textit{5.} Mutatis mutandis, on détermine la loi $D_{d_t}$ en calculant $\tilde{L}_\beta$ la quantité définie par
\[\tilde{L}_\beta(x) = \E\left[\int_0^{+\infty} \d t e^{-\beta t} \ind{D(d_t)\geq x}\right].\]
De la même façon que précédemment, on a
\begin{align*}
  \beta \tilde{L}_\beta(x) = & \beta \E\left[ \sum_{\lambda \geq 0}  \int_{\tau_{\lambda -}}^{\tau_\lambda} \d t e^{-\beta t} \ind{D_{d_t} \geq x}\right]\\
  = & \beta \E\left[ \sum_{\lambda \geq 0}  e^{-\beta \tau_l} \left(\ind{D(\tau_{\lambda-})\geq x}\ind{V(e_\lambda) \leq x} +\ind{V(e_\lambda) \geq x}\right) \int_0^{V(e_\lambda)} \d t e^{-\beta t} \right]\\
  = & \E\left[ \int_0^{+\infty} \d l e^{-\beta \tau_l} \ind{D_{\tau_l} \geq x} \right] \int_{\Omega^*} \n(\d \epsilon) (1-e^{-\beta V(\epsilon)})\ind{V(\epsilon) \leq x}\\
  & \qquad + \E\left[ \int_0^{+\infty} \d l e^{-\beta \tau_l}\right]\int_{\Omega^*} \n(\d \epsilon) (1-e^{-\beta V(\epsilon)})\ind{V(\epsilon) \geq x}.
\end{align*}
Calculons entre autres :
\[ \int_{\Omega^*} \n(\d \epsilon) (1-e^{-\beta V(\epsilon)})\ind{V(\epsilon) \geq x}  = \sqrt{\frac{2}{\pi x}}-f(x) \text{  et  } \]
\[ \int_{\Omega^*} \n(\d \epsilon) (1-e^{-\beta V(\epsilon)})\ind{V(\epsilon) \leq x} = \sqrt{2\beta} + f(x) - \sqrt{\frac{2}{\pi x}}. \]

On obtient donc
\[ \beta \tilde{L}_\beta (x) = \left(\sqrt{2\beta} + f(x) - \sqrt{\frac{2}{\pi x}}\right) \int_0^{+\infty} \d l \Phi_l(x,\beta) + \frac{\sqrt{\frac{2}{\pi x}}-f(x)}{\sqrt{2\beta}}, \]
ce qui permet de réécrire
\begin{align*}
  \beta \tilde{L}_\beta(x) = & \dfrac{f(x)}{f(x)+\sqrt{2\beta}}\left(1 -  \frac{\sqrt{\frac{2}{\pi x}}-f(x)}{\sqrt{2\beta}}\right)+ \frac{\sqrt{\frac{2}{\pi x}}-f(x)}{\sqrt{2\beta}}\\
  = & \dfrac{f(x)}{f(x)+\sqrt{2\beta}} + \frac{\sqrt{\frac{2}{\pi x}}-f(x)}{f(x)+\sqrt{2\beta}}\\
  = & \dfrac{\sqrt{\frac{2}{\pi x}}}{f(x)+\sqrt{2\beta}}.
\end{align*}  
On obtient donc une expression assez semblable à celle associée à $D_{g_t}$.

\textit{6.} Observons maintenant que par scaling brownien, on a 
\[ D_{g_t} \egaldistr t D_{g_1} \text{   et   } D_{d_t} \egaldistr t D_{d_1}. \]
On peut donc réécrire par exemple $L_\beta(x)$ comme :
\[ \int_0^{+\infty} \d t e^{-\beta t} \P\left(D_{g_1} \geq \frac{x}{t}\right) =  \int_0^{+\infty} \d t e^{-\beta t} \P\left( \frac{x}{D_{g_1}} \leq t \right). \]
On obtient ainsi
\[ L_1(x) = \E\left( e^{-\frac{x}{D_{g_1}}} \right),\]
et par conséquent,
\[ \E\left( \exp\left( -\frac{\lambda}{D_{g_1}} \right)\right) = \dfrac{f(\lambda)}{f(\lambda)+\sqrt{2}}. \]
De la même manière
\[ \E\left( \exp\left( -\frac{\lambda}{D_{d_1}} \right)\right) = \dfrac{\sqrt{\frac{2}{\pi \lambda}}}{f(\lambda)+\sqrt{2}}.\]
\end{proof}

\begin{exercice}
\label{exo_radonnikodym}
Soit $\P^{(3)}_0$ la loi d'un processus de Bessel de dimension 3 issu de 0, $\W^{(1)}_{0,0}$ la loi d'un pont brownien de longueur 1 et ${}^{(3)}\Pi^{(1)}$ celle d'une excursion brownienne normalisée, i.e. $(\frac{1}{\sqrt{d_t-g_t}}|B_{g_t+u(d_t-g_t)}|, u \leq 1)$.

Trouver $\phi$ et $\psi$ tels que, pour tout $u < 1$, on a :
\[ \left\{
\begin{array}{l}
  {\W^{(1)}_{0,0}}_{|\F_u} = \phi(u,X_u) \cdot \W_{|\F_u}\\
  {}^{(3)}\Pi^{(1)}_{|\F_u} = \psi(u, X_u) \cdot {\P^{(3)}_0}_{|\F_u}
\end{array}
\right.\]
\end{exercice}

\begin{proof}
\textbf{Nous donnerons deux démonstrations de ce résultat, qui présentent chacune leur intérêt propre}

\textit{Cette première démonstration peut facilement se généraliser à des processus de Markov satisfaisant des hypothèses assez générales.}

Commençons par étudier le pont brownien  $\beta$. C'est un mouvement brownien conditionné à revenir en 0 à l'instant 1. Soit $B$ un mouvement brownien, $F$ mesurable pour la tribu engendrée par $\mathcal{C}([0,u],\R)$, et $g : \R \to \R^+$, on calcule de deux manières différentes
\begin{multline*}
  \E\left( F(B_s,s \leq u) g(B_1) \right) = \E\left( \E(F(B_s,s\leq u)|B_1)g(B_1)\right)\\
   = \int_\R \dfrac{\d x}{\sqrt{2\pi}}e^{-\frac{x^2}{2}} g(x) \E(F(B_s,s\leq u)|B_1=x)
\end{multline*}
par conditionnement par rapport à $B_1$ d'une part, et
\begin{multline*}
  \E\left( F(B_s,s \leq u) g(B_1) \right)  = \E\left[ F(B_s, s\leq u) \E(g(B_1)|\F_u)\right]\\
   = \int_\R \dfrac{\d y}{\sqrt{2\pi (1-u)}} g(y) \E\left[\exp\left(-\frac{(y-B_u)^2}{2(1-u)}\right) F(B_s,s\leq u)\right]
\end{multline*}
par propriété de Markov d'autre part. Par égalité de ces deux quantités, on tire,
\[\E(F(\beta_s,s\leq u))= \E(F(B_s,s\leq u)|B_1=0) = \E\left[ \dfrac{\exp\left( - \frac{B_u^2}{2(1-u)} \right)}{\sqrt{1-u}} F(B_s,s\leq u) \right].\]
On en conclut 
\[{\W^{(1)}_{0,0}}_{|\F_u} = \dfrac{\exp\left( - \frac{B_u^2}{2(1-u)} \right)}{\sqrt{1-u}} \cdot W_{|\F_u}.\]

On s'intéresse maintenant au pont de Bessel de dimension 3 $\rho$. Soit $B$ un mouvement brownien de dimension 3, on peut écrire, de deux manières différentes :
\begin{multline*}
 \E\left[ F(\norme{B_s},s\leq u) g(\norme{B_1}) \right] = \int_{\R^3} \dfrac{\d x}{\sqrt{2\pi}^3}e^{-\frac{\norme{x}^2}{2}} g(\norme{x}) \E(F(\norme{B_s},s\leq u)|B_1=x)\\
  = \int_{\R^3} \dfrac{\d x}{\sqrt{2\pi (1-u)}^3} g(\norme{x}) \E\left[\exp\left(-\frac{\norme{x-B_u}^2}{2(1-u)}\right) F(\norme{B_s},s\leq u)\right].
\end{multline*}
Ceci peut encore \^etre réécrit :
\begin{align*}
  &  \E\left[ F(\norme{B_s},s\leq u) g(\norme{B_1}) \right]\\
  = & \int_{\R^+}  \d r \sqrt{\dfrac{2}{\pi}}r^2 e^{-\frac{r^2}{2}} g(r) \E^{(3)}_{0,r} (F(X_s,s\leq u))\\
  = & \int_{\R^+} \d r \sqrt{\dfrac{2}{\pi(1-u)^3}}r^2 g(r) \int_{[0,2\pi]^2} \d \theta \sin(\theta)\d \phi\\
  &\qquad \qquad \E\left[\exp\left(-\frac{\norme{x(r,\theta,\phi)-B_u}^2}{2(1-u)}\right) F(\norme{B_s},s\leq u)\right].
\end{align*}
On obtient dès lors,
\[ \E (F(\rho_s,s\leq u)) = \E\left[ \dfrac{\exp\left(-\frac{R_u^2}{2(1-u)}\right)}{\sqrt{(1-u)^3}} F(R_s, s \leq u)\right].\]
On en conclut ${}^{(3)}\Pi^{(1)}_{|\F_u} = \dfrac{\exp\left(-\frac{R_u^2}{2(1-u)}\right)}{\sqrt{(1-u)^3}} \cdot {P^{(3)}_0}_{|\F_u}.$

\textit{Nous pouvons donner une autre preuve de ce résultat, qui utilise fortement le caractère gaussien des processus considérés.}

Commençons par étudier le pont brownien $\beta$. On sait que pour tout $u<1$, $\beta_u$ est une variable aléatoire gaussienne centrée de variance $u(1-u)$. On tente de déterminer la dérivée de Radon-Nikodym de $\beta_u$ par rapport à $B_u$, une variable aléatoire gaussienne centrée de variance $u$. On obtient $\phi(u,X_u) = \sqrt{1-u} \exp(-\frac{X_u^2}{2(1-u)})$.

Montrons maintenant que la dérivée de Radon-Nikodym ainsi obtenue est celle du pont brownien jusqu'en l'instant $u$ par rapport à celle du mouvement brownien jusqu'en l'instant $u$. Soit $t_1 < t_2 < \cdots < t_n < u$, on calcule, en conditionnant par rapport à $\beta_u$ :
\[\E\left[f(\beta_{t_1},\beta_{t_2},\ldots,\beta_{t_n},\beta_u), \phi(u, \beta_u) \right] = \E\left[ \widehat{f}(\beta_u) \phi(u,\beta_u) \right] = \E(\widehat{f}(B_u)).\]

Or conditionnellement à $\beta_u=x$ ou à $B_u=x$, le processus jusqu'en l'instant $u$ est un pont brownien issu de $0$ et arrivant en $x$ à l'instant $u$. Par conséquent
\[\E(\widehat{f}(B_u)) = \E\left[ f(B_{t_1},\ldots, B_{t_n}, B_u) \right],\]
$\phi(u,X_u)$ est donc bien la dérivée de Radon-Nikodym espérée, et nous avons pour $u <1$ :
\[{\W^{(1)}_{0,0}}_{|\F_u} = \frac{1}{\sqrt{1-u}} \exp \left(-\frac{{X_u}^2}{2(1-u)}\right) \cdot \W_{|\F_u}.\]

De la même manière, un processus de Bessel ou un pont de Bessel sur $[0,1]$, tous deux conditionnés à arriver en $x$ à l'instant $u$ sont des ponts de Bessel de longueur $u$ arrivant en $x$. Par conséquent, il suffit de calculer la dérivée de Radon-Nikodym de $X_u$ sous les lois $\Pi^{(1)}$ et $\P^{(3)}_0$. Sous ${}^{(3)}\Pi^{(1)}$, on a :
\[\E(f(X_u)) = \dfrac{1}{(\sqrt{2\pi u(1-u)})^3}\int_{\R^3} \d x \d y \d z f(\sqrt{x^2+y^2+z^2}) \exp\left( -\frac{x^2+y^2+z^2}{2u(1-u)}\right),\]
et sous $\P^{(3)}_0$, on a :
\[\E^{(3)}_0(f(X_u)) = \dfrac{1}{(\sqrt{2\pi u})^3}\int_{\R^3} \d x \d y \d z f(\sqrt{x^2+y^2+z^2}) \exp\left( -\frac{x^2+y^2+z^2}{2u}\right).\]
Par conséquent on obtient la dérivée de Radon-Nikodym suivante :
\[{}^{(3)}\Pi^{(1)}_{|\F_u} = (1-u)^{-\frac{3}{2}}\exp\left( -\frac{X_u^2}{2(1-u)} \right) \cdot {\P^{(3)}_0}_{|\F_u}.\]
\end{proof}

\begin{exercice}[Une autre construction de l'excursion individuelle]
\label{exo_excursionindividuelle}
Soit $c \in \R^+$, on note :
\[s_c : \epsilon \mapsto \left( \frac{1}{\sqrt{c}}\epsilon(ct)\right)_{t \geq 0}.\]
\begin{enumerate}
  \item Montrer que pour tout $A$ borélien de $\Omega^*$, on a $\n(s_c^{-1}(A)) = \frac{1}{\sqrt{c}} \n(A)$.
  \item On pose $\nu(\epsilon) = s_{V(\epsilon)}(\epsilon)$ la \og forme \fg{} de l'excursion. Montrer qu'il existe une mesure $\gamma$ sur $\{\epsilon \in \Omega^* : V(\epsilon)=1\}$ telle que pour tout borelien $\Gamma$ et réel $c \geq 0$
\[ \gamma(\Gamma) = \dfrac{\n^+(\nu^{-1}(\Gamma) \cap \{V \geq c\})}{\n^+(V>c)} .\]
La mesure $\gamma$ est considérée comme la loi de l'excursion brownienne individuelle de longueur 1.
  \item Montrer que pour tout $S$ borélien de $\R^+$, on a :
\[ \n^+(\nu^{-1}(\Gamma) \cap \{V \in S\} ) = \gamma(\Gamma) \n^+(V \in S).\]
(Indépendance de la forme et de la longueur de l'excursion.)
  \item Montrer que
\[ \gamma(\Gamma) = \P(\nu(e^c) \in \Gamma),\]
où $e^c$ est la première excursion positive de longueur au moins $c$.
\end{enumerate}
\end{exercice}

\begin{proof}
\textit{1.} On remarque pour commencer que $\n(A)$ est la paramètre du processus de Poisson comptant le nombre d'excursions browniennes appartenant à $A$. On utilise alors la propriété de scaling du mouvement brownien, c'est-à-dire $(\frac{1}{\sqrt{c}} B_{ct}) \egaldistr (B_t)$. Les deux processus de Poissons 
\[N_l = \sum_{\lambda \leq l} \ind{e_\lambda \in A} \text{   et   }  N'_l = \sum_{\lambda \leq \frac{l}{\sqrt{c}}} \ind{s_c(e_\lambda) \in A}\]
sont de même paramètre, d'où $\n(s_c^{-1}(A)) = \frac{1}{\sqrt{c}}\n(A).$

\textit{2.} Montrons que $\dfrac{\n^+(\nu^{-1}(\Gamma) \cap \{V \geq c\})}{\n^+(V \geq c)}$ est indépendant du paramètre $c \geq 0$. Pour cela on utilise le résultat précédent, on a :
\[s_c^{-1}(\nu^{-1}(\Gamma) \cap \{V \geq 1\})=\nu^{-1}(\Gamma) \cap \{V \geq c\},\]
on en déduit :
\[\dfrac{\n^+(\nu^{-1}(\Gamma) \cap \{V \geq c\})}{\n^+(V \geq c)}  = \dfrac{\n^+(\nu^{-1}(\Gamma) \cap \{V \geq 1\})}{\n^+(V \geq 1)},\]
que l'on pose par définition comme égal à $\gamma(\Gamma)$.

\textit{3.} Par définition
\[\n^+(\nu^{-1}(\Gamma) \cap \{V \geq c\}) = \gamma(\Gamma)\n^+(V \geq c),\]
et par lemme des classes monotones, pour tout $S$ borélien de $\R^+$ :
\[\n^+(\nu^{-1}(\Gamma) \cap \{V \in S\}) = \gamma(\Gamma)\n^+(V \in S).\]

\textit{4.} On utilise maintenant le Lemme des réveils. Soit $N$ le processus de Poisson comptant les excursions de longueur au moins $c$ qui ont une forme dans $\Gamma$, et $N'$ celui comptant les excursions de longueur au moins $c$ qui ont une forme différente de celles dans $\Gamma$. Ces deux processus de Poissons sont indépendants car relatifs à des ensembles disjoints, par conséquent la probabilité que la première excursion de longueur $c$ ait une forme dans $\Gamma$ vaut
\[\frac{\n^+(\nu^{-1}(\Gamma) \cap \{V \geq c\})}{\n^+(\nu^{-1}(\Gamma) \cap \{V \geq c\}) + \n^+(\nu^{-1}(\Gamma)^c \cap \{V \geq c\})} = \gamma(\Gamma).\]
On a donc bien
\[\gamma(\Gamma) = \P(\nu(e^c) \in \Gamma). \]
\end{proof}

\begin{exercice}[Autour du méandre brownien]
\begin{enumerate}
  \item Montrer que le processus $(B_t,t\leq u)$ conditionné à $g_1=u$, le processus est un pont brownien de longueur $u$. En déduire que $B^{[0,g_1]}$ est un pont brownien indépendant de $g_1$ et de $(B_{g_1+u},u \geq 0)$ (utiliser $tB_{\frac{1}{t}}$ mouvement brownien).
  \item Soit $l^a$ les temps locaux en 1 du pont brownien, montrer que $L^a_{g_1}$ à la même loi que $\sqrt{g_1}l^{\frac{a}{g_1}}$.
  \item Montrer que $|B^{[g_1,1]}|$ est indépendant de $\F_{g_1}$ et que $\frac{|B_1|}{\sqrt{1-g_1}}$ a pour loi $\sqrt{2\e}$.
  \item Calculer la loi jointe de $(g_t,L_t,B_t).$
\end{enumerate}
\end{exercice}

\begin{proof}
\textit{1.} On sait que $\widehat{B}_t = t B_\frac{1}{t}$ est un autre mouvement brownien. Nous allons donc étudier le processus $(B_t, t \leq g_1)$ comme processus défini à partir de $\widehat{B}$. Remarquons pour commencer que $g_1= \frac{1}{\widehat{d}_1}$. Nous pouvons donc réécrire le processus étudié comme
\[ \left(t \widehat{B}_{\frac{1}{t}}, t \leq \frac{1}{\widehat{d}_1} \right). \]

Or, le processus $\tilde{B}_v = \widehat{B}_{d_1+ v}$ est un mouvement brownien indépendant de $\widehat{\F}_{\widehat{d}_1}$, donc en particulier de $\widehat{d}_1$. En d'autre termes, conditionnellement à $g_1=u$ (i.e. $\widehat{d}_1= \frac{1}{u}$), le processus
\[ \left(t \tilde{B}_{\frac{1}{t}-\frac{1}{u}}, t \leq u\right), \]
est un pont brownien standard indépendant de $\widehat{\F}_\frac{1}{u}$, ce que l'on peut vérifier en calculant les covariances impliquées.

Par conséquent, conditionnellement à $g_1=u$, $(B_t, t \leq u)$ est un pont brownien de longueur $u$. On en déduit sans peine que conditionnellement à $g_1$, $B^{[0,g_1]}$ est un pont brownien de longueur 1, donc en particulier est indépendant de $g_1$.

\textit{2.} Il suffit d'appliquer le résultat précédent pour obtenir l'identité en loi escomptée. Le processus $B^{[0,g_1]}$ suit la loi d'un pont brownien indépendant de $g_1$, par conséquent
\begin{align*}
  l^a  \egaldistr & \lim_{\epsilon \to 0} \int_0^1 \frac{1}{\epsilon} \ind{0<B_{ug_1}-a\sqrt{g_1}\leq \sqrt{g_1} \epsilon} \d u\\
   = & \lim_{\eta \to 0} \frac{1}{\sqrt{g_1}} \int_0^{g_1} \frac{1}{\eta} \ind{0<B_u-a\sqrt{g_1} \leq \eta} \d u = \frac{1}{\sqrt{g_1}} L^{a\sqrt{g_1}}_{g_1}.
\end{align*}
Par indépendance de $g_1$ et de $B^{[0,g_1]}$, on obtient l'égalité souhaitée.

\textit{3.} On s'intéresse maintenant à la loi du méandre brownien $|B^{[g_1,1]}|$. On souhaite montrer que ce processus est indépendant de $\F_{g_1}$, pour cela, nous allons à nouveau écrire ce processus à partir du mouvement brownien $\widehat{B}$. Observons en effet que $\F_{g_1}$ est la tribu du passé de $B$ avant l'instant $g_1$, donc également la tribu du futur de $\widehat{B}$ après l'instant $\widehat{d}_1$. Cette tribu est, en utilisant la propriété de Markov forte, indépendante du futur de $B$, entre les instants $g_1$ et $t$. Il reste donc à montrer que, convenablement renormalisé, le processus $B^{[g_1,1]}$ est bien indépendant de $g_1$.

Conditionnons par rapport à l'évènement $\{g_1=u\}$, nous allons discuter de la loi de $B^{[u,1]}$.  Observons que ce processus s'écrit 
\[(\frac{1}{\sqrt{1-u}} B_{(1-u)t+u}-B_u, 0 \leq t \leq 1) \egaldistr (\beta_t,0 \leq t \leq 1),\]
par propriété de scaling du mouvement brownien. Par conséquent, la loi de ce processus est indépendante de $g_1$ (la Proposition \ref{pro_decomposition-MBt} montre que ce processus est un méandre brownien).

On peut maintenant écrire :
\[|B_1|^2 = (1-g_1) \left(B^{[g_1,1]}_1\right)^2.\]
Or $B_1^2$ est le carré d'une gaussienne, et $(1-g_1)$ est une variable aléatoire suivant la loi de l'arcsinus, indépendante (par ce qui précède) de $B^{[g_1,1]}_1$. Par conséquent, on peut sans difficultés, par exemple par calcul de la transformée de Legendre, ou en utilisant les propriétés de l'algèbre beta-gamma vérifier que :
\[\left(B^{[g_1,1]}_1\right)^2 = 2\e.\]
Par conséquent, $|B^{[g_1,1]}_1| \egaldistr \sqrt{2\e}$ suit la loi de Rayleigh. Le signe de $B_1$ étant indépendant de toutes les autres variables aléatoires considérées, on peut sans peine vérifier que $B^{[g_1,1]}_1$ suit une loi de Rayleigh symétrique.

\textit{4.} La loi jointe de $(g_t, L_t,B_t)$ s'obtient grâce à la propriété de scaling : \[(g_t,L_t,B_t) \egaldistr (tg_1, \sqrt{t} L_1, \sqrt{t} B_1).\]
De plus, grâce aux calculs que nous avons effectué jusqu'ici, décomposant le mouvement brownien jusqu'en l'instant 1 en trois quantités indépendantes : un pont brownien de longueur 1, un méandre brownien de longueur 1, et une loi de l'arcsinus indiquant comment les deux processus précédents devaient être renormalisés, et mis bout-à-bout. C'est ainsi que pour calculer la loi jointe de $(g_1,L_1,|B_1|)$, il suffit de connaitre la loi de $l_1$ et de $\frac{1}{\sqrt{1-g_1}}|B_1|$, et de réaliser les bonnes renormalisations.

On a alors :
\[\P(g_1 \in \d s, L_1 \in \d l, |B_1| \in \d x) = \P(g_1 \in \d s, \sqrt{s}l_1\in \d l, \sqrt{1-s}\sqrt{2\e} \in \d x),\]
et ces trois quantités sont indépendantes. Dès lors, en observant que $\sgn(B_1)$ est indépendant des autres variables aléatoires, la densité jointe de $(g_t, L_t, B_t)$ est :
  \[\ind{l \geq 0, s \leq t} \frac{l}{\sqrt{2\pi s^3}} \exp\left( -\frac{l^2}{2s}\right) \frac{|x|}{\sqrt{2\pi (t-s)^3}} \exp\left( - \frac{x^2}{2(t-s)}\right) .\]
\end{proof}

\begin{exercice}[Lien entre le méandre brownien et l'excursion individuelle standard]
Soit $M$ un méandre brownien, montrer que conditionnellement à $M_1=0$, le méandre suit la loi d'une excursion brownienne standard.
\end{exercice}

\begin{proof}
Soit $B$ un mouvement brownien, on pose $M=B^{[g_1,1]}$, et on utilise le mouvement brownien $\beta_u = u B_\frac{1}{u}$. Par propriété de Markov, conditionnellement à $B_1=\beta_1=0$, le processus $\left(\beta_{1+s}, s \geq 0\right)$ est un mouvement brownien issu de 0. Dès lors le processus entre 1 et $g_1$ est une excursion de $\beta$, que l'on peut renormaliser en une excursion standard de longueur 1. Par retournement du temps, c'est également une excursion de $B$ donc $M$, conditionné à $M_1=0$ est bien une excursion brownienne standard.
\end{proof}

\begin{exercice}
Soit $T_c = \inf\{t \geq 0:B_t = c\}$ le premier temps d'atteinte de $c$ par le mouvement brownien $B$, $g_t = \sup\{s \leq t : B_s = 0\}$ le dernier zéro avant l'instant $t$, et $d_t = \inf\{s \geq t : B_s = 0\}$ le premier zéro après l'instant $t$.
\begin{enumerate}
  \item Calculer la loi de $M_c = \max_{g_{T_c} \leq t \leq d_{T_c}} B_t$.
  \item Calculer la loi de $S_{g_{T_c}}$.
\end{enumerate}
\end{exercice}

\begin{proof}
\textit{1.} Observons pour commencer que $M_c$ le maximum de $B_t$ entre les instants $g_{T_c}$ et $d_{T_c}$ vaut également $S_{d_{T_c}}$ le maximum de $B$ avant l'instant $d_{T_c}$. En effet, on a clairement $S_{g_{T_c}} < c$  et $M_c \geq c$.

On calcule alors, pour tout $x \geq c$ :
\[\P(M_c \geq x) = \P(T_x \leq d_{T_c}).\]
On utilise ensuite la propriété de Markov à l'instant $T_c$,
\[P(T_x \leq d_{T_c}) = \P_c(T_x \leq T_0) = \dfrac{c}{x}.\]
$M_c$ est donc une variable aléatoire à densité par rapport à la mesure de Lebesgue de densité $\frac{c}{x^2} \ind{x \geq c}$, c'est-à-dire $M_c \egaldistr \frac{c}{U}$.

\begin{remarque}
Ce résultat peut être obtenu simplement en réutilisant l'Exercice \ref{exo_supmartingale} : on a vu que le maximum d'une martingale locale positive issue de $c$ tendant vers 0 était égal en distribution à $\frac{c}{U}$. On applique ceci au mouvement brownien issu de $c$ et absorbé en 0, qui reprend la partie de l'excursion enjambant $T_c$ qui atteint le maximum sur l'excursion.
\end{remarque}

\textit{2.} Trouvons maintenant la loi de $S_{g_{T_c}}$ ; pour tout $x \leq c$, on calcule
\[
  \P(S_{g_{T_c}} \leq x)
   =  \P(g_{T_c} \leq T_x)
   =  \P(T_c \leq d_{T_x}).
\]
Par propriété de Markov forte appliquée en $T_x$, on en déduit $\P(S_{g_{T_c}} \leq x) = \frac{c-x}{c}$, par conséquent, $S_{g_{T_c}}$ est distribuée de manière uniforme sur $[0,c]$.
\end{proof}

\chapter[Grossissements de filtrations]{Grossissements de filtrations}

De façon générale, quand un espace de probabilité filtré est donné, il existe trois grandes façons de modifier les semi-martingales. On peut réaliser un changement de probabilité (théorème de Girsanov), un changement de temps (Doob-Meyer), ou, comme nous le verrons dans ce chapitre, un changement de filtration. Nous allons voir que sous certaines hypothèses, il est possible de grossir une filtration donnée avec une certaine quantité d'information de manière à ce que les martingales de la filtration d'origine demeurent des semi-martingales pour la nouvelle filtration. Il est même possible de donner explicitement leur décomposition en somme d'une martingale et d'un processus à variations finies prévisible, ce qui donne un résultat équivalent à la formule de Girsanov pour les changements de probabilité. Pour un panorama assez complet sur le changement de temps et le changement de probabilité, le lecteur pourra se référer à \cite{BNS}.

Cette théorie a commencé à être développée de façon indépendante par Jeulin et Yor \cite{JeY1978} et M. Barlow \cite{Bar1978} ; citons ensuite le volume d'application \cite{JeY1985} puis d'essai de popularisation \cite{MaY2006}. Cette théorie nécessite une bonne familiarité avec les principaux résultats de la théorie générale des processus ; on peut néanmoins reverser la situation en disant que la théorie du grossissement offre une famille variée  d'illustrations de cette théorie. Malgré de nombreux efforts, on peut dire que cette théorie n'est pas encore rentrée dans les mœurs, en tout cas pas au m\^eme degré que le calcul d'It\^o. Mais ne désespérons pas, ce petit chapitre représente un effort supplémentaire de popularisation !

\section{Théorème principal}

Citons pour commencer le théorème principal de ce chapitre, qui montre que tout temps aléatoire $\Lambda$ bien choisi peut \^etre utilisé pour enrichir une filtration sans modifier les caractères de semi-martingales, ce qui transforme $\Lambda$ en temps d'arr\^et.

\begin{theoreme}[Formule de grossissement]
\label{thm:grossissement}
Soit $(\F_t)$ une filtration et $\Gamma$ un ensemble $(\F_t)$-prévisible. On pose $\Lambda(\omega) = \sup \{t \geq 0 : (t,\omega) \in \Gamma\}$ un temps aléatoire, supposé fini presque sûrement.

Pour $t \geq 0$, on pose $\F^{\Lambda}_t = \F_t \vee \sigma(\Lambda\wedge s, s \leq t)$. On suppose que pour tout $(\F_t)$-temps d'arrêt $T$, on a $\P(\Lambda=T)=0$. Alors toute $(\F_t)$-martingale $M$ est une $(\F^\Lambda_t)$-semi-martingale dont on peut donner la décomposition suivante :
\[M_t = M^\Lambda_t + \left( \int_0^{t\wedge \Lambda} \dfrac{\d \crochet{M,Z}_s}{Z_s} + \int_{\Lambda \wedge t}^t \dfrac{\d \crochet{M,1-Z}_s}{1-Z_s}\right),\]
où on a posé $Z_t = \P(\Lambda>t|\F_t)$.
\end{theoreme}

\begin{remarque}
Si la martingale considérée est un mouvement brownien, alors la partie martingale du processus modifié est encore un mouvement brownien, en utilisant le théorème de Lévy.
\end{remarque}

La preuve de ce théorème s'obtient à l'aide des trois lemmes suivants.
\begin{lemme}
Soit $J$ processus $(\F^\Lambda_t)$-prévisible, il existe $J^+$ et $J^-$ deux processus $(\F_t)$-prévisibles tels que pour tout $u \geq 0$ :
\[ J_u = J^-_u \ind{u \leq \Lambda} + J^+_u \ind{u \geq \Lambda}.\]
\end{lemme}

\begin{lemme}
Soit $Z_t = \P(\Lambda > t|\F_t)$, on écrit la décomposition de $Z_t=\mu_t-A_t$ où $\mu$ est une martingale locale et $A$ un processus continu croissant.

Pour tout processus $\F$-prévisible $k$, la formule suivante est satisfaite :
\[\E(k_\Lambda) = \E\left[ \int_0^{+\infty} k_s \d A_s \right].\]
\end{lemme}

\begin{lemme}
Si $k$ est un processus continu à variations bornées $(\F_t)$-adapté, on a
\[\E(k_\Lambda) = \E\left( \int_0^{+\infty} Z_u dk_u \right) \text{   et}\]
\[\E(\int_\Lambda^{+\infty} dk_u = \E\left[\int_0^{+\infty} (1-Z_u)dk_u\right].\]
\end{lemme}

On insère pour finir un résultat concernant les diffusions (i.e. processus de Markov sur $\R$= continues. On a montré dans l'Exercice \ref{exo_radonnikodym} que dans le cas du mouvement Brownien ou du processus de Bessel, la loi d'un pont issu de $x$ terminant en $y$ admet une dérivée de Radon-Nikodym par rapport à la mesure d'origine du processus. Le résultat suivant en est une extension.

\begin{proposition}
Soit $X_t$ une diffusion sur $\R$, $\P_x(X_t \in \d y) =p_t(x,y)\d y$ son semi-groupe, que l'on suppose continu et strictement positif. On note $P^v_{x \to y}$ la loi de $(X_t, t\leq v)$ sous $\P_x$, conditionné par $X_v=y$. Si cette famille peut être choisie étroitement continue en $y$, alors on a pour tout $u \leq v$,
\[{P^v_{x \to y}}_{|\F_u} = \dfrac{p_{v-u}(X_u,y)}{p_v(x,y)}{\P_x}_{|\F_u}.\] 
\end{proposition}

\section{Exercices}

\begin{exercice}
Soit $B$ un mouvement brownien. Calculer $Z^\Lambda$ et décrivez $B$ dans la filtration $(\F^\Lambda_t)$ dans les cas suivants :
\begin{itemize}
 \item $\Lambda = g_t = \sup \{u \leq t : B_u = 0\}$,
 \item $\Lambda = g^a_t = \sup\{u \leq t : B_u=a\}$,
 \item $\Lambda = g_{T_a}$ avec $T_a = \inf\{t \geq 0:B_t=a\}$,
 \item $\Lambda = g_{T^*_a}$ avec $T^*_a = \inf\{t \geq 0:|B_t|=a\}$.
\end{itemize}
\end{exercice}

\begin{proof}
On vérifie pour commencer que tous les temps aléatoires cités précédemment sont bien les derniers instants d'ensembles prévisibles. On a par exemple dans le troisième cas : $g_{T_a} = \sup\{t \geq 0 : B_t = 0 \text{  et  } S_t < a\}$.

On calcule alors dans chacun des cas précédents la valeur de $Z^\Lambda$, on a, dans le premier cas
\begin{align*}
  Z^{g_t}_s & = \P(g_t > s | \F_s)\\
  & = \P(t>d_s|\F_s)\\
  & = \P_{B_s}(t-s>T_0)\\
  & = \Phi\left(\frac{|B_s|}{\sqrt{t-s}}\right)\ind{s <t},
\end{align*}
où on a posé $\Phi(x) = \sqrt{\frac{2}{\pi}}\int_x^{+\infty} e^{-\frac{u^2}{2}}\d u$. En effet, on a pour tout $x \geq 0$
\[\P_x(T_0<u) = \P(I_u<-x) = \P\left(|N|<\frac{-x}{\sqrt{u}}\right).\]
De la même manière, on a également avec des notations évidentes
\begin{align*}
  Z^{g^a_t}_s & = \P(t>d^a_s|F_s)\\
  & = \P_{B_s} (t-s>T_a)\\
  & = \Phi\left(\frac{|B_s-a|}{\sqrt{t-s}}\right)\ind{s <t}.
\end{align*}

On a également par des calculs similaires les égalités
\begin{align*}
  Z^{g_{T_a}}_s  = & \P(T_a>d_s|\F_s)\\
   = & \ind{T_a \geq s, B_s \geq 0} \P_{B_s}(T_a>d_s) + \ind{B_s\leq 0}\\
   = & 1 - \frac{B^+_{s\wedge T_a}}{a},
\end{align*}
et de la même manière :
\begin{align*}
  Z^{g_{T^*_a}}_s = & \P(T^*_a>d_s|\F_s)\\
   = & \ind{T^*_a \geq s} \P_{B_s}(T^*_a>d_s)\\
   = & 1 - \frac{|B|_{s\wedge T^*_a}}{a}.
\end{align*}

Pour déterminer la loi du mouvement brownien dans la filtration enrichie, il suffit alors d'appliquer le Théorème \ref{thm:grossissement}. On a en particulier, en notant $\beta$ un mouvement brownien bien choisi dans la filtration $(\F^\Lambda_t)$ :
\begin{itemize}
 \item pour $\Lambda = g_t$,
\begin{multline*}
  B_s = \beta_s  + \left( \int_0^{s\wedge g_t} \dfrac{\d u}{\sqrt{t-u}} \sgn(B_u) \sqrt{\frac{2}{\pi}} \frac{e^{-\frac{B_u^2}{2(t-u)}}}{Z^{g_t}_u}\right.\\
   \left. + \sgn(B_t)\int_{g_t \wedge s}^{t\wedge s} \dfrac{\d u}{\sqrt{t-u}} \sqrt{\frac{2}{\pi}} \frac{e^{-\frac{B_u^2}{2(t-u)}}}{1-Z^{g_t}_u}  \right),
\end{multline*}
 \item pour $\Lambda = g^a_t$, une expression similaire, en remplaçant $B_t$ par $B_t-a$
 \item pour $\Lambda=g_{T_a}$, en utilisant la formule de Tanaka, on observe que 
\[B_s = \beta_s + \left( \int_0^{s\wedge g_{T_a}} \d u\dfrac{\ind{B_u>0}}{a-B_u} + \int_{g_{T_a} \wedge s}^{T_a\wedge s}\dfrac{\d u}{B_u} \right),\]
 \item et de manière semblable, pour $\Lambda=g_{T^*_a}$, on a
\[B_s = \beta_s + \left( \int_0^{s\wedge g_{T^*_a}} \d u\dfrac{\sgn(B_u)}{a-|B_u|} + \int_{g_{T^*_a} \wedge s}^{T^*_a\wedge s}\dfrac{\d u}{B_u}\right).\]
\end{itemize}
\end{proof}

\begin{exercice}
Déduire de l'exercice précédent la loi de $(B_{g_{T_c}+u}, u \leq d_{T_c}-g_{T_c})$ qui est la première excursion qui enjambe l'instant $T_c$.
\end{exercice}

\begin{proof}
On s'intéresse au mouvement brownien dans la filtration $\F^{g_{T_c}}$, on sait que ce processus s'écrit :
\[B_s = \beta_s + \left( \int_0^{s\wedge g_{T_c}} \d u\dfrac{\ind{B_u>0}}{c-B_u} + \int_{g_{T_c} \wedge s}^{T_c \wedge s}\dfrac{\d u}{B_u}\right),\]
où $\beta$ est un mouvement brownien dans la nouvelle filtration.
Par conséquent, sachant que $B_{g_{T_c}}=0$, on a pour tout $u \in [0,T_c-g_{T_c}]$ :
\[B_{g_{T_c}+u} = \beta_{g_{T_c}+u}-\beta_{g_{T_c}} + \int_0^u \dfrac{\d v}{B_{g_{T_c}+v}}.\]

Par conséquent, $B_{g_{T_c}+.}$ est solution de l'équation différentielle stochastique
\[ \d X_t = \d B_t + \frac{\d t}{X_t} \]
entre les instants $0$ et $T_c-g_{T_c}$. C'est donc un processus de Bessel de dimension 3, issu de 0, jusqu'en son premier temps d'atteinte de $c$. De plus, entre les instants $T_c-g_{T_c}$ et $d_{T_c}-g_{T_c}$, le processus n'est pas modifié par l'agrandissement de filtration, donc reste un mouvement brownien.

En appliquant le théorème de retournement du temps de Nagasawa, le mouvement brownien issu de $c$ jusqu'en son premier temps d'atteinte de 0 peut être vu comme un processus de Bessel de dimension 3 issu de 0 jusqu'en son dernier temps d'atteinte de $c$. Par conséquent, l'excursion qui enjambe l'instant $T_c$ a la loi de deux processus de Bessel de dimension 3 mis dos-à-dos.
\end{proof}

\begin{exercice}
Soit $N$ une martingale positive issue de 1 tendant p.s. vers 0 en $+\infty$, étudier la martingale dans la filtration enrichie par $\Lambda = \sup\{t \geq 0 : N_t = \sup_{s \geq 0} N_s\}$ le temps d'atteinte du point maximal de la martingale.
\end{exercice}

\begin{proof}
Rappelons pour commencer que le maximum d'une martingale positive issue de $c>0$ tendant vers 0 en $+\infty$ est égal en loi à $\frac{c}{U}$.

On calcule alors
\[Z_t = \P(\Lambda > t | \F_t) = \P(\sup_{u \geq t} N_u > \sup_{s\leq t} N_u|\F_t),\]
et en décalant le temps d'un facteur $t$, c'est la probabilité pour que le maximum d'une martingale issue de $N_t$, tendant vers 0 en $+\infty$ reste plus grand que $S_t = \sup_{s \leq t}N_s$ qui est, conditionnellement à $\F_t$, une constante. On obtient donc
\[Z_t = \P(\frac{N_t}{U}>S_t) = \frac{N_t}{S_t}.\]

Par conséquent, en utilisant la formule de grossissement, avec la filtration $\F^\Lambda$, la martingale $N$ s'écrit
\[N_t = \tilde{N}_t + \int_0^{t\wedge \Lambda} \frac{\d \crochet{N}_s}{N_s} -  \int_{\Lambda \wedge t}^t \dfrac{\d \crochet{N}_s}{S_s-N_s}.\]
\end{proof}

\begin{exercice}
Soit $R_t$ un processus de Bessel de dimension 3. On pose $\gamma_a= \sup\{t \geq 0 : R_t=a\}$ et $(\mathcal{R}_t, t \geq 0)$ la filtration naturelle de $R$.

\begin{enumerate}
  \item Montrer la formule 
  \[\P(\gamma_a > t | \mathcal{R}_t) = \left( \frac{a}{R_t} \right)\wedge 1.\]
  \item Donner la décomposition d'Itô de cette sur-martingale.
  \item En déduire la loi du processus $(R_{\gamma_a+t}-a, t \geq 0).$
\end{enumerate}
\end{exercice}

\begin{proof}
\textit{1.} Rappelons pour commencer que $R_t$ est la solution de l'EDS
\[R_t = \int_0^t \dfrac{\d s}{R_s} + B_t\]

On utilise le fait que $(\frac{1}{R_t})$ est une martingale locale positive. Pour tout $t \geq 0$, $(\frac{1}{R_{t+s}}, s \geq 0)$ est une martingale locale positive issue de $\frac{1}{R_t}$. On connaît par conséquent la loi de son maximum,
\[\P(\gamma_a>t|\mathcal{R}_t) = \P\left(\left. \sup_{s \geq t} \frac{1}{R_s} > \frac{1}{a}\right|\mathcal{R}_t\right) = \left(\frac{a}{R_t}\right)\wedge 1.\]

\textit{2.} On va maintenant décomposer la sur-martingale $Z^{(a)}_t$ grâce à la formule d'Itô-Tanaka :
\[Z^{(a)}_t = \int_0^t \frac{-a}{R_s^2}\ind{R_s>a} \d R_s + \int_0^t \frac{a}{R_s^3}d \crochet{R}_s - \frac{1}{2a}L^a_t(R).\]

On obtient donc, en utilisant l'EDS satisfaite par $R_t$ :
\[Z^{(a)}_t = \int_0^t \frac{-a}{R_s^2}\ind{R_s>a}\d B_s - \frac{1}{2a}L^a_t(R).\]

\textit{3.} On note $\mathcal{R}^{(a)}_t=\mathcal{R}^{(\gamma_a)}_t$ la filtration qui transforme $\gamma_a$ en un temps d'arrêt. Dans cette filtration, le mouvement brownien $B$ s'écrit :
\begin{align*}
  B_t &= B^{(a)}_t + \int_0^{t\wedge \gamma_a} \dfrac{\d \crochet{B,Z^{(a)}}_s}{Z^{(a)}_s} + \int_{t\wedge \gamma_a}^t \dfrac{\d \crochet{B,1-Z^{(a)}}_s}{1-Z^{(a)}_s},\\
  & = B^{(a)}_t - \int_0^{t \wedge \gamma_a}  \frac{a\ind{R_s > a}\d s}{R_s^2 Z^{(a)}_s} + \int_{t \wedge \gamma_a}^t \frac{a \ind{R_s > a} \d s}{(1-Z^{(a)}_s)R_s^2}.
\end{align*}
Or, rappelons que $Z^{(a)}_t = \frac{a}{R_t} \wedge 1$, et que pour tout $t > \gamma_a$, $R_t > a$. Par conséquent
\[ B_t = B^{(a)}_t - \int_0^{t\wedge \gamma_a} \frac{\ind{R_s > a} \d s}{R_s} + \int_{t\wedge \gamma_a}^t \dfrac{a\d s}{R_s (R_s - a)}.\]

On pose alors $X_t = R_{t + \gamma_a}-a$, on observe que $X$ vérifie
\begin{align*}
  X_t & = R_{t+\gamma_a} - R_{\gamma_a}\\
  & = B_{t+\gamma_a}-B_{\gamma_a} + \int_{\gamma_a}^{t+\gamma_a} \dfrac{\d s}{R_s}\\
  & = B^{(a)}_{t+\gamma_a}-B^{(a)}_{\gamma_a} + \int_{\gamma_a}^{t+\gamma_a} \dfrac{a}{R_s (R_s-a)} + \dfrac{1}{R_s} \d s.
\end{align*}
On utilise maintenant que $\gamma_a$ est un temps d'arrêt pour le mouvement brownien $B^{(a)}$. On note $\widehat{B}_t = B^{(a)}_{\gamma_a+t} - B^{(a)}_{\gamma_a}$, qui est un mouvement brownien. On a alors
\begin{align*}
  X_t
  & = \widehat{B}_t + \int_0^t \dfrac{a}{R_{s+\gamma_a} (R_{s+\gamma_a}-a)} + \dfrac{1}{R_{s+\gamma_a}} \d s\\
  & = \widehat{B}_t + \int_0^t \d s \dfrac{a}{X_s(X_s + a)} + \frac{1}{X_s+a}\\
  & = & \widehat{B}_t + \int_0^t \frac{\d s}{X_s}.
\end{align*}
Le processus $X$ satisfaisant donc à cette équation différentielle stochastique, c'est un processus de Bessel de dimension 3 issu de 0.
\end{proof}

\begin{exercice}[Preuve de la formule de grossissement]
On n'utilisera bien entendu aucun des résultats listés dans ce chapitre pour résoudre cet exercice. On se place sous les hypothèses du Théorème \ref{thm:grossissement}
\begin{enumerate}
  \item On pose $J$ un processus $\left(\F^{\Lambda}_t\right)$-prévisible, montrer que l'on peut trouver $J^+$ et $J^-$ deux processus $(\F_t)$-prévisibles tels que pour tout $u \geq 0$ :
\[ J_u = J^-_u \ind{u \leq \Lambda} + J^+_u \ind{u \geq \Lambda}.\]
  \item Soit $Z_t = \P(\Lambda > t|\F_t)$, on écrit la décomposition de $Z_t=\mu_t-A_t$ où $\mu$ est une martingale locale et $A$ un processus continu croissant. Montrer que pour tout processus $(\F_t)$-prévisible $k$, la formule suivante est satisfaite :
\[\E(k_\Lambda) = \E\left[ \int_0^{+\infty} k_s \d A_s \right].\]
  \item Soit $k$ est un processus continu à variations bornées $(\F_t)$-adapté issu de 0, montrer que
\[\E(k_\Lambda) = \E\left( \int_0^{+\infty} Z_u dk_u \right) \text{   et}\]
\[\E\left(\int_\Lambda^{+\infty} dk_u\right) = \E\left[\int_0^{+\infty} (1-Z_u)dk_u\right].\]
  \item En déduire la formule de grossissement.
\end{enumerate}
\end{exercice}

\begin{proof}
\textit{1.} On prend pour $J_u$ un processus $\left(\F^{\Lambda}_t\right)$-prévisible simple de la forme 
\[J_u= H_s \ind{u \in (s, t]},\]
avec $H_s$ une variable aléatoire $\F^\Lambda_s$-mesurable de la forme
\[ H_s = K_s f(\Lambda \wedge s), \]
où $K_s$ est $\F_s$-mesurable, et $f$ est une fonction mesurable $\R_+ \to \R$.

On observe alors que $J^-_u$ vérifie :
\[J^-_u \ind{u \leq \Lambda}= J_u \ind{u \leq \Lambda} = H_s \ind{u \in (s,t]} \ind{u \geq \Lambda}.\]
Or lorsque ceci n'est pas nul, on sait que $\Lambda$ est au moins égal à $s$. Par conséquent, on peut réécrire $J^-_u = K_s f(s) \ind{u \in (s,t]}$.

On peut de la même manière écrire
\[J^+_u = K_s f(\sup\{u \leq s | (u, \omega) \in \Gamma\}) \ind{u \in (s,t]},\]
qui est bien $\F_s$-mesurable, car $\Gamma$ est un ensemble prévisible. On conclut par lemme des classes monotones.

\textit{2.} On choisit maintenant un processus $(\F_t)$-prévisible $k_u = H_s \ind{u \in (s,t]}$, avec $H_s$ une variable aléatoire $\F_s$-mesurable. On peut écrire
\[\E(k_\Lambda) = \E(H_s \ind{\Lambda \in (s,t]}) = \E(H_s (Z_s- Z_t)) = \E(H_s(-A_s+A_t))\]
par propriété de martingale de $\mu$. On peut donc bien écrire
\[ \E(k_\Lambda) = \E \left( \int_0^{+\infty} k_u \d A_u \right).\]
On conclut encore par lemme des classes monotones.

\textit{3.} On utilise le résultat précédent, le processus $k$ est bien prévisible, on a
\begin{align*}
  \E(k_\Lambda) = & \E\left( \int_0^{+\infty} k_u \d A_u \right)\\
   = & \E\left( - \int_0^{+\infty} k_u dZ_u \right)\\
   = & \E\left(\int_0^{+\infty} Z_u dk_u - k_\infty Z_\infty \right)\\
   = & \E\left(\int_0^{+\infty} Z_u dk_u\right).
\end{align*}
par intégration par parties, en utilisant le fait que $Z_\infty = 0$ p.s. De la même manière,
\begin{align*}
  \E\left(\int_\Lambda^{+\infty} dk_u\right) = & \E(k_\infty - k_\Lambda)\\
  = & \E\left( \int_0^{+\infty} (1-Z_u)dk_u \right).
\end{align*}

\textit{4.} Ces trois résultats nous permettent de déduire la formule de grossissement. Soit $(M_t)$ une $\F_t$-martingale, on souhaite calculer $\E(M_t-M_s|\F^\Lambda_s)$, pour calculer le \og défaut de martingale \fg{}  dans la nouvelle filtration.  Pour cela, soit $J$ un processus $(\F^\Lambda_t)$-prévisible borné, que l'on décompose grâce à la question \textit{1}, on calcule
\[\E\left( \int_0^{+\infty} J_u \d M_u \right) = \E\left( \int_0^{+\infty} J^+_u \ind{u \geq \Lambda}+J^-_u \ind{u \leq \Lambda} \d M_u \right).\]
En utilisant le fait que $(\int_0^t J^+_u \d M_u)$ est une $(\F_t)$-martingale, on peut réécrire :
\[\E\left( \int_0^{+\infty} J_u \d M_u \right)  = \E \left( \int_0^\Lambda (J^-_u-J^+_u) \d M_u\right).\]

On applique alors la question \textit{2} au processus $k_t = \int_0^t(J^-_u-J^+_u) \d M_u$, on obtient :
\begin{align*}
  \E\left( \int_0^{+\infty} J_u \d M_u \right) = & \E \left( \int_0^{+\infty} \d A_t  \int_0^t  (J^-_u-J^+_u) \d M_u\right)\\
  = & \E \left( A_\infty \int_0^{+\infty} (J^-_u-J^+_u) \d M_u \right) \text{  par intégration par parties.}
\end{align*}
Comme $Z_\infty=0$, on a immédiatement $\mu_\infty=A_\infty$, on peut donc réécrire cette dernière espérance comme l'espérance du produit de deux martingales locales, donc égale à l'espérance de leur crochet :
\[ \E \left( \mu_\infty \int_0^{+\infty} (J^-_u-J^+_u) \d M_u\right) = \E \left( \int_0^{+\infty} J^-_u-J^+_u \d \crochet{M,Z}_u \right).\]
On réécrit alors cette dernière quantité, en appliquant les résultats de la question \textit{3} aux processus prévisible à variations finies $\int_0^t \frac{J^\pm_u}{Z_u} \d \crochet{M,Z}_u$. On a
\begin{align*}
  \E \left( \int_0^{+\infty} J^-_u \d \crochet{M,Z}_u \right) &= \E\left( \int_0^{\Lambda} \frac{J^-_u}{Z_u}\d \crochet{M,Z}_u \right) \text{  et}\\
\E \left( \int_0^{+\infty} J^+_u \d \crochet{M,Z}_u \right) &= -\E\left( \int_{\Lambda}^{+\infty} \frac{J^+_u}{1-Z_u}\d \crochet{M,1-Z}_u \right).
\end{align*}

On en déduit
\[\E \left(\int_0^{+\infty} J_u \d M_u \right) = \E\left( \int_0^{\Lambda} \frac{J_u}{Z_u}\d \crochet{M,Z}_u + \int_\Lambda^{+\infty} \frac{J_u}{1-Z_u}\d \crochet{M,1-Z}_u\right),\]
par conséquent, le processus $\tilde{M}$ défini par
\[\tilde{M}_t = M_t - \int_0^{t \wedge \Lambda} \frac{\d \crochet{M,Z}_u}{Z_u} - \int_{t \wedge \Lambda}^\Lambda \frac{\d\crochet{M,1-Z}_u}{1-Z_u}\]
est une $\F^{\Lambda}$-martingale.
\end{proof}

\chapter[Décompositions de la mesure d'Itô]{Décomposition d'Itô et de Williams de la mesure d'Itô}

La mesure d'Itô étant une mesure $\sigma$-finie, il est naturel que les probabilistes aient cherché à la représenter de plusieurs façons, aussi naturelles que possible, comme intégrales ($\sigma$-finies) de probabilités sur l'espace des trajectoires. On a ainsi obtenu des décompositions par rapport à la longueur $V$, ou au maximum $M$ de l'excursion. Nous avons déjà déterminé précédemment la mesure image de $\n$ par les applications $V$ et $M$, il s'agit alors de déterminer la loi $n_V$ conditionnellement à la variable fixée.

On obtient ainsi deux décompositions différentes de la mesure d'Itô. Dans ces deux décompositions interviennent des processus dont la loi est connue ; ce sont précisément le pont de Bessel de dimension 3, ou deux processus de Bessel de dimension 3 mis dos à dos. Il est alors d'intérêt de déterminer une formule de \og changement de variables \fg{} entre ces deux décompositions, autrement dit de déterminer la dérivée de Radon-Nikod\'ym entre les lois de ces deux processus.

\section{Théorèmes principaux}

On donne pour commencer la décomposition d'Itô de la mesure $\n$. Celle-ci est réalisée par rapport à la longueur de l'excursion brownienne.

\begin{theoreme}[Décomposition d'Itô]
On note $\Pi^{(v)}$ la loi du pont de processus de Bessel de dimension 3 de longueur $v$, on a alors :
\[\n^+(\d \epsilon) = \int_0^{+\infty} \dfrac{\d v}{2\sqrt{2\pi v^3}} \Pi^{(v)}(\d \epsilon).\]
\end{theoreme}

\begin{remarque}
Afin de réaliser un processus de loi $\Pi^{(v)}$, on peut considérer un mouvement brownien $C$ de dimension 3 ; alors $(C_t - \frac{t}{v}C_v)_{t \leq v}$ est un pont brownien dans $\R^3$ issu de $0$. Enfin, le processus $(\norme{C_t - \frac{t}{v} C_v})_{t \leq v}$ suit la loi $\Pi^{(v)}$.
\end{remarque}

\begin{theoreme}[Décomposition de Williams]
Notons $\widehat{\Pi}^{(m)}$ la loi de deux processus de Bessel de dimension 3 considérés jusqu'en leurs premiers temps respectifs d'atteinte du niveau $m$ et mis dos-à-dos. On a alors :
\[\n^+(\d \epsilon) = \int_{\R^+} \dfrac{\d m}{2m^2} \widehat{\Pi}^{(m)}(\d \epsilon).\]
\end{theoreme}

\begin{remarque}
Afin de réaliser un processus de loi $\widehat{\Pi}^{(m)}$, on considère $R$ et $R'$ deux processus de Bessel jusqu'en leur premier temps d'atteinte de $m$, notés respectivement $T_m$ et $T'_m$. Le processus défini par :
\[X_t = \left\{
\begin{array}{ll}
  R_t & \text{ pour } t \leq T_m\\
  R_{T_m+T'_m-t} & \text{ pour } T_m \leq t \leq T_m +T'_m
\end{array}
\right.\]
a pour loi $\widehat{\Pi}^{(m)}$.
\end{remarque}

Le résultat suivant a pour objet de permettre de passer de la représentation d'Itô à celle de Williams, et réciproquement.

\begin{theoreme}[Formule de concordance]
Soit $(\rho_t, t \leq V)$ un processus suivant la loi $\widehat{\Pi}^{(1)}$, on note 
\[(\widehat{\rho}_v,v \leq 1) = \left(\dfrac{1}{\sqrt{V}} \rho_{vV}, v \leq 1\right).\]

Le processus $\widehat{\rho}$ a une loi équivalente à $\Pi^{(1)}$ la loi du pont de Bessel de dimension 3, que l'on note $(r_v ,v \leq 1)$. De plus, pour toute fonction $F$ mesurable positive
\[\E\left( F(r_v,v \leq 1) \right) = \sqrt{\frac{\pi}{2}}\E\left(\frac{1}{\max_{v \leq 1} \widehat{\rho}_v} F(\widehat{\rho}_v,v \leq 1)\right).\]
\end{theoreme}

\begin{remarque}
On rappelle que $\rho$ processus suivant la loi $\widehat{\Pi}^{(1)}$ peut être obtenu en collant dos-à-dos deux processus de Bessel de dimension 3 indépendants jusqu'en leurs premiers temps d'atteinte respectifs, notés $T_1$ et $T'_1$, de 1. Dans ce cas, $T_1+T'_1$ représente le temps de vie de cette excursion et $\widehat{\rho}$ est défini de la manière suivante :
\[\widehat{\rho}_v = \frac{1}{\sqrt{T_1+T'_1}} \rho((T_1+T'_1)v), v \leq 1.\]

En particulier, la variable aléatoire $\sqrt{T_1+T'_1}$ est l'inverse du maximum de $\widehat{\rho}$ sur $[0,1]$, 
la dérivée de Radon-Nikodym de la loi de $\widehat{\rho}$ par rapport à $\Pi^{(1)}$ peut donc se réécrire $C \sqrt{T_1+T'_1}$.
\end{remarque}

\section{Exercices}

\begin{exercice}
Retrouver la constante $C = \sqrt{\frac{2}{\pi}}$ de la dérivée de Radon-Nikodym de la loi de $\widehat{\rho}$ par rapport à $\Pi^{(1)}$. Pour ce faire, on calculera $\E(\sqrt{T_1+T'_1})$, où $T_1$ et $T'_1$ sont les deux temps d'atteinte de 1 par des processus de Bessel de dimension 3 indépendants, grâce à l'égalité suivante
\[\E\left(\exp\left(-\frac{\lambda^2}{2} T_1\right)\right) = \frac{\lambda}{\sh \lambda}\]
\end{exercice}

\begin{proof}
On observe pour commencer que 
\[\sqrt{\lambda} = \int_0^{+\infty} \d x \frac{1-e^{-\lambda x}}{2\sqrt{\pi x^3}}.\]
Par conséquent, en utilisant le théorème de Fubini, on a
\begin{align*}
  \E\left(\sqrt{T_1+T'_1}\right)  = & \int_0^{+\infty} \frac{1-\E(\exp(-(T_1+T'_1)x)}{2\sqrt{\pi x^3}} \d x\\
   = & \sqrt{\frac{2}{\pi}}\int_0^{+\infty} \frac{1-\E\left[\exp\left(-(T_1+T'_1)\frac{u^2}{2}\right)\right]}{u^2}\d u\\
   = & \sqrt{\frac{2}{\pi}}\int_0^{+\infty} \frac{\left(1-\frac{u^2}{(\sh u)^2}\right)}{u^2}\d u.
\end{align*}

Il suffit alors de calculer
\[ \int_0^{+\infty} \left(\frac{1}{u^2} - \frac{1}{(\sh u)^2}\right)\d u = \left[ \frac{-1}{x} + \coth(x) \right]_0^{+\infty} = 1, \]
dès lors $\E(\sqrt{T_1+T'_1}) = \sqrt{\frac{2}{\pi}}$. En particulier, lors des calculs de dérivée de Radon-Nikod\'ym :
\[1 = C \E\left( \frac{1}{\max_{v \leq 1} \widehat{\rho}_v}\right) = C \E (\sqrt{T_1+T'_1}),\]
d'où l'on tire $C = \sqrt{\frac{\pi}{2}}$
\end{proof}

\begin{exercice}[Lemme de la décomposition de Williams]
Soit $c \geq 0$, montrer que :
\[\E\left[ F(i_{g_{T_c}}) | \F_{g_{T_c}}, d_{T_c}-g_{T_c} = v \right] = \dfrac{\n(F \ind{V>T_c}|V=v)}{\n(V>T_c|V=v)}.\]

En déduire, en utilisant la Proposition \ref{pro_decomposition-MBTc} la décomposition de Williams.
\end{exercice}

\begin{proof}
Soit $H$ un processus prévisible et $f,F$ deux fonctions mesurables bornées, on calcule l'espérance conditionnelle suivante, grâce à la la formule clé additive modifiée :
\begin{align*}
   \E \left[ H_{g_{T_c}} F(i_{g_{T_c}}) f(d_{T_c}-g_{T_c}) \right]
   = & \E \left[ \sum_{s \in G} \ind{s \leq T_c} H_s F(i_s) f(V(i_s)) \ind{M(i_s) \geq c} \right]\\
   = & \E \left[ \int_0^{T_c} \d L_u H_u \int_{\Omega^*} \n(\d\epsilon) F(\epsilon)f(V(\epsilon)) \ind{M(\epsilon) \geq c} \right].
\end{align*}
Dès lors, e n utilisant la décomposition de la mesure d'It\^o selon la longueur des excursions, on a
\begin{multline*}
   \E \left[ H_{g_{T_c}} F(i_{g_{T_c}}) f(d_{T_c}-g_{T_c}) \right]\\
   = \E \left[ \int_0^{T_c} \d L_u H_u \int_0^{+\infty} \dfrac{\d v}{\sqrt{2\pi v^3}} f(v)\int_{\Omega^*} \n(\d \epsilon|V=v) F(\epsilon) \ind{M(\epsilon) \geq c} \right].
\end{multline*}
On pose $\mathbf{\nu}(v,G) = \int_{\Omega^*} G(\epsilon) \n(\d \epsilon|V=v)$, on a alors :
\[\E \left[ H_{g_{T_c}} F(i_{g_{T_c}}) f(d_{T_c}-g_{T_c}) \right] = \E \left[ \int_0^{T_c} \d L_u H_u \int_0^{+\infty} \dfrac{\d v}{\sqrt{2\pi v^3}} f(v) \mathbf{\nu}(v,F\ind{M \geq c})\right].\]

En particulier, on observe que $\frac{\mathbf{\nu}(d_{T_c}-g_{T_c},F\ind{M \geq c})}{\mathbf{\nu}(d_{T_c}-g_{T_c},\ind{M>c})}$ est une variable aléatoire mesurable par rapport à $d_{T_c}-g_{T_c}$ ; par conséquent, en posant $F=1$, on a
\begin{multline*}
   \E\left[ H_{g_{T_c}} f(g_{T_c}-d_{T_c}) \frac{\mathbf{\nu}(d_{T_c}-g_{T_c},F\ind{M \geq c})}{\mathbf{\nu}(d_{T_c}-g_{T_c},\ind{M>c})} \right]\\
  = \E \left[ \int_0^{T_c} \d L_u H_u \int_{0}^{+\infty} \dfrac{\d v}{\sqrt{2\pi v^3}} f(\mathbf{\nu}(v,\ind{M>c})  \dfrac{\mathbf{\nu}(v,F\ind{M \geq c})}{\mathbf{\nu}(v,\ind{M>c})} \right],
\end{multline*}
on obtient donc bien l'espérance conditionnelle escomptée.

On démontre maintenant la décomposition de Williams pour la mesure d'Itô, c'est-à-dire
\[\n( \cdot | M=m) = \widehat{\Pi}^{(m)}.\]
On note $e^{(c)}$ l'excursion individuelle qui enjambe $T_c$, soit $\Gamma$ un ensemble mesurable de $\Omega^*$, on a :
\[\P(e^{(c)} \in \Gamma) = \dfrac{\n(\Gamma \cap \{M \geq c\})}{\n(M \geq c)},\]
en effet $\{ V \geq T_c \} = \{ M \geq c \}$. Or $\n(M\geq c) = \int_c^{+\infty} \frac{\d m}{2m^2} = \frac{1}{2c}$, par conséquent
\begin{align*}
  \n(\Gamma \cap \{M \geq c\})  = & \frac{1}{2c} \P(e^{(c)} \in \Gamma)\\
   = & \frac{1}{2c} \int_c^{+\infty} \dfrac{c \d m}{m^2} \P(e^{(c)} \in \Gamma | M^{(c)}=m)\\
   = & \int_c^{+\infty} \dfrac{\d m}{2m^2} \widehat{\Pi}^{(m)}(\Gamma)
\end{align*}

On fait maintenant tendre $c$ vers 0, ce qui permet d'obtenir 
\[\n^+(\Gamma) = \int_0^{+\infty} \dfrac{\d m}{2m^2} \widehat{\Pi}^{(m)}(\Gamma),\]
qui est bien la décomposition de Williams de la mesure d'Itô.
\end{proof}

\begin{exercice}
\label{exo_tempslocalexcursion}
Grâce à la décomposition, on peut définir la notion de temps local en $x$ pour une excursion $\epsilon$ fixée de longueur $v$. On note cette quantité $l^x_v$.

Déterminer la mesure $\nu_x$ telle que pour toute fonction mesurable positive $f$ vérifiant $f(0)=0$ on ait :
\[\int_{\Omega^*} \n(\d \epsilon) f(l^x_{V(\epsilon)}) = \int_\R \nu_x(\d y) f(y).\]
\end{exercice}

\begin{proof}
Soit $\mu \geq 0$, on a
\begin{align*}
  \E(\exp(-\mu L^x_{\tau_l})  = & \E\left[ \exp\left( - \mu \sum_{\lambda \leq l} L^x_{\tau_l}-L^x_{\tau_{l-}} \right) \right]\\
   = & \E\left[\exp\left(- \mu \sum_{\lambda \leq l} l^x_{V(e_\lambda)}(e_\lambda)\right) \right]
\end{align*}
Il suffit d'appliquer la formule multiplicative pour obtenir :
\[\E\left[\exp\left(- \mu \sum_{\lambda \leq l} l^x_{V(e_\lambda)}(e_\lambda)\right) \right] = \exp\left( - l \int_{\Omega^*} \n(\d \epsilon) \left(1-e^{-\mu l^x_{V(\epsilon)}} \right) \right).\]

On utilise le théorème de Ray-Knight d'autre part, on a
\[\E(\exp(-\mu L^x_{\tau_l})) = Q^0_l(\exp(-\mu X_x)) = \exp\left(\frac{l}{2}\Phi'_\mu(0)\right),\]
par formule de Ricatti, où $\Phi_\mu$ est la seule solution décroissante minorée issue de 1 de l'équation de Sturm-Liouville associée à $\mu \delta_x(\d y)$
\[u'' = 2\mu \delta_x u.\]

En intégrant, on trouve, pour $t < x$ :
\[u'(t) = u'(0) \text{   et   } u(t) = 1 +t u'(0),\]
et, pour $t \geq x$ :
\[u'(t) = u'(0) + 2 \mu (1 + x u'(0)).\]
Comme $\Phi_\mu$ est souhaitée décroissante et minorée,
\[\Phi_\mu'(0) + 2 \mu (1 + x \Phi'_\mu(0)) = 0,\]
donc $\Phi'_\mu(0) = \dfrac{-2\mu}{1 + 2 x \mu}$.

On obtient alors, par identification :
\[\int_{\Omega^*}\n(\d \epsilon) \left(1-e^{-\mu l^x_{V(\epsilon)}} \right) = \frac{\mu}{1+2\mu x} = \int_{\R^+} \nu_x(\d y) (1-e^{-\mu y}).\]
Il reste donc à déterminer $\nu_x(\d y)$. On obtient, en dérivant l'égalité ci-dessus
où on obtient en dérivant
\begin{multline*}
  \int_{\R^+} \nu_x(\d y) y e^{-\mu y}  =  \dfrac{1}{(1+2\mu x)^2}\\
   =  \int_{\R^+} \d y y \exp(-(1 + 2 \mu x) y) = \frac{1}{4x^2} \int_{\R^+} y e^{-\frac{y}{2x}} e^{-\mu y}\d y,
\end{multline*}
on en conclut $\nu_x(\d y) = \frac{1}{4x^2} e^{-\frac{y}{2x}}$.

\begin{remarque}
On peut de la même manière calculer toutes les marginales finies-dimensionnelles de $(l^x_{V(\epsilon)}, x \geq 0)$, et donc obtenir un certain type de théorème de Ray-Knight pour les excursions sous la mesure d'Itô. La mesure obtenue pour ce \og processus \fg{} serait toutefois $\sigma$-finie.
\end{remarque}
\end{proof}

\begin{exercice}
Soit $\delta >0$, on pose :
\[ \eta_t(\delta) = \sum_{\lambda \geq 0 : \tau_\lambda \leq t} \ind{V(\epsilon) \geq \delta} \text{  et}\]
\[ \xi_t(\delta) = \sum_{\lambda \geq 0 : \tau_\lambda \leq t} V(\epsilon)\ind{V(\epsilon) \leq \delta}.\]

Montrer que :
\[\begin{array}{l}
  \left(\sqrt{\frac{\pi\delta}{2}} \eta_t(\delta), t \geq 0\right) \convd (L_t, t \geq 0) \text{  p.s.}\\
  \left(\sqrt{\frac{\pi}{2\delta}} \xi_t(\delta), t \geq 0\right) \convd (L_t,t\geq 0) \text{  p.s.}
\end{array}\]
\end{exercice}

\begin{proof}
Notons pour commencer que $\eta_{\tau_l}(\delta)$ est un processus de Poisson de paramètre $\n(V\geq \delta) = \sqrt{\frac{2}{\pi \delta}}$. Par continuité, pour obtenir la convergence en loi du processus, il suffit de montrer la convergence sur les instants rationnels, et grâce à la propriété de Markov, on peut se contenter d'une convergence à $t$ fixé.

Grâce à la loi des grands nombres, on a, pour tout $l \geq 0$, $\sqrt{\frac{\pi \delta}{2}} \eta_{\tau_l} \conv l$ p.s. La convergence a donc lieu en temps que processus, donc par composition des processus on obtient $\sqrt{\frac{\pi \delta}{2}}\eta_t \convd L_t$ p.s.

Intéressons-nous maintenant à $\xi$, on observe que
\[\xi_t(\delta) = - \int_0^{\delta} x \eta_t(\d x),\]
par conséquent, par convergence dominée, on obtient
\[ \sqrt{\frac{\pi}{2 \delta}} \xi_t(\delta) = - \frac{1}{\delta}\int_0^\delta x \sqrt{\frac{\pi \delta}{2}} \eta_t(\d x) \convd L_t.\]
\end{proof}

\begin{exercice}
\begin{enumerate}
  \item Calculer $\E\left( \exp \left(-\frac{\lambda^2 A^{(+)}_{\tau_l}}{2} \ind{S_{\tau_l} \leq a}\right)\right).$
  \item En déduire $\E\left( \exp \left(-\frac{\lambda^2 A^{(+)}_{\tau_l}}{2S_{\tau_l}^2} \right)\right) = \frac{2\lambda}{\sh(2\lambda)}$.
  \item En déduire l'identité de Knight $\frac{A^{(+)}_{\tau_l}}{2S_{\tau_l}^2} \egaldistr \inf\{s \geq 0 : R^{(3)}_s=2\}.$
\end{enumerate}
\end{exercice}

\begin{proof}
\textit{1.}
On utilise la formule multiplicative des excursions
\begin{align*}
  &\E\left( \exp \left(-\frac{\mu^2 A^{(+)}_{\tau_l}}{2} \ind{S_{\tau_l} \leq a}\right)\right)\\
  & \qquad \qquad \qquad = \E\left( \exp \left(- \sum_{\lambda \leq l} \frac{\mu^2}{2} V(e_\lambda) \ind{e_\lambda \geq 0} - \ln(\ind{M(e_\lambda)\leq a} ) \right) \right)\\
   & \qquad \qquad \qquad = \exp\left[ - l \int_{\Omega^+} \n(\d \epsilon) \left( 1- \exp\left( - \frac{\mu^2}{2}V(\epsilon) \right)\ind{M(\epsilon) \leq a} \right) \right]\\
   & \qquad \qquad \qquad = \exp\left[ - l \n(M(\epsilon) \geq a) - l \int_0^a \dfrac{\d m}{2m^2} \widehat{\Pi}^{(m)} \left( 1- \exp\left( - \frac{\mu^2}{2}V(\epsilon)\right)\right) \right],
\end{align*}
en utilisant la décomposition de Williams. On calcule alors séparément chacun des termes. On a pour commencer $\n(M(\epsilon) \geq a) = \frac{1}{2a}$. Ensuite
\begin{align*}
  \widehat{\Pi}^{(m)} \left( 1- \exp\left( - \frac{\mu^2}{2}V(\epsilon)\right)\right)
  & =  1 - \E\left(\exp\left(- \frac{\mu^2}{2} (T_m^{(3)}+{T'}_m^{(3)})\right) \right)\\
  & =  1 - \E\left(\exp\left(- \frac{\mu^2 m^2}{2} T_1^{(3)}\right) \right)^2\\
  & =  1 - \left(\dfrac{\mu m}{\sh(\mu m)}\right)^2.
\end{align*}
Par conséquent
\begin{align*}
  \int_0^a \dfrac{\d m}{2m^2} \widehat{\Pi}^{(m)} \left( 1- \exp\left( - \frac{\mu^2}{2}V(\epsilon)\right)\right)
   = & \int_0^a \d m \left(\dfrac{1}{2m^2} - \frac{\mu^2}{2\sh(\mu m)^2}\right)\\
   = & \left[ \mu \frac{\ch(\mu m)}{2\sh(\mu m)} - \frac{1}{2m} \right]_0^a\\
   = & \frac{\mu \ch(\mu a)}{2\sh(\mu a)} - \frac{1}{2a}.
\end{align*}

On en déduit la formule suivante
\[  \E\left( \exp \left(-\frac{\mu^2 A^{(+)}_{\tau_l}}{2} \ind{S_{\tau_l} \leq a}\right)\right) = \exp\left( - l \mu \coth(\mu a) \right).\]

\textit{2.}
Cette formule nous permet d'accéder à la loi jointe de $S_{\tau_l}$ et $A^{(+)}_{\tau_l}$. En effet, en dérivant par rapport à $a$ l'expression, on obtient la mesure de la transformée de Laplace de $A^{(+)}_{\tau_l}$ restreinte à l'ensemble $\{S_{\tau_l} \in da\}$. On peut alors calculer l'espérance de $\exp \left(-\frac{\mu^2 A^{(+)}_{\tau_l}}{2S_{\tau_l}^2}\right)$ de la façon suivante. Posons
$\phi(\mu,a) = \exp\left( - l \mu \coth(\mu a) \right)$,
on a
\begin{align*}
  \E\left[ \exp \left(-\frac{\mu^2 A^{(+)}_{\tau_l}}{2S_{\tau_l}^2}\right) \right]
   = & \int_{\R^+} \P(S_{\tau_l} \in da) \E\left[\left.\exp \left(-\frac{\mu^2 A^{(+)}_{\tau_l}}{2a^2}\right) \right| S_{\tau_l} \in \d a \right]\\
   = & -\int_{\R^+} \d a \partial_2\phi \left( \frac{\mu}{a},a \right)\\
   = & \frac{\mu^2}{2\sh(\mu)^2}\int_{\R^+} \dfrac{da}{a^2}\exp\left( -l \mu \frac{1}{2a}\coth(\mu) \right)\\
   = & \frac{\mu}{\sh(\mu)\ch(\mu)}\\
   = & \frac{2\mu}{\sh(2\mu)}.
\end{align*}

\textit{3.}
On tire de cette égalité l'identité de Knight, en effet, soit $T^{(3)}_2$ le premier temps d'atteinte de $2$ par un processus de Bessel de dimension 3 issu de 0, on a
\[ \E\left[e^{- \mu^2 T^{(3)}_2}\right] = \frac{2 \mu}{\sh(2\mu)}. \]
Par conséquent, $\frac{\mu^2 A^{(+)}_{\tau_l}}{2S_{\tau_l}^2} \egaldistr T^{(3)}_2$.

%
\end{proof}

\begin{exercice}[Retour sur une représentation de Skorokhod]
Donner à l'aide de la théorie des excursions une démonstration du fait que $B_{T_\mu}$ a pour loi $\mu$, où on a posé
\[T_\mu = \inf \{t \geq 0 : S_t \geq \Psi_\mu(B_t)\},\]
et $\Psi_\mu(t) = \frac{1}{\mu([t, +\infty))}\int_t^{+\infty} x \mu(\d x)$ est la fonction de Hardy-Littlewood associée à $\mu$.
\end{exercice}

\begin{proof}
On réécrit $T_\mu$ à l'aide de $\Phi_\mu$, l'inverse continue à droite de la fonction de Hardy-Littlewood $\Psi_\mu$. On déterminera tout d'abord la loi de $S_{T_\mu}$, puis on en déduira la loi de $B_{T_\mu}$. On a
\begin{multline*}
  T_\mu  =  \inf\{t \geq 0 : S_t \geq \Psi_\mu(B_t)\}\\
   =  \inf\{t \geq 0 : \Phi_\mu(S_t) \geq B_t\}
   =  \inf\{t \geq 0 : S_t - B_t \geq S_t - \Phi_\mu(S_t)\}
\end{multline*}

Utilisons alors le théorème de Lévy
\[(S_t-B_t, S_t, t \geq 0 ) \egaldistr (|B_t|,L_t, t \geq 0),\]
on peut remplacer $T_\mu$ par
\[T'_\mu = \inf\{t \geq 0 : |B_t| \geq L_t-\Phi_\mu(L_t)\},\]
et on cherche à caractériser la loi de $L_{T'_\mu}$. Or
\[\P(L_{T'_\mu} \geq x) = \P(T'_\mu \geq \tau_x).\]
On étudie alors le nombre $N_x$ d'excursions avant l'instant $\tau_x$ sur lesquelles la condition $|B_t| \geq L_t-\Phi_\mu(L_t)$ est vérifiée à un moment, plus précisément :
\[N_x = \sum_{l \leq x} \ind{M^*(e_l) \geq l -\Phi_\mu(l)},\]
où on a posé $M^*(\epsilon) = \max_{0 \leq s \leq V(\epsilon)}|\epsilon(s)|$. C'est, grâce à la théorie des excursions une variable aléatoire de Poisson dont on peut calculer l'espérance :
\begin{align*}
 \E(N_x)  = & \int_0^x \d l \int \n(\d \epsilon) \ind{M^*(\epsilon)>l - \Phi_\mu(l)}\\
  = & \int_0^x \d l \int_0^{+\infty} \dfrac{\d m}{m^2} \ind{m>l-\Phi_\mu(l)}\\
  = & \int_0^x \dfrac{\d l}{l - \Phi_\mu(l)}.
\end{align*}
On a donc
\[\P(S_{T_\mu} \geq x) = \P(L_{T'_\mu} \geq x) = \P(N_x =0) = \exp\left( - \int_0^x \dfrac{\d l}{l-\Phi_\mu(l)}\right).\]

On peut maintenant calculer la fonction de répartition de $B_{T_\mu}$
\[
  \P(B_{T_\mu} \geq x) = \P(\Psi_\mu(B_{T_\mu} \geq \Psi_\mu(x))= \P(S_{T_\mu} \geq \Psi_\mu(x)).
\]
A partir de ce point, les calculs sont identiques à ceux réalisés dans l'Exercice~\ref{exo_representationSkorokhod}.
On note $\nu$ la loi de $B_{T_\mu}$, on a alors par différentiation
\[\nu(\d x) = \dfrac{\d \Psi_\mu(x)}{\Psi_\mu(x)-x} \nu([x,+\infty)),\]
soit
\[\nu(\d x)  \Psi_\mu(x) + \nu([x,+\infty))\d \Psi_\mu(x) = x \nu(\d x),\]
et cette équation différentielle sur $\Psi$ se résout de la façon suivante
\[\Psi_\mu(x) = \dfrac{1}{\nu([x,+\infty))} \int_x^{+\infty} t \nu(\d t) = \Psi_\nu(x).\]
Par injectivité de la transformation de Hardy-Littlewood, on en déduit $\nu=\mu$, et $B_{T_\mu} \sim \mu$.
\end{proof}

\begin{exercice}[Araignée brownienne]
On admet l'existence d'un processus $W$ appelé araignée brownienne (ou processus de Walsh) vivant sur la réunion de $N$ demi-droites possédant la même origine. Celui-ci se comporte comme un mouvement brownien sur chaque demi-droite, et à chaque retour en 0, choisit instantanément la nouvelle branche uniformément au hasard, indépendamment du passé.

Plus précisément, la mesure d'Itô des excursions de ce processus est donnée par :
\[\mathbf{w} = \dfrac{1}{N} \sum_{i=1}^N \n_i,\]
où $\n_i$ est l'image de $\n^+$ par l'application qui à une excursion positive associe l'excursion semblable sur la $i^\text{ième}$ branche de la toile.

Prouver l'extension suivante de la loi de l'arcsinus : la loi conjointe des temps passés dans chacune des branches jusqu'à l'instant 1 est donnée par  :
\[(A^1_1, \ldots, A^N_1) \egaldistr \left( \frac{T_1}{\sum_{i=1}^N T_i}, \cdots , \frac{T_N}{\sum_{i=1}^N T_i} \right),\]
où $T_1, \ldots, T_N$ sont des variables aléatoires i.i.d. stables de paramètre $\frac{1}{2}$.

Pour ce faire, on calculera :
\[ \Psi(\mu,(\lambda_i)) = \int_0^{+\infty} \d t e^{-\mu t} \E\left[\exp\left( - \dfrac{1}{{L^0_t}^2} \sum_{i=1}^N \lambda_i A^i_t\right)\right].\]
\end{exercice}

\begin{remarque}
Le processus $(|W_t|)_{t \geq 0}$ de la distance à l'origine de $W$ est un mouvement brownien réfléchi.
\end{remarque}

\begin{proof}
On observe pour commencer que pour tout entier $i \leq N$ et pour tout $l \geq 0$, on peut écrire le temps passé dans la $i^\text{ième}$ demi-droite $D_i$ jusqu'en l'instant $\tau_l$ est donné par $A^i_{\tau_l} = \sum_{\lambda \leq l} V(e_\lambda) \ind{e_\lambda \in D_i}$. De plus en utilisant les propriétés d'invariance d'échelle du mouvement brownien, on a, pour tout $t > 0$ :
\[(A^1_t, \ldots, A^N_t) \egaldistr t(A^1_1, \ldots, A^N_1).\]

On calcule alors 
\begin{align*}
 \Psi(\mu,(\lambda_i)) = & \int_0^{+\infty}\d t e^{-\mu t} \E\left[\exp\left( - \dfrac{1}{{L^0_t}^2} \sum_{i=1}^N \lambda_i A^i_t\right)\right]\\
  = & \E\left[\int_0^{+\infty} \d t \exp\left(- \mu t - \dfrac{1}{{L^0_t}^2} \sum_{i=1}^N \lambda_i \int_0^t \d s \ind{W_s \in D_i} \right) \right]\\
  = & \E\left[ \sum_l \int_{\tau_{l-}}^{\tau_l} \d t \exp\left(- \mu t - \dfrac{1}{l^2} \int_0^t \d s \sum_{i=1}^N \lambda_i \ind{W_s \in D_i} \right)\right]\\
  = & \E \left[ \sum_l \exp\left( - \mu \tau_{l-} - \frac{1}{l^2}\int_0^{\tau_{l-}} \d s \sum_{i=1}^N \lambda_i \ind{W_s \in D_i} \right) f(l,e_l)\right]
\end{align*}
où on a posé $f(\lambda,\epsilon) =  \int_0^{V(\epsilon)} \d t \exp\left( - \mu t - \frac{t}{\lambda^2} \sum_{i=1}^N \lambda_i \ind{\epsilon \in D_i} \right)$. Grâce à la formule additive des excursions, et par le théorème de Fubini, on obtient
\begin{equation*}
  \Psi(\mu,(\lambda_i))
  =\int_0^{+\infty} \d l \phi(l) \E\left[\exp\left( - \mu \tau_l - \frac{1}{l^2}\int_0^{\tau_l}\d s \sum_{i=1}^N \lambda_i \ind{W_s \in D_i}\right)\right],
\end{equation*}
où on a posé, $ \phi(l) =  \dfrac{1}{N}\sum_{i=1}^N \int_{\Omega^+}\n(\d \epsilon) \int_0^{V(\epsilon)}\d t \exp\left(- \left(\mu + \frac{\lambda_i}{l^2}\right) t\right)$.

On calcule d'une part, par décomposition de la mesure d'Itô,
\begin{align*}
  \phi(l)
  & =  \dfrac{1}{N} \sum_{i=1}^N \int_0^{+\infty} \dfrac{\d v}{2\sqrt{2\pi v^3}}\int_0^v \d t \exp\left(- \left(\mu +\frac{\lambda_i}{l^2}\right) t\right)\\
  & =  \dfrac{1}{N} \sum_{i=1}^N \dfrac{1}{\mu +\frac{\lambda_i}{l^2}}\int_0^{+\infty} \dfrac{\d v}{2 \sqrt{2\pi v^3}} \left( 1- \exp\left(- \left(\mu +\frac{\lambda_i}{l^2}\right) v\right) \right)\\
  & =  \dfrac{1}{N} \sum_{i=1}^N \dfrac{1}{\sqrt{2\left( \mu + \frac{\lambda_i}{l^2} \right)}} =  \dfrac{1}{\sqrt{2}N} \sum_{i=1}^N \dfrac{l}{\sqrt{l^2 \mu +  \lambda_i}}.
\end{align*}
D'autre part, on a par formule multiplicative,
\begin{align*}
   & \E\left[\exp\left( - \mu \tau_l - \frac{1}{l^2}\int_0^{\tau_l} \sum_{i=1}^N \lambda_i \ind{W_s \in D_i} \d s \right)\right]\\
   = & \exp\left[ - \frac{l}{N}\sum_{i=1}^N \int \n^+(\d \epsilon) \left[1- \exp\left(-\left(\mu + \frac{\lambda_i}{l^2}\right)V(\epsilon)\right) \right] \right]\\
   = & \exp\left[ -\frac{l}{N} \sum_{i=1}^N \int_0^{+\infty} \dfrac{\d v}{2\sqrt{2\pi v^3}} \left( 1-\exp\left( -\left(\mu +\frac{\lambda_i}{l^2}\right)v\right) \right) \right]\\
   = & \exp\left[ -\frac{1}{\sqrt{2}N} \sum_{i=1}^N \sqrt{l^2\mu +\lambda_i} \right]
   =  \exp[-\Phi(l)],
\end{align*}
où on a posé $\Phi(l) = \frac{1}{\sqrt{2}N} \sum_{i=1}^N \sqrt{l^2\mu +\lambda_i}$. On observe que $\Phi'(l) = \mu \phi(l)$.
Par conséquent, on obtient
\[  \Psi(\mu,(\lambda_i))= \frac{1}{\mu} \int_0^{+\infty} \mu \phi(l) \exp[-\Phi(l)]\d l = \frac{1}{\mu}\exp[-\Phi(0)].
\]
On en déduit
\[\int_0^{+\infty} e^{-\mu t} \E\left[\exp\left( - \dfrac{1}{{L^0_t}^2} \sum_{i=1}^N \lambda_i A^i_t\right)\right]\d t = \frac{1}{\mu} \exp\left[ - \frac{1}{\sqrt{2}\mu N} \sum_{i=1}^N \sqrt{\lambda_i}\right].\]

D'autre part, lorsqu'on prend $T_1, \ldots, T_N$ des variables aléatoires stables i.i.d. de paramètre $\frac{1}{2}$, on a
\[\int_0^{+\infty} e^{-\mu t} \E\left[ \exp\left( - \sum_{i=1}^N \lambda_i c^2 T_i\right) \right]\d t = \dfrac{1}{\mu} \exp\left(- c \sum_{i=1}^N \sqrt{\lambda_i}\right).\]
Par conséquent, en utilisant deux fois de suite l'injectivité de la transformée de Fourier et la propriété de scaling du mouvement brownien, on obtient
\[\left( \dfrac{A^1_t}{{L^0_t}^2}, \cdots, \dfrac{A^N_t}{{L^0_t}^2} \right) \egaldistr \dfrac{1}{2N^2} (T_1, \ldots, T_N).\]
Il suffit maintenant d'utiliser l'égalité triviale $\sum_{i=1}^N A^i_t = t$ pour conclure. En effet on obtient alors
\[\left( \dfrac{A^1_t}{{L^0_t}^2}, \cdots, \dfrac{A^N_t}{{L^0_t}^2}, \dfrac{t}{{L^0_t}^2} \right) \egaldistr \dfrac{1}{2N^2} \left(T_1, \ldots, T_N, \sum_{i=1}^N T_i\right),\]
d'où on conclut en particulier $\left( A^1_t, \ldots, A^N_t\right) \egaldistr t\left( \dfrac{T_1}{\sum_{i=1}^N T_i}, \cdots, \dfrac{T_N}{\sum_{i=1}^N T_i} \right)$.

\begin{remarque}
Il a été montré par Tsirel'son que pour tout $N \geq 3$, la filtration de l'araignée à $N$ branche est une filtration brownienne faible, mais pas forte. En d'autres termes,
\begin{itemize}
  \item toute martingale par rapport à la filtration $(\mathcal{A}^N_t)$ de l'araignée brownienne peut être représentée comme une intégrale stochastique par rapport au mouvement brownien $\beta$ vérifiant $|W_t| = \beta_t+L_t$,
  \item néanmoins, il n'existe pas de mouvement brownien $B$ tel que $\mathcal{A}^N$ est la filtration canonique de $B$.
\end{itemize}

Bien entendu, pour $N=2$, l'araignée à deux branches est un mouvement brownien, donc la filtration associée est fortement brownienne.
\end{remarque}
\end{proof}

\begin{exercice}[Représentation de Lévy-Khintchine des carrés de Bessel]
On rappelle que les lois de carrés de Bessel de dimension $\delta$ issu de $x$ vérifient la relation d'additivité
\[Q^{\delta+\delta'}_{x+x'} = Q^\delta_x \ast Q^{\delta'}_{x'}.\]
Par conséquent, il est naturel de dire que $Q^\delta_x$ est indéfiniment divisible.

Exprimer en fonction de $\n$, la mesure d'Itô du mouvement brownien, la mesure de Lévy associée à la loi $Q^\delta_x$. Plus précisément, on recherche deux mesures positives $M$ et $N$, sigma-finies sur $\mathcal{C}(\R^+,\R^+)$, vérifiant pour tout $\delta,x \geq 0$ et toute mesure $\mu$ sigma-finie sur $\R^+$ :
\begin{multline*}
  \qquad Q^\delta_x \left( \exp\left( -\int X_t \mu(\d t) \right) \right)\\
  = \exp\left[ -\int(x M+\delta N)(\d \omega)\left( 1-\exp\left(-\int \omega(t)\mu(\d t)\right)\right)\right].\qquad
\end{multline*}
\end{exercice}

\begin{proof}
On utilise la relation d'additivité : $Q^\delta_x = Q^\delta_0 \ast Q^0_x$ pour accéder séparément aux mesures $M$ et $N$.
Commençons par déterminer la mesure $M$ ; on s'intéresse à la loi $Q^0_x$. On utilise le théorème de Ray-Knight $(L^y_{\tau_l})_{y \geq 0}$ suit la loi de $Q^0_l$. Soit $f : \R^+ \to \R^+$ fonction continue à support compact, on calcule alors :
\[ \E\left(\exp\left(-\int_0^{+\infty} L^x_{\tau_l} f(x) \d x \right)\right) = \E\left( \exp\left( - \int_0^{\tau_l} f(B_s) \d s\right) \right).\]
Utilisant alors la formule multiplicative, on obtient :
\[\E\left( \exp\left( - \sum_{\lambda \leq l} \int_0^{V(e_\lambda)} f(e_\lambda(s)) \d s\right) \right) = \exp\left[ -l \int \n^+(\d \epsilon) \left(1- e^{-\int_0^{V(\epsilon)} f(\epsilon(s))\d s}\right) \right].\]

On définit alors le temps local d'une excursion $l^x_t$, qui est une semi-martingale --car une excursion normalisée suit la loi d'un pont de Bessel de dimension 3 de longueur 1-- on a alors :
\begin{multline*}
  \qquad\E\left(\exp\left(-\int_0^{+\infty} L^x_{\tau_l} f(x) \d x \right)\right)\\
  =  \exp\left( - l \int \n^+(\d \epsilon) \left(1 - \exp \left( - \int_0^{+\infty} l^x_{V(\epsilon)} f(x)\d x \right) \right)\right).\qquad
\end{multline*}

On peut étendre ce résultat à toute mesure positive $\mu(\d t)$ sur $\R^+$ ; on a :
\[Q^0_l\left( \exp\left( - \int X_t \mu(\d t) \right) \right) = \exp\left[ - l \int \n^+(\d \epsilon) \left( 1 - \exp \left( - \int_0^{+\infty} l^x_{V(\epsilon)} \mu(\d x) \right) \right) \right].\]
Par conséquent, $M$ est la mesure image de $\n^+$ par la fonction $\epsilon \mapsto (l^y_{V(\epsilon)}(\epsilon))_{y \geq 0}$, dont on connaît les marginales finies-dimensionnelles grâce  à l'Exercice \ref{exo_tempslocalexcursion}.

Pour la mesure $N$, on cherche une représentation de $Q^\delta_0$ à l'aide des théorèmes de Ray-Knight. L'Exercice \ref{exo_representationbessel} nous garantit que $\left(L^y_{\infty}\left(|B| +\frac{2}{\delta} L\right)\right)_{y \geq 0}$ suit la loi $Q^\delta_0$. On calcule alors, pour toute fonction $f : \R^+ \to \R^+$ continue à support compact :
\begin{align*}
  & \E\left[ \exp\left( - \int f(y) L^y_{+\infty} \left(|B|+\frac{2}{\delta}L\right)\d y \right) \right]\\
  = & \E\left[ \exp\left( - \int_0^{+\infty} f\left(|B_t| +\frac{2}{\delta} L_t\right) \d t \right) \right]\\
  = & \E\left[ \exp\left( - \sum_{\lambda \geq 0} \int_0^{V(e_\lambda)} f\left(|e_\lambda(t)| +\frac{2}{\delta} \lambda \right) \d t \right) \right]\\
  = & \exp\left[ - \int_0^{+\infty} \d l \int \n(\d \epsilon) \left[1- \exp\left( - \int_0^{V(\epsilon)} f\left(|\epsilon(t)| + \frac{2}{\delta}l\right)\d t \right)\right]\right].
\end{align*}
On obtient, après changement de variables, le résultat suivant :
\begin{multline*}
  \qquad\exp\left[ - \delta \int N(\d \omega) \left( 1 - e^{-\int \omega(t)f(t)\d t}\right) \right]\\
  = \exp\left[ - \delta \int_{\R^+} \d l \int \n^+ (\d \epsilon) \left( 1 - e^{-\int_0^{V(\epsilon)} f(\epsilon(t)+l) \d t} \right) \right].\qquad
\end{multline*}
Par conséquent, utilisant à nouveau le temps local des excursions, $N$ est la mesure image de $\d l \otimes \n^+$ par l'application $(l, \epsilon) \mapsto (l^y_{V(\epsilon)}(\epsilon + l))_{y \geq 0}$.

\begin{remarque}
On peut réécrire $N$ de la façon suivante : on pose $M_x$ la mesure obtenue de la même façon que $M$ à partir de la loi des excursions partant de $x$ au lieu de $0$. On a alors $N = \int_0^{+\infty} \d x M_x$.
\end{remarque}
\end{proof}

\begin{exercice}[Lien avec les processus de Lévy]
Soit $f \in L^1_\text{loc}(\R)$ et $g \in L^2_\text{loc}(\R)$.
\begin{enumerate}
  \item Montrer que les processus
\[\left( \int_0^{\tau_l} \d s f(B_s), l \geq 0\right)  \text{  et  } \left( \int_0^{\tau_l} \d s g(B_s)\d B_s, l \geq 0 \right)\]
sont des processus de Lévy.
  \item Exprimer leurs mesures de Lévy $m^{(f)}$ et $\tilde{m}^{(g)}$ comme images de la mesure d'It\^o.
  \item Montrer que tout subordinateur stable peut \^etre représenté sous la forme $\int_0^{\tau_l} \d s f(B_s)$ pour une certaine fonction $f$.
  \item Montrer que tout processus de Lévy stable symétrique peut être représenté sous la forme $\int_0^{\tau_l} g(B_s) \d B_s$.
\end{enumerate}
\end{exercice}

\begin{proof}
\textit{1.} Le fait que $X_l = \int_0^{\tau_l} \d s f(B_s)$ soit un processus de Lévy s'obtient en observant que ce processus est càdlàg, et à accroissements indépendants et stationnaires grâce aux propriétés du mouvement brownien. Le raisonnement est le m\^eme pour $Y_l = \int_0^{\tau_l} g(B_s) \d B_s$, gr\^ace aux propriétés de l'intégrale stochastique.

\textit{2.} Pour calculer la mesure de Lévy de ces processus, il suffit de calculer leur transformée de Laplace
\begin{align*}
  \E \left( \exp(-\mu X_l \right)
   = & \E\left[\exp\left( - \mu \sum_{\lambda \leq l} \int_{\tau_{\lambda-}}^{\tau_\lambda} \d s f(B_s) \right) \right]\\
   = & \E\left[\exp\left( - \mu \sum_{\lambda \leq l} \int_{0}^{V(e_\lambda)} \d s f(e_\lambda(s)) \right) \right]\\
   = & \exp\left( - \int l \n(\d \epsilon) \left( 1 - \exp\left( -\lambda \int_{0}^{V(\epsilon)} \d s f(\epsilon(s)) \right) \right) \right).
\end{align*}
Par conséquent, $m^{(f)}$ est l'image de $\n$ par l'application $\epsilon \mapsto \int_0^{V(\epsilon)} \d s f(\epsilon(s))$.

De la m\^eme manière, on obtient également que $\tilde{m}^{(g)}$ est l'image par $\n$ de l'application $\epsilon \mapsto \int_0^{V(\epsilon)} \d \epsilon(s) f(\epsilon(s))$, et cette application est bien définie, car $\epsilon$ est une semi-martingale.

\textit{3.} Un subordinateur stable $(S_l)$ d'exposant $\alpha$ est un processus de Lévy qui vérifie 
\[ \E(e^{-\lambda S_l}) = e^{-l \lambda^\alpha}.\]
La mesure de Lévy de ce subordinateur est donc
\[ \nu(\d x) = c_\alpha \dfrac{\d x}{x^{\alpha+1}},\]
où $c_\alpha$ est une constante à préciser.

On cherche $f$ sous la forme $c |x|^\gamma$. Le processus $(X_l)$ est alors bien un subordinateur (un processus de Lévy croissant), et pour vérifier que ce subordinateur est stable d'indice $\alpha$, il suffit de vérifier que $X_{al} \egaldistr a^\frac{1}{\alpha} X_l$.

Pour cela, on utilise la propriété de scaling du mouvement brownien, on a
\[  \int_0^{\tau_{al}} B_s^\gamma \d s
  \egaldistr  a^\gamma \int_0^{a^2 \tau_l} B_{\frac{s}{a^2}}^\gamma \d s   \egaldistr  a^{2+\gamma} \int_0^{\tau_l} B_s^\gamma \d s.
\]
Par conséquent, pour tout $\gamma >-1$, $(\int_0^{\tau_l} |B_s|^\gamma \d s)$ est un subordinateur stable d'indice $\frac{1}{2+\gamma}$, on obtient donc un subordinateur stable d'indice $1<\alpha<2$ en utilisant la fonction $f \mapsto |x|^{\frac{1}{\alpha}-2}$.

%
%

\textit{4.} On s'intéresse maintenant aux processus de Lévy stables symétriques de paramètres $\alpha$, c'est-à-dire vérifiant :
\[ \E\left( e^{i \lambda S_t} \right) = \exp(-tc|\lambda|^\alpha), \]
de mesure de Lévy $\frac{\d x}{|x|^{\alpha+1}}$. On conjecture que ceux-ci sont construit en utilisant une fonction $g(x) = |x|^\gamma$ pour $\gamma<\frac{1}{2}$. On utilise à nouveau la propriété de scaling du mouvement brownien pour conclure :
\begin{align*}
  \int_0^{\tau_{al}} |B_s|^\gamma \d B_s
  \egaldistr & a^{\gamma+1} \int_0^{a^2 \tau_l} |B_{\frac{s}{a^2}}|^\gamma \d B_{\frac{s}{a^2}}\\
  \egaldistr & a^{\gamma+1} \int_0^{\tau_l} |B_s|^\gamma \d B_s.
\end{align*}
Par conséquent, un processus de Lévy stable symétrique d'indice $0<\alpha<2$ est obtenu en utilisant la fonction $g : x \mapsto \frac{1}{\alpha}-1$.

On peut ainsi obtenir de très nombreux processus de Lévy, en utilisant le mouvement brownien. Tous ne sont néanmoins pas représentables de cette manière.

\begin{remarque}
La construction des processus de Lévy stables symétriques peut être obtenue en étendant l'intégrale d'une fonction du mouvement brownien en utilisant la notion de valeur principale. Formellement, en utilisant la formule d'Itô, on a
\[ \int_0^{\tau_l} |B_s|^\gamma \d B_s = -\frac{\gamma}{2} \int_0^{\tau_l} \sgn(B_s) |B_s|^{\gamma-1} \d s, \]
mais cette dernière intégrale n'est pas bien définie car la fonction $x \mapsto x^{\gamma-1}$ n'est pas nécessairement dans $L^1_\text{loc}$. En revanche, si on la définit comme :
\[ \lim_{\epsilon \to 0} \frac{-\gamma}{2} \int_0^t \sgn(B_s) |B_s|^{\gamma-1} \ind{|B_s|> \epsilon} \d s,\]
l'égalité précédente tient toujours. La limite existe, comme on peut s'en convaincre en transformant l'égalité grâce aux temps locaux du mouvement brownien :
\[ \frac{\gamma}{2}\int_\epsilon^{+\infty} x^{\gamma-1} (L^{-x}_t-L^x_t).\]
\end{remarque}
\end{proof}

\chapter[Formule de Feynman-Kac]{Représentations de l'intégrale de la mesure de Wiener et formule de Feynman-Kac}

La formule de Feynman-Kac, découverte de façon concomitante par chacun des deux auteurs R. Feynman et M. Kac est quelquefois appelée méthode de la transformée de Laplace double : en effet, elle donne un accès analytique à la transformée de Laplace d'une fonctionnelle additive du mouvement brownien, prise en un temps exponentiel indépendant. Nous montrerons ici, en nous inspirant fortement de \cite{JPY1997} comment les principales formules de la théorie des excursions browniennes permettent de comprendre en profondeur cette formule de Feynman-Kac

\section{Théorèmes principaux}

Nous allons commencer par définir quelques notations relatives aux loi de trajectoires de longueur finie. Soit $P$ la loi de $(X_u,u\leq T)$ et $P'$ celle de $(\tilde{X}_u,u \leq S)$ :
\begin{itemize}
  \item pour $A$ temps aléatoire, on note $P^{(A)}$ la loi de $(X_u,u\leq A \wedge T)$ ;
  \item $P\cdot P'$ est la loi de $(X_u, u \leq T ; \tilde{X}_{u-T}, T \leq u \leq T+S)$ ;
  \item pour finir, $\check{P}$ est la loi de $(X_{T-u},u \leq T)$.
\end{itemize}

Nous allons maintenant donner plusieurs représentations de l'intégrale sur $\R^+$ de la mesure de Wiener stoppée en chaque instant $t$, grâce à plusieurs mesures de trajectoires que nous avons déjà obtenu.

\begin{theoreme}
\label{the_representationIntegrale}
Les formules intégrales suivantes sont satisfaites :
\begin{align*}
  \int_0^{+\infty} \d t \W^{(t)} &= \left(\int_0^{+\infty} \d l \W^{(\tau_l)}\right) \cdot \left( \int_0^{+\infty} \d v \n^{(v)}(.,v <V(\epsilon)) \right) ,\\
  \int_0^{+\infty} \d l \W^{(\tau_l)} &= \int_0^{+\infty} \dfrac{\d u}{\sqrt{2\pi u}}\W_{0,0}^{(u)} ,
\end{align*}
De la même façon, on a
\begin{align*}
  \int_0^{+\infty} \d v \n^{+,(v)}(.,v <V(\epsilon)) &= \int_0^{+\infty} \d a {P^{(3)}_0}^{(\gamma_a)},\\
  \int_0^{+\infty} \d v \n^{+,(v)}(.,v <V(\epsilon)) &= \int_{-\infty}^{+\infty} \d a \check{\W}^{(T_0)}_a.
\end{align*}
\end{theoreme}

On va maintenant utiliser ces représentations intégrales pour obtenir la formule de Feynman-Kac.

\begin{theoreme}[Formule de Feynman-Kac]
Soit $q \in \mathcal{C}_c(\R)$, $f \geq 0$ mesurable. On pose :
\[K = \int_0^{+\infty} \d t e^{-kt} \E\left[q(B_t) \exp\left( - \int_0^t f(B_s) \d s \right) \right].\]

On a alors $K = \int_\R \d x q(x) U^{(f)}(x)$, où $U^{(f)}(x)$ est la seule solution de l'équation de Sturm-Liouville :
\[ u''(x) = \left( k + f(x) \right) u(x),\]
avec les conditions :
\begin{itemize}
  \item $u'$ existe et est uniformément bornée pour $x \neq 0$,
  \item $u$ s'annule en $\pm \infty$,
  \item $u'(0^+)-u'(0^-)=-2$
\end{itemize}
\end{theoreme}

\section{Exercices}

\begin{exercice}
\label{exo_integration}
Soit $R$ un processus de Bessel de dimension 3 et $J_t= \inf_{s \geq t} R_s$, montrer que :
\[\E \left[ \int_0^{+\infty} dJ_u F(R_s, s \leq u) \right] = \E\left[ \int_0^{+\infty} \dfrac{\d u}{2R_u}F(R_s,s \leq u) \right],\]
en utilisant les résultats sur le maximum de semi-martingales positives.
\end{exercice}

\begin{proof}
On montre l'égalité pour des processus simples $H_u=A \ind{u \in (t,t+s]}$, avec $A$ une variable aléatoire $\mathcal{R}_t$-mesurable. L'égalité à démontrer devient :
\[\E\left[A(J_{t+s}-J_t)\right] = \E \left[A \int_t^{t+s} \dfrac{\d u}{2R_u}\right].\]

Grâce à la propriété de Markov, on a $\E\left[ J_{t+s}-J_t|\F_t \right] = \E_{R_t}(J_s-J_0)$. De plus,
\[J_t = \inf_{s \geq t} R_s = \dfrac{1}{\sup_{s\geq t} \frac{1}{R_s}} \egaldistr U R_t.\]
Dès lors, pour tout $r \geq 0$,
\begin{align*}
\E_r\left[ J_{s}-J_0|\F_t \right] &= \frac{1}{2}\E_r(R_t-r)
= \frac{1}{2} \E_r\left[ \int_0^s \frac{\d u}{R_u} \right]
\end{align*}
grâce à l'équation différentielle stochastique satisfaite par $R$.

On en déduit
\[\E\left[ J_{t+s}-J_t|\F_t \right] = \E\left[\left. \int_0^s \frac{\d u}{2R_u}\right|\F_t \right], \]
donc, pour $H_u = A \ind_{u \in (t,t+s]}$, on a
\[\E \left[ \int_0^{+\infty} dJ_u F(R_s, s \leq u) \right] = \E\left[ \int_0^{+\infty} \dfrac{\d u}{2R_u}F(R_s,s \leq u) \right],\]
et on conclut en appliquant le lemme des classes monotones.
\end{proof}

\begin{exercice}
Montrer la troisième égalité du Théorème \ref{the_representationIntegrale} :
\[\int_0^{+\infty} \d v {\n^+}^{(v)}(.,v <V(\epsilon)) = \int_0^{+\infty} \d a {P^{(3)}_0}^{(\gamma_a)}(.).\]
\end{exercice}

\begin{proof}
Soit $F$ une fonctionnelle mesurable positive, on commence par utiliser la décomposition de $\n^+$ en fonction de la longueur d'une excursion, pour calculer :
\[\int_0^{+\infty} \d v {\n^+}^{(v)}\left(F(\epsilon_u,u\leq v)\ind{v <V(\epsilon)}\right) = \int_0^{+\infty} \d v \int_v^{+\infty} \dfrac{\d t}{2\sqrt{2\pi t^3}} {\Pi^t}^{(v)}(F(R_u, u \leq v)),\]
où $\Pi^t$ est la loi du pont de Bessel de longueur $t$. Grâce à la relation d'absolue continuité entre le pont de Bessel et le processus de Bessel, on peut encore réécrire cette intégrale :
\[\int_0^{+\infty} \d v \int_v^{+\infty} \dfrac{\d t}{2\sqrt{2\pi t^3}} P^{(3)}_0\left[\left(\frac{t}{t-v}\right)^{\frac{3}{2}} \exp\left( -\frac{R_v^2}{2(t-v)} \right) F(R_u, u \leq v) \right].\]

On utilise alors la formule de Fubini-Tonelli, et par changement de variables :
\[\int_0^{+\infty} \d v P^{(3)}_0\left[ F(R_u, u \leq v) \int_0^{+\infty} \dfrac{\d t}{2\sqrt{2\pi t^3}} \exp\left( -\frac{R_v^2}{2t} \right) \right].\]
on obtient alors, en intégrant par rapport à $\d t$
\[\int_0^{+\infty} \d v {\n^+}^{(v)}\left(F(\epsilon_u,u\leq v)\ind{v <V(\epsilon)}\right) = P^{(3)}_0\left[\int_0^{+\infty} \dfrac{\d v}{2R_v} F(R_u, u \leq v) \right].\]

Par la suite, nous utilisons le résultat de l'Exercice \ref{exo_integration}, cette intégrale peut encore se réécrire en l'intégrale par rapport à $J_u$ suivante :
\[P^{(3)}_0\left[\int_0^{+\infty}\dfrac{\d v}{2R_v} F(R_u, u \leq v) \right] = P^{(3)}_0\left[ \int_0^{+\infty} dJ_v F(R_u, u \leq v) \right],\]
il ne reste plus qu'à utiliser un changement de variables pour obtenir :
\[\int_0^{+\infty} \d v {\n^+}^{(v)}\left(F(\epsilon_u,u\leq v)\ind{v <V(\epsilon)}\right) = \int_0^{+\infty} \d a P^{(3)}_0\left[ F(R_u, u \leq \gamma_a) \right],\]
ce qui nous permet d'obtenir l'égalité.
\end{proof}

\begin{exercice}
\begin{enumerate}
  \item Montrer que 
\[\E(H_{\gamma_a}|\gamma_a=t) = \E(H_t|R_t=a) \text{  }\d tda-\text{p.p.}\]
pour tout processus $(\mathcal{R}_u)$-prévisible borné $H$.
  \item Donner une expression semblable de $\E(H_{T_a}|T_a=t)$ dans le cas du mouvement brownien.
\end{enumerate}
\end{exercice}

\begin{proof}
\textit{1.} Soit $H$ un processus $(\mathcal{R}_u)$-prévisible borné, et $\phi(t,x)$ une fonction continue bornée. On calcule pour commencer :
\[\E\left(\int_0^{+\infty} \d a \phi(\gamma_a, R_{\gamma_a}) H_{\gamma_a}\right) = \E\left(\int_0^{+\infty} dJ_u \phi(u, R_u) H_u\right),\]
par changement de variables $\gamma_a=u$. Ensuite, en utilisant l'Exercice \ref{exo_integration}, on obtient
\[\E\left(\int_0^{+\infty} \d a \phi(\gamma_a, a) H_{\gamma_a}\right) = \E\left( \int_0^{+\infty} \dfrac{\d u}{2R_u} \phi(u,R_u) H_u \right).\]

D'autre part, en conditionnant par rapport à $\gamma_a$ ou $R_u$ on obtient, par formule de Fubini :
\[\E\left[ \int_0^{+\infty} \d a \phi(\gamma_a,a) \E(H_{\gamma_a}|\gamma_a)\right] = \E\left[ \int_0^{+\infty} \dfrac{\d u}{2R_u} \phi(u,R_u) \E(H_u|R_u) \right].\]
On réalise alors un nouveau changement de variables $a=J_u$, on obtient alors :
\[\int_0^{+\infty} \d a \E\left[ \phi(\gamma_a,a) \E(H_{\gamma_a}|\gamma_a)\right] = \int_0^{+\infty} \d a \E\left[ \phi(\gamma_a,a) \E(H_u|R_u=a) \right].\]
Cette égalité étant valable pour toute fonction $\phi(\gamma_a,a)$, on a pour $\lambda$-presque tout couple $(t,a)$ :
\[\E(H_{\gamma_a}|\gamma_a=t) = \E(H_t|R_t=a).\]

\textit{2.} Nous nous intéressons maintenant à $\E(H_{T_a}|T_a=t)$ pour le mouvement brownien. On réalise alors les même calculs que précédemment en utilisant que $(S_t)$ est l'inverse continu à droite de $(T_a)$. Pour toute fonction $\phi$ continue bornée, on a
\[\E\left(\int_0^{+\infty} \d a \phi(T_a, a) H_{T_a}\right) = \E\left(\int_0^{+\infty} dS_t \phi(t, S_t) H_t\right).\]

On utilise alors le théorème d'équivalence de Lévy, et l'Exercice \ref{exo_integrationtempslocal}, 
\[ \E\left(\int_0^{+\infty} dS_t \phi(t,S_t)H_t\right) = \int_0^{+\infty} \dfrac{\d u}{\sqrt{2\pi u}} \E(\phi(u,S_u) H_u|S_u-B_u=0).\]
On peut alors, en réalisant à nouveau le changement de variables $t = T_a$, par conditionnement par rapport à $S_t$
\[\E\left(\int_0^{+\infty} \d a \phi(T_a, a) H_{T_a}\right) = \E\left(\int_0^{+\infty} \d a \phi(T_a, a) \E(H_t|S_t-B_t=0,B_t=a)\right).\]
Cette égalité étant valable pour toute fonction $\phi(T_a,a)$, on a pour $\lambda$-presque tout couple $(t,a)$ :
\[\E(H_{T_a}|T_a=t) = \E(H_t|S_t=B_t=a).\]
\end{proof}

\begin{exercice}
Soit $f$ fonction mesurable positive, on note $X_t$ le processus de Markov tué en $T^{(f)}$, vérifiant pour tout $q \in \mathcal{C}_c$ :
\[\E(q(X_t) \ind{T^{(f)}>t} ) = \E\left[q(B_t) \exp\left(- \int_0^t f(B_s)\d s\right) \right].\]
\begin{enumerate}
  \item Donner une expression de la transformée de Laplace de $T^{(f)}$.
  \item Calculer cette expression de façon explicite lorsque $f(x) = c\ind{x \geq 0}$.
\end{enumerate}
\end{exercice}

\begin{proof}
\textit{1.} On souhaite calculer, pour $\lambda \geq 0$
\[\E(\exp(-\lambda T^{(f)})) = \lambda \int_0^{+\infty} \d t \exp(-\lambda t) \P(T^{(f)} \leq t).\]
On observe que :
\[\E(\exp(-\lambda T^{(f)})) = 1 - \int_0^{+\infty} \d t \lambda e^{-\lambda_t} \P(T^{(f)} > t).\]
Soit $q^n$ une suites de fonctions continues à support compact convergeant en croissant vers 1. Par convergence monotone, on a :
\[\int_0^{+\infty} \d t \lambda e^{-\lambda_t} \P(T^{(f)} > t) = \int_\R \d x u^{(\lambda)}(x),\]
où $u^{(\lambda)}$ est la solution de l'équation de Sturm-Liouville $u''= (\lambda + f) u$.
En particulier, on obtient :
\[ \P(T^{(f)}=+\infty) = \int_\R \d x u^{(0)}(x)\]

\textit{2.} Il reste donc maintenant à déterminer les solutions de l'équation de Sturm-Liouville associée à $f = c \ind{x\geq 0}$. Soit $\lambda \geq 0$, les solutions de 
\[u'' = \lambda u\]
sont de la forme $u(x) = A e^{\sqrt{\lambda}x} + B e^{-\sqrt{\lambda} x}$. Étant données les contraintes recherchées, on obtient, pour $f(x) = c\ind{x \geq 0}$ :
\[u^{(\lambda)}(x) = \dfrac{2}{\sqrt{\lambda+c}+\sqrt{\lambda}} \left[ e^{-\sqrt{\lambda +c} x}\ind{x\geq 0} + e^{\sqrt{\lambda}x}\ind{x \leq 0} \right],\]
d'où on tire :
\[\E(\exp(-\lambda T^{(f)})) = \dfrac{2}{\sqrt{\lambda(\lambda+c)}}.\]
\end{proof}

\chapter[Diffusions linéaires]{Temps locaux et Excursions de diffusions linéaires}

Dans ce dernier chapitre, on étend les définitions des temps locaux, et on établit une théorie des excursions pour les diffusions réelles. Cette classe de processus de Markov nous semble représenter une famille adéquate pour de telles généralisations, par rapport au cadre des semi-martingales d'une part, et du mouvement brownien d'autre part. Il existe toutefois de nombreuses études de temps locaux et d'excursions dans des cadres Markoviens beaucoup plus généraux, mais ceci dépasse les limites de notre propos...

On a déjà vu que les temps locaux de semi-martingales peuvent se comporter de manière contre-intuitive dans certains cas : ainsi le temps local en chaque point d'un processus à variations finies est nul. Nous définissons ici les temps locaux de diffusion, par certains côtés plus intuitifs, qui reflétera mieux certains aspects du temps infinitésimal passé à un niveau donné. De façon analogue nous pourrons étendre la théorie des excursions 
aux solution d'EDS à \og bons coefficients \fg{} .

\section{Théorèmes principaux}

Soit $X^x$ la solution de l'équation différentielle stochastique
\[X^x_t = x + \int_0^t \sigma(X_s) \d B_s + \int_0^t b(X_s)\d s. \]

Nous allons pour commencer définir la notion de fonction d'échelle.

\begin{definition}
On note $T_a = \inf\{t \geq 0 : X_t = a\}$, pour $a \in \R$. On dit que $h$ est une fonction d'échelle associée à $X$ si pour tous $a < x < b$, on a
\[\P_x(T_a \leq T_b) = \frac{h(b)-h(x)}{h(b)-h(a)}.\]

Si une diffusion admet l'identité pour fonction d'échelle, on dit que cette diffusion est \og en échelle normale \fg{}.
\end{definition}

\begin{proposition}
Les fonctions d'échelles sont définies à affinité près.

Une fonction d'échelle est monotone, peut donc être choisie croissante.

Une fonction d'échelle associée à $X$ est donnée par :
\[ h(x) = \int_0^x \d y \exp\left( - \int_0^y  \d z \frac{2b(z)}{\sigma^2(z)} \right).\]
\end{proposition}

\begin{theoreme}
La fonction $h$ est une fonction d'échelle associée à $X$ si et seulement si $h(X)$ est une martingale locale. En particulier $h(X)$ est une diffusion en échelle normale.
\end{theoreme}

\begin{definition}
Soit $Y$ une diffusion en échelle normale et $W$ un mouvement brownien de temps local $L$. Pour $\mu$ mesure de Radon et $t \geq 0$, on pose $A ^\mu_t = \int_\R \mu (\d x) L^x_t$. L'inverse continu à droite de ce processus croissant est noté $\gamma^\mu$.

La mesure de vitesse associée à $Y$ est la mesure de Radon $m$ qui vérifie l'égalité en loi suivante :
\[(Y_t, t \geq 0) \egaldistr (W_{\gamma^m_t}, t \geq 0).\]

Si $X$ est une diffusion réelle, et $h$ sa fonction d'échelle, la mesure de vitesse associée à $X$ est l'image par $h$ de la mesure de vitesse associée à $h(X)$.
\end{definition}

\begin{propriete}
Le semi-groupe de $X$ est absolument continu par rapport à la mesure de vitesse $m$, et la dérivée de Radon-Nikodym est symétrique. En d'autres termes :
\[ \P_x(X_t \in \d y) = p_t(x,y)m(\d y) \mathrm{  avec  } p_t(x,y)=p_t(y,x). \]
\end{propriete}

\begin{theoreme}
Il existe une famille $(l^x_t,x \in \R, t \geq 0)$ conjointement continue en $t$ et en $x$, telle que pour toute fonction $f$ positive mesurable :
\[ \int_0^t f(X_s) \d s = \int_\R m(\d x) l^x_t f(x). \]

Cette famille est appelée les temps locaux de diffusion de $X$.
\end{theoreme}

Nous allons maintenant nous intéresser à une théorie des excursions pour les diffusions réelles. Désormais, les diffusions considérées seront récurrentes, c'est-à-dire que $\P_x(T_y<+\infty)=1$, pour tous $x,y$. On note $(\tau_l,l \geq 0)$ l'inverse continu à droite de $(l^0_t,t \geq 0)$. On pose alors $e_l^X(s) = X_{(\tau_{l-}+s)\wedge \tau_l}$.

\begin{theoreme}[Théorème d'It\^o]
Le processus $(e_l)_{l \geq 0}$ est un processus de Poisson ponctuel, i.e. pour tout choix d'ensembles mesurables disjoints $\Gamma_1, \ldots, \Gamma_k$ de $\Omega^*$, les processus à valeurs entières 
\[N^{\Gamma_j}_l = \sum_{\lambda \leq l} \ind{e_\lambda \in \Gamma_j}\]
sont soit infinis, soit des processus de Poisson indépendants de paramètre $\n^X(\Gamma_j)$.

La mesure $\n^X$ est sigma-finie sur $\Omega^*$. On l'appelle la mesure d'Itô de la diffusion~$X$. Pour tout $l \geq 0$, on écrit $X_l = (B_t,t \leq \tau_l)$. Le processus de Poisson pour la concaténation $(X_l, l \geq 0)$ a pour mesure de Lévy $\n^X$.
\end{theoreme}

Nous allons maintenant donner deux décompositions possible de la mesure d'Itô, en fonction de la longueur ou de la hauteur des excursions, comme dans le cas du mouvement brownien. Pour cela, il est nécessaire de définir un certain nombre de mesures sur les trajectoires.

On note $\P_{x,y}$ la loi de $(X_u,u \leq T_y)$ la loi de la diffusion issue de $x$ tuée en son premier temps d'atteinte de $y$. Pour $u \geq 0$, on note $\P_{x,y,u}$ la loi $\P_{x,y}$ conditionnée à $T_y=u$. C'est donc une diffusion issue de $x$, qui, après $u$, atteint $y$ pour la première fois.

On note en particulier $\P_{x,0,u} = \lim_{y \to 0} \P_{x,y,u}$, et pour finir $\Pi_u = \lim_{x \to 0} \P_{x,0,u}$. De manière heuristique, cette dernière loi est celle de la diffusion $X$ issue de 0, conditionnée à retourner en 0 exactement à l'instant $u$. C'est l'exact analogue du pont de Bessel 3 dans le cas du mouvement brownien.

De la même manière, on va construire une mesure sur les excursions de hauteur $m$. Pour cela, on note $\widehat{P}_x$ la loi de la diffusion issue de $x$ conditionnée à ne jamais retourner en 0. Cette loi est un bon analogue du processus de Bessel 3. On note enfin $\widehat{\Pi}_m$ la loi de deux processus suivant la loi $\widehat{P}_0$, arrêtés en leurs premiers temps d'atteinte de $m$ respectifs, et mis dos-à-dos.

On obtient, gr\^ace à ces lois, les décompositions d'Itô et de Williams de la mesure d'excursion de $X$.
\begin{theoreme}[Décomposition d'It\^o et de Williams]
La mesure d'It\^o de $X$ peut se décomposer des deux manières suivantes :
\[\n^X(\d \epsilon) = \int_0^{+\infty} n_V(\d v) \Pi_v(\d \epsilon) = \int_0^{+\infty} n_M(\d m) \widehat{\Pi}_m(\d \epsilon),\]
où $n_V$ et $n_m$ sont les mesures images de $\n$ par $\epsilon \mapsto V(\epsilon)$ et $\epsilon \mapsto M(\epsilon)$.
\end{theoreme}

\section{Exercices}

\begin{exercice}
Soit $X$ une solution d'équation différentielle stochastique.
\begin{enumerate}
  \item Trouver $\phi$ telle que $(\phi(X_t), t \geq 0)$ est une martingale locale.
  \item Trouver alors une fonction d'échelle $h$ associée à $X$.
  \item Réciproquement, montrer que si $h$ est une fonction d'échelle associée à $X$, alors $h(X_t)$ est une martingale locale.
\end{enumerate}
\end{exercice}

\begin{proof}
\textit{1.} On va utiliser la formule d'Itô pour déterminer une fonction $\phi$ de classe $\mathcal{C}^2$ telle que $\phi(X)$ est une martingale locale. On a :
\[ \phi(X_t) = \phi(x) +  \int_0^t \phi'(X_s) \sigma(X_s) \d B_s + \int_0^t \phi'(X_s) b(X_s) + \phi''(X_s) \frac{\sigma^2(X_s)}{2} \d s. \]
Par conséquent, on choisit pour $\phi$ une solution de l'équation différentielle :
\[\frac{\sigma^2}{2}\phi'' + b\phi' = 0.\]

\textit{2.} Sur $[0,T_a \wedge T_b]$, $X$ est borné donc $\phi(X)$ l'est également. On peut appliquer le théorème d'arrêt à cette martingale, on obtient :
\[\phi(a) \P_x(T_a \leq T_b) + \phi(b)(1-\P_x(T_a \leq T_b) = \E( \phi(X_{T_a \wedge T_b} ) = \phi(x),\]
d'où on obtient :
\[\P_x(T_a \leq T_b) = \frac{\phi(b)-\phi(x)}{\phi(b)-\phi(a)},\]
la fonction $\phi$ est donc bien une fonction d'échelle associée à $X$.

\textit{3.} Soit $h$ une fonction d'échelle associée à $X$, on pose $Y_t=h(X_t)$, et pour $a \in \R$, on pose :
\[S_a = \inf\{t \geq 0 : Y_t = a\} \text{  et  } S_{a,b} = S_a \wedge S_b.\]

Soit $a<y<b$, montrons que $(Y_{t \wedge S_{a,b}}, t \geq 0)$ est une martingale. On utilise la propriété de Markov pour calculer l'espérance conditionnelle de ce processus :
\[ \E_y(Y_{t+s \wedge S_{a,b}}|\F_t) = \E_{Y_{t\wedge S_{a,b}}}(Y_{s\wedge S_{a,b}}),\]
il suffit donc de prouver que l'espérance de ce processus est constante, i.e. que pour tout $t \geq 0$, on a $\E_y(Y_{t\wedge S_{a,b}})=y$. On pose alors $a = h(\alpha)$, $b=h(\beta)$ et $y=h(x)$, on a alors :
\[ \E_y(Y_{t\wedge S_{a,b}}) = \E_x(h(X_{t\wedge T_{\alpha,\beta}})).\]

Lorsqu'on fait tendre $t$ vers $+\infty$, on obtient :
\[\E_x(h(X_{T_\alpha,\beta}) = \frac{h(\beta)-h(x)}{h(\beta)-h(\alpha)} h(\alpha) + \frac{h(x)-h(\alpha)}{h(\beta)-h(\alpha)} h(\beta) = h(x) = y.\]
De plus, cette espérance est décroissante au cours du temps, ceci indique qu'elle est constante, donc que le processus $(Y_t, t \geq 0)$ est une martingale locale.
\end{proof}

\begin{exercice}
Soit $X$ la solution de l'équation différentielle stochastique :
\[ X_t = \int_0^t \sigma(X_s) \d B_s + \int_0^t b(X_s) \d s, \]
on note $h$ une fonction d'échelle associée à $X$ et $Y_t=h(X_t)$.

\begin{enumerate}
  \item Déterminer la mesure de vitesse associée à $Y$.
  \item En déduire la mesure de vitesse associée à $X$. En particulier, donner la mesure de vitesse associée au processus de Bessel de dimension $d$.
\end{enumerate}
\end{exercice}

\begin{proof}
\textit{1.} Soit $W$ un mouvement brownien, et $\rho$ une fonction mesurable positive, on calcule $A^\mu_t$ pour $\mu(\d x)= \rho(x)\d x$. On a :
\[A^\mu_t = \int_0^t \rho(W_s) \d s.\]

En utilisant la formule d'Itô, $Y$ satisfait l'équation différentielle stochastique
\[ Y_t = \int_0^t h'(h^{-1}(Y_s)) \sigma(h^{-1}(Y_s)) \d B_s = \int_0^t \tilde{\sigma}(Y_s)\d B_s.\]
En particulier,
\[\crochet{Y}_t = \int_0^t \tilde{\sigma}(Y_s)^2 \d s.\]

D'autre part, on a $\crochet{W}_{\gamma^\mu_t}=\gamma^\mu_t$, on souhaite donc trouver $\mu$ tel que :
\[ \gamma^\mu_t = \int_0^t \tilde{\sigma}(Y_s)^2 \d s,\]
autrement dit, $W$ est le mouvement brownien de Dubins-Schwarz associé à $Y$. Supposons que $m$ existe, et que $m(\d x)=\rho(x)\d x$, on obtient :
\[t = \int_0^{A_t} \tilde{\sigma}(W_{\gamma_s})^2\d s = \int_0^t \tilde{\sigma}(W_u)^2 \rho(W_u) \d u,\]
d'où on déduit $\rho = \frac{1}{\tilde{\sigma}^2}$.
Il suffit maintenant de vérifier que pour ce $\rho$, l'égalité en loi est bien satisfaite, ce qui est automatique (c'est la transformée réciproque de celle de Dubins-Schwarz).

\textit{2.} La mesure de vitesse associée à $X$ est alors la mesure image de $m$ par $h$, on en conclut que cette mesure est également à densité par rapport à la mesure de Lebesgue, et de densité $\frac{1}{h'(x)\sigma^2(x)}$.
En particulier, un processus de Bessel de dimension $d$ est la solution de l'équation différentielles stochastique suivante
\[X_t = B_t + \int_0^t \dfrac{d-1}{2 X_s} \d s.\]

Une fonction d'échelle associée à $X$ est donc de la forme :
\[ h(x) = \left\{
\begin{array}{ll}
  x & \text{  si } d = 1\\
  \log(x) & \text{  si } d = 2 \\
  \frac{1}{x^{d-2}} & \text{  sinon.}
\end{array}
\right.\]
Par conséquent, la mesure de vitesse associée au processus de Bessel de dimension $d$ est donnée par $m(\d x) = x^{d-1}\d x$.
\end{proof}

\begin{exercice}[Lien entre les différentes notions de temps locaux]
Soit $X, h, Y$ définis comme précédemment.

\begin{enumerate}
  \item Montrer l'existence de temps locaux de diffusions pour $Y$, et les exprimer en fonction des temps locaux de semi-martingales associés à $Y$.
  \item Déterminer une expression semblable pour les temps locaux de diffusion associés à $X$.
  \item La théorie des processus de Markov permet de prouver qu'il existe, à une constante multiplicative près, une unique fonctionnelle additive continue $A$ associée à un processus de Markov $X$ dont l'ensemble des instants de croissance est p.s. inclus dans $\{t \geq 0 : X_t=x\}$. $A$ est appelé temps local markovien de $X$ en $x$ si et seulement si il vérifie pour tout $y \in \R$ :
\[ \E_y \left( \int_0^{+\infty} \d A_s e^{-s} \right) = \E_y(e^{-T_x}).\]
Donner un lien entre les temps locaux de processus de Markov, les temps locaux de diffusions et les temps locaux de semi-martingales de la diffusion $X$.
\end{enumerate}
\end{exercice}

\begin{proof}
\textit{1.} On observe dans un premier temps que pour toute fonction $f$ mesurable positive, on a
\[
  \int_0^t f(Y_s) \d s =  \int_0^t \frac{f(Y_s)}{\tilde{\sigma}(Y_s)^2} \d \crochet{Y}_s =  \int_\R \d x L^x_t \frac{f(x)}{\tilde{\sigma}(x)^2} =  \int_\R m(\d x) L^x_t f(x).
\]
Par conséquent, pour les diffusions en échelle naturelle, on a $L^x_t=l^x_t$, les temps locaux de semi-martingales et de diffusions sont les mêmes.

\textit{2.} On s'intéresse maintenant à une diffusion arbitraire $X$ :
\[  \int_0^t f(X_s) \d s  =  \int_0^t \frac{f(X_s)}{\sigma(X_s)^2} \d \crochet{X}_s\\
   =  \int_\R \d x L^x_t \frac{f(x)}{\sigma(x)^2} \\
   =  \int_\R m(\d x) L^x_t f(x) h'(x).
\]
Donc $l^x_t = h'(x)L^x_t$ pour une solution d'équations différentielle stochastique générale. En d'autres termes les temps locaux de diffusion prennent en compte le temps passé au voisinage de $x$ en échelle naturelle.

\textit{3.} Pour étudier le lien entre temps local de diffusion et temps local de processus de Markov, il suffit d'utiliser le résultat suivant : pour toute fonction mesurable positive $f$, on a :
\[ \int_0^t f(s, X_s) \d s = \int m(\d x) \int_0^t \d_sl^x_s f(s,x). \]

En prenant l'espérance de cette égalité pour des fonctions $f$ bien choisies, on obtient :
\[ \int_\R \d s e^{-s} p_s(x,y) = \E_y\left(\int_0^t e^{-s} \d_sl^x_s\right).\]

On pose $R(x,y) = \int_0^{+\infty} \d s e^{-s} p_s(x,y)$, il reste donc à rappeler que 
\[\E_y(e^{-T_x}) = \dfrac{R(x,y)}{R(x,x)},\]
pour déduire $l^x_t = R(x,x) L^{(M) x}_t$. Un résultat similaire peut s'obtenir avec les temps locaux de semi-martingales et la résolvante associée au processus de Markov.
\end{proof}

\begin{exercice}
Soit $u$ et $v$ deux fonctions de classe $\mathcal{C}^1$, $v$ strictement croissante, et $B$ un mouvement brownien.
\begin{enumerate}
  \item Montrer que $X = (u(t)B_{v(t)}, t \geq 0)$ est une semi-martingale.
  \item En déduire une relation entre $L^0_t(X)$ et $L^0_t(B)$.
  \item Développer une théorie des excursions pour $X$.
  \item Application au pont brownien $\beta_t = (1-t)B_\frac{t}{1-t}$.
  \item Étudier de la m\^eme manière le processus $Y_t = u(g_t)B_t$.
\end{enumerate}
\end{exercice}

\begin{proof}
\textit{1.} On s'intéresse dans un premier temps à $Z_t = X_{v^{-1}(t)} = u(v^{-1}(t))B_t$. Par intégration par parties, on obtient :
\[ Z_t = \int_0^t  u(v^{-1}(s))\d B_s + \int_0^t B_s \frac{u(v^{-1}(s))}{v'(v^{-1}(s))} \d s. \]
Par changement de temps déterministe,
\[ X_t = \int_0^t u(s) \d B_{v(s)} + \int_0^t B_{v(s)} u'(s) \d s \]
donc $X$ est une semi-martingale, dont la partie martingale est $\int_0^t u(s)\d B_{v(s)}$, de crochet $\int_0^t u(s)^2v'(s)\d s$.

\textit{2.} On calcule ensuite le temps local en 0 de $Z$. On a :
\[|Z_t| = \int_0^t \sgn(Z_s)dZ_s + L^0_t(Z).\]
D'un autre côté, en intégrant par parties, on a :
\begin{multline*}
 |u(v^{-1}(t))| |B_t|  =  \int_0^t |u(v^{-1}(s))| \sgn(B_s) \d B_s + \int_0^t |u(v^{-1}(s))| \d L^0_s(B)\\
 + \int_0^t |B_s| \sgn(u(v^{-1}(s))) \d u(v^{-1}(s)),
\end{multline*}
Dès lors,
\[L^0_t(Z) = \int_0^t \d_sL^0_s(B)u(v^{-1}(s)).\]
Par changement de variables déterministe, on a également
\[ L^0_t(X) = \int_0^{v(t)} \d_sL^0_s(B) u(v^{-1}(s)). \]

\textit{3.} Pour développer une théorie des excursions, on étudie pour commencer la transformation d'une excursion du mouvement brownien $B$ commençant en l'instant $v(s)$. L'excursion correspondante du processus $X$ commence en $s$, et est donnée par :
\[i^X_s(u) = u(s+u)i^B_{v(s)}(v(s+u)).\]
On note $f(s,\epsilon)$ la fonction qui à une excursion de $B$ à l'instant $v(s)$ associe l'excursion de $X$ à l'instant $s$. Soit $F(s, \omega, \epsilon)$ une fonctionnelle mesurable positive, on calcule grâce à la formule additive modifiée des excursions :
\begin{align*}
  \E\left[ \sum_{s \leq t} F(s, \omega, i_s^X) \right] & =  \E\left[ \sum_{s \leq t} F(s, \omega,f(s, i_{v(s)}^B)) \right]\\
  & =  \E\left[ \sum_{u \leq v(t)} F(v^{-1}(u), \omega, f(v^{-1}(u), i_u^B))\right]\\
  & =  \E\left[ \int_0^{v(t)} \d_sL^0_s(B) \int \n(\d \epsilon) F(v^{-1}(s), \omega, f(v^{-1}(s),\epsilon)) \right]\\
  & =  \E\left[ \int_0^t \frac{\d_sL^0_s(X)}{u(v^{-1}(s))}\int \n(\d \epsilon) F(s,\omega, f(s,\epsilon)) \right]
\end{align*}

Cette formule additive modifiée obtenue pour le processus $X$ permet également d'écrire une formule additive plus \og classique \fg{} par simple changement de variables :
\[ \E \left[ \sum_{\lambda \leq l} F(\lambda, \omega, e^X_\lambda) \right] = \E\left[ \int_0^l \d u \int \n(\d \epsilon) \frac{F(u, \omega, f(\tau^X_u,\epsilon))}{u(v^{-1}(\tau^X_u))} \right], \]
ce qui nous indique en particulier que la \og mesure de Lévy instantanée \fg{} de $X$ en $s$ est la mesure image de $\n$ par $f(s,.)$.

On pourrait également développer une formule multiplicative adaptée au processus $X_t$, pour cela, on pose $\Phi$ une fonctionnelle positive mesurable, et on note
\[F_l = \exp\left(-\sum_{\lambda  \leq l} \Phi(e^X_\lambda)\right).\]
On a
\begin{align*}
 \E(F_l)  = & 1 + \E\left[ \sum_{\lambda \leq l} F_{\lambda-}\left( 1 - e^{-\Phi(e^X_\lambda)}\right)\right]\\
  = & 1 + \E\left[ \int_0^l \d u F_u \int \n(\d \epsilon) \frac{1-e^{-\Phi(f(\tau^X_u,\epsilon))}}{u(v^{-1}(\tau^X_u))} \right],
\end{align*}
et les calculs ne peuvent être terminés dans le cas général.

\textit{4.} Dans le cas du pont brownien $\beta_t = (1-t)B_{\frac{t}{1-t}}$, on a 
\[u(t)=1-t, \text{  } v(t)=\frac{t}{1-t} \text{  et  } v^{-1}(t) = \frac{t}{1+t}.\]
En particulier, $u(v^{-1}(t)) = \frac{1}{1+t}$. La formule additive du pont s'écrit donc
\[ \E\left[ \sum_{\lambda \leq l} F(\lambda, \omega, e_\lambda^\beta) \right] = \int_0^l \d u \int \n(\d \epsilon) \E\left[ (1+\tau^\beta_u) F(u,\omega, f(\tau^\beta_u,\epsilon)) \right]. \]

On peut également écrire, de façon synthétique la formule additive modifiée :
\begin{align*}
   \E\left[ \sum_{s \leq t} F(s, \omega, i_s^\beta) \right]
  = & \E\left[ \int_0^{\frac{t}{1-t}} \d_sL^0_s(B) \int \n(\d \epsilon) F\left(\frac{s}{1+s}, \omega, f\left(\frac{s}{1+s},\epsilon\right) \right) \right]  \\
  = & \int_0^\frac{t}{1-t} \dfrac{\d u}{\sqrt{2\pi u}} \int \n (\d \epsilon) \E\left[ \left. F\left( \frac{s}{1+s}, \omega, f\left(\frac{s}{1+s}, \epsilon\right) \right)\right| B_s = 0\right]\\
  = & \int_0^t \dfrac{\d u}{\sqrt{2\pi u(1-u)^3}} \int \n(\d \epsilon) \E\left[ \left. F\left( u, \omega, f(u, \epsilon) \right)\right| \beta_u = 0\right],\\
\end{align*}
en observant que, conditionnellement à $\beta_u=0$, le processus $(\beta_v,v \leq u)$ est un pont brownien de longueur $u$.

\textit{5.} Intéressons-nous maintenant au processus $Y_t$. En utilisant la formule de balayage, on a $Y_t = \int_0^t u(g_s) \d B_{s}$, donc $\crochet{Y}_t = \int_0^t u(g_s)^2 \d s$.

On calcule ensuite le temps local de $Y$ en 0. Par formule d'Itô-Tanaka, on a :
\[|Y_t| = \int_0^t \sgn(Y_s)u(g_s)\d B_{s} + L^0_t(Y),\]
que l'on peut également écrire :
\[|u(g_t)| |B_{t}| = \int_0^t |u(g_s)| \sgn(B_{s}) \d B_{s} + \int_0^t |u(g_s)|\d L^0_{s}(B),\]
d'où on tire immédiatement $L^0_t(Y) = \int_0^t |u(g_s)|\d L^0_{s}(B)$.

La théorie des excursions se développe de la même manière, on observe que l'excursion commençant à l'instant $s$ est dilatée d'un facteur $u(s)$.
\end{proof}

\begin{exercice}
Soit $\xi$ un processus de Lévy, et $\theta$ une fonction mesurable positive rendant intégrable la fonction définie par la suite.
\begin{enumerate}
  \item[1.] Déterminer $\lambda \geq 0$ tel que $M_t = \exp\left(\lambda t - \sum_{s \leq t} \theta(\Delta \xi_s)\right)$ est une martingale.
  \item[2.] Calculer la loi de $\xi$ sous $\P^\theta_{|\F_t} = M_t.P_{|\F_t}$.
\end{enumerate}

Soit $b$ une fonction mesurable et $(X_t)_{t \geq 0}$ la solution de l'équation différentielle stochastique :
\[X_t = X_0+ B_t + \int_0^t b(X_s)\d s.\]

\begin{enumerate}
  \item[3.] Trouver $\phi$ tel que $M_t = \phi(X_t) \exp(- \frac{1}{2}\int_0^t \frac{\phi''(X_s)}{\phi(X_s)}\d s)$ est la dérivée de Radon-Nikodym de la loi de $X$ par rapport à la mesure de Wiener.
  \item[4.] En déduire la mesure d'Itô des excursions du processus $(X_t)_{t \geq 0}$.
  \item[5.] En déduire la loi de $(\tau_l, l \geq 0)$ pour le processus d'Ornstein-Uhlenbeck de paramètre $\lambda$ (i.e. $b(x)=\lambda x$).
\end{enumerate}
\end{exercice}

\begin{proof}
\textit{1.} Soit $\xi$ un processus de Lévy, dont on note $\nu$ la mesure de Lévy associée. On calcule pour commencer :
\[ \E\left( \exp\left( - \sum_{s \leq t} \theta(\Delta \xi_s ) \right) \right) = \exp\left[ -t \int \nu(\d x) \left(1 - e^{-\theta(x)}\right) \right].\]
On observe alors que le processus 
\[M_t = \exp\left[ t \int \nu(\d x) \left(1 - e^{-\theta(x)}\right)  - \sum_{s \leq t} \theta(\Delta \xi_s)  \right],\]
est un processus à espérance constante, et de plus, on voit aisément que :
\[\E(M_{t+s}|\F_t) = M_t \E(M_s) = M_t.\]
Par conséquent, $M$ est une martingale.

\textit{2.} Observons pour commencer que la partie continue de $\xi$ est un processus stochastique indépendant de $(M_t)$, donc n'est pas modifiée par le changement de loi. On calcule la transformée de Fourier, sous $\P^\theta$ de $\xi^{(d)}$. On a :
\begin{align*}
  \E^\theta(\exp(i \lambda \xi^{(d)}_t))
  = & \E\left(M_t \exp(i \lambda \xi^{(d)}_t) \right)\\
  = & \exp\left[ t \int \nu(\d x) \left(1 - e^{-\theta(x)}\right) \right] \E\left[ \exp\left( - \sum_{s \leq t} \theta(\Delta \xi_s) - i \lambda \Delta \xi_s \right)  \right] \\
  = & \exp\left[ t \int \nu(\d x) \left(1 - e^{-\theta(x)}\right) - t \int \nu(\d x) \left(1 - e^{-\theta(x) + i \lambda x}\right)  \right].
\end{align*}
Par conséquent, sous $\P^\theta$, $\xi$ est encore un processus de Lévy, de mêmes dispersions et dérive que sous $\P$, et de mesure de Lévy associée donnée par $\n(\d x)e^{-\theta (x)}$.

\textit{3.} Nous allons raisonner de la même manière dans le cas du mouvement brownien, vu comme un processus de Lévy. Soit $\phi$ une fonction de classe $\mathcal{C}^2$, on a par formule d'Itô, sous $\W$ :
\[ M_t = \phi(X_t) \exp\left( -\frac{1}{2}\int_0^t \frac{\phi''(X_s)}{\phi(X_s)} \d s\right) = \phi(0) + \int_0^t \frac{\phi'(X_s)}{\phi(X_s)} M_s \d X_s. \]

En particulier, $(M_t, t \geq 0)$ est une martingale. Par conséquent, sous $\P^\phi = M . \W$, par formule de Girsanov, $X$ satisfait l'équation différentielle stochastique suivante :
\[ X_t = B_t + \int_0^t \frac{\phi'(X_s)}{\phi(X_s)} \d s.\]
Par conséquent, $X$ satisfait l'équation différentielle recherchée si $\phi$ satisfait l'équation différentielle :
\[ \phi' - b \phi = 0. \]
On note par conséquent $\phi(x) = \exp\left( \int_0^x b(y)\d y\right)$.

\textit{4.} L'existence d'une théorie des excursions de $X$ sous $\P^\phi$ est assez simple à vérifier. On note $(L_t)_{t \geq 0}$ le temps local en $0$ de $X$ et $(\tau_l)_{l \geq 0}$ l'inverse continu à droite de $L$. Pour tout $l \geq 0$, on pose $e_l(s) = B_{(\tau_{l-}+s)\wedge \tau_l}$. Remarquons que pour tout $0 \leq l \leq L_\infty$, on a :
\[X_{\tau_l}  = 0 = B_{\tau_l} + \int_0^{\tau_l} b(X_s) \d s,\]
par conséquent on a également :
\[X_{\tau_l+t} = B_{\tau_l+t}-B_{\tau_l} + \int_0^t b(X_{\tau_l + s}) \d s.\]

On obtient, par propriété de Markov forte, que conditionnellement à $L_\infty \geq l$ le processus $(X_{\tau_l +t})_{t \geq 0}$ est de même loi que $X$ et indépendant de $\F_{\tau_l}$. Il est aisé d'en déduire que le processus $(e_l)_{l \leq L_\infty}$ est un processus de Poisson ponctuel sur $\Omega^*$, de mesure de Lévy $\n_b$.

Afin d'obtenir cette mesure de Lévy, on utilise la formule multiplicative des excursions de deux manières différentes. Pour commencer, on a :
\[ \E^\phi\left( \exp\left( - \sum_{\lambda \leq l} \psi(e_\lambda) \right) \right) = \exp\left( - l \int \n_b(\d \epsilon) \left( 1 - e^{-\psi(\epsilon)} \right) \right).\]
D'un autre côté, on a également :
\begin{align*}
  &\E^\phi\left( \exp\left( - \sum_{\lambda \leq l} \psi(e_\lambda) \right) \right)\\
  = & \E\left( \exp\left( - \sum_{\lambda \leq l} \psi(e_\lambda) \right) \phi(X_{\tau_l}) \exp\left( - \frac{1}{2}\int_0^{\tau_l} \frac{\phi''(X_s)}{\phi(X_s)}\d s \right) \right)\\
  = & \E\left( \exp\left( - \sum_{\lambda \leq l} \psi(e_\lambda) + \int_0^{V(e_\lambda)} \frac{\phi''(e_\lambda(s))}{2\phi(e_\lambda(s))}\d s \right) \right)\\
  = & \exp\left( - l \int \n(\d \epsilon) \left( 1 - e^{-\psi(\epsilon) - \int_0^{V(\epsilon)} \frac{\phi''(\epsilon_s)}{2\phi(\epsilon_s)}\d s} \right)\right).
\end{align*}

Étant donné que l'on a en particulier 
\[\int \n(\d \epsilon) \left( 1 - e^{- \int_0^{V(\epsilon)} \frac{\phi''(\epsilon_s}{2\phi(\epsilon_s)}\d s}\right)=0,\]
on peut en déduire que $\n_b$ est à densité par rapport à $\n$, et que l'on a :
\[ \dfrac{\d \n_b}{\d \n}(\epsilon) = \exp\left( - \frac{1}{2}\int_0^{V(\epsilon)} b(\epsilon_s)^2 + b'(\epsilon_s) \d s \right). \]

\begin{remarque}
Ce calcul est en réalité exactement le pendant de celui qui a été réalisé précédemment.
\end{remarque}

\textit{5.} Intéressons-nous maintenant à la loi du processus $(\tau_l,l \geq 0)$ pour le processus d'Ornstein- Uhlenbeck. Grâce à la propriété de Markov forte, on observe immédiatement que $(\tau_l, l \geq 0)$ est un subordinateur.
Il reste donc à calculer la mesure de Lévy de $\tau_l$. Pour cela, on va calculer sa transformée de Laplace :
\begin{align*}
  \E(e^{-\mu \tau_l})
  = & \exp\left( - l \int \n_b(\d \epsilon) \left(1 - e^{-\mu V(\epsilon)} \right) \right)\\
  = & \exp\left[ - l \int \n(\d \epsilon) \left[1 - \exp\left(-\left(\mu-\frac{\lambda}{2}\right) V(\epsilon)- \int_0^{V(\epsilon)} \frac{\lambda^2}{2} \epsilon_s^2 + \lambda \d s \right)\right] \right]\\
  = & \exp\left( -l \int_0^{+\infty} \dfrac{e^{-\frac{\lambda}{2}v}\d v}{\sqrt{2\pi v^3}} \left( 1 -e^{-\mu v} \right) \right).
\end{align*}
La mesure de Lévy de $\tau_l$ est donc $\dfrac{e^{-\frac{\lambda}{2}v}\d v}{\sqrt{2\pi v^3}}$.
\end{proof}

\begin{exercice}
Montrer que pour une diffusion à valeurs positives, la mesure $n_M$ image de la mesure de Lévy par $\epsilon \mapsto M(\epsilon)$ est donnée par :
\[n_M([a,+\infty)) = \n(M \geq a) = \frac{1}{h(a)},\]
où $h$ est une bonne fonction d'échelle associée à $X$.
\end{exercice}

\begin{proof}
Nous allons démontrer ce résultat en nous basant sur la démonstration réalisée dans le cas du mouvement brownien. En effet, on a, pour tout $a \geq 0$ :
\[\P(X_{\tau_l^X} \leq a) = \P(h(X_{\tau_l^X} \leq h(a)) = \P(T_{h(a)}(h(X)) \geq \tau_l^{h(X)}),\]
en choisissant judicieusement la fonction d'échelle comme étant nulle et de dérivée 1 en 0. On utilise alors pour $Y_t = h(X_t)$ une martingale locale positive, la propriété suivante :
\[\P(T_{h(a)} \geq \tau_l) = \P(L_{T_{h(a)}} \geq l) = \exp\left(- \frac{l}{h(a)}\right). \]
En effet, par utilisation de la propriété de Markov, on observe que :
\[\P(L_{T_{(h(a)}} \geq l+l' | L_{T_{(h(a)}} \leq l) = \P(L_{T_{(h(a)}} \geq l'), \]
donc la loi de cette variable aléatoire est bien exponentielle, on utilise alors le calcul :
\[\E(L_{T_{h(a)}} ) = \E(|Y_{T_{h(a)}}|) = h(a),\]
pour trouver le paramètre de cette variable aléatoire exponentielle.
On en déduit, en utilisant la formule multiplicative des excursions que l'on a bien
\[\n(M \geq a) = \frac{1}{h(a)}.\]
\end{proof}

\end{document}